\documentclass[11pt,a4paper]{article}
\usepackage[T1]{fontenc}
\usepackage[utf8]{inputenc}
\usepackage{lmodern}
\usepackage[margin=28mm,headheight=14pt]{geometry}
\usepackage{amsmath,amssymb,amsthm,mathtools}
\usepackage{microtype,booktabs,longtable,array,tabularx,enumitem}
\usepackage{xcolor,xurl}
\usepackage[pdfusetitle,colorlinks=true,linkcolor=blue!40!black,citecolor=blue!40!black,urlcolor=blue!40!black]{hyperref}
\hypersetup{pdftitle={The shortest figure-eight geodesics on random hyperbolic surfaces},pdfsubject={Poisson critical-window statistics, cuff-shape marks and exponential moments of the non-simple systole}}
\numberwithin{equation}{section}
\newtheorem{theorem}{Theorem}[section]
\newtheorem{proposition}[theorem]{Proposition}
\newtheorem{lemma}[theorem]{Lemma}
\newtheorem{corollary}[theorem]{Corollary}
\theoremstyle{definition}
\newtheorem{definition}[theorem]{Definition}

\theoremstyle{remark}
\newtheorem{remark}[theorem]{Remark}
\allowdisplaybreaks[1]
\setlist{itemsep=2pt,topsep=4pt}
\title{The shortest figure-eight geodesics\\on random hyperbolic surfaces}
\author{Guan HUANG}
\date{}
\begin{document}
\maketitle
\begin{abstract}
We consider primitive unoriented figure-eight geodesics on a closed hyperbolic surface of genus $g$ sampled from Weil--Petersson measure. The shortest of these geodesics realizes the non-simple systole. He, Shen, Wu and Xue proved that its length is tight after centering at $a_g=\log g-\log\log g$. We prove that, as $g\to\infty$, the point process of centered lengths converges to a Poisson process on $\mathbb R$ with intensity $e^s\,ds/(8\pi^2)$. Its first arrival gives a reflected Gumbel limit for the centered systole $M_g$.

We also prove convergence of the marked process in which each figure-eight carries the three cuff lengths of its supporting pair of pants, normalized by their sum. In the limiting process, the third cuff accounts for half the perimeter. The two winding cuffs, ordered at random, divide the other half uniformly. These splits are independent across points and independent of the limiting lengths. Uniform exponential bounds on both tails of $M_g$ give convergence of all fixed polynomial moments. At the critical exponent $-1$, we obtain $\mathbb E_g e^{-M_g}\sim(\log g)/(16\pi^2)$.
\end{abstract}
\tableofcontents

% BEGIN sections/01_introduction.tex
\section{Introduction and main results}\label{sec:introduction}

On a closed hyperbolic surface, a shortest primitive closed geodesic with a self-intersection is a figure-eight~\cite[Theorem~4.2.4]{Buser2010}. On a Weil--Petersson random surface of genus $g$, He, Shen, Wu and Xue~\cite{HSWX} located its length near $\log g-\log\log g$ and computed the expected number of figure-eights up to a given length. Their estimate suggests a question about dependence. How are the lengths of distinct figure-eights distributed within a fixed distance of that scale? We prove that the centered lengths form a Poisson process in the limit. We also determine the limiting cuff proportions of the pants supporting the first arrivals and use tail estimates to obtain the moments of the centered minimum.

\subsection{Counting conventions and the point process}
Let $\mathcal M_g$ denote the moduli space of closed oriented hyperbolic surfaces of genus $g\ge3$. Write $V_g$ for its Weil--Petersson volume, $\mathbb P_g$ for the corresponding probability measure, and $\mathbb E_g$ for expectation. We sample $X$ according to $\mathbb P_g$ and denote the length of a geodesic $\gamma$ on $X$ by $\ell_X(\gamma)$. We count primitive unoriented closed geodesics, identifying the two orientations and excluding iterates. A \emph{figure-eight} has exactly one transverse self-intersection. Distinct geodesics of the same length contribute separately.

Put
\begin{equation}\label{i:center}
H=\log g,\qquad a_g=H-\log H,\qquad
\Xi_g^{(8)}=\sum_{\gamma\text{ a figure-eight on }X}
\delta_{\ell_X(\gamma)-a_g}.
\end{equation}
The random measure $\Xi_g^{(8)}$ takes values in $\mathcal M_p(\mathbb R)$, the space of locally finite counting measures equipped with the vague topology; repeated atoms are allowed. Local finiteness follows from finiteness of the length spectrum below each fixed bound on a compact surface. To see measurability, enumerate the primitive unoriented free-homotopy classes with self-intersection number one on a reference topological surface, and use their continuous length functions on Teichm\"uller space. The resulting point measure descends to $\mathcal M_g$ because the mapping class group permutes these classes. Section~\ref{mk3:subsec:marking} gives the corresponding construction with persistent marks.

The first-moment estimate~\eqref{i:external-mean} gives $\mathbb E_g\Xi_g^{(8)}((-\infty,s])\to e^s/(8\pi^2)$ for each fixed $s$. This identifies the intensity in the proposed Poisson law.
\begin{equation}\label{i:intensity}
\lambda(ds)=\frac{e^s}{8\pi^2}\,ds,
\qquad \Lambda(s)=\lambda((-\infty,s])=\frac{e^s}{8\pi^2},
\qquad \Pi\sim\operatorname{PPP}(\lambda).
\end{equation}
Here $\operatorname{PPP}(\lambda)$ denotes the law of a Poisson point process with intensity $\lambda$, and $\mathcal L_g$ denotes distribution under $\mathbb P_g$. The symbol $\Longrightarrow$ denotes convergence in distribution in the indicated space. For probability measures $P,Q$ on the same measurable space, we use
$d_{\rm TV}(P,Q)=\sup_A|P(A)-Q(A)|$, with the supremum over measurable sets.

\begin{theorem}[Critical-window process]\label{thm:main}
As $g\to\infty$,
\begin{equation}\label{i:ppp}
\Xi_g^{(8)}\Longrightarrow\Pi
\quad\text{in }\mathcal M_p(\mathbb R)\text{ with the vague topology}.
\end{equation}
For every fixed finite family of pairwise disjoint bounded half-open intervals $I_1,\ldots,I_d$,
\begin{equation}\label{i:tv}
d_{\rm TV}\left(
\mathcal L_g\bigl((\Xi_g^{(8)}(I_i))_{i=1}^d\bigr),
\bigotimes_{i=1}^d\operatorname{Poisson}(\lambda(I_i))
\right)\longrightarrow0.
\end{equation}
For every fixed finite set of real thresholds $s_1,\ldots,s_d$, the vector
$(\Xi_g^{(8)}((-\infty,s_i]))_{i=1}^d$ also converges in total variation to the corresponding count vector of $\Pi$.
\end{theorem}

For overlapping intervals, refine their union into disjoint cells and sum the corresponding counts. Thus Theorem~\ref{thm:main} also gives total variation convergence for any fixed finite vector of bounded interval counts, with the dependence prescribed by $\Pi$.

Write $\ell^{(8)}_{g,1}\le\ell^{(8)}_{g,2}\le\cdots$ for the ordered figure-eight lengths, counted with multiplicity whenever the entries exist. For each fixed $k$, every surface has at least $k$ such curves when $g\ge k+2$; see Section~\ref{sec:probability}. Let $\ell_{\rm sys}^{\rm ns}$ denote the length of the shortest genuinely non-simple closed geodesic, excluding repeated traversals of simple geodesics.

\begin{theorem}[First arrivals and the non-simple systole]\label{thm:arrivals}
For each fixed $k\ge1$,
\begin{equation}\label{i:arrivals}
(\ell^{(8)}_{g,j}-a_g)_{j=1}^k
\Longrightarrow (\log(8\pi^2\Gamma_j))_{j=1}^k,
\qquad \Gamma_j=E_1+\cdots+E_j,
\end{equation}
where the $E_i$ are independent exponential variables of rate one. In particular, if $M_g=\ell_{\rm sys}^{\rm ns}-a_g$, then
\begin{equation}\label{i:minimum}
\lim_{g\to\infty}\mathbb P_g(M_g>s)
=\exp\left(-\frac{e^s}{8\pi^2}\right),\qquad s\in\mathbb R.
\end{equation}
Thus $\log(8\pi^2)-M_g$ converges to the standard Gumbel law with distribution function $x\mapsto\exp(-e^{-x})$.
\end{theorem}

The first arrival is the non-simple systole. The subsequent arrivals in~\eqref{i:arrivals} are ordered within the figure-eight spectrum. Section~\ref{i:subsec:moments} states the corresponding moment results.

\subsection{Cuff shapes and marked extreme statistics}

The two loops of a figure-eight distinguish two of the three boundary components of its intrinsic supporting pair of pants. We call these the winding cuffs and keep all three cuffs distinct even if two have the same image on $X$ (Section~\ref{sec:foundations}). To record a shape, we label the winding cuffs in a uniformly random order, independently for each curve conditional on $X$. The order is chosen once for each curve and used at every length cutoff. If $x_\gamma,y_\gamma$ are the winding cuff lengths in this order and $z_\gamma$ is the third cuff length, the mark records their proportions of the total cuff length.
\begin{equation}\label{mk:full}
 \begin{gathered}
 \Delta=\{(z_1,z_2,z_3)\in[0,1]^3:z_1+z_2+z_3=1\},
 \qquad E=\mathbb R\times\Delta,\\
 Z_\gamma=\frac{(x_\gamma,y_\gamma,z_\gamma)}
                    {x_\gamma+y_\gamma+z_\gamma},\qquad
 \mathsf N_g=\sum_{\gamma\text{ a figure-eight}}
       \delta_{(\ell_X(\gamma)-a_g,Z_\gamma)}.
 \end{gathered}
\end{equation}
Section~\ref{mk3:subsec:marking} constructs these marks measurably on the whole countable curve set. Write $\mathbb P_g^\circ$, $\mathbb E_g^\circ$ and $\mathcal L_g^\circ$ for probability, expectation and distribution under the enlarged law. Projection onto the length coordinate recovers the unmarked process:
$\mathsf N_g(B\times\Delta)=\Xi_g^{(8)}(B)$ for every Borel $B\subset\mathbb R$.

Let
\begin{equation}\label{mk:intensity}
 \begin{gathered}
 \nu=\mathcal L\left(U/2,(1-U)/2,1/2\right),
 \quad U\sim\operatorname{Unif}[0,1],\\
 \mu(ds,dz)=\lambda(ds)\nu(dz),\qquad
 \Pi^{\mathrm{mark}}\sim\operatorname{PPP}(\mu).
 \end{gathered}
\end{equation}
In the limiting shape law, the third cuff contributes half the total cuff length, while the winding cuffs divide the other half uniformly. In the next theorem, $\partial_\Delta A$ denotes the boundary of $A$ relative to $\Delta$, and $C_c(E)$ is the space of continuous compactly supported functions on $E$.

\begin{theorem}[Marked critical-window process]\label{mk:main}
As $g\to\infty$,
\begin{equation}\label{mk:vague}
 \mathsf N_g\Longrightarrow\Pi^{\mathrm{mark}}
 \quad\text{in }\mathcal M_p(E)\text{ with the vague topology}.
\end{equation}
For every fixed finite family $C_j=I_j\times A_j$, where
$A_j\subset\Delta$ is Borel with
$\nu(\partial_\Delta A_j)=0$ and either $I_j=(u_j,v_j]$ with finite
$u_j<v_j$ or $I_j=(-\infty,v_j]$, one has
\begin{equation}\label{mk:count-tv}
 d_{\mathrm{TV}}\!\left(
  \mathcal L_g^\circ((\mathsf N_g(C_j))_{j=1}^d),
  \mathcal L((\Pi^{\mathrm{mark}}(C_j))_{j=1}^d)\right)
       \longrightarrow0.
\end{equation}
For pairwise disjoint cells the limiting coordinates are independent
Poisson variables with means
$\mu(C_j)=(e^{v_j}-e^{u_j})\nu(A_j)/(8\pi^2)$, where
$e^{-\infty}=0$; zero means are allowed. For overlapping cells the
target vector has the dependence induced by their common parts.
For every nonnegative $f\in C_c(E)$,
\begin{equation}\label{mk:laplace}
 \mathbb E_g^\circ e^{-\mathsf N_g(f)}
 \longrightarrow
 \exp\!\left[-\int_{\mathbb R}\int_\Delta
       (1-e^{-f(s,z)})\,\nu(dz)\,\frac{e^s}{8\pi^2}\,ds\right].
\end{equation}
\end{theorem}

Attach the corresponding marks $Z_{g,j}$ to the lengths
$\ell^{(8)}_{g,j}$ ordered with multiplicities. At length ties use any
fixed measurable rule, for example lexicographic order of the mark
coordinates; identical marked atoms require no distinction. Put
$T_{g,j}=\ell^{(8)}_{g,j}-a_g$.

\begin{theorem}[Marked arrivals]\label{mk:arrivals}
Let $\Gamma_j=E_1+\cdots+E_j$ as in Theorem~\ref{thm:arrivals}, and
let $U_1,U_2,\ldots$ be independent uniform variables on $[0,1]$,
independent of the exponential sequence. For every fixed $k\ge1$,
\begin{equation}\label{mk:arrival-vector}
 ((T_{g,j},Z_{g,j}))_{j=1}^k
 \Longrightarrow
 \left(\left(\log(8\pi^2\Gamma_j),
      \left(\frac{U_j}{2},\frac{1-U_j}{2},\frac12\right)\right)\right)_{j=1}^k
\end{equation}
weakly in $E^k$, independently of the chosen tie rule. In particular,
for every $s\in\mathbb R$ and every Borel $A\subset\Delta$ with
$\nu(\partial_\Delta A)=0$,
\begin{equation}\label{mk:shortest}
 \mathbb P_g^\circ\{T_{g,1}>s,\ Z_{g,1}\in A\}
       \longrightarrow e^{-\Lambda(s)}\nu(A).
\end{equation}
\end{theorem}

The random choice of order can be removed by sorting the winding cuffs. The continuous map
\begin{equation}\label{mk:quotient}
 Q(z_1,z_2,z_3)=(\min(z_1,z_2),\max(z_1,z_2),z_3)
\end{equation}
leaves the third cuff in its distinguished position. For each fixed number of arrivals, the sorted marks converge jointly to independent copies of $(V/2,(1-V)/2,1/2)$, where $V=\min(U,1-U)$ has density $2$ on $[0,1/2]$. These limiting marks are independent of the limiting arrival lengths.

The total variation statement applies to counts in fixed continuity cells. The ordered marked arrivals converge weakly, and the marked process converges in distribution for the vague topology. The full marked-process laws have total variation distance one at every $g\ge3$ (Proposition~\ref{mk3:prop:tv-fails}).

\subsection{Exponential and polynomial moments of the minimum}\label{i:subsec:moments}

To pass from the limiting distribution of the minimum to its moments, we need uniform control of both tails. Put $q=(8\pi^2)^{-1}$ and $M=\log(\mathcal E/q)$, where $\mathcal E$ has the rate-one exponential law.

For $p\ge1$, write
\[
 W_p(P,Q)=\inf_{\pi\in\mathcal C(P,Q)}
       \left(\int_{\mathbb R^2}|x-y|^p\,d\pi(x,y)\right)^{1/p},
\]
where $\mathcal C(P,Q)$ is the set of probability couplings with marginals
$P,Q$.

\begin{theorem}[Exponential moments, polynomial moments and Wasserstein limits]\label{thm:all-moments}
There is a universal $\kappa>0$ such that
$\sup_{g\ge3}\mathbb E_g e^{\theta|M_g|}<\infty$ for every
$0<\theta<\kappa$. For real $s<\kappa$,
\[
 \lim_{g\to\infty}\mathbb E_g e^{sM_g}
  =\begin{cases}
    (8\pi^2)^s\Gamma(1+s),&-1<s<\kappa,\\
    +\infty,&s\le-1.
  \end{cases}
\]
At the endpoint $s=-1$, the rate of divergence is
$\mathbb E_g e^{-M_g}\sim(\log g)/(16\pi^2)$ as $g\to\infty$.
For every fixed real $r>0$, there is a finite constant $C_r$ such that
\begin{equation}\label{i:all-moment-bound}
 \sup_{g\ge3}\mathbb E_g|M_g|^r\le C_r.
\end{equation}
For every fixed real $p>0$ and every $R>0$,
\begin{equation}\label{i:all-moment-ui}
 \sup_{g\ge3}\mathbb E_g\!\left[
 |M_g|^p\mathbf1_{\{|M_g|>R\}}\right]\le \frac{C_{p+1}}R.
\end{equation}
In particular, $\{|M_g|^p:g\ge3\}$ is uniformly integrable, and
\begin{equation}\label{i:all-moment-limit}
 \mathbb E_g|M_g|^p\longrightarrow
 \mathbb E|M|^p=\int_0^\infty |\log(u/q)|^p e^{-u}\,du.
\end{equation}
More generally, $\mathbb E_g f(M_g)\to\mathbb E f(M)$ for every
continuous $f:\mathbb R\to\mathbb R$ satisfying
$|f(x)|\le C_f(1+|x|^p)$ for some fixed finite $p>0$.
For every fixed $1\le p<\infty$,
\begin{equation}\label{i:all-wasserstein}
 W_p(\mathcal L_g(M_g),\mathcal L(M))\longrightarrow0.
\end{equation}
In particular,
\begin{equation}\label{i:second-variance}
 \begin{gathered}
 \mathbb E_g M_g^2\longrightarrow
       [\log(8\pi^2)-\gamma_{\!E}]^2+\frac{\pi^2}{6},\\
 \operatorname{Var}_g(\ell_{\rm sys}^{\rm ns})\longrightarrow\frac{\pi^2}{6}.
 \end{gathered}
\end{equation}
The constants $C_r$ may depend on the fixed polynomial order $r$.
\end{theorem}

The expectation and stop-loss limits follow from the case $p=1$. For $x>0$, define
\begin{equation}\label{i:exponential-integrals}
 \operatorname{Ein}(x)=\int_0^x\frac{1-e^{-u}}u\,du,
 \qquad \operatorname E_1(x)=\int_x^\infty\frac{e^{-u}}u\,du.
\end{equation}
The integrand defining $\operatorname{Ein}$ is extended continuously at zero. The function $\operatorname E_1$ is distinct from the exponential random variable $E_1$ in~\eqref{i:arrivals}. Write $x_+=\max\{x,0\}$. The quantities $\mathbb E_g[(-M_g-R)_+]$ and $\mathbb E_g[(M_g-R)_+]$ below are the integrated negative and positive tails, also called stop-loss functions.

\begin{corollary}[Expected non-simple systole]\label{thm:expected-minimum}
\begin{equation}\label{i:expected-minimum}
 \mathbb E_g\ell_{\rm sys}^{\rm ns}
 =\log g-\log\log g+\log(8\pi^2)-\gamma_{\!E}+o(1),
\end{equation}
where $\gamma_{\!E}$ denotes Euler's constant. The two stop-loss functions satisfy
\begin{align}
 \sup_{R\ge0}\left|\mathbb E_g[(-M_g-R)_+]-\operatorname{Ein}(qe^{-R})\right|&\longrightarrow0,\label{i:negative-stop}\\
 \sup_{R\ge0}\left|\mathbb E_g[(M_g-R)_+]-\operatorname E_1(qe^R)\right|&\longrightarrow0.\label{i:positive-stop}
\end{align}
\end{corollary}

The constant term in~\eqref{i:expected-minimum} is the mean of the limiting reflected Gumbel law. The tail estimates in Section~\ref{sec:tails} justify taking expectations in the distributional limit.

\subsection{Relation to previous work}
Nie, Wu and Xue~\cite[Theorem~4]{NWX} showed that $(1-\varepsilon)\log g<\ell_{\rm sys}^{\rm ns}<2\log g$ with probability tending to one for every fixed $\varepsilon>0$. He, Shen, Wu and Xue~\cite[Theorems 1.1--1.2]{HSWX} then located the non-simple systole in $a_g\pm\omega(g)$ with probability tending to one whenever $\omega(g)\to\infty$ and $\omega(g)=o(\log\log g)$. Since $\omega(g)$ may grow arbitrarily slowly, this gives tightness of $M_g$. They also proved $\mathbb E_g\ell_{\rm sys}^{\rm ns}/\log g\to1$.

Let $N_{\rm f8}(X,L)$ count primitive unoriented figure-eights on $X$ of length at most $L$. With $c_0=3+4\log2$, their first-moment formula is
\begin{equation}\label{i:external-mean}
\mathbb E_g N_{\rm f8}(X,L)
=\frac{(L-c_0)e^L}{8\pi^2g}
+O\left(\frac{L^2e^{L/2}}g+\frac{L^4e^L}{g^2}\right),
\end{equation}
uniformly for $g>2$ and $L\ge2\operatorname{arccosh}3$, with an absolute implied constant~\cite[Proposition 3.1]{HSWX}. At $L=a_g+s$ it gives the limiting mean $\Lambda(s)$ in~\eqref{i:intensity}. Our contribution at this scale is the joint counting law. For the minimum, the first moment bounds the negative tail. To obtain the centered expectation and the exponential moments we also estimate the positive tail in Section~\ref{sec:tails}.

He, Shen, Wu and Xue use a lower cutoff of $10\log L$ on the three cuffs in their second-moment argument. They distinguish coincident supports, disjoint closures, distinct supports with overlapping interiors, and contacts along cuffs alone~\cite[Definition~4.2 and equation~(1.8)]{HSWX}. We use the relations between a support and each member of a tuple to estimate how many disjoint tuples of each fixed order are affected when a curve with a neighbor is removed (Theorem~\ref{thm:weighted-loss}).

Mirzakhani's integration formula~\cite{Mirzakhani} expresses curve counting integrals in terms of Weil--Petersson volumes. We use the large-genus asymptotics of Mirzakhani and Zograf~\cite{MZ} and the boundary-length estimates of Nie, Wu and Xue~\cite[Lemma 22]{NWX}. Anantharaman and Monk~\cite[Theorem 5.7 and Section 5.5]{AMI} treat ordered multi-loops, including joint length tests on disjoint filled supports. We use their integration conventions with the multiplicities for unoriented figure-eights. Their cuff coordinates also yield the leading one-point winding split (Section~\ref{sec:marked-intensity}). Our marked process identifies the cuff shapes of several curves selected by their order in the figure-eight length spectrum. Anantharaman and Monk's second paper~\cite[Theorems 1.4, 1.15 and 1.6]{AMII} establishes the weak Friedman--Ramanujan property for general non-simple local types and applies it to the spectral gap.

Interactions between short curves also arise in estimates for tangles. Monk and Thomas~\cite{MT} prove that a Weil--Petersson random surface is $(a\log g)$-tangle-free with high probability for every fixed $0<a<1$. Anantharaman and Monk express the indicator of the complement of their tangle-free set as a weighted sum over subsurfaces and simple closed geodesics~\cite[Theorem 1.2 and Corollary 1.3]{AMMobius}. For length tests supported on $[0,6\log g]$, they also estimate a weighted geodesic count multiplied by the number of unordered pairs of distinct short tangles~\cite[Lemma~4.1]{AMGap}. Here a surface with intersecting supports can carry many disjoint tuples, so the factorial-moment argument estimates their expected loss directly.

Other counting results use different orders of limits. Le Guilloux~\cite{LG25,LG26} integrates length tests over mapping-class orbits, first for figure-eights and then for general closed loops. The long-length asymptotics in~\cite[Theorem 1.1 and Corollary 1.2]{LG26} fix the genus, the number of boundaries and their lengths. Liu, Rafi, Souto and Trin~\cite[equation (1.1)]{LRST} study large-genus asymptotics of orbit frequencies defined by first sending the length bound to infinity at fixed genus. Our length bound $a_g+s$ grows with the genus.

For closed Weil--Petersson surfaces, Mirzakhani and Petri~\cite[Theorem 4.1]{MP19} prove joint Poisson limits for primitive geodesic counts in fixed bounded length intervals. Rudnick~\cite[Theorem 1.1]{RudnickSI} obtains short-interval second-moment asymptotics for weighted geodesic counts, and Leung~\cite[Theorem 1.1]{Leung} proves a central limit theorem. Both authors count orientations and iterates with a primitive-length weight. They first take the large-genus limit and then let the norm cutoff grow, where the norm is the exponential of the geodesic length.

For surfaces with a growing number of cusps, Budd and Lions~\cite[Theorem 1.2]{BL} obtain a Poisson limit for rescaled primitive tight geodesics in a Boltzmann ensemble. In their proof, surfaces with intersecting or self-intersecting short tight geodesics contribute $o(1)$ to the expected number of tuples of disjoint simple tight geodesics~\cite[Proposition~5.3]{BL}. For Weil--Petersson surfaces with $n=n(g)$ cusps, Gatignol and Thomas~\cite[Theorems 1.1--1.3]{GT26} obtain Poisson laws in fixed or rescaled windows. When $n/\sqrt g\to c\in(0,\infty)$, their Theorem~1.3 gives independent Poisson counts on disjoint intervals with upper endpoints below $2\operatorname{arcsinh}(1)$. In the same regime Wang et al.~\cite[Theorem 14.10]{WangEtAl26} identify the limiting primitive spectrum on bounded windows as the sum of a Poisson background and an independent cluster process. Each cluster is the full primitive spectrum of a two-cusp pair of pants.

\subsection{Outline of the proof}
The proof begins with tuples of figure-eights whose supports have disjoint closures. Their expected counts admit a joint volume integral. To compare these counts with those of the full process, fix a centered length ceiling $S\in\mathbb R$ and set
\begin{equation}\label{i:cutoff}
L=a_g+S,\qquad b=10\log L,
\end{equation}
with $g$ sufficiently large. Retain curves of length at most $L$ whose supporting pants have embedded ambient closures, three distinct cuff geodesics, connected individual complements, and all cuff lengths at least $b$. Distinct cuffs may have equal lengths.

Form the conflict graph on the retained curves, joining two vertices when their ambient closed supports meet. Let $R_{g,S}$ count the discarded figure-eights of length at most $L$, let $B_{g,S}$ count the unordered conflict edges, and let $\widehat\Xi_{g,S}$ place atoms at the centered lengths of isolated vertices. Proposition~\ref{prop:trim} and Proposition~\ref{prop:conflicts} give
\begin{equation}\label{i:coupling-errors}
\mathbb E_gR_{g,S}=O_S(\log H/H),\qquad
\mathbb E_gB_{g,S}=O_S(H^{-2}).
\end{equation}
Isolation is tested against every retained curve at or below the ceiling, including those outside the observation windows. Thus full and isolated counts agree with probability tending to one. A single curve removed at the isolation step can belong to many disjoint tuples, so we also estimate how many disjoint tuples isolation removes.

Consider ordered retained $m$-tuples $\tau$ whose ambient closed supports are pairwise disjoint, allowing every topology of the joint complement. Let $\mathcal W(\tau)$ consist of retained curves outside the tuple whose closed supports meet a selected support; call these curves witnesses, and put $d(\tau)=|\mathcal W(\tau)|$. Denote by $Q_m$ the number of such tuples with $d(\tau)>0$. Theorem~\ref{thm:weighted-loss} proves
\begin{equation}\label{i:local-loss}
\mathbb E_gQ_m=O_{m,S}(H^{-2}+g^{-1})
\quad\text{for every fixed }m\ge1.
\end{equation}
An affected tuple may have several witnesses. Giving each one an equal share counts the tuple exactly once, as expressed by the identity
\begin{equation}\label{i:normalized-witness}
Q_m=\sum_\tau\sum_{\beta\in\mathcal W(\tau)}\frac1{d(\tau)},
\end{equation}
where the outer sum is over disjoint tuples and the inner sum is zero when $d(\tau)=0$. Section~\ref{sec:disjoint} bounds the contribution from tuples with disconnected complement. For a connected complement, Sections~\ref{sec:witnesses}--\ref{sec:higher} group witnesses according to which supports they meet and how they meet them. For a crossing witness, the subsurface filled by the witness and the supports it crosses constrains the positions of the other selected pants. Cores of bounded complexity can then be summed type by type. A separate subsurface estimate handles growing complexity. The small-core pair bounds used here come from Lemmas~\ref{tail:planar}--\ref{tail:small}. Those lemmas use only local geometry and volume estimates, so they can be read before the moment argument in Section~\ref{sec:tails}.

Theorem~\ref{thm:disjoint} compares disjoint tuple measures with products of the exact retained one-point measure and yields the limiting product moments determined by $\lambda$. Subtracting the expected loss~\eqref{i:local-loss} gives the mixed factorial-moment limits for the isolated process. Section~\ref{sec:probability} then combines these limits with the coupling~\eqref{i:coupling-errors} and a left-tail truncation to prove the unmarked process and arrival theorems.

Section~\ref{sec:marked-intensity} determines the limiting cuff shape under the normalized expected counting measure. Section~\ref{sec:marked-factorials} incorporates cuff weights in the tuple comparison. In Section~\ref{sec:marked-process} we keep the winding order of each curve fixed across cutoffs, pass to the marked process on sets whose boundaries have zero limiting mass, and order the curves by length to obtain the marked arrivals.

The two tails of the minimum enter differently in Section~\ref{sec:tails}. On the short side, a pair bound complements the first-moment upper bound and gives $\mathbb P_g(\ell_{\rm sys}^{\rm ns}\le uH)\sim uHe^{uH}/(8\pi^2g)$ for fixed $0<u<1$. Integration yields the coefficient $1/(16\pi^2)$ at the critical exponent. On the long side, a surface with $M_g>t$ has no figure-eight below $L=a_g+t$. We bound the probability of this event through the second moment of the count on good supports with all cuffs at least $L/5$. This gives $\mathbb P_g(M_g>t)\le Ce^{-t}$ for $0\le t\le H/20$ and all sufficiently large $g$. Monotonicity and the deterministic $O(H)$ bound on the minimum extend the exponential estimate to all $t\ge0$ with a smaller exponent.

Appendices~\ref{app:genus}--\ref{app:planar} supply the complementary-volume sums and the integrals for attached pants and planar cores. Appendix~\ref{app:sources} gives the conventions for the imported estimates. Unless stated otherwise, implicit constants may depend on the fixed parameters in their subscripts and are uniform in $g$.

% END sections/01_introduction.tex

% BEGIN sections/02_foundations.tex
\section{Trimming and pair interactions}
\label{sec:foundations}

Below a fixed centered length ceiling, we delete curves in two steps.
First we discard figure-eights with exceptional supports or short cuffs.
Among those retained, we then delete every curve whose closed support
meets that of another retained curve. The survivors have pairwise
disjoint closed supports. To compare their point measure with the
original process, we bound the expected number discarded at the first
step and the expected number of pairs whose supports meet at the second.

Throughout this section, $X$ is a closed oriented hyperbolic surface of
genus $g\geq3$, with normalized Weil--Petersson law $\mathbb P_g$.
Write $\mathbb E_g$ for expectation, $V_g$ for the unnormalized volume,
and $V_{h,n}(\mathbf x)$ for the volume with labeled boundary lengths
$\mathbf x$. We abbreviate $V_{h,n}=V_{h,n}(\mathbf0)$ and write
$S_{h,n}$ for the orientable surface of genus $h$ with $n$ boundary
components. Figure-eights are counted as primitive unoriented geodesics
with exactly one transverse self-intersection; distinct geodesics of
equal length are counted separately.

We write $L$ for the upper length cutoff and $b$ for the lower cuff
cutoff. The estimates below give their admissible ranges. For the
critical window we take
\begin{equation}
 H=\log g,\qquad a_g=H-\log H,\qquad
 L=a_g+S,\qquad b=10\log L,
 \label{f:critical-parameters}
\end{equation}
where $S\in\mathbb R$ is fixed and $g$ is sufficiently large. All
constants implicit in $O_S(\cdot)$ may depend on $S$, but not on $g$.

\subsection{Supports and the one-pants integration formula}

Let $P_\gamma$ be the interior of the geodesic pants support of a
figure-eight $\gamma$. Its intrinsic metric completion is a pair of
pants; its closure in $X$ may identify two boundary copies. We use
$\overline P_\gamma$ for the latter closure when discussing intersections
in $X$. A support is called \emph{good} if the natural map from its
intrinsic completion to $X$ is an embedding, its three cuffs are distinct
ambient geodesics, and the complement of its image is connected. Distinct cuffs need not have different numerical lengths.

The support description in
\cite[Section 2.4 and Section 3]{HSWX} gives the following facts. Each
figure-eight has a unique pants support. On an intrinsic pair of pants
there are exactly three primitive unoriented figure-eights, one for
each unordered pair of winding cuffs. If the winding cuffs have lengths
$x,y$ and the remaining cuff has length $z$, its length $\ell(x,y,z)$
satisfies Rivin's identity~\cite[Theorem~1.3]{Rivin}
\begin{equation}
 \cosh\frac{\ell(x,y,z)}2
 =\cosh\frac z2+2\cosh\frac x2\cosh\frac y2.
 \label{f:pants-length}
\end{equation}
In particular, every figure-eight has length at least
$L_*=2\operatorname{arccosh}3$.

\begin{lemma}[Support types]
\label{lem:f-support-types}
The following four classes partition the figure-eights on $X$:
\begin{enumerate}
 \item good supports, with complement of type $S_{g-2,3}$;
 \item supports whose ambient closure is a one-holed torus obtained
 by identifying two cuffs, with complement $S_{g-1,1}$;
 \item embedded supports with complement
 $S_{h_1,1}\sqcup S_{h_2,2}$, where $h_1,h_2\geq1$ and
 $h_1+h_2=g-1$;
 \item embedded supports with complement
 $S_{h_1,1}\sqcup S_{h_2,1}\sqcup S_{h_3,1}$, where each $h_i\geq1$
 and $h_1+h_2+h_3=g$.
\end{enumerate}
Complementary components are understood with their metric completions.
\end{lemma}

\begin{proof}
The support construction in \cite[Section 3, equation (3.1)]{HSWX}
allows precisely the embedded-pants closure and the identified-cuff
one-holed-torus closure. In the embedded case, distributing the three
distinct cuffs among complementary components gives the boundary
partitions $3$, $1+2$, and $1+1+1$. No complementary component is a
geodesic-boundary disk or annulus: Gauss--Bonnet would give nonpositive
area. Thus a component with one or two boundaries has positive genus.
Euler characteristic gives the displayed genus relations. The
one-holed-torus closure has one boundary and its complement is connected
of genus $g-1$. These alternatives are mutually exclusive and exhaustive.
\end{proof}

Let $N_T(X,L),N_{D_2}(X,L),N_{D_3}(X,L)$ count the actual curves of
length at most $L$ in the last three classes. The counts include every
admissible complementary genus; the three genera in the last class can
be listed in nonincreasing order so that its topological types are not
repeated. We use the following estimates of He, Shen, Wu and Xue.
\begin{equation}
 \begin{aligned}
 \mathbb E_gN_T(X,L)&\leq C_Tg^{-1}L^2e^{L/2},\\
 \mathbb E_gN_{D_2}(X,L)&\leq C_2g^{-2}L^2e^L,\\
 \mathbb E_gN_{D_3}(X,L)&\leq C_3g^{-3}L^2e^L.
 \end{aligned}
 \label{f:exceptional-supports}
\end{equation}
Here $g\geq3$, $L\geq L_*$, and $C_T,C_2,C_3$ are absolute
constants. This is \cite[Lemma 3.4]{HSWX}, with the primitive unoriented
counting convention stated in its Section 2.5.

The one-pants integral uses two volume estimates. First, the upper
boundary-profile estimate
\cite[Lemma 22]{NWX} states that, for $h,n\geq1$ and all $x_i\geq0$,
\begin{equation}
 \frac{V_{h,n}(x_1,\ldots,x_n)}{V_{h,n}}
 \leq\prod_{i=1}^n\frac{\sinh(x_i/2)}{x_i/2}.
 \label{f:upper-profile}
\end{equation}
Quotients at zero are defined by continuity. In particular this is a
global upper bound, with no restriction of the form $x_i=o(\sqrt h)$.
For each fixed $n\geq0$, the volume asymptotic as $h\to\infty$
\cite[Theorem 1.2]{MZ} has a universal positive leading constant $c$:
\begin{equation}
 V_{h,n}
 =c\frac{(2h-3+n)!(4\pi^2)^{2h-3+n}}{\sqrt h}
   \bigl(1+O_n(h^{-1})\bigr).
 \label{f:volume-asymptotic}
\end{equation}
Consequently, as also recorded in \cite[equation (3.3)]{HSWX},
\begin{equation}
 \rho_g:=\frac{V_{g-2,3}}{V_g}
 =\frac1{8\pi^2g}\bigl(1+O(g^{-1})\bigr),\qquad
 C_V:=\sup_{g\geq3}g\rho_g<\infty.
 \label{f:one-pants-volume-ratio}
\end{equation}
Indeed, the factorial gap in the ratio is one, and its remaining
square-root factor is $\sqrt{g/(g-2)}=1+O(g^{-1})$. This bounds
$g\rho_g$ for large $g$; the remaining finitely many values are finite,
so the supremum defining $C_V$ is finite.

\begin{lemma}[Restricted one-point integration]
\label{lem:f-one-point-integral}
Let $w:\mathbb R_{\geq0}^3\longrightarrow[0,\infty]$ be measurable and
symmetric under every permutation. Define
\[
 N_w(X,L)=\sum_{\substack{\gamma\text{ a figure-eight},\ \ell_X(\gamma)\leq L\\
                         P_\gamma\text{ good}}}
 w(\text{the three cuff lengths of }P_\gamma).
\]
Then, for $g\geq3$ and $L\geq0$,
\begin{equation}
 \mathbb E_gN_w(X,L)
 =\frac1{2V_g}\int_{\ell(x,y,z)\leq L}
 w(x,y,z)V_{g-2,3}(x,y,z)xyz\,dx\,dy\,dz.
 \label{f:one-point-integral}
\end{equation}
Moreover, with
$\Phi(x,y,z)=\sinh(x/2)\sinh(y/2)\sinh(z/2)$, the density satisfies
\begin{equation}
 \frac1{2V_g}V_{g-2,3}(x,y,z)xyz
 \leq4\rho_g\Phi(x,y,z)
 \leq\frac{4C_V}{g}\Phi(x,y,z).
 \label{f:one-point-density}
\end{equation}
\end{lemma}

\begin{proof}
Mirzakhani integration in the ordered-cuff orbit that cuts $X$ into
$S_{0,3}$ and $S_{g-2,3}$ has coefficient one when $g>2$; see
\cite[Theorem 2.2 and the derivation of equation (3.2)]{HSWX}.
The pants volume is $V_{0,3}=1$. A good pair of pants with unlabeled
cuffs is represented by six ordered triples of its cuff geodesics. Its three actual
figure-eights contribute the indicators
\[
 \mathbf1_{\{\ell(x,y,z)\leq L\}},\qquad
 \mathbf1_{\{\ell(z,x,y)\leq L\}},\qquad
 \mathbf1_{\{\ell(y,z,x)\leq L\}}.
\]
After integration, these three terms are equal by symmetry of $w$, of
$V_{g-2,3}$, and of $xyz$. Their coefficient is therefore $3/6=1/2$.
This proves \eqref{f:one-point-integral}, initially for bounded
nonnegative weights and then by monotone convergence. The labelings
concern distinct cuff geodesics even if their lengths coincide, so no
metric-symmetry factor is introduced. Finally apply
\eqref{f:upper-profile} at genus $g-2\geq1$; multiplying the three
profile factors by $xyz/2$ gives $4\Phi$.
\end{proof}

\subsection{Uniform trimming}

We discard curves with exceptional supports or short cuffs and estimate
the expected number lost. For $L,b\geq0$, let $\mathcal F(X,L)$ be the figure-eights of length
at most $L$, and let $\mathcal F^{\mathrm{tr}}(X;L,b)$ consist of those
with good supports and all three cuffs at least $b$. Put
\[
 R_g(X;L,b)=\#\mathcal F(X,L)-\#\mathcal F^{\mathrm{tr}}(X;L,b).
\]
These are finite counts on every $X$.

\begin{proposition}[Trimming estimate]
\label{prop:trim}
There is an absolute constant $C$ such that, for all $g\geq3$,
$L\geq L_*$, and $0\leq b\leq L$,
\begin{equation}
 \mathbb E_gR_g(X;L,b)
 \leq C\left[
 \frac{be^L+Le^{(L+b)/2}+L^2e^{L/2}}g
 +\frac{L^2e^L}{g^2}\right].
 \label{f:trimming-uniform}
\end{equation}
For the parameters \eqref{f:critical-parameters}, write
$R_{g,S}=R_g(X;L,b)$. Then
\begin{equation}
 \mathbb E_gR_{g,S}
 \leq C_S\left(\frac{\log H}{H}
          +g^{-1/2}H^{11/2}+\frac Hg\right)
 =O_S\left(\frac{\log H}{H}\right)=o(1).
 \label{f:trimming-critical}
\end{equation}
\end{proposition}

\begin{proof}
Let $D_L=\{(x,y,z)\geq0:\ell(x,y,z)\leq L\}$. The identity
\[
 \cosh\frac\ell2
 =\cosh\frac z2+\cosh\frac{x+y}2+\cosh\frac{x-y}2
\]
implies
\begin{equation}
 D_L\subseteq\{x,y,z\geq0:x+y\leq L,\ z\leq L\}.
 \label{f:one-point-relaxation}
\end{equation}
Set
\[
 J(t)=\int_0^t\sinh(u/2)\,du
     =2(\cosh(t/2)-1)\leq e^{t/2}\qquad(t\geq0).
\]
Using also $\sinh(t/2)\leq e^{t/2}/2$, a short winding cuff gives
\begin{equation}
 \begin{aligned}
 \int_{D_L\cap\{x<b\}}\Phi\,dx\,dy\,dz
 &\leq J(L)\int_0^b\sinh(x/2)J(L-x)\,dx\\
 &\leq e^{L/2}\int_0^b\tfrac12e^{x/2}e^{(L-x)/2}\,dx
 =\frac b2e^L.
 \end{aligned}
 \label{f:short-winding-cuff}
\end{equation}
The same estimate applies to $y<b$. For the third cuff,
\begin{equation}
 \begin{aligned}
 \int_{D_L\cap\{z<b\}}\Phi\,dx\,dy\,dz
 &\leq J(b)\int_0^L\sinh(x/2)J(L-x)\,dx\\
 &\leq e^{b/2}\int_0^L\tfrac12e^{x/2}e^{(L-x)/2}\,dx
 =\frac L2e^{(L+b)/2}.
 \end{aligned}
 \label{f:short-third-cuff}
\end{equation}
Taking the union of these conditions yields
\begin{equation}
 \int_{D_L\cap\{\min(x,y,z)<b\}}\Phi\,dx\,dy\,dz
 \leq be^L+\frac L2e^{(L+b)/2}.
 \label{f:short-cuff-integral}
\end{equation}
If $N_{\mathrm{short}}(X;L,b)$ counts good-support curves with at least
one cuff shorter than $b$, the weight
$\mathbf1_{\{\min(x,y,z)<b\}}$ is symmetric. Hence
Lemma~\ref{lem:f-one-point-integral} gives
\begin{equation}
 \mathbb E_gN_{\mathrm{short}}(X;L,b)
 \leq\frac{4C_V}{g}
       \left(be^L+\frac L2e^{(L+b)/2}\right).
 \label{f:short-cuff-count}
\end{equation}
The factor $1/2$ in the exact integral already accounts for the curve
choices; there is no additional factor three in this step.

By Lemma~\ref{lem:f-support-types}, the exact deletion identity is
\begin{equation}
 R_g(X;L,b)=N_T(X,L)+N_{D_2}(X,L)+N_{D_3}(X,L)
                     +N_{\mathrm{short}}(X;L,b).
 \label{f:deletion-partition}
\end{equation}
The four counts are disjoint, and a curve with several short cuffs is
counted only once in the last term. Combine
\eqref{f:exceptional-supports} and \eqref{f:short-cuff-count}, and use
$g^{-3}\leq g^{-2}$, to prove \eqref{f:trimming-uniform}.

For fixed $S$, we have $L/H\to1$, $L\to\infty$, and
$10\log L/L\to0$. Thus the hypotheses hold eventually, in fact with
$1\leq b\leq L/3$. The exact substitutions are
\begin{equation}
 \begin{aligned}
 \frac{be^L}{g}&=10e^S\frac{\log L}{H},&
 \frac{Le^{(L+b)/2}}g&=e^{S/2}g^{-1/2}L^6H^{-1/2},\\
 \frac{L^2e^{L/2}}g&=e^{S/2}g^{-1/2}L^2H^{-1/2},&
 \frac{L^2e^L}{g^2}&=e^S\frac{L^2}{gH}.
 \end{aligned}
 \label{f:trimming-substitution}
\end{equation}
They follow from $e^L=e^Sg/H$ and $e^{b/2}=L^5$. Since $L\asymp_S H$,
they imply \eqref{f:trimming-critical}.
\end{proof}

\subsection{Different curves on the same pants}

We begin the pair estimates with curves that share a support.
Let $B_g^{\mathrm{same}}(X;L,b)$ count unordered pairs of distinct
members of $\mathcal F^{\mathrm{tr}}(X;L,b)$ with the same support,
and set $Q_g^{\mathrm{same}}=2B_g^{\mathrm{same}}$ for ordered pairs.
Pair order is independent of geodesic orientation, which is never
counted here. At the parameters \eqref{f:critical-parameters}, we
abbreviate $B_{g,S}^{\mathrm{same}}=B_g^{\mathrm{same}}(X;L,b)$.

\begin{proposition}[Same-pants pairs]
\label{prop:f-same-pants}
For every $g\geq3$ and $L,b\geq0$,
\begin{equation}
 \begin{aligned}
 \mathbb E_gQ_g^{\mathrm{same}}(X;L,b)
 &\leq8\rho_g e^{L-b/2}\leq\frac{8C_V}{g}e^{L-b/2},\\
 \mathbb E_gB_g^{\mathrm{same}}(X;L,b)
 &\leq4\rho_g e^{L-b/2}\leq\frac{4C_V}{g}e^{L-b/2}.
 \end{aligned}
 \label{f:same-pants-bound}
\end{equation}
Both counts vanish identically when $2b\geq L$. In the critical
regime \eqref{f:critical-parameters},
\begin{equation}
 \mathbb E_gB_{g,S}^{\mathrm{same}}
 \leq\frac{4C_Ve^S}{HL^5}=O_S(H^{-6}).
 \label{f:same-pants-critical}
\end{equation}
\end{proposition}

\begin{proof}
On a pair of pants with cuff labels $x,y,z$, write $\gamma_{xy}$ for the
figure-eight winding around the first two cuffs, and put
$\ell_{xy}(x,y,z)=\ell(x,y,z)$ and
$\ell_{xz}(x,y,z)=\ell(x,z,y)$. An ordered pair of distinct curves
uniquely labels the cuffs: $x$ is the cuff common to their winding
pairs, $y$ is the other winding cuff of the first curve, and $z$ is
the other winding cuff of the second. Thus the ordered pair is
$(\gamma_{xy},\gamma_{xz})$ in exactly one cuff labeling.
Length coincidences do not affect this correspondence.

Define
\[
 D^{(2)}_{L,b}
 =\{(x,y,z)\in[b,\infty)^3:
              \ell_{xy}(x,y,z)\leq L,\ \ell_{xz}(x,y,z)\leq L\}.
\]
The ordered-cuff integration formula used in
Lemma~\ref{lem:f-one-point-integral} now gives
\begin{equation}
 \mathbb E_gQ_g^{\mathrm{same}}(X;L,b)
 =\frac1{V_g}\int_{D^{(2)}_{L,b}}
            V_{g-2,3}(x,y,z)xyz\,dx\,dy\,dz.
 \label{f:same-pants-integral}
\end{equation}
The unique labeling gives coefficient one in
\eqref{f:same-pants-integral}. Passing to unordered pairs divides it by two.

Equation~\eqref{f:pants-length} gives the strict inequalities
$\ell_{xy}>x+y$ and $\ell_{xz}>x+z$. Therefore
\begin{equation}
 D^{(2)}_{L,b}\subseteq
 \Omega_{L,b}:=\{x,y,z\geq b:x+y\leq L,\ x+z\leq L\}.
 \label{f:same-pants-relaxation}
\end{equation}
The exact domain is empty for $2b\geq L$. If $0\leq b\leq L/2$, put
\[
 J_b(L)=\int_{\Omega_{L,b}}\Phi(x,y,z)\,dx\,dy\,dz
       =\int_b^{L-b}\sinh(x/2)
          \left(\int_b^{L-x}\sinh(t/2)\,dt\right)^2dx.
\]
Then
\begin{equation}
 \begin{aligned}
 J_b(L)
 &\leq\frac18\int_b^{L-b}e^{x/2}
             \left(\int_b^{L-x}e^{t/2}\,dt\right)^2dx\\
 &\leq\frac12e^L\int_b^{L-b}e^{-x/2}\,dx\\
 &=e^L\bigl(e^{-b/2}-e^{-(L-b)/2}\bigr)
 \leq e^{L-b/2}.
 \end{aligned}
 \label{f:same-pants-elementary-integral}
\end{equation}
The relaxed domain is empty for $b>L/2$ and has zero volume at equality.
By \eqref{f:upper-profile}, the density in
\eqref{f:same-pants-integral} is at most $8\rho_g\Phi$. This proves
\eqref{f:same-pants-bound}. Finally
$e^{L-b/2}/g=e^S/(HL^5)$ proves \eqref{f:same-pants-critical}.
\end{proof}

At $b=0$, Proposition~\ref{prop:f-same-pants} yields $O(e^L/g)$ for
good supports. No estimate for equal-support pairs in the exceptional
topologies of Lemma~\ref{lem:f-support-types} is needed. If
$U_{g,S}^{\mathrm{same}}$ counts retained curves having a retained
neighbor on the same support, the pointwise incidence bound and its
expectation and probability consequences are
\begin{equation}
 U_{g,S}^{\mathrm{same}}\leq2B_{g,S}^{\mathrm{same}},\qquad
 \mathbb E_gU_{g,S}^{\mathrm{same}}=O_S(H^{-6}),\qquad
 \mathbb P_g(B_{g,S}^{\mathrm{same}}>0)=O_S(H^{-6}).
 \label{f:same-pants-deletion}
\end{equation}
Indeed, a pair of pants carrying $n\in\{0,1,2,3\}$ retained curves
contributes $n\mathbf1_{\{n\geq2\}}\leq n(n-1)$ to this comparison.

\subsection{The geometry and labeling of different supports}

From now on we use the critical parameters \eqref{f:critical-parameters}.
Join two distinct retained curves by an edge when their closed supports
meet. Denote the number of unordered edges by $B_{g,S}$. For different
supports, separate touching pairs, with disjoint interiors and a common
cuff, from pairs with overlapping interiors. In the latter case form
the intrinsic geodesic filled union $Y$: fill contractible boundary
components of the union by disks and straighten the essential boundary
components to geodesics. Boundary copies in the intrinsic completion
may represent the same ambient geodesic. The classification uses the
intrinsic topological type of $Y$. Throughout this section,
$p(Y)=\ell(\partial Y)$ denotes the sum of the lengths of its intrinsic
boundary circles, counting both copies when their ambient images
coincide.

\begin{lemma}[Geometry of a pair core]
\label{lem:f-pair-core}
For two distinct retained supports with overlapping interiors, their
filled union $Y$ is connected and satisfies
\begin{equation}
 p(Y):=\ell(\partial Y)\leq4L,\qquad
 \operatorname{Area}(Y)\leq4L+4\pi.
 \label{f:pair-core-geometry}
\end{equation}
The two actual figure-eights fill $Y$ and have total length at most
$2L$. For sufficiently large $g$, $Y$ is proper, has nonempty boundary,
and has $m=|\chi(Y)|\geq2$. The possibilities with $m=2$ are
$S_{1,2}$ and $S_{0,4}$.
\end{lemma}

\begin{proof}
If $u,v$ are the winding cuffs of a retained curve and $w$ its third
cuff, \eqref{f:pants-length} implies
\begin{equation}
 \ell_X(w)<L,\qquad \ell_X(u)+\ell_X(v)<L.
 \label{f:actual-to-relaxed-lengths}
\end{equation}
Thus each supporting pair of pants has perimeter at most $2L$. The geodesic-union
construction has boundary length at most the sum of the two original
perimeters, and area at most their total area plus that perimeter sum;
these are the geometric inequalities
\cite[equations (4.17)--(4.18)]{HSWX}. They follow by straightening
essential boundary pieces and applying the hyperbolic isoperimetric
inequality to the disks and annular strips added in the construction.
Each pair of pants has area $2\pi$, giving \eqref{f:pair-core-geometry}.

We justify that the actual curves, without any increase of the length
cutoff, fill this same intrinsic surface. A figure-eight fills its
pants; the complementary components there are peripheral annuli.
An arc crossing a pair of pants while avoiding its figure-eight must therefore
be boundary-parallel. A geodesic crossing segment of this kind would
form a bigon with a geodesic cuff, which is impossible. Consequently a
simple closed geodesic avoiding both figure-eights cannot cross either
support. After disks are filled and essential boundary components
straightened, the only remaining components disjoint from the curves
are disks and annuli parallel to the boundary of $Y$. This is the
filling condition. The argument is carried out with distinct intrinsic
boundary copies, so it also applies when some copies have the same
ambient image. For a direct proof that the curve graph has this filled
interior, see the argument following Lemma~\ref{lem:barrier}.

Neither support is strictly contained in the other, because they have
equal area and are distinct geodesic subsurfaces. The filled union
contains both and has area strictly greater than $2\pi$. Hence
Gauss--Bonnet gives $m\geq2$. Its area is $O_S(\log g)$, whereas
$\operatorname{Area}(X)=4\pi(g-1)$, so for large $g$ it is a proper
subsurface with nonempty boundary. Solving $2h-2+k=2$ with $k\geq1$
gives precisely $(h,k)=(0,4),(1,2)$.
\end{proof}

Disjoint interiors can meet only along common entire cuffs. Indeed,
distinct geodesic boundaries meet transversely unless they coincide,
and a transverse crossing forces the two chosen interior sides to
overlap locally. Touching pants share one or two cuffs: sharing all
three would give a closed genus-two component, incompatible with
connected $X$ of genus $g>2$. Thus the conflict classes are exhaustive.
For unordered pairs of retained curves on distinct supports, let
$B_{g,S}^{\mathrm{touch}}$ count touching pairs and
$B_{1,2},B_{0,4},B_{\geq3}$ count overlapping pairs, classified by the
intrinsic core in Lemma~\ref{lem:f-pair-core}.

For comparison with the HSWX counts, let $\mathcal N_L(X)$ be
the set of ordered triples $\Gamma=(c_1,c_2,c_3)$ of distinct simple
geodesics bounding a good pair of pants such that
\begin{equation}
 \ell_X(c_1)\leq L,\qquad
 \ell_X(c_2)+\ell_X(c_3)\leq L,\qquad
 \min_i\ell_X(c_i)\geq10\log L.
 \label{f:relaxed-triples}
\end{equation}
This is \cite[Definition 4.2]{HSWX}. Denote by
$C_{1,2}^{\mathrm{rel}},C_{0,4}^{\mathrm{rel}},C_{\geq3}^{\mathrm{rel}}$
its ordered pair counts for distinct overlapping supports, and by
$D^{\mathrm{rel}}$ its ordered pair count for touching supports. The
core classification is exactly the intrinsic construction above.

\begin{lemma}[Transfer of actual pairs]
\label{lem:f-eight-transfer}
Pointwise on every surface in this regime,
\begin{equation}
 8B_{\tau}\leq C_{\tau}^{\mathrm{rel}}
 \quad(\tau=(1,2),(0,4),\geq3),\qquad
 8B_{g,S}^{\mathrm{touch}}\leq D^{\mathrm{rel}}.
 \label{f:eight-transfer}
\end{equation}
\end{lemma}

\begin{proof}
For a retained curve with winding cuffs $u,v$ and third cuff $w$,
\eqref{f:actual-to-relaxed-lengths} places both $(w,u,v)$ and $(w,v,u)$
in $\mathcal N_L(X)$. These are two distinct cuff-geodesic triples,
even when lengths agree. Either one recovers its support and the unique
unoriented figure-eight winding around its second and third cuffs.
An unordered actual pair on distinct supports therefore produces
$2\cdot2\cdot2=8$ ordered source pairs, from pair order and the two
winding orders for each member. Different actual pairs have disjoint
image sets, and all eight images have the same support type. This proves
\eqref{f:eight-transfer}. The source's reverse construction of a curve
from a relaxed triple uses an additive length allowance; that direction
is not used here.
\end{proof}

\subsection{Touching and the two smallest overlapping types}

For overlapping supports, Lemma~\ref{lem:f-eight-transfer} compares
actual pairs with relaxed pairs. For $L>1$ and sufficiently large $g$,
He, Shen, Wu and Xue~\cite[Propositions 4.7 and 4.11]{HSWX} give the bounds
\begin{equation}
 \mathbb E_g C_{1,2}^{\mathrm{rel}}(X,L)
 \leq C\frac{e^{2L}}{g^2},\qquad
 \mathbb E_g C_{0,4}^{\mathrm{rel}}(X,L)
 \leq C\frac{Le^{2L}}{g^2}.
 \label{f:imported-small-overlaps}
\end{equation}
Here $C$ is absolute, and the estimates apply to the full counts just
defined. In particular,
the complementary-type sums in the source's equations (4.41) and
(4.57) include disconnected complements. Its Proposition 4.11 also
includes the intrinsic $S_{0,4}$ case with an identified pair of
boundary copies. No assertion that the joint complement of two good
pants is connected is needed.

For touching pairs, either curve can be chosen as the single mark,
and its support has connected complement by goodness. Thus
$2B_{g,S}^{\mathrm{touch}}=I_{T,1}^{\rm conn}(L,b)$.
At the retained cutoff, Proposition~\ref{prop:touching} bounds its
expectation by $O_S(\log H/H^{10})$. The corresponding comparison
for relaxed cuff triples is given in~\eqref{f:touching-source-estimate}.

At \eqref{f:critical-parameters},
\begin{equation}
 \frac{e^{2L}}{g^2}=\frac{e^{2S}}{H^2}.
 \label{f:small-overlap-substitution}
\end{equation}
For intrinsic $S_{1,2}$ pairs, apply \eqref{f:eight-transfer} and
\eqref{f:imported-small-overlaps}. Splitting the intrinsic $S_{0,4}$
pairs according to whether two boundary copies are identified gives
$2B_{0,4}=Q_{04^0}(L,b)+Q_{04^1}(L,b)$, where the right-hand counts
are ordered. At $b=10\log L$, Lemmas~\ref{tail:planar}
and~\ref{tail:small} bound each expectation by $Ce^{2L}/g^2$.
Consequently,
\begin{equation}
 \mathbb E_gB_{g,S}^{\mathrm{touch}}=O_S(H^{-8}),\qquad
 \mathbb E_gB_{1,2}=O_S(H^{-2}),\qquad
 \mathbb E_gB_{0,4}=O_S(H^{-2}).
 \label{f:small-interaction-rates}
\end{equation}

\subsection{Higher-complexity cores}

For each fixed core type, we combine a count of filling pairs with a
volume estimate. The constant in the pair count depends on the
complexity. To sum over the remaining types, we double the core along
its boundary and use a bound for primitive curves on the resulting
closed surface.

For a compact hyperbolic surface $Y$ with
nonempty geodesic boundary, let $N_k^{\mathrm{fill}}(Y,R)$ count
ordered $k$-tuples of nonperipheral primitive unoriented geodesics
which fill $Y$ and whose \emph{sum} of lengths is at most $R$.
Here filling means that every complementary component is a disk or an
annulus parallel to a boundary component of $Y$. By
\cite[Definition 17 and Theorem 18]{WXPGT}, for fixed $k\geq1$,
$0<\delta<1/2$, and $m=2h-2+n\geq1$, there is a constant
$C_{k,\delta,m}$ such that for every such metric of type $S_{h,n}$ and
every $R>0$,
\begin{equation}
 N_k^{\mathrm{fill}}(Y,R)
 \leq C_{k,\delta,m}(1+R)^{k-1}
       e^{R-(1-\delta)p(Y)/2}.
 \label{f:filling-tuple-count}
\end{equation}
The constant is uniform in the boundary lengths and the metric. We
also write $N_{\mathrm{prim}}(Z,R)$ for the number of primitive
unoriented closed geodesics on $Z$ of length at most $R$. The
metric-uniform estimate \cite[Theorem 2.8]{HSWX} states that, on every
closed genus-$G$ hyperbolic surface $Z$, for $G\geq2$ and $R>0$,
\begin{equation}
 N_{\mathrm{prim}}(Z,R)\leq(G-1)e^{R+7}.
 \label{f:primitive-count}
\end{equation}

Fix $a>0$ and put $T=a\log g$. Let $\mathcal S_T(X)$ be the family
of connected geodesic-boundary subsurfaces with
\begin{equation}
 p(Y)\leq2T,\qquad \operatorname{Area}(Y)\leq4T.
 \label{f:subsurface-family}
\end{equation}
As in \cite[Section 7, equation (21) and the subsequent definition]{WX},
distinct intrinsic boundary circles are allowed to represent the same
ambient geodesic. Complements may be disconnected. For every fixed
integer $m_0\geq1$, Proposition 33 of that paper gives, for sufficiently
large $g$,
\begin{equation}
 \mathbb E_g\sum_{\substack{Y\in\mathcal S_T(X)\\|\chi(Y)|\geq m_0}}
 e^{-p(Y)/4}
 \leq C_{a,m_0}\frac{T e^{7T/2}}{g^{m_0}}.
 \label{f:subsurface-tail}
\end{equation}
The upper complexity cutoff in the source is automatic from the area
condition: $|\chi(Y)|\leq\lfloor4T/(2\pi)\rfloor$.
For each fixed type $S_{h,n}$ with $m=2h-2+n\in[1,16]$ and each fixed
$\delta>0$, Proposition 34 gives
\begin{equation}
 \mathbb E_g\sum_{\substack{Y\in\mathcal S_T(X)\\Y\cong S_{h,n}}}
 e^{-(1-\delta)p(Y)/2}
 \leq C_{a,\delta,h,n}\frac{T^{66}e^{\delta T}}{g^m}.
 \label{f:subsurface-fixed-type}
\end{equation}
This follows from the stated unnormalized estimate by division by $V_g$
and cancellation of its common factor $e^{T/2}$. Both estimates hold
for sufficiently large $g$, with the threshold allowed to depend on
the displayed fixed parameters. They sum all complementary
realizations, including boundary identifications.

\begin{lemma}[High-complexity pair estimate]
\label{lem:f-high-complexity}
For the retained actual pairs with intrinsic core complexity at least
three,
\begin{equation}
 \mathbb E_g B_{\geq3}
 =O_S\bigl(g^{-3/4}H^{65}+g^{-1}\bigr).
 \label{f:high-complexity-bound}
\end{equation}
\end{lemma}

\begin{proof}
Take $T=2H$, a permitted choice in \eqref{f:subsurface-family}.
Eventually $L\leq H$ and $4L+4\pi\leq8H$. Thus every core arising
from an overlapping retained pair belongs to $\mathcal S_{2H}(X)$,
by Lemma~\ref{lem:f-pair-core}.

First suppose $3\leq m\leq10$. Each ordered actual pair fills its
unique intrinsic core with total length at most $2L$. Removing the
retention conditions only increases the filling-tuple count. There are
finitely many intrinsic types in this range, all within the range of
\eqref{f:subsurface-fixed-type}. Apply
\eqref{f:filling-tuple-count} with $k=2$ and $R=2L$, then sum
\eqref{f:subsurface-fixed-type}, to obtain
\begin{equation}
 \begin{aligned}
 \mathbb E_gB_{3\leq m\leq10}
 &\leq C_\delta Le^{2L}
       \sum_{\substack{h\geq0,\ n\geq1\\3\leq2h-2+n\leq10}}
       \frac{T^{66}e^{\delta T}}{g^{2h-2+n}}\\
 &\leq C_\delta\frac{L H^{66}e^{2L+2\delta H}}{g^3}.
 \end{aligned}
 \label{f:fixed-complexity-pairs}
\end{equation}
The ordered filling-pair count bounds the desired unordered count
from above. With $\delta=1/8$,
$e^{2L}=e^{2S}g^2/H^2$, and $L=O_S(H)$, the last expression is
$O_S(g^{-3/4}H^{65})$.

For $m\geq11$ we use a bound whose constant does not depend on $m$.
Double the intrinsic bordered surface $Y\cong S_{h,n}$ along its
geodesic boundary. The resulting closed hyperbolic surface has genus
$2h+n-1=m+1$. An interior primitive geodesic remains primitive in the
double, and distinct unoriented classes remain distinct. To see this,
fold the double continuously onto the chosen copy of $Y$. This is a
retraction: a purported root or conjugacy in the double would retract
to a root or conjugacy in $Y$. Actual primitive curves in $X$ are
primitive in $Y$ as well, since an intrinsic proper power would remain
a proper power in $X$. Their geodesic lengths are unchanged on doubling.
Equation~\eqref{f:primitive-count} therefore bounds the number of
eligible curves in $Y$ by $m e^{L+7}$.

Let $M_L(Y)$ count ordered primitive pairs, each member of length at
most $L$, whose filled core is $Y$. When $Y$ arises from an overlapping
retained pair, \eqref{f:pair-core-geometry} and Gauss--Bonnet give
$2\pi m\leq4L+4\pi$, hence $m=O(L)$ in the present regime. Consequently
\begin{equation}
 M_L(Y)\leq C L^2e^{2L}
          \leq C L^2e^{3L}e^{-p(Y)/4}.
 \label{f:growing-complexity-pointwise}
\end{equation}
The second inequality uses $p(Y)\leq4L$. Crucially, the first bound
retains the individual cutoff $L$ on each curve.

We now sum \eqref{f:growing-complexity-pointwise} over cores arising
from overlapping retained pairs with $m\geq11$, and enlarge the resulting
nonnegative weighted subsurface sum to $\mathcal S_{2H}(X)$. Applying \eqref{f:subsurface-tail} with $m_0=11$
gives
\begin{equation}
 \begin{aligned}
 \mathbb E_g B_{m\geq11}
 &\leq C L^2e^{3L}\frac{H e^{7H}}{g^{11}}\\
 &=C e^{3S}\frac{L^2}{gH^2}=O_S(g^{-1}).
 \end{aligned}
 \label{f:growing-complexity-pairs}
\end{equation}
This enlargement applies to the weighted sum; the pointwise bounds
$m=O(L)$ and $p(Y)\leq4L$ were used only for cores arising from
retained pairs.
Combining the two complexity ranges proves the lemma.
\end{proof}

\begin{remark}
The perimeter cutoff used throughout this argument is $4L$. The
displayed union-perimeter inequality in \cite[equation (4.17)]{HSWX},
together with the two individual perimeter bounds $2L$, does not by
itself yield the $2L$ cutoff used in its Lemma 4.5. Our proof uses
\eqref{f:pair-core-geometry} and $T=2\log g$ and does not invoke that
paper's Proposition 4.6. The estimate above is asserted at the fixed
critical cutoff; no extension to an arbitrary uniform $L$ range is
needed.
\end{remark}

\subsection{Conflict bounds and the isolation coupling}

\begin{proposition}[Conflict estimate]
\label{prop:conflicts}
For each fixed $S\in\mathbb R$, the retained family satisfies
\begin{equation}
 \mathbb E_gB_{g,S}=O_S(H^{-2}),\qquad
 \mathbb P_g(B_{g,S}>0)=O_S(H^{-2}).
 \label{f:conflicts-total}
\end{equation}
The same $O_S(H^{-2})$ expectation bound holds for all distinct-support
conflicts alone.
\end{proposition}

\begin{proof}
The geometric classification gives the pointwise identity
\[
 B_{g,S}=B_{g,S}^{\mathrm{same}}+B_{g,S}^{\mathrm{touch}}
          +B_{1,2}+B_{0,4}+B_{\geq3}.
\]
Proposition~\ref{prop:f-same-pants}, \eqref{f:small-interaction-rates},
and Lemma~\ref{lem:f-high-complexity} bound these terms. Since
$g^{-3/4}H^{65}$ and $g^{-1}$ are $o(H^{-2})$, summation yields the
expectation claim, and Markov's inequality gives its probability
consequence.
\end{proof}

Define the retained point measure and its isolated submeasure by
\[
 \Xi_{g,S}^{\mathrm{tr}}
   =\sum_{\gamma\in\mathcal F^{\mathrm{tr}}(X;L,b)}
           \delta_{\ell_X(\gamma)-a_g},\qquad
 \widehat\Xi_{g,S}
   =\sum_{\substack{\gamma\in\mathcal F^{\mathrm{tr}}(X;L,b)\\
                   \gamma\text{ isolated in the retained conflict graph}}}
           \delta_{\ell_X(\gamma)-a_g}.
\]
Isolation is tested against all retained curves with centered length at
most $S$.
Write $\Xi_g^{(8)}$ for the full centered figure-eight point process.

\begin{corollary}[Same-surface comparison]
\label{cor:f-isolation-coupling}
On the same sampled surface,
\begin{equation}
 \widehat\Xi_{g,S}\leq\Xi_{g,S}^{\mathrm{tr}}
       \leq\Xi_g^{(8)}\big|_{(-\infty,S]},\qquad
 0\leq\Xi_{g,S}^{\mathrm{tr}}(\mathbb R)
              -\widehat\Xi_{g,S}(\mathbb R)\leq2B_{g,S}.
 \label{f:submeasure-coupling}
\end{equation}
Moreover,
\begin{equation}
 \begin{aligned}
 \mathbb P_g\left\{\Xi_g^{(8)}\big|_{(-\infty,S]}
                       \ne\Xi_{g,S}^{\mathrm{tr}}\right\}
 &\leq\mathbb E_gR_{g,S},\\
 \mathbb P_g\left\{\Xi_g^{(8)}\big|_{(-\infty,S]}
                       \ne\widehat\Xi_{g,S}\right\}
 &\leq\mathbb E_gR_{g,S}+\mathbb E_gB_{g,S}
  =O_S\left(\frac{\log H}{H}\right).
 \end{aligned}
 \label{f:full-to-isolated-probability}
\end{equation}
Each right-hand side also bounds the total variation distance between
the corresponding two laws, with
$d_{\mathrm{TV}}(P,Q)=\sup_A|P(A)-Q(A)|$.
\end{corollary}

\begin{proof}
Trimming removes exactly $R_{g,S}$ atoms. Isolation then removes at
most $2B_{g,S}$, since every nonisolated vertex is an endpoint of an
edge. This proves \eqref{f:submeasure-coupling}, including when several
atoms have the same length. A deletion occurs in the first step exactly
when $R_{g,S}\geq1$, and in the second exactly when $B_{g,S}\geq1$.
The union bound and Markov's inequality give
\eqref{f:full-to-isolated-probability}. Coupling the measures on the
same surface gives the asserted total variation bounds.
\end{proof}

In particular, with probability
$1-O_S(\log\log g/\log g)$, every figure-eight of length at most
$a_g+S$ has a good support with all cuffs at least $b$, and these closed
supports are pairwise disjoint. At each fixed ceiling $S$, isolation
deletes $O_S(H^{-2})$ curves in expectation. This is the deletion
estimate at order one. To compare factorial moments at higher orders,
we count disjoint tuples in Section~\ref{sec:disjoint} and bound the
number of those tuples affected by isolation in Theorem~\ref{thm:weighted-loss}.

% END sections/02_foundations.tex

% BEGIN sections/03_disjoint.tex
\section{Disjoint supports and factorial moments}\label{sec:disjoint}

Consider an ordered tuple of retained figure-eights with pairwise
disjoint closed supports. Cutting out a single support leaves a connected
surface. Cutting out the whole tuple can leave several components.
The exact integration formula below includes each possible
complement. Its connected term is asymptotic to a product of one-point
measures, and each additional complementary component costs a power
of the genus. At the critical cutoff this gives the mixed moments of
disjoint tuples at every fixed order $m$. Such a tuple can contain
a curve whose support meets that of a curve outside the tuple.
Theorem~\ref{thm:weighted-loss} will control the loss when we pass
to factorial moments of isolated curves. All asymptotic estimates in
this section hold at fixed $m$.

\subsection{Counting conventions and complementary realizations}

We use the good-support convention of Section~\ref{sec:foundations}:
the supporting pair of pants is embedded, its three cuffs are distinct
geodesics, and its complement is connected. Write
$\ell_\alpha=\ell_X(\alpha)$. For $L,b\geq0$ and a nonnegative Borel
function $F:[0,\infty)^m\to[0,\infty]$, let
\begin{equation}\label{a:disjoint-count}
 D_{g,m}^{L,b}(F;X)
 =\sum_{(\alpha_1,\ldots,\alpha_m)}
 F(\ell_{\alpha_1},\ldots,\ell_{\alpha_m}),
\end{equation}
where the sum runs over ordered tuples of primitive unoriented figure-eights
$\alpha_i$ with good supports $P_i$ whose closures $\overline P_i$ are
pairwise disjoint. Every curve has length at most $L$, and every support
cuff has length at least $b$. Curves and cuffs retain their identities
when lengths coincide.

With the first two cuffs designated as winding cuffs, put
\begin{equation}\label{a:length-domain}
 \begin{split}
 \ell(x,y,z)&=2\operatorname{arccosh}
 \bigl(\cosh(z/2)+2\cosh(x/2)\cosh(y/2)\bigr),\\
 \mathcal D_{L,b}&=\{(x,y,z)\in(0,\infty)^3:
                  \ell(x,y,z)\leq L,\ \min(x,y,z)\geq b\}.
 \end{split}
\end{equation}
The length identity is \cite[Eq.~(2.3)]{HSWX}.  The identity
$2\cosh(x/2)\cosh(y/2)=\cosh((x+y)/2)+\cosh((x-y)/2)$ gives
\begin{equation}\label{a:domain-inclusion}
 \mathcal D_{L,b}\subseteq
 \{(x,y,z)\in[0,\infty)^3:x+y\leq L,\ z\leq L\}.
\end{equation}
In particular every cuff coordinate is at most $L$.  We retain the exact
domain in all identities and leading-order comparisons.

Label the cuffs of the pants in slot $i$ by $(i,1),(i,2),(i,3)$, with
labels $1,2$ for the winding cuffs, and set
$E_m=\{1,\ldots,m\}\times\{1,2,3\}$.  A set partition $\pi$
of $E_m$ records the complementary components. Each block $C\in\pi$
contains the cuff labels incident to one component, of genus $h_C$.
The blocks are unordered. Write $q=|\pi|$. To read connectedness from
these data, form the incidence multigraph $G_\pi$ with a vertex $p_i$
for each selected pair of pants and a vertex for each block $C\in\pi$.
Each cuff $(i,a)$ gives an edge from $p_i$ to its block, so two vertices
may be joined by several edges.

\begin{definition}\label{a:def-admissible}
The datum $(\pi,(h_C)_{C\in\pi})$ is admissible at genus $g$ if
$G_\pi$ is connected, deleting any one pants vertex and all its incident
edges leaves a connected graph, every $h_C$ is a nonnegative integer
with $2h_C-2+|C|>0$, and
\begin{equation}\label{a:genus-identity}
 \sum_{C\in\pi}h_C=g-2m+q-1.
\end{equation}
A graph with one vertex is connected.  Denote the finite set of
admissible data by $\mathcal R_{g,m}$.
\end{definition}

\begin{lemma}\label{a:topological-classification}
Admissible data classify the mapping-class orbits of ordered
configurations of good pants with individually labeled cuffs and
pairwise disjoint closures.
\end{lemma}

\begin{proof}
Cut along all $3m$ cuffs.  The pieces are the selected pants and the
complementary components $Y_C$.  Every piece has negative Euler
characteristic: in particular, the complementary components are compact
hyperbolic surfaces with geodesic boundary.  Connectivity of the original
surface is equivalent to connectivity of $G_\pi$.  Removing $P_i$ while
retaining the other pants corresponds to deleting $p_i$ and its incident
edges.  Hence connectedness of the individual complement is exactly
the second graph condition.  Euler characteristic gives
\[
 2-2g=-m+\sum_C(2-2h_C-|C|)
      =2q-2\sum_Ch_C-4m,
\]
which is \eqref{a:genus-identity}.

Conversely, take $m$ ordered labeled oriented pants and, for each block
$C$, an oriented surface of type $S_{h_C,|C|}$.  Glue matching cuff labels
by orientation-reversing boundary identifications.  The graph conditions
give a connected surface with connected individual pants complements,
and \eqref{a:genus-identity} gives its genus.  Stability excludes disks
and annuli among the pieces.  More generally, an annulus between two
cut curves, or a disk bounded by one cut curve, would be a union of cut
pieces and would have negative Euler characteristic, a contradiction.
Thus all cuffs are essential and pairwise nonisotopic.  Choose positive
matching boundary lengths, hyperbolic metrics on the pieces, and arbitrary
gluing twists.  This realizes the datum by geodesic pants with disjoint
closures.  Finally, classification of oriented bordered surfaces and
boundary gluing shows that two realizations of the same datum are
related by an orientation-preserving homeomorphism respecting all labels.
The twists do not supply additional mapping-class orbits.
\end{proof}

For $m\geq2$, every complementary block meets at least two different
pants.  Otherwise deleting its only neighboring pants vertex isolates
that block from another remaining pants vertex.  For $m=1$, the deletion
condition forces the unique partition with one block of size three.
Consequently every complementary block has size at least two, and
\begin{equation}\label{a:component-bound}
 q\leq\lfloor3m/2\rfloor\quad(m\geq2).
\end{equation}
There are no once-holed-torus complementary factors $V_{1,1}$ in this
classification.  Disk, annulus, and cuff-identification realizations do
not belong to these configurations of embedded pants with disjoint
closures.

\subsection{Exact integration and multiplicity}

The topological classification allows us to sum the integration formula
over all complementary realizations. We use bordered Weil--Petersson
volumes with individually labeled boundary components and $V_{0,3}=1$.
For $\mathbf x=(x_{i,a})_{(i,a)\in E_m}$, write
$\mathbf x_i=(x_{i,1},x_{i,2},x_{i,3})$ and
$\mathbf x_C=(x_{i,a})_{(i,a)\in C}$. Define
\begin{equation}\label{a:volume-kernel}
 K_{g,m}(\mathbf x)
 =\sum_{(\pi,\mathbf h)\in\mathcal R_{g,m}}
   \prod_{C\in\pi}V_{h_C,|C|}(\mathbf x_C),
 \qquad \mathbf x=(x_{i,a})_{(i,a)\in E_m}.
\end{equation}
The order of entries within $\mathbf x_C$ is immaterial, since a bordered
volume is symmetric in its boundary lengths.  An empty sum is zero.

\begin{proposition}[Exact disjoint-pants integration]
\label{prop:exact-integration}
For every $g\geq3$, $m\geq1$, $L,b\geq0$, and nonnegative Borel $F$,
\begin{equation}\label{a:exact-integration}
 \begin{split}
 \mathbb E_g D_{g,m}^{L,b}(F)
 =\frac{1}{2^mV_g}
 \int_{\mathcal D_{L,b}^{\,m}}
 &F\bigl(\ell(\mathbf x_1),\ldots,\ell(\mathbf x_m)\bigr)
 K_{g,m}(\mathbf x)\\[-2mm]
 &{}\times\prod_{i=1}^m\prod_{a=1}^3x_{i,a}\,d\mathbf x.
 \end{split}
\end{equation}
The identity holds in $[0,\infty]$.  It also holds with finite absolute
integrals for bounded signed or complex Borel functions.
\end{proposition}

\begin{proof}
For a fixed admissible realization, orient each labeled cuff with its
selected pants on its left.  These orientations are prescribed by the
configuration, rather than additional binary choices.  By
Lemma~\ref{a:topological-classification}, the resulting ordered oriented
multicurves form one mapping-class orbit $\mathcal O_{\pi,\mathbf h}$.
The oriented ordered convention of \cite[Theorem~2.5 and
Remark~2.7]{AMI} gives, initially for a bounded compactly supported
nonnegative Borel function $W$,
\begin{equation}\label{a:oriented-orbit}
 \int_{\mathcal M_g}\sum_{\Gamma\in\mathcal O_{\pi,\mathbf h}}
 W(\ell_\Gamma(X))\,dX
 =\int W(\mathbf x)
   \prod_{C\in\pi}V_{h_C,|C|}(\mathbf x_C)
   \prod_{i,a}x_{i,a}\,d\mathbf x.
\end{equation}
The selected pants contribute $V_{0,3}^m=1$.  There is no extra symmetry
coefficient in this oriented ordered formula.  Its special convention
for $V_{1,1}$ has no effect here, since no such factor occurs; the ambient
genus assumption $g\geq3$ also agrees with that source.

Each curve in the tuple determines its third cuff and allows two orders
for its winding cuffs, giving exactly $2^m$ labelings. Conversely,
the oriented labeled cuffs recover the selected pair of pants on
their prescribed side and the unique unoriented figure-eight winding
around its first two cuffs, as in \cite[Section~2.4 and Eq.~(3.2)]{HSWX}.
The tuple positions are already ordered. A homeomorphism fixing the
ordered oriented cuffs preserves their left sides and each complementary
component's set of labels. Thus even equal genera or equal cuff lengths
introduce no additional $m!$ or $q!$ factor. The moduli-space integral
uses the prescribed orbifold measure.

In \eqref{a:oriented-orbit}, take
$W(\mathbf x)=\mathbf1_{\mathcal D_{L,b}^{\,m}}(\mathbf x)
F(\ell(\mathbf x_1),\ldots,\ell(\mathbf x_m))$, sum over all admissible
data, and divide by $2^mV_g$.  For bounded $F$, the cuff range lies in $[0,L]^{3m}$, where the finite
volume-polynomial sum is bounded. This gives absolute integrability.
Monotone convergence proves the nonnegative unbounded case and linearity
proves the signed and complex cases.
\end{proof}

For $m=1$, the formula reduces to
\begin{equation}\label{a:one-pants}
 \mathbb E_gD_{g,1}^{L,b}(F)
 =\frac{1}{2V_g}\int_{\mathcal D_{L,b}}
 F(\ell(x,y,z))V_{g-2,3}(x,y,z)xyz\,dx\,dy\,dz,
\end{equation}
agreeing with \cite[Eq.~(3.2)]{HSWX}.  Appendix~\ref{app:genus} gives
the complete two-pants formula, including its genus-three endpoint.
The same labeling proof permits any nonnegative cuff kernel invariant
under exchanging the two winding labels in each slot; we use this
extension for short-cuff indicators.

\subsection{Comparison with the exact one-point measure}

When the joint complement is connected, its volume is the only factor
in the cuff density that couples the different slots. To compare the tuple
integral with a product, define the one-point measures on cuff triples
and on lengths by
\begin{equation}\label{a:one-point-measures}
 d\eta_g^{L,b}(\mathbf x)
 =\frac{\mathbf1_{\mathcal D_{L,b}}(\mathbf x)}{2V_g}
 V_{g-2,3}(\mathbf x)xyz\,d\mathbf x,
 \qquad \nu_g^{L,b}=\ell_*\eta_g^{L,b}.
\end{equation}
Thus $\nu_g^{L,b}$ is the exact expected length measure of the retained
good-support curves.  Let $\mathcal C_{g,m}^{L,b}(F)$ denote the
expectation in \eqref{a:disjoint-count} restricted to connected joint
complement.  Such a complement has type $S_{g-2m,3m}$. Restoring all selected pants
except $P_i$ to that connected component shows that
$X\setminus\overline P_i$ is connected, for every $i$.  The single-block term in
Proposition~\ref{prop:exact-integration}, used here for $g\geq2m+1$
so that the complementary genus is positive, is therefore
\begin{equation}\label{a:connected-integral}
 \mathcal C_{g,m}^{L,b}(F)
 =\frac1{2^mV_g}\int_{\mathcal D_{L,b}^{\,m}}
 F\bigl(\ell(\mathbf x_1),\ldots,\ell(\mathbf x_m)\bigr)
 V_{g-2m,3m}(\mathbf x)\prod_{i,a}x_{i,a}\,d\mathbf x.
\end{equation}

\begin{proposition}[Uniform connected-complement factorization]
\label{prop:factorization}
Fix $m\geq1$ and $A>0$.  There are $C_m<\infty$ and $g_0(m,A)$ such
that, whenever $g\geq g_0(m,A)$ and $0\leq b\leq L\leq A\log g$,
\begin{equation}\label{a:factorization}
 \left|\mathcal C_{g,m}^{L,b}(F)
       -(\nu_g^{L,b})^{\otimes m}(F)\right|
 \leq \varepsilon_{g,m,L}(\nu_g^{L,b})^{\otimes m}(|F|),
 \qquad \varepsilon_{g,m,L}=C_m\frac{1+L^2}{g},
\end{equation}
for every real or complex Borel $F$ integrable with respect to
$(\nu_g^{L,b})^{\otimes m}$.  The expected sum of $|F|$ over tuples with connected joint complement
is finite.
The same estimate holds directly on cuff space with
$(\eta_g^{L,b})^{\otimes m}$, for any integrable cuff kernel invariant
under winding-cuff exchange in each slot.  Constants are independent of
$b,L,F$ in the stated range; no uniformity in growing $m$ is asserted.
\end{proposition}

\begin{proof}
Put $p(x)=\sinh(x/2)/(x/2)$, with $p(0)=1$, and
$P_n(\mathbf x)=\prod_{j=1}^np(x_j)$.  For positive genus $h$ and
$n\geq1$, \cite[Lemma~22]{NWX} gives
\begin{equation}\label{a:boundary-profile}
 q_{h,n}(\mathbf x)
 :=\frac{V_{h,n}(\mathbf x)}{V_{h,n}P_n(\mathbf x)},
 \qquad 0<q_{h,n}\leq1,
 \qquad 0\leq1-q_{h,n}(\mathbf x)
       \leq\frac{c(n)}h\sum_jx_j^2.
\end{equation}
This estimate holds for all nonnegative lengths, with $c(n)$ independent
of $h$ and of the lengths. In the logarithmic range used below, the
deficit tends to zero, so the lower bound also controls the denominator
in the density comparison.

Substituting $(h,n)=(g-2m,3m)$ and $(g,0)$ into the fixed-boundary-number
asymptotic \eqref{f:volume-asymptotic} gives
\begin{equation}\label{a:volume-ratio}
 \begin{split}
 r_{g,m}:=\frac{V_{g-2m,3m}}{V_g}
 &=(4\pi^2)^{-m}\sqrt{\frac g{g-2m}}
   \frac{(2g-m-3)!}{(2g-3)!}\bigl(1+O_m(g^{-1})\bigr)\\
 &=(8\pi^2g)^{-m}\bigl(1+O_m(g^{-1})\bigr).
 \end{split}
\end{equation}
Indeed, the factorial quotient is
$\prod_{j=0}^{m-1}(2g-3-j)^{-1}=(2g)^{-m}(1+O_m(g^{-1}))$, and the
square-root factor is $1+O_m(g^{-1})$.  In particular
\begin{equation}\label{a:ratio-of-ratios}
 R_{g,m}:=\frac{r_{g,m}}{r_{g,1}^m}=1+O_m(g^{-1}).
\end{equation}

For $\mathbf x\in[0,L]^{3m}$, write
$q=q_{g-2m,3m}(\mathbf x)$,
$q_i=q_{g-2,3}(\mathbf x_i)$, and $T=\prod_iq_i$.
Because $q_i\in(0,1]$, \eqref{a:boundary-profile} implies
\begin{equation}\label{a:product-deficit}
 1-T\leq\sum_i(1-q_i)
       \leq\frac{c(3)}{g-2}\sum_{i,a}x_{i,a}^2,
 \qquad
 1-q\leq\frac{c(3m)}{g-2m}\sum_{i,a}x_{i,a}^2.
\end{equation}
Choose $g_0(m,A)$ so that $g\geq4m+2$ and
$3mc(3)A^2(\log g)^2/(g-2)\leq1/2$ for all $g\geq g_0(m,A)$.
Then $T\geq1/2$ and
\[
 \left|\frac qT-1\right|
 \leq2\bigl((1-q)+(1-T)\bigr)
 \leq\frac{C_m}{g}\sum_{i,a}x_{i,a}^2.
\]
On the common support $\mathcal D_{L,b}^{\,m}$, the ratio of the
density in \eqref{a:connected-integral} to the product density from
\eqref{a:one-point-measures} is
\begin{equation}\label{a:pointwise-density}
 \Theta_{g,m}(\mathbf x)
 =\frac{V_{g-2m,3m}(\mathbf x)/V_g}
        {\prod_i[V_{g-2,3}(\mathbf x_i)/V_g]}
 =R_{g,m}\frac qT,
 \qquad
 |\Theta_{g,m}(\mathbf x)-1|
 \leq C_m\frac{1+\sum_{i,a}x_{i,a}^2}{g}
 \leq C'_m\frac{1+L^2}{g}.
\end{equation}
The factors $2^{-m}$, the cuff Jacobians, and the boundary profiles
cancel exactly.  Integrate this pointwise comparison against the product
one-point density and $|F(\ell_1,\ldots,\ell_m)|$ over the unchanged
domain.  This proves \eqref{a:factorization} and absolute integrability.
The same argument with a cuff kernel proves the asserted extension.
\end{proof}

For nonnegative slot functions $f_i$ with finite one-point integrals,
Tonelli's theorem yields the particularly useful consequence
\begin{equation}\label{a:slot-factorization}
 \left|\mathcal C_{g,m}^{L,b}\Bigl(\prod_i f_i\Bigr)
         -\prod_i\nu_g^{L,b}(f_i)\right|
 \leq\varepsilon_{g,m,L}\prod_i\nu_g^{L,b}(f_i).
\end{equation}
If a factor vanishes, both sides being compared vanish.  Signed slot
functions instead require the product of their absolute integrals in
the error bound.

Under the hypotheses of Proposition~\ref{prop:factorization}, the
corresponding leading one-point profile is
\begin{equation}\label{a:leading-profile}
 d\overline\eta_g^{L,b}(\mathbf x)
 =\frac{\mathbf1_{\mathcal D_{L,b}}(\mathbf x)}{2\pi^2g}
   \sinh(x/2)\sinh(y/2)\sinh(z/2)\,d\mathbf x.
\end{equation}
Equations~\eqref{a:boundary-profile} and \eqref{a:volume-ratio} also
compare the connected density with
$(\overline\eta_g^{L,b})^{\otimes m}$ with relative error
$O_m((1+L^2)/g)$.  In particular its coefficient is
$2^{-m}8^m(8\pi^2g)^{-m}=(2\pi^2g)^{-m}$.
The exact length domain is retained in this comparison and in the
derivation of the limiting intensity.

\subsection{Disconnected complementary components}

It remains to estimate the terms with disconnected joint complement.
Each additional complementary component contributes an extra inverse
power of $g$ through the genus sum. Let $D_{g,m;q}^{L,b}(F)$ denote \eqref{a:disjoint-count} restricted to
exactly $q$ joint complementary components, and let
$D_{g,m}^{\mathrm{disc},L,b}(F)=\sum_{q\geq2}D_{g,m;q}^{L,b}(F)$.
For $m=1$ this disconnected count is identically zero.

\begin{proposition}[Suppression of disconnected complements]
\label{prop:disconnected}
For every fixed $m\geq2$, there are $C_m<\infty$ and $g_0(m)$ such that
for all $g\geq g_0(m)$, $L\geq0$, $0\leq b\leq L$, integers $q\geq2$,
and bounded real or complex Borel $F$,
\begin{equation}\label{a:disconnected-bound}
 \mathbb E_g|D_{g,m;q}^{L,b}(F)|
 \leq C_m\|F\|_\infty g^{-(q-1)}
                  \left(\frac{Le^L}{g}\right)^m,
 \qquad
 \mathbb E_g|D_{g,m}^{\mathrm{disc},L,b}(F)|
 \leq\frac{C_m\|F\|_\infty}{g}
                  \left(\frac{Le^L}{g}\right)^m.
\end{equation}
The constants are uniform in $L,b,q,F$ in this range.  In particular
$b=0$ is included and no logarithmic restriction on $L$ is needed.
At $L=a_g+S$, for fixed $S$, these bounds become respectively
$O_{m,S}(\|F\|_\infty g^{-(q-1)})$ and
$O_{m,S}(\|F\|_\infty g^{-1})$.
\end{proposition}

\begin{proof}
For an admissible partition with $q$ blocks, write $n_C=|C|$.  Then
$\sum_Cn_C=3m$, all $n_C\geq2$, and the genus condition is equivalently
\begin{equation}\label{a:complexity-sum}
 \sum_C(2h_C-2+n_C)=2g-2-m.
\end{equation}
Lemma~\ref{lem:genus-convolution}, with $a=m$ and $N=3m$, gives
the full stable-genus sum
\begin{equation}\label{a:genus-bound-used}
 \frac1{V_g}\sum_{\substack{h_C\geq0,\ 2h_C-2+n_C\geq1\\
                         \sum_C(2h_C-2+n_C)=2g-2-m}}
          \prod_CV_{h_C,n_C}
 \leq C_mg^{-m-q+1}.
\end{equation}
That proof includes the genus-zero endpoints.  By
Lemma~\ref{a:global-boundary-bound}, also in the appendix, every stable
type with $2\leq n\leq3m$ satisfies
$V_{h,n}(\mathbf x)\leq B_mV_{h,n}P_n(\mathbf x)$ for all nonnegative
lengths.  Since every cuff belongs to precisely one complementary block,
\begin{equation}\label{a:profile-product}
 \prod_CP_{n_C}(\mathbf x_C)\prod_{i,a}x_{i,a}
 =2^{3m}\prod_{i=1}^m\Phi(\mathbf x_i),
 \qquad \Phi(x,y,z)=\sinh(x/2)\sinh(y/2)\sinh(z/2).
\end{equation}
Put $J(t)=\int_0^t\sinh(u/2)\,du
=2(\cosh(t/2)-1)\leq e^{t/2}$.  By
\eqref{a:domain-inclusion},
\begin{equation}\label{a:profile-integral-bound}
 \int_{\mathcal D_{L,b}}\Phi
 \leq J(L)\int_0^L\sinh(x/2)J(L-x)\,dx
 \leq\frac L2e^L.
\end{equation}
Here the enlarged domain is used only for an upper bound.

Apply Proposition~\ref{prop:exact-integration} with $F=1$, restricting
its partition sum to $q$ blocks.  Bound the volumes by their profiles,
use \eqref{a:genus-bound-used}, and integrate using
\eqref{a:profile-product}--\eqref{a:profile-integral-bound}.  The number of partitions is at most the Bell number for a set of $3m$
labels.  The partition count and
$B_m^q$ can be absorbed into $C_m$, because
$q\leq\lfloor3m/2\rfloor$.  Consequently
\[
 \mathbb E_gD_{g,m;q}^{L,b}(1)
 \leq C_mg^{-m-q+1}2^{-m}2^{3m}
                    (Le^L/2)^m
 \leq C'_mg^{-m-q+1}(Le^L)^m.
\]
The pointwise bound
$|D_{g,m;q}^{L,b}(F)|\leq\|F\|_\infty D_{g,m;q}^{L,b}(1)$ proves the
first assertion.  Summing over the finite $q$ range proves the second.
Finally, writing $H=\log g$ and $L=H-\log H+S$ gives
$Le^L/g=e^SL/H=O_S(1)$, proving the critical-cutoff statements.
\end{proof}

\subsection{Critical one-point measures}

We now pass to centered lengths at the critical cutoff and record the retained
one-point measure. Fix $S\in\mathbb R$ and
henceforth take
$H=\log g$, $a_g=H-\log H$, $L=a_g+S$, and $b=10\log L$,
for sufficiently large $g$.  Let $\widehat\nu_{g,S}$ be the pushforward
of $\nu_g^{L,b}$ by $t\mapsto t-a_g$.
These measures have uniformly bounded total mass. Using the finite
constant $C_V$ from \eqref{f:one-pants-volume-ratio}, the upper bound in
\eqref{a:boundary-profile}, and \eqref{a:profile-integral-bound}, we obtain
\begin{equation}\label{a:mass-bound}
 M_b:=\eta_g^{L,b}((0,\infty)^3)
 \leq M_0:=\eta_g^{L,0}((0,\infty)^3)
 \leq2C_V\frac{Le^L}{g}
 =2C_Ve^S\frac LH=O_S(1).
\end{equation}

Let $N_g(s)=\Xi_g^{(8)}(( -\infty,s])$ be the full figure-eight
count with centered length at most $s$. For each fixed finite $s\leq S$,
\cite[Proposition~3.1]{HSWX}, applied at $a_g+s$, gives
\begin{equation}\label{a:full-one-point-expansion}
 \mathbb E_gN_g(s)
 =\frac{e^s}{8\pi^2}
     \left(1-\frac{\log H+c_0-s}{H}\right)
  +O_s(g^{-1/2}H^{3/2}+g^{-1}H^3),
 \qquad c_0=3+4\log2.
\end{equation}
For $s\leq S$, the number deleted below $a_g+s$ is bounded pointwise
by $R_{g,S}$. Thus Proposition~\ref{prop:trim} applies with the common
cutoffs $L,b$; the cuff cutoff remains fixed as $s$ varies. Therefore
\begin{equation}\label{a:centered-one-point}
 \widehat\nu_{g,S}(( -\infty,s])
 =\frac{e^s}{8\pi^2}+O_{S,s}(\log H/H).
\end{equation}
Subtracting at finite endpoints gives, for each fixed $u<v\leq S$,
\begin{equation}\label{a:interval-intensity}
 \widehat\nu_{g,S}((u,v])
 =\lambda((u,v])+O_{S,u,v}(\log H/H),
 \qquad \lambda(ds)=\frac{e^s}{8\pi^2}\,ds.
\end{equation}
Equation~\eqref{a:centered-one-point} handles lower-unbounded intervals
directly. The intensity is obtained from counts of the actual
unoriented figure-eights with their original length cutoff.

\subsection{All fixed disjoint mixed moments}

We now combine the product approximation for the connected complement
with the bound for all other complements. For a bounded
Borel function $F:\mathbb R^m\to\mathbb C$, let $D_{g,m}(F)$ be the
retained disjoint tuple count with weight
$F(\ell_{\alpha_1}-a_g,\ldots,\ell_{\alpha_m}-a_g)$.
In particular, $D_{g,m}(J_1,\ldots,J_m)$ uses the product of the
indicators of slot intervals $J_i\subseteq(-\infty,S]$.

\begin{theorem}[Disjoint mixed moments]\label{thm:disjoint}
For every fixed $m\geq1$, $S\in\mathbb R$, and bounded real or complex
Borel $F$,
\begin{equation}\label{a:full-disjoint-measure}
 \left|\mathbb E_gD_{g,m}(F)
             -\widehat\nu_{g,S}^{\otimes m}(F)\right|
 \leq C_{m,S}\left[
     \frac{H^2}{g}\widehat\nu_{g,S}^{\otimes m}(|F|)
     +\frac{\|F\|_\infty}{g}\right].
\end{equation}
In particular, for fixed slot intervals $J_1,\ldots,J_m$, each of the
form $(u_i,v_i]$ or $(-\infty,v_i]$ with $v_i\leq S$,
\begin{equation}\label{a:disjoint-mixed-limit}
 \mathbb E_gD_{g,m}(J_1,\ldots,J_m)
 =\prod_{i=1}^m\lambda(J_i)
   +O_{m,S,\mathbf J}(\log H/H),
 \qquad \lambda(ds)=\frac{e^s}{8\pi^2}\,ds.
\end{equation}
Slot intervals may repeat.  Repeating a window $I_j$ in $r_j$ specified slots,
with $m=\sum_jr_j$, gives the mixed disjoint-selection count
$D_{g,\mathbf r}$ and
\begin{equation}\label{a:disjoint-multiindex}
 \mathbb E_gD_{g,\mathbf r}\longrightarrow
                   \prod_j\lambda(I_j)^{r_j}.
\end{equation}
For total order zero the conventional empty count and empty product are
both one.
\end{theorem}

\begin{proof}
For $m\geq2$, split the tuple measure into connected and disconnected
joint complements.  Proposition~\ref{prop:factorization} bounds the
connected difference by
$C_{m,S}H^2\widehat\nu_{g,S}^{\otimes m}(|F|)/g$.
Proposition~\ref{prop:disconnected} bounds the absolute disconnected
contribution by $C_{m,S}\|F\|_\infty/g$.  Their sum proves
\eqref{a:full-disjoint-measure}.  For $m=1$ the comparison is exact by
definition.  The total mass bound \eqref{a:mass-bound} makes the right
side $O_{m,S}(H^2/g)$ uniformly over $\|F\|_\infty\leq1$.

For a product of interval indicators, use
\eqref{a:centered-one-point}--\eqref{a:interval-intensity} in each slot.
The factors are bounded by \eqref{a:mass-bound}, so telescoping their
product gives an error $O_{m,S,\mathbf J}(\log H/H)$.
Combining with \eqref{a:full-disjoint-measure} proves
\eqref{a:disjoint-mixed-limit}.  Repeated intervals require no
multinomial coefficient: the slots already encode the specified order.
This proves \eqref{a:disjoint-multiindex}.
\end{proof}

% END sections/03_disjoint.tex

% BEGIN sections/04_witnesses.tex
\section{Witnesses, touching, and one crossed support}\label{sec:witnesses}

A disjoint tuple is lost under isolation if some retained curve outside
the tuple has closed support meeting a selected support. We call each
such curve a \emph{witness}. Even a single witness can meet several
selected supports. To count each lost tuple once, we divide its weight
equally among all its witnesses before grouping them by how they meet
the tuple.

We first treat witnesses whose supports coincide with a selected support
or touch one or more supports only along cuffs. For a crossing witness,
we choose one crossed support and estimate the pair it forms with the
witness. A joint cutting integral includes the other slots when the pair
has bounded complexity. At higher complexity we bound their multiplicity
on each surface. The localization result in Section~\ref{sec:localization}
then gives the estimate for witnesses crossing exactly one support.
Sections~\ref{sec:localization} and~\ref{sec:higher} treat witnesses
crossing several supports.

\subsection{Witnesses and their normalized shares}

Fix an integer \(m\geq1\).  Write
\[
 H=\log g,\qquad a_g=H-\log H,\qquad
 L=L_g=a_g+S,\qquad b=b_g=10\log L ,
\]
where \(S\) is fixed and \(g\) is sufficiently large.  When an estimate is
stated for arbitrary \(L,b\), these symbols instead denote its indicated
independent parameters.  We use the retention rule of Section~\ref{sec:foundations}: a curve is
a primitive unoriented figure-eight of length at most \(L\), its support
is good, and all three cuffs have length at least \(b\).
Let \(\mathcal R(X)\) be the finite set of retained curves.

A \emph{mark} in this section is a retained curve selected for a
specified tuple slot. Let \(\mathcal T_m(X)\) be the set of ordered tuples of such curves
\(\tau=(\alpha_1,\ldots,\alpha_m)\) with pairwise disjoint closed supports.
Each slot may additionally be restricted to any prescribed Borel window
for \(\ell(\alpha_i)-a_g\) contained in \((-\infty,S]\).
Put \(P_i=P_{\alpha_i}\) and define
\begin{equation}\label{w:witness-set}
 \mathcal W(\tau)=
 \{\beta\in\mathcal R(X)\setminus\{\alpha_1,\ldots,\alpha_m\}:
       \overline P_\beta\cap\overline P_i\ne\varnothing
       \text{ for some }i\},
 \qquad d(\tau)=|\mathcal W(\tau)|.
\end{equation}
A witness may lie outside every slot window, provided it satisfies the
common retention rule up to the ceiling \(S\). Define the local loss by
\begin{equation}\label{w:local-loss}
 Q_m(X)=\sum_{\tau\in\mathcal T_m(X)}
                \mathbf1_{\{d(\tau)>0\}}.
\end{equation}
Thus \(Q_m\) counts exactly the disjoint marked tuples lost when all
nonisolated retained curves are deleted. Superscripts \({\rm conn}\)
and \({\rm disc}\) refer to the complement of the union of the closed
marked supports; the witness is not removed when forming this complement.

For each mark, the witness support has exactly one of four relations:
closed-disjointness, equality, disjoint interiors with nonempty boundary
contact, or distinctness with intersecting interiors.  Denote these relations by
\(0,E,T,X\), respectively, with \(X\) denoting a crossing relation here.  Two distinct geodesics cannot be tangent.
Thus a \(T\)-contact consists precisely of common cuffs: a transverse
boundary crossing would force intersection of the chosen interior
half-disks.  Moreover, two distinct geodesic pants with intersecting
interiors have transversely crossing boundaries.  Indeed, without such
crossings, no nonshared boundary geodesic of one pair of pants could lie inside
the other, since a geodesic pair of pants has no essential simple closed
geodesics other than its cuffs.  An interior intersection would consequently put one support
inside the other.  Their equal area \(2\pi\) would force equality.

An equality at one mark forces closed-disjointness from all other marks.
We may therefore partition the witnesses into equal-support witnesses,
pure-touching witnesses, and crossing witnesses, with respective counts
\(d_E,d_T,d_X\). A pure-touching witness touches at least one mark
and crosses none. A witness touching one mark and crossing another belongs
to the last class. Give each witness the share \(1/d(\tau)\) when
\(d(\tau)>0\). With the convention \(0/0=0\), set
\begin{equation}\label{w:normalized-partition}
 Q_{A,m}=\sum_{\tau\in\mathcal T_m(X)}\frac{d_A(\tau)}{d(\tau)},
 \qquad A\in\{E,T,X\}.
 \qquad
 Q_m=Q_{E,m}+Q_{T,m}+Q_{X,m}.
\end{equation}
Each affected tuple contributes one to this decomposition. A witness
meeting several marks is counted once in its class. In every subsequent
subdivision or removal, the denominator remains the original
\(d(\tau)\) from \eqref{w:witness-set}.  In particular, for any witness subclass \(A\),
\begin{equation}\label{w:raw-domination}
 0\leq Q_{A,m}\leq I_{A,m}:=\sum_{\tau\in\mathcal T_m(X)}d_A(\tau),
 \qquad
 Q_{A,m}^{\rm disc}\leq
       |\mathcal T_m^{\rm disc}(X)|.
\end{equation}
The disconnected estimate uses the fact that the normalized share of
each tuple is at most one; the raw incidence can be larger.

These are measurable counts.  On a reference marked surface there are
countably many primitive unoriented free-homotopy classes.  Their length
conditions are Borel, and support relations are determined by the
geodesic subsurface classes.  A compact hyperbolic surface has only
finitely many length-bounded closed geodesics.  The resulting counts and
ratios are Borel and invariant under mapping-class reindexing.
Nonnegative sums below are consequently subject to Tonelli's theorem.

Different witnesses on the same surface can have the same filled
subsurface type. To recover a witness, cut along the \(3m\) marked cuffs.
Record the coincident cuffs and their sides, the closed witness cuffs
in complementary pieces, and the simultaneous systems of proper arcs
formed by crossing cuffs. An essential closed witness cuff cannot lie
strictly inside a marked pair of pants, and a boundary-parallel
geodesic arc would form a bigon. The endpoint orders, gluing matches
and twist data specify how to reassemble the cuffs. The pants side and
the unordered winding pair then specify the figure-eight. These
relative isotopy data are countable. The local kernels below sum over
all witness configurations of the given type.

\begin{lemma}[Equal-support loss]\label{w:equal-support}
For every fixed \(m\) and \(S\),
\[
 \mathbb E_g I_{E,m}=O_{m,S}(H^{-6}+g^{-1}).
\]
The assertion includes disconnected marked complements and arbitrary
slot windows.
\end{lemma}

\begin{proof}
Use the cuff measure \(\eta_g^{L,b}\) from
\eqref{a:one-point-measures} and its mass
\(M_b=\eta_g^{L,b}(\mathbb R_+^3)\), bounded in \eqref{a:mass-bound}.  The first two coordinates of a
cuff triple are its winding cuffs. Define
\[
 u_L(x,y,z)=
  \mathbf1_{\{\ell(x,z,y)\leq L\}}+
  \mathbf1_{\{\ell(y,z,x)\leq L\}}.
\]
The kernel counts the other two figure-eights on the same support and is
invariant under exchanging \(x,y\).  Their cuff restrictions and
connected-complement condition are inherited from the mark.
The ordered same-pants estimate in
Proposition~\ref{prop:f-same-pants} gives
\begin{equation}\label{w:equal-kernel}
 \int u_L\,d\eta_g^{L,b}
 =\mathbb E_g Q_g^{\rm same}(L,b)
 \leq 8\frac{V_{g-2,3}}{V_g}e^{L-b/2}.
\end{equation}
There is no additional factor two: the right-hand count is already an
ordered count of two distinct actual curves.

Drop the marked windows.  Apply the pointwise cuff-density comparison of
Proposition~\ref{prop:factorization} to the nonnegative kernel
\(\sum_{i=1}^m u_L(\mathbf x_i)\).  It gives
\[
 \mathbb E_g I_{E,m}^{\rm conn}
 \leq
 \left(1+O_m\!\left(\frac{1+L^2}{g}\right)\right)
 m M_b^{m-1}\mathbb E_g Q_g^{\rm same}(L,b).
\]
At the retained cutoff \(M_b=O_S(1)\), and
\eqref{w:equal-kernel} is \(O_S(H^{-6})\).
Each pair of pants supports exactly three primitive unoriented
figure-eights.
Consequently a tuple has at most \(2m\) equal-support witnesses, so
\[
 I_{E,m}^{\rm disc}\leq
 2m|\mathcal T_m^{\rm disc}(X)|.
\]
Proposition~\ref{prop:disconnected} proves the remaining
\(O_{m,S}(g^{-1})\) bound.  At \(m=1\) this disconnected term is zero.
\end{proof}

\subsection{Pure touching at every fixed marked order}

For pure touching, the marked supports and witness support have disjoint
interiors and meet along common cuffs. We integrate over their union,
including every topology of its possibly disconnected complement.
For independent \(L\geq1\) and \(0\leq b\leq L\), write
\[
 \mathcal A_g(L)=\frac{Le^L}{g},
 \qquad \Theta(b)=(1+b)e^{-b}.
\]
Let \(I_{T,m}^{\rm conn}(L,b)\) count ordered marked tuples with connected
joint complement and one distinguished pure-touching witness; all
curves have the good-support and cutoff conditions just described.
This count has no slot-window restrictions. Inserting such restrictions
only decreases it.

\begin{proposition}[Pure-touching incidence]\label{prop:touching}
For fixed \(m\), sufficiently large \(g\), all \(L\geq1\), and
\(0\leq b\leq L\),
\begin{equation}\label{w:touching-bound}
 \mathbb E_g I_{T,m}^{\rm conn}(L,b)
 \leq C_m\mathcal A_g(L)^{m+1}
 \left((1+g^{-1})\Theta(b)+L^{-1}\Theta(b)^2\right).
\end{equation}
In particular, at the retained cutoff,
\begin{equation}\label{w:touching-critical}
 \mathbb E_g I_{T,m}^{\rm conn}
 =O_{m,S}\!\left(\frac{\log H}{H^{10}}\right),
 \qquad
 \mathbb E_g Q_{T,m}
 =O_{m,S}\!\left(\frac{\log H}{H^{10}}+g^{-1}\right).
\end{equation}
\end{proposition}

\begin{proof}
Put \(p=m+1\), with witness slot \(0\) and marked slots
\(1,\ldots,m\).  Label the winding cuffs of each selected pair of pants by \(1,2\),
and its third cuff by \(3\). A nonempty matching \(M\) identifies \(k\) witness cuffs
with \(k\) distinct marked cuffs.  Let \(\pi\) be the partition of the \(3p-2k\) unmatched labels into
complementary blocks \(A\), with \(n_A=|A|\) and stable genus \(h_A\).
Write \(q=|\pi|\). The incidence multigraph \(G\) has the selected pants and the complementary
blocks as vertices; matched cuffs give witness--mark edges and every
unmatched cuff joins its pants to its complementary block.

The exact topological conditions are
\begin{equation}\label{w:touching-admissibility}
 \begin{gathered}
  G-\{1,\ldots,m\}\text{ is connected},\qquad
  G-\{0\}\text{ is connected},\\
  h_A\geq0,\qquad 2h_A-2+n_A>0,\qquad
  \sum_A h_A=g+q+k-2m-3 .
 \end{gathered}
\end{equation}
The first graph condition is connectedness of the marked complement,
and the second is connectedness of the witness complement.  Every
complementary block must meet both a witness cuff and a marked cuff.
For \(m\geq2\) some marked cuff is unmatched.  For \(m=1\), the only case
with no complementary piece glues two pants along all three cuffs and
has genus two.  Excluding this small genus gives
\begin{equation}\label{w:touching-topology}
 1\leq k\leq2,\qquad 1\leq q\leq3-k,\qquad n_A\geq2.
\end{equation}
The Euler equation in \eqref{w:touching-admissibility} follows from
\(\sum_A n_A=3(m+1)-2k\).

Conversely, glue the labeled pants and stable complementary pieces
according to these data.  The graph conditions give the stated
connected complements.  For each mark, restoring all the other marked pants to the connected
marked complement gives its connected individual complement.  All pieces have negative Euler characteristic;
a disk or annulus bounded by one or two distinct gluing curves would
contradict additivity of Euler characteristic.  Thus the cut curves are
essential and distinct except for the prescribed matched labels.
No edge joins two marks.  The marked closed supports are therefore
disjoint, and classification of labeled oriented bordered surfaces
gives one orbit per datum.  Positive matching boundary lengths and
arbitrary gluing twists realize the topological gluing.

Let \(\mathcal E_M\) be the \(3p-k\) distinct cuff classes after matching,
with lengths \(t_e\). Let \(\mathbf t_v\) be the incident cuff triple
at pants \(v\), and \(\mathbf t_A\) the boundary-length vector of
complementary block \(A\). Recall the exact length domain from
\eqref{a:length-domain}:
\begin{equation}\label{w:true-pants-domain}
 \begin{gathered}
 \ell(x,y,z)=2\operatorname{arccosh}
 \bigl(\cosh(z/2)+2\cosh(x/2)\cosh(y/2)\bigr),\\
 \mathcal D_{L,b}=
 \{(x,y,z)>0:\min(x,y,z)\geq b,\quad\ell(x,y,z)\leq L\}.
 \end{gathered}
\end{equation}
The ordered oriented multicurve formula
\cite[Theorem~2.5 and Remark~2.7]{AMI} gives
\begin{equation}\label{w:touching-exact-integral}
 \mathbb E_g I_{T,m}^{\rm conn}(L,b)
 =\frac{1}{2^pV_g}
 \sum_{M,\pi,\mathbf h\ {\rm admissible}}
 \int
 \prod_{v=0}^m\mathbf1_{\mathcal D_{L,b}}(\mathbf t_v)
 \prod_{A\in\pi}V_{h_A,n_A}(\mathbf t_A)
 \prod_{e\in\mathcal E_M}t_e\,d\mathbf t .
\end{equation}
One can order the cut curves by their least cuff label, orient a shared
curve with the marked pants on the left, and orient an unshared curve
with its selected pants on the left.  This recovers the labeled selected
pieces and makes the coefficient of each labeled orbit one.
Every actual marked tuple and distinguished witness has exactly \(2^p\)
winding-cuff labelings.  There is no tuple factorial, complementary-block
factorial, or contact-number factor.  A shared cuff has a single length
variable and a single twist integral, hence a single factor \(t_e\).
It occurs in no complementary boundary vector.  Since \(n_A\geq2\),
there is no exceptional one-holed-torus convention in this exact
formula.  Additional slot indicators may be inserted without changing
the coefficient.

For fixed \(m\) there are finitely many matching and partition data.
Their complementary complexities sum to \(2g-2-(m+1)\), and
\(\sum_A n_A=3(m+1)-2k\leq3m+1\).  Lemma~\ref{lem:genus-convolution},
with \(a=m+1\) and \(N=3m+1\), therefore gives
\begin{equation}\label{w:touching-genus}
 \frac1{V_g}\sum_{\mathbf h\ {\rm admissible}}
       \prod_A V_{h_A,n_A}
 \leq C_m g^{-m-q}.
\end{equation}
Put \(\psi(t)=\sinh(t/2)/(t/2)\), with \(\psi(0)=1\).
Since \(2\leq n_A\leq3m+1\), Lemma~\ref{a:global-boundary-bound},
applied with boundary parameter \(m+1\), gives
\begin{equation}\label{w:touching-profile}
 V_{h,n}(\mathbf t)\leq C_m V_{h,n}\prod_i\psi(t_i).
\end{equation}
This includes the stable genus-zero pieces.  An unmatched cuff has
\(t\psi(t)=2\sinh(t/2)\leq e^{t/2}\), whereas a shared cuff keeps
just \(t\).

The true domain is contained in
\(\{x,y,z\geq b:x+y\leq L,\ z\leq L\}\).
For \(s\in\{0,1,2\}\) prescribed cuffs of total length \(T\), integrate
\(\exp(\frac12\sum u)\) in the remaining variables on this enlarged
domain and call the result \(B_s\).  Direct integration gives
\begin{equation}\label{w:touching-slices}
 B_0\leq4Le^L,\qquad
 B_1\leq4Le^{L-T/2},\qquad
 B_2\leq2e^{L-T/2}.
\end{equation}
For \(B_0\), the winding triangle
integrates to at most \(2Le^{L/2}\), and the third interval to at most
\(2e^{L/2}\).  If a prescribed cuff \(t\) is winding, the free winding
interval ends at \(L-t\); if it is the third cuff, the winding triangle
gives \(2Le^{L/2}\leq2Le^{L-t/2}\), because \(t\leq L\).
If two winding cuffs are prescribed, \(T\leq L\), and the free third
interval gives \(2e^{L/2}\leq2e^{L-T/2}\).  If a winding cuff \(t_w\)
and the third cuff \(t_z\) are prescribed, the free integral is at most
\(2e^{(L-t_w)/2}\leq2e^{L-(t_w+t_z)/2}\).
Infeasible prescribed data contribute zero.

Let \(r_2(M)\) count the selected pants incident to two shared cuffs.
It is zero for \(k=1\), two when two cuffs join the same mark, and one
when they join two different marks.  Holding the shared lengths fixed,
all unmatched variables separate between the selected pants.  Each
shared length enters the prescribed sum at two pants.  The resulting
boundary integral \(J_M\) is therefore bounded by
\begin{equation}\label{w:touching-length-integral}
 \begin{split}
 J_M&\leq4^p L^{p-r_2(M)}e^{pL}
       \prod_{j=1}^k\int_b^\infty t_j e^{-t_j}\,dt_j\\
    &=4^p L^{p-r_2(M)}e^{pL}\Theta(b)^k .
 \end{split}
\end{equation}
The extension of the shared-length domain is made after the sliced
bounds have been applied on feasible data.
Combining \eqref{w:touching-genus} and
\eqref{w:touching-length-integral} in
\eqref{w:touching-exact-integral} bounds each family by
\[
 C_m g^{-(q-1)}
       \mathcal A_g(L)^{m+1}L^{-r_2(M)}\Theta(b)^k.
\]
The possibilities \((k,q)=(1,1),(1,2),(2,1)\) give
\eqref{w:touching-bound}.

At the retained cutoff,
\(\mathcal A_g(L)=e^S L/H=O_S(1)\) and
\(\Theta(b)=(1+10\log L)L^{-10}=O_S(\log H/H^{10})\).
This proves the raw connected estimate.  On disconnected marked
complements use only
\(Q_{T,m}^{\rm disc}\leq|\mathcal T_m^{\rm disc}(X)|\) and
Proposition~\ref{prop:disconnected} for \(m\geq2\); for \(m=1\),
the disconnected term is zero. Together with \eqref{w:raw-domination},
this proves \eqref{w:touching-critical}.
\end{proof}

At order three one can enumerate the matching and partition data explicitly.
For \(k=1\) there are \(3\cdot9=27\) matchings.
For \(q=2\), the two unmatched witness labels distinguish the blocks;
of the \(2^8\) assignments of the remaining marked labels, precisely
\(2^3\) fail connectivity after the witness is deleted.  Hence the
families with \((k,q)=(1,1),(1,2),(2,1)\) have sizes \(27\),
\(27\cdot248=6696\), and \(\binom32\,9\cdot8=216\), respectively.  The last family includes 54 matchings to one
mark twice and 162 matchings to two marks.  These multiplicities cannot
replace the labeled integrals when cuff roles have different length
constraints.

\subsection{Local masses for complexity-two crossing pairs}

Choose a marked curve crossed by a witness and call it the
\emph{anchor}. For this ordered anchor/witness pair, let \(Y\) be the
canonical intrinsic geodesic filled surface. Its completed image
\(Z\subset X\) is obtained by identifying boundary circles of \(Y\)
that represent the same ambient geodesic.  Lemma~\ref{lem:f-pair-core} gives
\(|\chi(Y)|\geq2\).  The generating supports lie in \(Z\):
the straightened boundary classes are disjoint from their original
cuffs and cannot run strictly inside a geodesic pants, so choosing
the filled side retains each generating support.  The placement of other marked supports will be determined by
localization in Section~\ref{sec:localization}.

For \(|\chi(Y)|=2\), the intrinsic types are \(S_{1,2}\) and \(S_{0,4}\).
At most two intrinsic boundary copies can represent one ambient
geodesic, since the interior is embedded and the geodesic has two
sides.  If every boundary disappeared under these identifications, the
image would be a closed genus-two surface, open and closed in \(X\).
For \(g>2\) the possible completions are consequently
\[
 S_{1,2}\longrightarrow S_{1,2},\qquad
 S_{0,4}\longrightarrow S_{0,4}\ \hbox{or}\ S_{1,2}.
\]
In each case, \(Z\) has distinct outer boundary geodesics. In the last case,
the identified geodesic is an internal seam; cutting it recovers the
intrinsic \(S_{0,4}\).  The seam is determined by the pair's canonical
filled surface, and is not a freely marked curve.

Call a supporting pair of pants \(P\subset Z\) \emph{extendible} if attaching
one connected stable exterior along all outer boundaries of \(Z\)
makes the complement of \(P\) connected.  Equivalently, every component of \(Z\setminus\overline P\), understood
with its metric completion, must meet an outer boundary of \(Z\) not
occupied by a cuff of \(P\).
This follows by contracting the connected exterior to a vertex in
the complementary incidence graph.  It also proves necessity whenever
the ambient complement of \(P\) is connected, even if the actual
exterior of \(Z\) is disconnected.  The test is made in completed
\(Z\), so a determined internal seam is already glued.

Use \(\zeta=12\) for intrinsic type \(S_{1,2}\), \(\zeta=04^0\)
for intrinsic \(S_{0,4}\) with no identification, and \(\zeta=04^1\)
for intrinsic \(S_{0,4}\) with one boundary pair identified. Put
\[
 (h_{12},n_{12})=(h_{04^1},n_{04^1})=(1,2),\qquad
 (h_{04^0},n_{04^0})=(0,4).
\]
On the bordered moduli space
\(\mathcal M_{h_\zeta,n_\zeta}(\mathbf t)\) with labeled outer
boundaries, let \(N^\zeta_L(Z)\) count ordered crossing pairs of primitive unoriented
figure-eights of type \(\zeta\), whose completed image is \(Z\), and
whose supports are embedded and extendible, each with three distinct
cuff geodesics in \(Z\).  Impose length at most \(L>1\) and cuff threshold
\(10\log L\).  Define the following nonnegative integrals, initially allowing the
value \(+\infty\):
\begin{equation}\label{w:low-kernels}
 K^\zeta_L(\mathbf t)=
 \int_{\mathcal M_{h_\zeta,n_\zeta}(\mathbf t)}
        N^\zeta_L(Z)\,dV_{\rm WP}(Z),
 \qquad
 \mathcal J^\zeta_L=
 \int_{\mathbb R_+^{n_\zeta}}
        K^\zeta_L(\mathbf t)\prod_i t_i\psi(t_i)\,d\mathbf t.
\end{equation}
The kernels retain all internal crossing classes, including the seam
length and twist in type \(04^1\). To bound their mass, we attach a
connected exterior of large genus. Cutting along the outer cuffs
then expresses an expected count of actual pairs in terms of the
same local kernel.

\begin{lemma}[Complexity-two local mass]\label{w:low-local-mass}
For every \(L>1\),
\begin{equation}\label{w:low-local-bound}
 \mathcal J^\zeta_L\leq C e^{2L}.
\end{equation}
The constant is independent of \(L\).
\end{lemma}

\begin{proof}
Fix \(L\), and let a separate auxiliary genus \(G\) tend to infinity.
Attach a connected exterior to the completed core.  Stability of this
exterior preserves the essential distinct outer boundary classes and
creates no additional disk or annulus in the filling construction.
Extendibility makes both individual supports good in the auxiliary
surface.  Conversely, every good pair of this type with connected
exterior gives the specified local data.  Let \(U_G^\zeta(L)\) count
these ordered pairs.

Cutting the outer cuffs gives an exterior of type
\(S_{G-h_\zeta-n_\zeta+1,n_\zeta}\), and
\begin{equation}\label{w:auxiliary-exact}
 \mathbb E_G U_G^\zeta(L)=
 \frac{1}{n_\zeta!V_G}
 \int_{\mathbb R_+^{n_\zeta}}
 K^\zeta_L(\mathbf t)
 V_{G-h_\zeta-n_\zeta+1,n_\zeta}(\mathbf t)
 \prod_i t_i\,d\mathbf t.
\end{equation}
There are exactly \(n_\zeta!\) outer-boundary labelings.  For large
\(G\), core and exterior have distinct types, so there is no exchange
of these sides; neither cut piece has type \((1,1)\).  The local
kernel already counts actual ordered curves and needs no winding-cuff
division.  To justify this formula for an internal moduli weight,
disintegrate the Weil--Petersson measure over the cut lengths, the two
bordered moduli spaces and the outer twists.  The local count is
independent of the latter twists, whose integration gives
\(\prod_i t_i\).  Tonelli then produces \(K^\zeta_L\).

Put \(b=10\log L\). If \(b>L\), all the local counts vanish by
\eqref{f:actual-to-relaxed-lengths}. Otherwise, for \(\zeta=12,04^1\),
\(U_G^\zeta(L)\leq Q_{12}(L,b)+Q_{04^1}(L,b)\), since the latter count
includes every exterior. By \eqref{tail:other-small}, its expectation is
\(O(G^{-2}e^{2L})\), since \(L^5e^{-b/2}=1\), and both
\(L^6e^{-L/2}\) and \(L^5e^{-L}\) are bounded for \(L>1\).

For type \(04^0\), we likewise have \(U_G^{04^0}(L)\leq Q_{04^0}(L,b)\)
when \(b\leq L\). Lemma~\ref{tail:planar} bounds its expectation by
\(Ce^{2L}/G^2\) for \(L\geq L_0\) with \(L_0\) absolute.
For \(1<L<L_0\), ordering the winding cuffs in both slots gives four
distinct relaxed pairs per ordered actual pair, as in
Lemma~\ref{lem:f-eight-transfer}. By \cite[Proposition~4.16]{HSWX}
the expected relaxed count is at most \(CL e^{2L}/G^2\), and
\(L\leq L_0\) absorbs the extra factor. Hence for every \(L>1\) and every subtype,
\begin{equation}\label{w:auxiliary-pair-bound}
 G^2\mathbb E_G U_G^\zeta(L)\leq C e^{2L}
\end{equation}
for all sufficiently large \(G\), with \(C\) independent of \(L\).

In all three cases the exterior factorial index is \(2G-5\).
The fixed-boundary volume asymptotic
\cite[Theorem~1.2]{MZ} and \cite[Lemma~22]{NWX} therefore imply, for
every fixed positive boundary vector,
\begin{equation}\label{w:auxiliary-volume-limit}
 G^2\frac{V_{G-h_\zeta-n_\zeta+1,n_\zeta}(\mathbf t)}{V_G}
 \longrightarrow
 (8\pi^2)^{-2}\prod_i\psi(t_i).
\end{equation}
Apply Fatou's lemma to \(G^2\) times
\eqref{w:auxiliary-exact}.  Equations
\eqref{w:auxiliary-pair-bound}--\eqref{w:auxiliary-volume-limit} give
\[
 \frac{\mathcal J^\zeta_L}{n_\zeta!(8\pi^2)^2}
 \leq\liminf_{G\to\infty}G^2\mathbb E_G U_G^\zeta(L)
 \leq C e^{2L}.
\]
This proves both integrability and the claimed bound.  One may first
truncate the local count and the boundary domain and then use monotone
convergence.  Uniform convergence in boundary lengths is unnecessary.
Here \(L\) is fixed during the auxiliary-genus limit. The resulting
constant is independent of \(L\), so the bound then applies at
\(L=L_g\).
\end{proof}

\subsection{Exterior marks and the pair-complexity tail}

Fix an integer \(r\geq0\). An exterior configuration consists of an
ordered crossing anchor/witness pair and \(r\) additional retained
marked curves in specified slots. The additional curves have pairwise
disjoint closed supports, all disjoint from the completed pair image
\(Z\).
No joint-complement condition is imposed.  Let \(I^{\rm ext}_r\)
count these ordered configurations.  For intrinsic type \(S_{h,n}\)
and \(j\) identified boundary pairs, write
\(I^{\rm ext}_{h,n,j;r}\) for the restricted count. At complexity two,
write \(I^{\zeta,\rm ext}_r\) for the corresponding subtype count.
Put \(c=2h-2+n\) and
\(\nu=n-2j\).  For fixed \(c\) and large genus, \(\nu\geq1\);
otherwise the completed image would be the entire closed surface
of genus \((c+2)/2\).

Appendix~\ref{app:attachments} supplies separate bounds for the core
and for the addition of exterior marks. Lemma~\ref{lem:fixed-core-kernel}
bounds the local mass \(\mathcal J_{L,b}^{h,n,j}\) defined in
\eqref{att:fixed-kernel-def}: for fixed \(c\geq3\), \(L\geq1\),
and \(b\geq0\),
\begin{equation}\label{w:fixed-exterior-kernel}
 \mathcal J_{L,b}^{h,n,j}
 \leq C_{c,\delta}(1+L)^{3c+1-2j}e^{(2+2\delta)L},
 \qquad 0<\delta<\tfrac12.
\end{equation}
It uses the ordered filling-pair theorem on intrinsic \(Y\), the
perimeter bound \(p(Y)\leq4L\), and the fixed-type volume
polynomial.  Lemma~\ref{lem:attachment-bounds}, with no attached
pants, inserts either this kernel or
\eqref{w:low-local-bound} into one joint cutting integral with the
exterior marks.  It yields
\begin{align}
 \mathbb E_g I^{\zeta,\rm ext}_r(L)
 &\leq C_r\mathcal A_g(L)^r\,g^{-2}e^{2L},
       \label{w:low-exterior}\\
 \mathbb E_g I^{\rm ext}_{h,n,j;r}(L,b)
 &\leq C_{c,r,\delta}\mathcal A_g(L)^r\,g^{-c}
       (1+L)^{3c+1-2j}e^{(2+2\delta)L}.
       \label{w:fixed-exterior}
\end{align}
The first estimate uses \(L>1\) and the cuff threshold \(10\log L\);
the second holds for \(L\geq1\) and \(b\geq0\). Both are raw incidence
bounds and include disconnected complements. The genus exponent follows
from the complementary topology: a partition of the \(\nu+3r\) exposed
cuffs into \(q\) complementary blocks has
\[
 \sum_A(2h_A-3+n_A)=2g-q-r-c-2.
\]
Its factorial gap against \(2g-3\) is \(c+r+q-1\), so its volume
ratio is \(O_{c,r}(g^{-c-r-q+1})\).  Each exterior pair of pants contributes at most \(4Le^L\).
Combining these factors inside the joint integral gives
\(\mathcal A_g(L)^r g^{-c}\) after summing the complementary types.

At the critical cutoff, \eqref{w:low-exterior} is \(O_{r,S}(H^{-2})\)
for all three types.
Taking \(\delta=1/8\) in \eqref{w:fixed-exterior} gives
\begin{equation}\label{w:fixed-exterior-critical}
 O_{c,r,S}\!\left(g^{-c+9/4}H^{3c-5/4-2j}\right).
\end{equation}
In a fixed finite band \(3\leq c\leq M\), take the maximum of the
finitely many constants and sum the displayed bounds. This gives
\(O_{M,r,S}(g^{-3/4}H^{31/4})\), with no uniformity in growing \(M\).

For growing pair complexity, we instead bound the number of additional
curves on each metric surface. The next estimate permits arbitrary
placements and increases the complexity threshold with the number of
additional slots. In the next proposition, \(K\geq3\) and \(r\geq0\)
are integers. Both members of each crossing pair have length at most
\(L\), and all their support cuffs have length at least \(b\).

\begin{proposition}[High pair complexity with additional marks]
\label{prop:witness-pair-tail}
Let \(P_{\geq K}^{\rm ord}(X;L,b)\) count ordered good crossing pairs
of intrinsic complexity at least \(K\).  Let \(J_{\geq K;r}(X;L,b)\)
count such pairs together with an ordered list of \(r\) additional primitive unoriented
geodesics of length at most \(L\); these additional curves may repeat
and have arbitrary placements.  For fixed \(K\geq3,r,S\), at
\(L=H-\log H+S\), uniformly in \(b\geq0\),
\begin{equation}\label{w:pair-tail-bound}
 \mathbb E_g J_{\geq K;r}
 \leq C_{K,r,S}g^{10+2r-K}H^{-r}.
\end{equation}
In particular \(K=11+2r\) makes this \(O_{r,S}(g^{-1}H^{-r})\).
\end{proposition}

\begin{proof}
Apply the general growing-core estimate of Lemma~\ref{lem:uniform-tail},
proved in Section~\ref{sec:higher}, with $q=2$.  Choose its tuple
subfamily to consist of the ordered good crossing pairs counted here.
Lemmas~\ref{lem:f-pair-core} and~\ref{l:entry-containment} show that
their actual curves have connected union and fill the canonical
intrinsic pair core.  The additional $r$ curves have exactly the
unrestricted placements and repetition convention in
\eqref{h:enlarged-tail-count}.  Hence \eqref{h:uniform-tail} gives
\[
 \mathbb E_g J_{\geq K;r}
 \leq A_{2,K,r,S}g^{10+2r-K}H^{-r}.
\]
This estimate is uniform in $b\geq0$.  Taking $K=11+2r$ proves the
last assertion.  The proof of Lemma~\ref{lem:uniform-tail} uses the
filling-face bound and the imported primitive-count and weighted-subsurface
estimates, independently of this proposition.
\end{proof}

Combining \eqref{w:low-exterior},
\eqref{w:fixed-exterior-critical}, and
Proposition~\ref{prop:witness-pair-tail} with the fixed splitting
threshold \(K=11+2r\) proves
\begin{equation}\label{w:all-exterior}
 \mathbb E_g I^{\rm ext}_r=O_{r,S}(H^{-2}).
\end{equation}
The low-complexity contribution uses the retained cuff threshold.
Only the tail and the complexity-at-least-three bounds were uniform
in an arbitrary \(b\geq0\).

\subsection{Boundary placements and exactly one crossed mark}

Call a tuple/witness configuration \emph{boundary-only} relative to
its crossed anchor if every other marked interior avoids the completed
pair image \(Z\), but at least one other closed marked support meets
\(Z\).  The contacts are outer common cuffs.  Pairwise closed-disjointness
of the marks makes their used cuff sets disjoint.  The witness lies in
\(Z\), so it crosses only this anchor, although it may touch other
marks.

Appendix~\ref{app:attachments} proves that the complexity-two
boundary-only patterns exhaust the following list:
one single-cuff attachment for intrinsic \(S_{1,2}\);
one single-cuff, one double-cuff, or two single-cuff attachments for
unidentified \(S_{0,4}\); and one single-cuff attachment for the
identified \(S_{0,4}\) completion. At higher fixed complexity, it suffices
to allow each attached mark to use one, two, or three outer cuffs, with
disjoint cuff sets for distinct marks. Patterns that cannot be realized
contribute zero. Let \(I_m^{\rm bd}\) count ordered
retained marked tuples with one distinguished boundary-only witness.
The shared-cuff calculation and the complete complementary genus sum
give
\begin{equation}\label{w:all-boundary-only}
 \mathbb E_g I_m^{\rm bd}
 =O_{m,S}\!\left(\frac{\log H}{H^{12}}\right).
\end{equation}
The fixed higher-complexity contribution is
\[
 \mathbb E_g I^{\rm bd}_{3\leq c\leq M,m}
 =O_{M,m,S}\!\left(g^{-3/4}H^{-9/4}\log H\right).
\]
Take the fixed band up to \(K_m-1\), where \(K_m=2m+9\), and apply
Proposition~\ref{prop:witness-pair-tail} with \(r=m-1\) to the
remaining tail. The constants in the fixed-complexity bounds are
therefore used only over a finite range.
The class is empty at \(m=1\).

\begin{proposition}[Exactly one crossed mark]\label{prop:single-crossing}
For every fixed marked order \(m\) and \(S\), let \(I_{1,m}(X)\)
count ordered retained closed-disjoint marked tuples and a
distinguished witness whose support crosses exactly one marked
support.  Then
\begin{equation}\label{w:single-crossing-bound}
 \mathbb E_g I_{1,m}=O_{m,S}(H^{-2}).
\end{equation}
This raw estimate includes all marked-complement types and arbitrary
slot windows.
\end{proposition}

\begin{proof}
Take the unique crossed mark as the anchor. Every other marked support
is closed-disjoint from it. If the interior of such a support met the
completed pair image, Proposition~\ref{l:localization} would force the
witness to cross this second mark. That proposition follows from the
intrinsic barrier in Lemma~\ref{lem:barrier}, independently of the
incidence estimates here. Thus all the other marked interiors avoid
the completed pair image.

Such a configuration is exterior if all the other closed supports
avoid that image; otherwise it is boundary-only.  These alternatives
are exhaustive, including all boundary identifications.  Sum
\eqref{w:all-exterior}, with \(r=m-1\), and
\eqref{w:all-boundary-only} over the finitely many possible ordered
anchor positions. This proves \eqref{w:single-crossing-bound}.
\end{proof}

The equal-support, pure-touching, and exactly-one-crossing classes
are now controlled. Sections~\ref{sec:localization} and~\ref{sec:higher}
treat witnesses crossing two or more marks; Theorem~\ref{thm:weighted-loss}
combines the class estimates with the original denominator \(d(\tau)\).
All conclusions above continue to hold after imposing marked-slot windows.  A bounded
real or complex joint weight in the marked lengths is handled by
its supremum norm after taking absolute values.

% END sections/04_witnesses.tex

% BEGIN sections/05_localization.tex
\section{Localization and two crossed supports}\label{sec:localization}

Suppose a witness crosses two marked supports. To count the resulting
triples, we first locate the other marked supports relative to the core
filled by the witness and the two crossed curves. The interior of a
further marked support overlaps the core interior precisely when the
witness crosses that support. For exactly two crossings and connected
joint complement, the core has complexity at least three. We estimate
the complexity-three case on a five-holed sphere, use filling-tuple
estimates for higher fixed complexities, and apply
Lemma~\ref{lem:uniform-tail} to the remaining range.

Throughout this section the marked curves are retained primitive unoriented
figure-eight geodesics.  Their ambient support closures are embedded pants
with three distinct cuffs and are pairwise disjoint.  The witness is also
retained.  A witness \emph{crosses} a marked support when their distinct
supports have overlapping interiors.  All curve lengths are at most
$L\geq1$, and all support cuffs have length at least $b\geq0$.
The local estimates keep $L$ and $b$ free; in the asymptotic conclusions
we take $H=\log g$ and $L=H-\log H+S$ for fixed $S$, and, when specified,
$b=b_L=10\log L$, as in \eqref{f:critical-parameters}.
Individual connected-complement conditions may be discarded in upper
counting kernels; joint connectedness is imposed only where stated.

We distinguish an intrinsic filled core $U\cong S_{h,n}$ from its completed
image $Z\cong S_{h+j,\nu}$, where $\nu=n-2j$.  The map
$\pi:U\longrightarrow Z\subset X$ embeds the interior and identifies $j$
disjoint pairs of boundary copies.  Set
\begin{equation}\label{l:core-notation}
 C=-\chi(U)=2h-2+n,\qquad
 \Omega=\pi(U^\circ).
\end{equation}
The image $\Omega$ is connected and dense in $Z$. Its complement
consists of the outer cuffs of $Z$ and the identified seams.
We test filling on $U$, where the complementary regions of the curve
family are disks and peripheral annuli, each annulus adjacent to one
intrinsic cuff. When $j>0$, gluing cuff pairs places seams in $Z^\circ$,
so filling $U$ need not imply filling $Z$.

\subsection{The intrinsic barrier and localization}

\begin{lemma}[Intrinsic geodesic barrier]\label{lem:barrier}
Let a finite geodesic family $\Gamma\subset\Omega$ fill a connected intrinsic
geodesic-boundary surface $U$ as above.  If an ambient simple closed geodesic
$\delta$ is disjoint from $\Gamma$, then $\delta\cap\Omega=\varnothing$.
The assertion also applies to a single figure-eight in its supporting pants.
\end{lemma}

\begin{proof}
If $\delta$ is an intrinsic boundary image, it misses $\Omega$.  Otherwise
its intersections with the finitely many boundary images are transverse and
finite: a tangency of two geodesics would force coincidence.  Cut $\delta$
at those intersections and lift the pieces entering $\Omega$ to $U$.
Such a piece is either a simple closed geodesic in $U^\circ$ or a properly
embedded geodesic arc with endpoints on $\partial U$.

A closed piece avoiding $\Gamma$ lies in a disk or a peripheral annulus.
A disk contains no closed geodesic.  In a peripheral annulus the piece is
freely homotopic to the cuff, so uniqueness of the geodesic representative
forces it to equal that cuff, contrary to its lying in the interior.

Under the filling convention, only peripheral annuli meet $\partial U$.
A proper piece therefore lies in a peripheral annulus
and has both endpoints on that annulus's one outer cuff.  A simple proper
arc in an annulus with endpoints on the same boundary component is
boundary-parallel.  It and a segment of the cuff bound a disk with two
geodesic sides, an impossible geodesic bigon.  This argument is made after
cutting the intrinsic surface.  At an identified seam the two boundary
copies belong to their respective peripheral annuli, so identification
introduces no additional case.
\end{proof}

The canonical core of a connected actual curve graph is obtained by taking
a small regular neighborhood, filling its contractible boundary components,
and straightening the remaining essential boundary classes.  Boundary
copies remain separate when they straighten to the same ambient geodesic.
The actual curves lie in the intrinsic interior and fill it by disks and
peripheral annuli.  Write $\Gamma$ for the graph and $A$ for the union of
its open supporting pants. By Lemma~\ref{lem:barrier} applied to each
figure-eight, a simple closed geodesic misses $\Gamma$ if and only if it
misses $A$. Thus every connected geodesic-boundary interior containing
$\Gamma$ also contains $A$, since its boundary images miss $A$ and each open
pant is connected and contains its figure-eight. The converse follows
from $\Gamma\subset A$.

Filling and straightening gives the smallest such interior in either
construction. To see this, first cut along the boundary images of any
containing interior. A disk to be filled cannot contain an essential cut
geodesic, and each essential boundary class straightens in the cut
surface, possibly to a cuff. Thus the curve graph $\Gamma$ and the union
$A$ of its open supporting pants have the same smallest containing
interior with geodesic boundary. In the
intrinsic metric completion each incident side of a boundary image gives
a cuff, so the two copies of an identified seam remain distinct. This is
the convention of \cite[Section~4.3]{HSWX}.

\begin{lemma}[Actual entry and support containment]\label{l:entry-containment}
If two distinct supporting pants have overlapping interiors, each actual
figure-eight enters the interior of the other support and the two actual
curves intersect.  Consequently, if a witness $\beta$ crosses each of the
pairwise closed-disjoint supports of $\alpha_1,\ldots,\alpha_k$, their actual
curve graph is connected.  Every one of their closed supporting pants is
contained in the completed canonical core of this graph.
\end{lemma}

\begin{proof}
Suppose $P_\alpha^\circ\cap P_\beta^\circ\ne\varnothing$ but $\beta$
avoids $P_\alpha^\circ$.  It then avoids $\partial P_\alpha$: a transverse
intersection forces entry, and tangency would force coincidence with a
simple cuff.  Apply Lemma~\ref{lem:barrier} to the single curve $\beta$
filling $P_\beta$.  The boundary of $P_\alpha$ misses $P_\beta^\circ$.
Connectedness and the interior overlap then give
$P_\beta^\circ\subset P_\alpha^\circ$, contradicting the assumed avoidance
by $\beta$.

If $\beta$ entered $P_\alpha^\circ$ but avoided $\alpha$, every excursion
would lie in a peripheral annulus of $P_\alpha\setminus\alpha$.  Such an
excursion need not be simple.  Lift the annulus and an excursion to
$\mathbb H^2$.  Its outer cuff generates the annulus's fundamental group,
so both endpoints of the lifted excursion lie on the same lift of that
geodesic cuff.  The unique geodesic segment between them lies on that line,
which contradicts genuine interior entry.  A closed geodesic wholly in the
annulus would instead be a power of the simple cuff, also impossible for a
figure-eight.  Hence $\alpha$ and $\beta$ intersect, proving graph
connectedness.

Every image of an intrinsic core boundary avoids the actual curves.
The single-curve case of Lemma~\ref{lem:barrier} makes it disjoint from
each support interior.  That connected support interior contains its
actual curve in $\Omega$, and so lies entirely in $\Omega$.  Taking
closures proves the last assertion, including at identified seams.
\end{proof}

\begin{proposition}[Localization]
Let\label{l:localization} the actual curves $\alpha_1,\ldots,\alpha_k,\beta$ fill the intrinsic
core $U$, and suppose the witness crosses every selected marked support.
If an embedded geodesic pair of pants $Q$ is closed-disjoint from every selected
marked support, then
\begin{equation}\label{l:localization-equivalence}
 Q^\circ\cap Z^\circ\ne\varnothing
 \quad\Longleftrightarrow\quad
 \beta\cap Q^\circ\ne\varnothing
 \quad\Longleftrightarrow\quad
 P_\beta^\circ\cap Q^\circ\ne\varnothing.
\end{equation}
If these conditions fail, $Q^\circ\cap Z=\varnothing$.
There is no restriction on core complexity, boundary identification, or
ambient complementary topology.
\end{proposition}

\begin{proof}
Suppose $\beta$ avoids $Q^\circ$.  It also avoids $\partial Q$, by
transversality and geodesic uniqueness.  The selected actual curves avoid
all of $Q$.  Lemma~\ref{lem:barrier} therefore gives
$\partial Q\cap\Omega=\varnothing$.  Connectedness puts $\Omega$ in one
component of $X\setminus\partial Q$.  This component is outside $Q$, since
it contains a selected actual curve.  Openness of $Q^\circ$ and density of
$\Omega$ in $Z$ give $Q^\circ\cap Z=\varnothing$.  Conversely,
$\beta\subset\Omega\subset Z^\circ$ proves that actual entry implies
interior contact.  The last equivalence is the single-support argument in
the proof of Lemma~\ref{l:entry-containment}, applied to $P_\beta$ and $Q$.
\end{proof}

We now apply localization to the witness decomposition.  For a marked
tuple $\tau=(\alpha_1,\ldots,\alpha_m)$, keep the original witness set
$\mathcal W(\tau)$ and denominator
\begin{equation}\label{l:original-denominator}
 d(\tau)=|\mathcal W(\tau)|.
\end{equation}
For a witness define the crossed slots and their number by
\begin{equation}\label{l:crossing-signature}
 \mathcal C(\tau,\beta)
 =\{i:P_\beta\ne P_{\alpha_i},\ 
       P_\beta^\circ\cap P_{\alpha_i}^\circ\ne\varnothing\},
 \qquad \kappa(\tau,\beta)=|\mathcal C(\tau,\beta)|.
\end{equation}
Let $Z_i$ be the completed core of $(\alpha_i,\beta)$ for a crossed slot
$i$.  Proposition~\ref{l:localization} gives the pointwise identity
\begin{equation}\label{l:slot-localization}
 \{j\ne i:P_{\alpha_j}^\circ\cap Z_i^\circ\ne\varnothing\}
 =\mathcal C(\tau,\beta)\setminus\{i\}.
\end{equation}
In particular, if the witness crosses exactly one slot, every other marked
support has interior disjoint from its pair core.  Either all these supports
are closed-disjoint from the core, or at least one attaches along an outer
cuff.  These are the exterior and attachment families considered earlier.
An identified seam is not a further attachment site, since the completed
core occupies both of its sides.

\subsection{Triple cores of least complexity}

Select two crossed marks, and write
$U=U(\alpha_1,\alpha_2,\beta)$.  Lemma~\ref{l:entry-containment} proves
that the actual triple graph is connected and that its completed core
contains all three closed supports.  The two marked pants have total area
$4\pi$.  The connected witness pants meets both and has a positive-area
open part outside their disconnected union.  Hence
\begin{equation}\label{l:triple-minimum-perimeter}
 C\ge3,\qquad p(U)\le6L.
\end{equation}
The complexity bound follows from Gauss--Bonnet:
$\operatorname{area}(U)=2\pi C>4\pi$ and $C$ is an integer.
For the perimeter inequality, each supporting pants has perimeter at most
$2L$ by the figure-eight length formula.  The boundary of the union of
the three pants consists of subarcs of their boundaries, with total length
at most the sum of the three perimeters.  Filling contractible boundary
components discards boundary length, and replacing each remaining boundary
walk by its geodesic representative does not increase length.  This counts
both intrinsic copies whenever two representatives coincide.  The same
argument gives $p(U(\alpha,\beta))\le4L$ for a pair core.

\begin{lemma}[Connected complexity-three cores]\label{l:planar-core}
Suppose $C=3$ and the ambient complement of the two selected closed marked
pants is connected.  Then $U=Z\cong S_{0,5}$ and $j=0$.  Each marked pair of pants
has two outer cuffs of $U$ and one internal cuff.  The residual component is
a pair of pants joining the two internal cuffs to the fifth outer cuff.
Each anchor--witness pair core is an embedded $S_{0,4}$ or all of $U$.
\end{lemma}

\begin{proof}
Remove the two closed marked pants from $Z$ and take the closure of the
remaining positive-area surface.  It has area $6\pi-4\pi=2\pi$.
Its boundaries are distinct ambient geodesics: disk and annulus components,
or zero-area regions between duplicate boundaries, are impossible.
There is therefore exactly one stable component $B$, with
$\chi(B)=-1$.  Since $Z$ is connected and the marked pants are
closed-disjoint, $B$ attaches to each mark along at least one cuff.  Thus
$B$ is a pair of pants, rather than a one-holed torus.  If $k_i\ge1$ cuffs attach
to the $i$th mark and $e$ remain exposed, then $k_1+k_2+e=3$.

If $e=0$, all of $\partial B$ belongs to the two marks, so $B^\circ$ is
an isolated component of $X$ minus their closed supports.  The remaining three
marked cuffs are distinct outer cuffs of $Z$, and an exterior collar gives
another complementary component.  This contradicts connectedness.  Hence
$e=1$ and $k_1=k_2=1$.  The three pants form a tree and give $S_{0,5}$.
Its completed genus is $h+j=0$, whence $h=j=0$.

Every outer cuff of $U$ avoids the actual anchor--witness pair.
Lemma~\ref{lem:barrier} makes it disjoint from the connected intrinsic
pair interior.  That interior contains its actual curves inside $U$, so
cannot leave $U$.  Its completed image is therefore planar, has area at
most $6\pi$, and has complexity at least two because its distinct
overlapping pants have union area greater than $2\pi$.  The only
possibilities are embedded $S_{0,4}$ and $S_{0,5}$, without identified
boundaries.
\end{proof}

The hypothesis holds for any selected pair from a marked tuple with
connected joint complement.  Restore the other closed-disjoint marked
pants to the complement.  Each restored pair of pants has exterior boundary
collars in the original connected complement, so the resulting complement
remains connected.

\subsection{The planar local mass}

At complexity three the general filling-triple bound is too weak. We
exploit the five-holed sphere in two ways. We count witness pants
having two outer cuffs, and we count triples in which the witness
and one anchor already fill the whole core. These two families cover
the triples of interest and can be bounded separately. Write
\begin{equation}\label{l:profiles}
 \psi(t)=\frac{\sinh(t/2)}{t/2},\qquad
 \rho(t)=\frac{t}{e^t-1},\qquad
 \psi(0)=\rho(0)=1,\qquad \mathcal A_g(L)=\frac{Le^L}{g}.
\end{equation}
For a five-holed sphere with labeled outer lengths $\mathbf t$, let
$N^\star_{L,b}(U)$ count ordered actual triples as above, filling $U$, with
the two marked closed pants disjoint and the witness crossing both.
Define
\begin{equation}\label{l:planar-kernel}
 K^\star_{L,b}(\mathbf t)
 =\int_{\mathcal M_{0,5}(\mathbf t)}N^\star_{L,b}(U)\,dV_{\rm WP},
 \qquad
 \mathcal J^\star_{05}(L,b)
 =\int_{\mathbb R_+^5}K^\star_{L,b}(\mathbf t)
       \prod_{a=1}^5t_a\psi(t_a)\,d\mathbf t.
\end{equation}
The count imposes the length and cuff cutoffs but no condition on the
ambient complement of an individual support.  It retains every internal
isotopy and twist configuration and both possible anchor--witness pair
types.

\begin{lemma}[Planar triple mass]\label{lem:planar-mass}
For $L\ge1$, $b\ge0$, and fixed $0<\delta<1/4$,
\begin{equation}\label{l:planar-mass}
 \mathcal J^\star_{05}(L,b)
 \le C_\delta\left[
 (1+L)^{11}e^{(5/2+2\delta)L}
 +(1+L)^9\bigl(e^{5L/2}+e^{3L-b}\bigr)\right].
\end{equation}
In particular,
\begin{equation}\label{l:planar-critical-mass}
 \mathcal J^\star_{05}(L,10\log L)\le C\frac{e^{3L}}L
 \qquad (L\ \hbox{sufficiently large}).
\end{equation}
\end{lemma}

\begin{proof}
We derive the two counting bounds here and integrate them in
Appendix~\ref{app:planar}.
An embedded pair of pants in $S_{0,5}$ has either two outer cuffs and one internal
cuff (an end pants), or one outer cuff and two internal cuffs (a central
pants).  Three disjoint nonparallel internal cuffs cannot occur in
$S_{0,5}$, and three outer cuffs would give no attachment to the rest of
the connected surface.  A central pair of pants has two complementary end pants.
Consequently two closed-disjoint marked pants must both be end pants:
if one were central, the other would lie in a complementary end pants.
A geodesic pair of pants contained in a geodesic pair of pants must equal it, and would
therefore share a cuff with the central pants, a contradiction.

If a central witness had both pair cores $Y_1,Y_2\cong S_{0,4}$, write its
outer cuff as $t$ and its internal cuffs as $u,v$.  Each $Y_i$ has three
outer cuffs of $U$ and one internal boundary $d_i$.  Both $u,v$ lie in
$Y_i$, but a four-holed sphere has at most one disjoint essential internal
curve.  At least one of $u,v$ is therefore peripheral in $Y_i$ and equals
$d_i$ by geodesic uniqueness.  If $d_1=d_2$, the side containing the
witness determines $Y_1=Y_2$.  This area-$4\pi$ surface could not contain
the two disjoint marked pants together with a connected witness having an
open part outside their union.  Thus, after renaming, $d_1=v$ and $d_2=u$.

Let $x,y$ be the internal cuffs of the two marked end pants.  The curve
$x$ avoids $v$ because its pants lies in $Y_1$; it also avoids $t$.
It must cross $u$, since distinct geodesic pants with noncrossing
boundaries cannot have overlapping interiors: cutting along their
boundaries would force nesting, and a geodesic pair of pants has no nonperipheral
simple closed geodesic.  Every component of $x$ in the witness pants is
therefore an essential $u$-to-$u$ wave.  Similarly $y$ supplies an essential
$v$-to-$v$ wave.  Such waves intersect.  Indeed, cutting along the first
wave produces two annuli, one containing $v$ and the other $t$; a disjoint
wave based at $v$ would lie in the former annulus and be boundary-parallel
there, hence also in the original pants.  This contradicts essentialness
and the closed-disjointness of the marked supports.  Consequently every
triple has either an end witness or an anchor--witness pair filling all
of $U$.  The two upper families may overlap.

For the end-pants count, fix two outer cuffs of positive lengths $s,t$.
The general McShane--Mirzakhani identity and its nonnegative summands give
\begin{equation}\label{l:mcshane-end}
 \#\{z\le L:s,t,z\text{ bound an end pants}\}
 \le\frac{s}{\mathcal R(s,t,L)}
 \le100(1+s)\left(1+e^{(L-s-t)/2}\right).
\end{equation}
The count is over simple internal geodesics $z$, with $z$ also denoting
their lengths.  Here $\mathcal R$ decreases in its third argument;
see \cite[Theorem~2.11 and Lemma~2.12]{HSWX}.  The identity applies to
general bordered surfaces, including $S_{0,5}$.  Zero boundary lengths
in the integrals are understood by limits.  A retained end figure-eight requires
$b\le s,t\le L$ and $z\le L$.  There are ten choices of the outer pair
and three actual unoriented figure-eights per pants.  Thus their total
number on any metric satisfies
\begin{equation}\label{l:end-count}
 E_{L,b}(U)\le C(1+L)(1+e^{L/2-b})
             \le C(1+L)e^{L/2}.
\end{equation}
This counts all simple internal geodesics appearing in the identity.

For the whole-core-pair family, the ordered filling-pair theorem
\cite[Definition~17 and Theorem~18]{WXPGT} gives
\begin{equation}\label{l:whole-pair-count}
 N_2^{\rm fill}(U,2L)
 \le C_\delta(1+L)e^{2L-(1-\delta)p(U)/2}.
\end{equation}
This estimate is uniform in the metric and the boundary lengths.
The remaining mark is an end figure-eight, so \eqref{l:end-count} bounds
its choices after the disjointness and crossing restrictions are discarded.  Also $p(U)\le4L$ in this family.  The resulting integral
is bounded in \eqref{pl:pair-mass}.

For the end-witness family write the marked decomposition as
\begin{equation}\label{l:marked-planar-decomposition}
 A=P(a_1,a_2,u),\qquad M=P(u,v,e),\qquad B=P(c_1,c_2,v).
\end{equation}
The figure-eight identity gives
\begin{equation}\label{l:marked-planar-lengths}
 a_1,a_2,c_1,c_2,u,v\le L,\qquad
 a_1+a_2+u\le2L,\qquad c_1+c_2+v\le2L.
\end{equation}
To control the fifth cuff $e$, follow the witness from one marked
interior to the other and then back. Between its last departure and
its first arrival on each passage, the witness travels through $M$
from $u$ to $v$ or from $v$ to $u$. The two subarcs occupy disjoint
intervals of the parameter circle, so their lengths may be added.
If $d_M(u,v)$ is the shortest such joining distance, then
$2d_M(u,v)\le\ell(\beta)\le L$.  Following the two cuffs and a shortest
joining arc in both directions, with the pants boundary orientations,
gives a loop in the class of $e$.  Geodesic minimization yields
\begin{equation}\label{l:fifth-cuff}
 e\le u+v+2d_M(u,v)\le u+v+L.
\end{equation}
The subarcs may intersect on the surface. The estimate uses only
their lengths along the witness. Unfold the two internal marked cuffs, retain
\eqref{l:marked-planar-lengths}--\eqref{l:fifth-cuff}, and use
\eqref{l:end-count} for all witness choices on each metric.
Equation~\eqref{pl:end-mass} gives the second term of
\eqref{l:planar-mass}.  The two upper families prove that equation.
Finally take $\delta=1/8$ and $b=10\log L$.  The term
$(1+L)^9e^{3L-b}$ is $O(e^{3L}/L)$, and the other terms are exponentially
smaller, proving \eqref{l:planar-critical-mass}.
\end{proof}

\subsection{Noncrossed marks at the planar core}

\begin{lemma}[Planar contacts and connectedness]\label{l:planar-contacts}
Suppose the intrinsic and completed triple core coincide as an embedded
$S_{0,5}$.  A remaining mark not crossed by the witness is either
closed-disjoint from the whole core or attached along its fifth outer
cuff $e$.  At most one mark attaches.  Such an attachment forces the
joint marked complement to be disconnected.  Therefore, when the joint
marked complement is connected, all noncrossed remaining marks are
closed-disjoint from the entire core.  Under this connected-complement
hypothesis, if exactly the two selected marks are crossed, every remaining
mark is exterior.
\end{lemma}

\begin{proof}
Proposition~\ref{l:localization} makes a noncrossed mark's interior
disjoint from the core.  Any contact of geodesic boundaries must be a whole
common outer cuff: transverse intersection would force interior overlap,
and tangency forces coincidence.  The four cuffs $a_1,a_2,c_1,c_2$ already belong
to the selected marked pants, so marked closed-disjointness excludes them.
Only $e$ is available.  Two attached marked pants would share $e$,
again violating closed-disjointness.

An attachment $Q=P(e,x,y)$ gives an enlarged core
$R=U\cup_e Q\cong S_{0,6}$.  Every cuff of the middle pants $M$ now
belongs to one of $A,B,Q$.  Any other closed marked support therefore
avoids $\partial M$.  Its connected interior, if it met $M^\circ$,
would lie in $M^\circ$.  A geodesic pair of pants contains no proper
geodesic pair of pants: every essential simple closed curve in it is
boundary-parallel, and its geodesic representative is the corresponding
cuff.
Equality with $M$ would also violate closed-disjointness from $A,B,Q$.
Thus no other mark, crossed or noncrossed, enters $M^\circ$, which is
an isolated component of $X$ minus the union of the closed marked supports.  An exterior collar
along an outer cuff belonging to $A$ lies outside all marks after
shrinking the collar if necessary: the other finitely many compact
supports have positive separation from $A$.
This gives a second complementary component and proves the obstruction.
\end{proof}

By Lemma~\ref{l:planar-contacts}, every remaining noncrossed mark is
exterior when the joint marked complement is connected. Let $I^\star_{05;r}$ count a
planar triple and an ordered list of $r$ further marked figure-eights
whose closed supports are disjoint from the core and from one another.
Summing over all stable complementary genera as in Appendix~\ref{app:planar} gives
\begin{equation}\label{l:planar-exterior-incidence}
 \mathbb E_g I^\star_{05;r}
 \le C_r\mathcal A_g(L)^r g^{-3}\mathcal J^\star_{05}(L,b).
\end{equation}
At the critical cutoff and retained threshold this is $O_{r,S}(H^{-4})$.

\subsection{Attachments at higher fixed triple complexities}

For a fixed intrinsic type $(h,n)$ and identification number $j$, define
$K^{(3),h,n,j}_{L,b}$ and $\mathcal J^{(3),h,n,j}_{L,b}$ as in
\eqref{l:planar-kernel}, replacing $\mathcal M_{0,5}$ by
$\mathcal M_{h+j,\nu}$ and using the $\nu=n-2j$ outer profiles.  The
integrand defining $K^{(3),h,n,j}_{L,b}$ counts ordered actual triples on
each completed metric with the prescribed intrinsic type and
identification number.  Each triple
determines its seam system; all such configurations remain in the count.  We only need $\nu\ge1$: at fixed $C$ a closed completed
core would equal $X$ and force $2g-2=C$.

The actual triples satisfy the hypotheses of
Lemma~\ref{lem:fixed-core-kernel} with \(q=3\).  Their curves fill
the cut intrinsic surface and determine its seam system.  Each remains
primitive after cutting, since an intrinsic root would give an ambient
root, and a figure-eight geodesic is nonperipheral.  The total length
is at most \(3L\), and \eqref{l:triple-minimum-perimeter} gives
\(p(U)\leq6L\), counting both copies of every seam.  Therefore
\eqref{att:tuple-kernel-bound} gives
\begin{equation}\label{l:triple-kernel-bound}
 \mathcal J^{(3),h,n,j}_{L,b}
 \le C_{C,\delta}(1+L)^{3C+2-2j}e^{(3+3\delta)L}.
\end{equation}

Now require every other marked support to be noncrossed.  By
Proposition~\ref{l:localization} it is exterior or attaches along a
nonempty subset of the $\nu$ outer cuffs.  Identified internal seams are
not attachment sites, because $Z$ occupies both sides locally.  Let
$p_s$ marks attach along exactly $s$ cuffs, $s=1,2,3$, and let $r$ be
exterior.  Set
\begin{equation}\label{l:attachment-counts}
 p=p_1+p_2+p_3,\quad k=p_1+2p_2+3p_3,\quad
 d_*=p_2+p_3,\quad m=2+p+r.
\end{equation}
The shared subsets are pairwise disjoint, so $p\le k\le\nu\le C+2$.
At fixed $C,m$ there are finitely many assignments of outer cuffs and
marked slots.  For each assignment, the local kernel still counts all
internal curve configurations.
The enlarged connected core $R$ has
\begin{equation}\label{l:enlarged-type}
 \chi(R)=-C-p,\quad n_R=\nu+3p-2k,\quad
 h_R=h+j+k-p.
\end{equation}
All exposed cuffs are distinct ambient geodesics.  Equality of a free
attached cuff with a core cuff would be another recorded contact;
equality with another marked cuff contradicts closed-disjointness.
If $n_R=0$, then $R=X$ and $2g-2=C+p$, which is absent for sufficiently
large $g$ at fixed $C,m$.

For example, at $C=4$ the complete list of intrinsic types and
identification numbers with positive completed boundary is
\begin{equation}\label{l:complexity-four-types}
 (h,n,j)\in
 \{(0,6,0),(0,6,1),(0,6,2),(1,4,0),(1,4,1),(2,2,0)\}.
\end{equation}
The completed types belong to $\{S_{0,6},S_{1,4},S_{2,2}\}$.
In the order listed in \eqref{l:complexity-four-types}, the local powers
$3C+2-2j$ are $14,12,10,14,12,14$.
The remaining formal types have closed completed genus three and
cannot occur for $g>3$.  At marked order three, the additional mark
may be exterior or attached along $s=1,2,3$ cuffs when $\nu\ge4$,
and along $s=1,2$ when $\nu=2$.  All these enlarged cores have positive
outer boundary.  This is a covering catalog; realizability with
connected marked complement is not asserted for every entry.

\begin{proposition}[Fixed triple bands]\label{prop:triple-fixed}
For fixed $C\ge4,m$, $0<\delta<1/2$, and every intrinsic type and contact
pattern just described, the raw ordered incidence satisfies, for all
sufficiently large $g$ and all $L\ge1,b\ge0$,
\begin{equation}\label{l:fixed-triple-incidence}
 \mathbb E_g I_{h,n,j;\mathbf p,r}(L,b)
 \le C_{C,m,\delta}\mathcal A_g(L)^{m-2}g^{-C}
 (1+L)^{3C+2-2j}e^{(3+3\delta)L}
 L^{-d_*}\rho(b)^k.
\end{equation}
All complementary topologies are included.  At the critical cutoff
$L=H-\log H+S$, for any fixed $M\ge4$, the raw count in which each
tuple--witness pair with exactly two crossed marks and $4\le C\le M$
is counted once satisfies
\begin{equation}\label{l:fixed-triple-band}
 \mathbb E_g I_{m,=2,[4,M]}
 =O_{m,M,S}\left(g^{-3/4}H^{43/4}\right).
\end{equation}
This last bound is uniform in $b\ge0$.
\end{proposition}

\begin{proof}
The completed core has type $S_{h+j,\nu}$, with $\nu=n-2j\geq1$
and complexity $C$.  Its local kernel counts all intrinsic configurations,
is independent of the outer gluing twists, and determines the selected
core.  Localization gives exactly the attached and exterior placements
in Lemma~\ref{lem:attachment-bounds}; their parameters satisfy
$p+r=m-2$ and $d_*=p_2+p_3$.  The marked closed supports are pairwise
disjoint, so distinct attached marks use disjoint sets of outer cuffs.
Internal seams have already been integrated in $\mathcal J^{(3),h,n,j}_{L,b}$.

Apply \eqref{att:joint-bound} with $c=C$.  For a complement with $q$
components this gives
\begin{equation}\label{l:triple-before-q-sum}
 \mathbb E_g I_{h,n,j;\mathbf p,r;q}
 \le C_{C,m}g^{-(q-1)}\mathcal A_g(L)^{m-2}g^{-C}
 \mathcal J^{(3),h,n,j}_{L,b}L^{-d_*}\rho(b)^k.
\end{equation}
The lemma includes the shared-cuff Jacobians, all cuff roles and all
stable complementary topologies, including genus zero and one-boundary
components.  There are finitely many cuff assignments and
$1\leq q\leq\nu+3(p+r)-2k$.  Summing them and using
\eqref{l:triple-kernel-bound} proves \eqref{l:fixed-triple-incidence}.

Take $\delta=1/12$ at the critical cutoff.  The resulting rate is
\begin{equation}\label{l:fixed-triple-critical}
 O_{C,m,S}\left(
 g^{-C+13/4}H^{3C-5/4-2j-d_*}\rho(b)^k\right).
\end{equation}
Relative to $C=4,j=d_*=0$, this contains the factor
$(H^3/g)^{C-4}H^{-2j-d_*}$.  Since $0\le\rho(b)\le1$, maximizing
constants over the finite set $4\le C\le M$ proves
\eqref{l:fixed-triple-band}.  Each actual exactly-two-crossed
configuration has a unique unordered crossed pair; ordering it by slot
labels and summing the $\binom m2$ choices costs only a fixed constant.
The argument makes no assertion uniform in growing $M$.
\end{proof}

\subsection{The bound for two crossed supports}

\begin{proposition}[Two crossed supports]\label{prop:two-crossing}
At every fixed marked order $m\ge2$, let $I_{m,=2}^{\rm conn}$ count each
ordered marked tuple and retained witness once, with connected joint
marked complement and exactly two crossed marked supports.  At the
critical cutoff and retained cuff threshold,
\begin{equation}\label{l:two-crossed-raw}
 \mathbb E_g I_{m,=2}^{\rm conn}=O_{m,S}(H^{-4}).
\end{equation}
Consequently the same bound holds for its expected normalized witness share with
the original denominator $d(\tau)$.  Without the connected-complement
restriction, its expected normalized share is
$O_{m,S}(H^{-4}+g^{-1})$.  Slot restrictions only reduce these
nonnegative counts, and bounded real or complex slot weights cost the
product of their absolute supremum norms.
\end{proposition}

\begin{proof}
The triple complexity is at least three by
\eqref{l:triple-minimum-perimeter}.  For $C=3$,
Lemma~\ref{l:planar-core} makes the core an embedded $S_{0,5}$, and
Lemma~\ref{l:planar-contacts} makes all $m-2$ remaining marks exterior.
Equation~\eqref{l:planar-exterior-incidence} and
Lemma~\ref{lem:planar-mass} bound the expectation of this raw contribution by
$O_{m,S}(H^{-4})$.

To separate fixed from growing complexity, choose the endpoints,
depending only on the fixed marked order $m$, as
\begin{equation}\label{l:triple-cutoffs}
 M_m=2m+11,\qquad K_m^{\rm tri}=M_m+1=2m+12.
\end{equation}
Proposition~\ref{prop:triple-fixed} handles $4\le C\le M_m$ with
$O_{m,S}(g^{-3/4}H^{43/4})=o(H^{-4})$.  For the entire remaining
range $C\ge K_m^{\rm tri}$ apply Lemma~\ref{lem:uniform-tail} with
three actual core curves and $r=m-2$ arbitrary additional primitive
geodesics.  Its enlarged count allows arbitrary placements and
repetitions of the latter, so contains the desired family.  It gives
\begin{equation}\label{l:triple-growing-tail}
 \mathbb E_g J_{\ge K_m^{\rm tri};m-2}^{(3)}
 \le C_{m,S}g^{15+2(m-2)-(2m+12)}H^{-(m-2)-1/2}
 =O_{m,S}\left(g^{-1}H^{-m+3/2}\right)=o(H^{-4}).
\end{equation}
Lemma~\ref{lem:uniform-tail} is proved independently of the present
two-crossed estimate.  It handles the whole growing range at once,
without summing fixed-complexity constants.  Combining the three
complexity ranges proves
\eqref{l:two-crossed-raw}.

For a nonempty witness class, division by the original integer
$d(\tau)\ge1$ can only reduce its raw incidence.  The disconnected
normalized share is bounded by the number of retained disjoint marked
tuples with disconnected complement, whose expected number is $O_{m,S}(g^{-1})$ by
Proposition~\ref{prop:disconnected}.  Finally take absolute values of
slot weights before applying the nonnegative upper bounds.
\end{proof}

This bounds the contribution from exactly two crossed supports to
Theorem~\ref{thm:weighted-loss}. The remaining witnesses cross at
least three supports. Their cores have a larger minimum complexity
when the joint marked complement is connected.

% END sections/05_localization.tex

% BEGIN sections/06_higher_cores.tex
\section{Three or more crossed supports and the deletion estimate}\label{sec:higher}

For a witness that crosses several marked supports, two features of its
core control the count. Its area and perimeter are bounded in terms
of the total length of the curves that fill it. This gives an estimate uniform over all
large core complexities, including those left over from
Sections~\ref{sec:witnesses} and~\ref{sec:localization}. At fixed
complexity, connectedness of the marked complement bounds how many
marked supports the witness can cross. We use the two estimates to
bound the witnesses with three or more crossings and then sum all
witness classes to prove Theorem~\ref{thm:weighted-loss}.

Throughout this section, $S\in\mathbb R$ is fixed and
\begin{equation}\label{h:critical-scale}
 H=\log g,\qquad L=H-\log H+S.
\end{equation}
All asymptotic assertions are for sufficiently large $g$ at this scale.
The fixed-core estimates are also stated for arbitrary $L\geq1$.
We allow any cuff threshold $b\geq0$ until the final deletion argument,
where we set $b=b_L=10\log L$.
All curves are actual primitive unoriented geodesics.  A good figure-eight
has an embedded supporting pair of pants with three distinct ambient
cuffs and connected individual complement, as in the preceding sections.
Local upper kernels may discard the individual-complement condition.
The ambient closures of marked supports are pairwise disjoint. The
joint marked complement means $X$ minus their union; the witness is
not removed.

\subsection{Intrinsic filling and a uniform geometric bound}

Let $\alpha_1,\ldots,\alpha_k$ be marked figure-eights and let a witness
$\beta$ cross each of their supports. Crossing means interior overlap
of distinct supports. By Lemma~\ref{l:entry-containment}, the witness
enters each marked interior and intersects each actual curve
$\alpha_i$. Consequently
\begin{equation}\label{h:connected-graph}
 \Gamma=\alpha_1\cup\cdots\cup\alpha_k\cup\beta
\end{equation}
is a connected finite geodesic graph.

Take the canonical intrinsic filled surface of this graph as in
Section~\ref{sec:localization}. Straightening the boundary can send
two intrinsic cuffs to the same ambient geodesic. With $j$ such pairs,
the core and its completed image are
\begin{equation}\label{h:intrinsic-completion}
 \pi:U\cong S_{h,n}\longrightarrow Z\cong S_{h+j,\nu},
 \qquad \nu=n-2j,\qquad C=2h-2+n=|\chi(Z)|.
\end{equation}
The map is an embedding on $U^\circ$, and its image $\Omega=\pi(U^\circ)$ is connected
and dense in the completed image $Z$. The actual graph lies in $\Omega$.
The complement of its lift in $U$ consists of disks and peripheral
annuli, each annulus adjacent to one intrinsic boundary cuff.  These are the
intrinsic filling conventions used in Lemma~\ref{lem:barrier}.

By Lemma~\ref{l:entry-containment}, $Z$ contains the ambient closures
of all the supporting pants. We estimate its area on $U$, where each
face of the graph is a disk or a peripheral annulus.

\begin{lemma}[Filling-face area and perimeter]\label{lem:filling-area}
Let $\Gamma\subset U^\circ$ be a finite geodesic graph filling a compact
intrinsic hyperbolic surface $U$ with geodesic boundary. Suppose that
every component of $U\setminus\Gamma$ is a disk or a peripheral annulus
containing one boundary cuff. Let $\operatorname{length}(\Gamma)$ be
the sum of its geometric edge lengths, each edge counted once. Then
\begin{equation}\label{h:face-bound}
 \operatorname{Area}(U)\leq2\operatorname{length}(\Gamma),
 \qquad p(U)\leq2\operatorname{length}(\Gamma).
\end{equation}
The assertion includes a closed $U$, when every face is a disk.
\end{lemma}

\begin{proof}
If $U$ has boundary, first double it to a closed hyperbolic surface.
A disk face lifts injectively into $\mathbb H^2$.  In upper-half-plane
coordinates the field $V=-y\partial_y$ has hyperbolic norm one and
divergence one.  The divergence theorem on a compact piecewise smooth
truncation $D$ of the lifted face gives
\begin{equation}\label{h:disk-flux}
 \operatorname{Area}(D)
 =\int_{\partial D}\langle V,n\rangle\,ds
 \leq\operatorname{length}(\partial D).
\end{equation}
Exhausting the face bounds its area by the length of its graph-side
boundary walk, with multiplicity.

A peripheral annulus lifts injectively into the cyclic cover associated
with its cuff.  This cover is a hyperbolic cylinder with Fermi metric
\[
 ds^2=dr^2+\cosh^2r\,dt^2,\qquad t\in\mathbb R/\ell\mathbb Z.
\]
The vector field $V=\tanh r\,\partial_r$ satisfies
\[
 \|V\|\leq1,\qquad
 \operatorname{div}V
 =\frac1{\cosh r}\frac{\partial}{\partial r}
       (\cosh r\tanh r)=1.
\]
It vanishes on the axis $r=0$, representing the geodesic cuff.  Its flux
through that cuff is thus zero, and the divergence theorem bounds the
annular area by its graph-side boundary length alone.

The asserted lifts are injective: their projections are the original
embedded faces, whose fundamental groups are respectively trivial or
generated by the cuff.  To handle repeated boundary vertices or edges,
remove a sufficiently small neighborhood of the finite graph, apply the
preceding argument to the resulting piecewise smooth domains, and let
the width tend to zero.  The graph-side portions of the truncated boundaries converge
in length to the corresponding boundary walks, counted with multiplicity.  Each edge of $\Gamma$ has two sides.  Summing over all
faces proves the area estimate.

Finally, the graph-side walk of a peripheral annulus is freely homotopic
to its cuff.  The geodesic representative minimizes length in this
class.  Summing this comparison over the peripheral annuli, whose
walks use at most the full two-sided graph perimeter, proves the
perimeter estimate.
\end{proof}

For $q$ distinct actual geodesics of lengths at most $L$, no two share
an open segment.  Thus their graph length is at most $qL$.
Lemma~\ref{lem:filling-area} and Gauss--Bonnet give
\begin{equation}\label{h:uniform-image-geometry}
 \operatorname{Area}(Z)=\operatorname{Area}(U)\leq2qL,
 \qquad p(Z)\leq p(U)\leq2qL,
 \qquad C\leq\frac{qL}{\pi}.
\end{equation}
Both copies of every identified cuff pair occur in $p(U)$.
Completion identifies boundary circles only, preserves area and Euler
characteristic, and removes their copies from the exposed perimeter.
The completed image is embedded with distinct outer cuffs.  If it had
no outer boundary, it would equal the connected closed surface $X$;
then $4\pi(g-1)\leq2qL$.  For fixed $q,S$ this is impossible for large
$g$.  Hence all such images are proper and bordered, uniformly over
the whole possible complexity range.

\subsection{A tail estimate uniform in the core complexity}

The constants in the fixed-type estimates can depend on the topology
of the core. To sum over growing complexities, we use the area and
perimeter bounds to place all completed images in one family of
subsurfaces, then count the possible curves on each image.

For a fixed integer $q\geq2$, consider ordered $q$-tuples of distinct
primitive unoriented geodesics of length at most $L$ whose union is
connected and has a canonical intrinsic filled core as above.  A
specified subfamily may impose goodness, supporting-pants, crossing,
or disjointness conditions.  Let $P_{\geq K}^{(q)}(X)$ count the tuples
in that subfamily whose intrinsic complexity is at least $K$.
In particular, this definition includes the anchor--witness families
used here, with $q-1$ anchors and one witness.  Adjoin an ordered list
of $r$ arbitrary primitive unoriented geodesics of length at most $L$,
allowing repetitions and all placements, and set
\begin{equation}\label{h:enlarged-tail-count}
 J_{\geq K;r}^{(q)}(X)
 =P_{\geq K}^{(q)}(X)\bigl(N_X^{\rm prim}(L)\bigr)^r.
\end{equation}
Here $N_X^{\rm prim}(L)$ counts all primitive unoriented closed
geodesics on $X$ of length at most $L$.

\begin{lemma}[Uniform growing-core tail]\label{lem:uniform-tail}
For fixed integers $q\geq2$, $K\geq1$, $r\geq0$, and fixed
$S\in\mathbb R$,
\begin{equation}\label{h:uniform-tail}
 \mathbb E_gJ_{\geq K;r}^{(q)}
 \leq A_{q,K,r,S}
 g^{5q+2r-K}H^{1-q/2-r}.
\end{equation}
The estimate is uniform in any additional cuff threshold $b\geq0$.
No connected-complement condition is required, and every intrinsic
boundary identification and complementary topology is included.
\end{lemma}

\begin{proof}
We use one bound for geodesics on a given surface and another for
subsurfaces of a random surface. By
\cite[Section~2.5 and Theorem~2.8]{HSWX}, on every closed genus-$G$
hyperbolic surface $D$ and for every $L>0$,
\begin{equation}\label{h:primitive-count}
 N_D^{\rm prim}(L)\leq(G-1)e^{L+7}.
\end{equation}
The convention is primitive unoriented geodesics.

For $T=a\log g$ with any fixed $a>0$, let $\operatorname{Sub}_T(X)$
be the family in \cite[Section~7.1, equation~(21) and the following definition]{WX}: connected
geodesic-boundary subsurfaces with perimeter at most $2T$ and area at
most $4T$, allowing two intrinsic boundary copies to represent one
ambient geodesic.  Proposition~33 of that paper, divided by $V_g$,
states that for every fixed integer $K\geq1$,
\begin{equation}\label{h:wx-subsurface-tail}
 \mathbb E_g\sum_{\substack{Y\in\operatorname{Sub}_T(X)\\
                            |\chi(Y)|\geq K}}
              e^{-p(Y)/4}
 \leq A_{K,a}\frac{T e^{7T/2}}{g^K}.
\end{equation}
The upper complexity cutoff in the source follows from the defining
area bound for $\operatorname{Sub}_T(X)$; its perimeter indicator is
also automatic. The estimate includes all complementary topologies.

Take $T=qH$. For sufficiently large $g$, we have $L\leq H$, so
\eqref{h:uniform-image-geometry} puts every actual completed image in
$\operatorname{Sub}_{qH}(X)$. We first count tuples on these embedded
images, using their exposed perimeter $p(Z)$. Internal seams do not
contribute to this perimeter. We will then enlarge the nonnegative
image sum to the full source family.

Suppose that $Z$ has genus $h'$ and $\nu\geq1$ outer cuffs.  Its
hyperbolic double has genus
\begin{equation}\label{h:double-genus}
 G=2h'+\nu-1=C+1.
\end{equation}
Each actual curve, lying in the interior of the first copy, remains a
geodesic of the same length.  It is primitive in $Z$, since a root in
$Z$ would map to an ambient root in $X$.  The folding map from the
double onto its first copy is a retraction.  Applying it to a proposed
root in the double gives a root in $Z$, which is impossible.  Likewise,
a conjugacy or inverse-conjugacy of two classes in the double retracts
to one in $Z$ and then maps to one in $X$.  Distinct actual unoriented
ambient curves therefore remain distinct in the double.

By \eqref{h:primitive-count}, the number of ordered actual $q$-tuples
with this completed image is at most
\begin{equation}\label{h:image-tuple-count}
 C^q e^{q(L+7)}\leq A_qL^qe^{qL},
 \qquad A_q=(q/\pi)^qe^{7q}.
\end{equation}
Every tuple determines its canonical image and its internal seams.
There is no further seam choice in the image sum: the right side of
\eqref{h:image-tuple-count} already counts all possible internal
configurations on that image.

For an actual image, $p(Z)\leq2qL$.  Inserting $e^{p(Z)/4}$ and its
reciprocal in \eqref{h:image-tuple-count}, and then enlarging the
weighted image sum, yields the pointwise inequality
\begin{equation}\label{h:weighted-image-sum}
 P_{\geq K}^{(q)}(X)
 \leq A_qL^qe^{3qL/2}
 \sum_{\substack{Y\in\operatorname{Sub}_{qH}(X)\\|\chi(Y)|\geq K}}
       e^{-p(Y)/4}.
\end{equation}
Only actual images enter the bound on the number of tuples per image.
Once each summand is a nonnegative perimeter weight, the sum may include
all subsurfaces in the indicated family. Taking expectation and applying
\eqref{h:wx-subsurface-tail} gives
\begin{align}\label{h:tail-arithmetic}
 \mathbb E_gP_{\geq K}^{(q)}
 &\leq A_{q,K}L^qe^{3qL/2}H e^{7qH/2}g^{-K}\notag\\
 &\leq A_{q,K,S}(L/H)^qg^{5q-K}H^{1-q/2}.
\end{align}
Here $e^{3qL/2}=e^{3qS/2}g^{3q/2}H^{-3q/2}$.
Finally, the ambient form of \eqref{h:primitive-count} implies
deterministically
\begin{equation}\label{h:extra-curves}
 \bigl(N_X^{\rm prim}(L)\bigr)^r
 \leq A_{r,S}g^{2r}H^{-r}.
\end{equation}
Multiply this inequality by $P_{\geq K}^{(q)}$ before taking
expectation.  Equations \eqref{h:enlarged-tail-count} and
\eqref{h:tail-arithmetic} prove \eqref{h:uniform-tail}; no independence
or factorization of expectations is used.
\end{proof}

For a witness crossing $k$ of the $m$ marks, the core is filled by
$q=k+1$ curves and there are $r=m-k$ remaining marks. Choosing $K$
so that $5q+2r-K=-1$ in \eqref{h:uniform-tail} gives the following
complexity cutoffs and tail bounds.
\begin{equation}\label{h:tail-cutoffs}
\begin{array}{c|c|c|c}
\text{crossed marks }k&q&r&
 K\ \text{and resulting bound for }\mathbb E_gJ_{\geq K;r}^{(q)}\\ \hline
1&2&m-1&2m+9:\quad O_{m,S}(g^{-1}H^{1-m})\\
2&3&m-2&2m+12:\quad O_{m,S}(g^{-1}H^{3/2-m})\\
3&4&m-3&2m+15:\quad O_{m,S}(g^{-1}H^{2-m})\\
4\leq k\leq m&k+1&m-k&
2m+3k+6:\quad O_{m,S}(g^{-1}H^{(k+1)/2-m}).
\end{array}
\end{equation}
Use the row with the actual number $k\leq m$ of crossed slots. Once
these slots are specified, their labels order the crossed curves, and
the remaining marks keep their relative order. This embeds the
tuple and its witness in the enlarged count. Summing over the
$\binom mk$ slot choices bounds the raw incidence in the tail at fixed
$m$, for any joint complement and any complexities of the individual
anchor--witness pairs.

\subsection{The connected-complement complexity gap}

Connected marked complement forces every residual component of the core
to have an exposed cuff. The resulting gluing graph gives the lower
complexity bound below and, at equality, excludes contacts with
noncrossed marks when the full marked complement is connected.

\begin{lemma}[Core topology and noncrossed marks]\label{lem:core-gap}
Let $k\geq2$ closed-disjoint marked pants be crossed by one witness,
and let $U\longrightarrow Z$ be the intrinsic core and completed image
of their $k+1$ actual curves.  Suppose $Z$ is proper and the complement
of these $k$ marked pants in $X$ is connected.  Then
\begin{equation}\label{h:core-gap}
 C=|\chi(U)|=|\chi(Z)|\geq2k-1.
\end{equation}
At equality, $U=Z\cong S_{0,2k+1}$ and there are no identified seams.

Without requiring connected selected-mark complement, any further
marked pants that is closed-disjoint from the selected
marks and is not crossed by the witness has interior disjoint from
$Z$.  Its contact with $Z$, if nonempty, consists of whole outer cuffs.
If the full marked complement is connected, the $k$ selected marks
are exactly the crossed marks, and $C=2k-1$, then every further
marked support is closed-disjoint from the entire core.
\end{lemma}

\begin{proof}
Remove the closed selected pants from $Z$ and let $B_1,\ldots,B_{q_0}$ be
the closures of the remaining positive-area components.  The witness
joins distinct closed-disjoint selected supports, so $q_0\geq1$.
All boundaries are geodesic.  Distinct ambient cuffs exclude
duplicate-boundary strips, while Gauss--Bonnet excludes positive-area
disks and annuli.  Thus each $B_a$ is stable.

Let $t$ be the sum of their genera, let $a$ be the total number of
cuffs attaching them to the selected pants, and let $e_a$ be the
number of outer core cuffs belonging to $B_a$.  Write $e=\sum_a e_a$
and $G=h+j$ for the completed genus.  Each $e_a\geq1$.  Otherwise
every boundary of $B_a$ would belong to a selected mark, making
$B_a^\circ$ an isolated component of $X$ minus the closed selected
supports. Properness of $Z$ supplies a nonempty exterior region in
that same complement, contradicting connectedness.

The gluing graph has $k$ selected-pants vertices, $q_0$ remaining
component vertices, and $a$ edges.  It is connected and bipartite:
selected pants do not meet, and distinct $B_a$ do not glue directly.
Its cycle rank is $a-k-q_0+1$, hence
\begin{equation}\label{h:gluing-genus}
 G=t+a-k-q_0+1.
\end{equation}
Additivity of Euler characteristic gives
\[
 C=k+2t+a+e-2q_0.
\]
Eliminating $a$ yields the exact identity
\begin{equation}\label{h:core-identity}
 C=2k-1+G+t+(e-q_0).
\end{equation}
All terms after $2k-1$ are nonnegative.  Equality forces $G=t=0$
and $e_a=1$ for every $a$.  The gluing graph is then a tree and
$a=k+q_0-1$.  Each planar stable $B_a$ has at least three cuffs,
only one of which is exposed; consequently $a\geq2q_0$ and
\begin{equation}\label{h:equality-component-range}
 1\leq q_0\leq k-1.
\end{equation}
The completed genus is zero and its complexity is $2k-1$, so it has
$2k+1$ outer boundaries.  Since $G=h+j=0$, both $h$ and $j$ vanish.

Now let $Q$ be a further noncrossed marked pants. It cannot have the
same support as $\beta$: that support crosses the selected marks,
whereas $Q$ is closed-disjoint from them.
Proposition~\ref{l:localization} therefore gives
\begin{equation}\label{h:noncrossed-exterior-interior}
 Q^\circ\cap Z=\varnothing.
\end{equation}
Two geodesic-boundary surfaces with disjoint interiors can meet only
along entire common cuffs: transverse boundary intersections would
force interior overlap, and tangent geodesics coincide.  An internal
identification seam cannot be a contact, since $Z$ occupies both sides
locally.  An anchor-owned outer cuff cannot be shared with $Q$, by
marked closed-disjointness.  This proves the asserted contact
classification.

Suppose finally that $C=2k-1$, the full marked complement is
connected, and the selected marks are exactly the crossed marks.
If a further mark $Q$ contacted $Z$, it would use
the unique exposed cuff of some $B_a$.  Every other boundary of
$B_a$ belongs to a selected mark.  All further noncrossed marks have
interior disjoint from $Z$ by \eqref{h:noncrossed-exterior-interior}.
It follows that $B_a^\circ$ is an isolated component of the full
marked complement.  There is another complementary region: the
number of outer core cuffs owned by anchors is
\begin{equation}\label{h:anchor-exposed-cuffs}
 3k-a=2k-q_0+1\geq k+2.
\end{equation}
A collar immediately outside $Z$ along any such cuff misses its
anchor on the exterior side and, after shrinking, misses every other
marked support.  The latter are finitely many compact sets
closed-disjoint from that anchor.  This collar is distinct from
$B_a^\circ$, contradicting connectedness.  Hence no contact occurs.
\end{proof}

If the full marked complement is connected, so is the complement of
any selected subfamily. Indeed, each omitted pants has exterior
boundary collars in the original connected complement, so restoring
it preserves connectedness. Thus the hypothesis in
Lemma~\ref{lem:core-gap} applies to the full crossed set and to any
selected subset of it.  Properness at the scale
\eqref{h:critical-scale} follows from
\eqref{h:uniform-image-geometry}.

For $k=3$, equality gives $S_{0,7}$.  The decomposition in the proof
has exactly two possible gluing patterns.  If $q_0=1$, the
remaining component is an $S_{0,4}$ attached once to each of the
three anchors and having one exposed cuff.  If $q_0=2$, the remaining
components are pants and the gluing tree is an alternating path
anchor--pants--anchor--pants--anchor.  Each remaining pants has one
exposed cuff. These two topological patterns are both included in
the kernels below. For each pattern the kernel counts all compatible
actual filling tuples; it makes no finite-orbit assumption on them.

\subsection{Fixed higher-core bands, including all attachments}

Fix $k\geq3$, an intrinsic type $(h,n)$, and an identification number
$j$, giving $U\cong S_{h,n}\longrightarrow Z\cong S_{h+j,\nu}$ with
$C=2h-2+n$ and $\nu=n-2j\geq1$. The cutoffs in this subsection may
be arbitrary $L\geq1$ and $b\geq0$ until we specialize to
\eqref{h:critical-scale}.
On a completed metric $Z\in\mathcal M_{h+j,\nu}(\mathbf x)$, let
$N_{L,b}^{(k+1),h,n,j}(Z)$ count ordered actual tuples
$(\alpha_1,\ldots,\alpha_k,\beta)$ with these canonical intrinsic
data. Every curve has length at most $L$ and an embedded support
with three distinct cuffs, each of length at least $b$.
The anchor closures are pairwise disjoint, the witness crosses every
anchor, and the completed image is all of $Z$. We omit individual
and selected-complement connectedness conditions in this upper
count. Define
\begin{align}\label{h:kernel-definition}
 K_{L,b}^{(k+1),h,n,j}(\mathbf x)
 &=\int_{\mathcal M_{h+j,\nu}(\mathbf x)}
        N_{L,b}^{(k+1),h,n,j}(Z)\,dV_{\rm WP},\notag\\
 \mathcal J_{L,b}^{(k+1),h,n,j}
 &=\int_{\mathbb R_+^\nu}
       K_{L,b}^{(k+1),h,n,j}(\mathbf x)
       \prod_{\ell=1}^{\nu}x_\ell\psi(x_\ell)\,d\mathbf x,
 \qquad \psi(x)=\frac{\sinh(x/2)}{x/2}.
\end{align}
We take $\psi(0)=1$.

\begin{lemma}[Fixed higher-core mass]\label{lem:higher-core-kernel}
For every fixed $C,k$, $0<\delta<1/2$, $L\geq1$, and $b\geq0$,
\begin{equation}\label{h:kernel-bound}
 \mathcal J_{L,b}^{(k+1),h,n,j}
 \leq A_{C,k,\delta}(1+L)^{3C+k-2j}
                     e^{(k+1)(1+\delta)L}.
\end{equation}
\end{lemma}

\begin{proof}
Apply Lemma~\ref{lem:fixed-core-kernel} with \(q=k+1\).
The actual family fills its cut intrinsic surface and determines the
identification seams.  Its total length is at most \((k+1)L\), and
\eqref{h:uniform-image-geometry} gives
\(p(U)\leq2(k+1)L\), including both intrinsic copies of every seam.
The curves remain primitive after cutting because a root would map to
an ambient root; each figure-eight is nonperipheral.  Thus all the
hypotheses of the common kernel bound hold.  Its polynomial power is
\(3C+(k+1)-1-2j=3C+k-2j\), and its exponential is
\(e^{(k+1)(1+\delta)L}\), proving \eqref{h:kernel-bound}.
\end{proof}

Let $p_s$ further marks attach along exactly $s$ outer core cuffs,
for $s=1,2,3$, and let $r$ further marks be closed-exterior to the
whole core.  Put
\begin{equation}\label{h:attachment-indices}
 p=p_1+p_2+p_3,\qquad
 a=p_1+2p_2+3p_3,\qquad d=p_2+p_3,\qquad m=k+p+r.
\end{equation}
These additional marks satisfy the same length and cuff bounds and
are required to be noncrossed.  Let
$I_{h,n,j;\mathbf p,r}^{(k+1)}(L,b)$ count the corresponding ordered
actual configurations, including all compatible slot, cuff, and
figure-eight-role assignments.  The count imposes no connectedness
on the joint marked complement or the enlarged-core complement.

\begin{proposition}[Fixed bands for higher crossing counts]
\label{prop:higher-fixed-bands}
For fixed $m,k,C$ with $3\leq k\leq m$ and $C\geq2k-1$, fixed
$0<\delta<1/2$, and sufficiently large $g$ depending on $C,m$,
\begin{align}\label{h:pattern-bound}
 \mathbb E_g I_{h,n,j;\mathbf p,r}^{(k+1)}(L,b)
 \leq {}&A_{C,m,\delta}\left(\frac{Le^L}{g}\right)^{m-k}g^{-C}
 (1+L)^{3C+k-2j}e^{(k+1)(1+\delta)L}\notag\\
 &\hspace{20mm}\cdot L^{-d}\rho(b)^a,
 \qquad \rho(t)=\frac{t}{e^t-1},\quad\rho(0)=1.
\end{align}
This holds for $L\geq1$ and $b\geq0$.
At the critical scale, every fixed band
$2k-1\leq C\leq M$ satisfies
\begin{equation}\label{h:fixed-band-rate}
 \mathbb E_g I_{m,=k,[2k-1,M]}
 =O_{m,M,S}\left(g^{-k+9/4}H^{6k-17/4}\right).
\end{equation}
Here each tuple and witness is counted once, with exactly $k$ crossed
marks and with the core of the full crossed set in the stated band.
All noncrossed placements and complementary topologies are included.
The lower complexity endpoint is imposed in this count, rather than
asserted for arbitrary disconnected marked complements.
\end{proposition}

\begin{proof}
The noncrossed-mark part of Lemma~\ref{lem:core-gap} did not require
connected marked complement.  It classifies every remaining mark
as exterior or attached along a nonempty subset of the outer cuffs.
Subsets assigned to distinct attached marks are disjoint, and
$p\leq a\leq\nu\leq C+2$.  At fixed $C,m$ there are finitely many
such subsets, slot assignments, and cuff roles. For each assignment,
the kernel retains all internal curve configurations.

Use the completed core of type $S_{h+j,\nu}$ and complexity $C$.
Each internal curve configuration determines its selected core.
The local count includes all such configurations and is independent of
the outer gluing twists.  Thus
Lemma~\ref{lem:attachment-bounds} applies with $c=C$, $p+r=m-k$,
attachment-cuff count $a$, and $d=p_2+p_3$.  The disjointness and
placement conditions were verified above.  Internal seams remain
inside the kernel and are not available as attachment cuffs.
Let $q_1$ be the number of connected components of $X$ after removing
$Z$ and all $p+r$ further closed marked supports.  Then
\eqref{att:joint-bound} gives
\begin{equation}\label{h:global-pattern-integration}
 \mathbb E_g I_{h,n,j;\mathbf p,r;q_1}^{(k+1)}
 \leq A_{C,m}g^{-(q_1-1)}
 \left(\frac{Le^L}{g}\right)^{m-k}g^{-C}
 \mathcal J_{L,b}^{(k+1),h,n,j}L^{-d}\rho(b)^a.
\end{equation}
This includes all cuff roles, all attachment sizes, and all stable
complementary topologies.  In particular, each three-cuff attachment
saves one power of $L$ relative to $Le^L$, as encoded by $d$.
Summing the finitely many assignments and
$1\leq q_1\leq\nu+3(p+r)-2a$, then using \eqref{h:kernel-bound},
proves \eqref{h:pattern-bound}.

For the critical rate choose $\delta=1/[4(k+1)]$, so that
$(k+1)(1+\delta)=k+5/4$.  Since $Le^L/g=O_S(1)$, the pattern
bound becomes
\begin{equation}\label{h:critical-pattern-rate}
 O_{C,m,S}\left(
 g^{-C+k+5/4}H^{3C-5/4-2j-d}\rho(b)^a\right).
\end{equation}
For $C\geq2k-1$ its unweighted factor is
\[
 g^{-k+9/4}H^{6k-17/4}
 \left(\frac{H^3}{g}\right)^{C-(2k-1)}H^{-2j-d}.
\]
The number of intrinsic types, seam orbit types, and outer contact
assignments in a fixed band is finite.  Every actual exact-$k$
configuration has the classification already proved, and its unique
crossed set is ordered by slots.  Sum the at most $\binom mk$ slot
choices and use $0\leq\rho(b)\leq1$. Maximizing constants over
this finite set before taking $g\to\infty$ proves
\eqref{h:fixed-band-rate}.
\end{proof}

\subsection{Exactly three crossings and the higher-crossing strata}

\begin{proposition}[The complete three-crossing class]
\label{prop:three-crossing}
For every fixed $m\geq3$, let $I_{m,=3}^{\rm conn}(L,b)$ count
ordered marked tuples and witnesses once, with pairwise closed-disjoint
marked supports, connected joint marked complement, and exactly three
crossed marks.  Then, at \eqref{h:critical-scale},
\begin{equation}\label{h:three-crossing}
 \mathbb E_g I_{m,=3}^{\rm conn}(L,b)
 =O_{m,S}\left(g^{-3/4}H^{55/4}\right)=o(1),
\end{equation}
uniformly for $b\geq0$.
\end{proposition}

\begin{proof}
Order the unique crossed triple by its marked slots and form the core
of these three curves and the witness.  Lemma~\ref{lem:core-gap} gives
$C\geq5$.  Set $K_m^{(4)}=2m+15$.  Proposition~\ref{prop:higher-fixed-bands}
with $k=3$ and $M=K_m^{(4)}-1$ controls the entire fixed band:
\begin{equation}\label{h:four-curve-fixed-band}
 \mathbb E_g I_{m,=3,[5,K_m^{(4)}-1]}
 =O_{m,S}(g^{-3/4}H^{55/4}).
\end{equation}
This includes every noncrossed placement and complementary topology;
restricting to connected marked complement only decreases the count.

The four-curve row of \eqref{h:tail-cutoffs}, obtained from
Lemma~\ref{lem:uniform-tail} with $q=4$ and $r=m-3$, gives
\begin{equation}\label{h:four-curve-tail}
 \mathbb E_g I_{m,=3,C\geq K_m^{(4)}}
 =O_{m,S}(g^{-1}H^{2-m}).
\end{equation}
These two complexity ranges exhaust the class.  The estimates count
each tuple and witness once, including the choices of crossed slots.
Consequently
\begin{equation}\label{h:three-crossing-assembly}
 \mathbb E_g I_{m,=3}^{\rm conn}
 \leq A_{m,S}\left(g^{-3/4}H^{55/4}+g^{-1}H^{2-m}\right).
\end{equation}
The second term is smaller than the first for each fixed $m$.
Both inputs are uniform in $b\geq0$, proving the proposition.
\end{proof}

\begin{corollary}[All higher crossing counts]\label{cor:higher-crossings}
For every fixed $m\geq4$, the once-per-tuple/witness incidence with
connected marked complement and at least four crossed marks satisfies
\begin{equation}\label{h:all-higher-crossings}
 \mathbb E_g I_{m,\kappa\geq4}^{\rm conn}(L,b)
 \leq A_{m,S}\left(
 g^{-7/4}H^{79/4}+g^{-1}H^{(1-m)/2}\right)
 =O_{m,S}(g^{-1}H^{(1-m)/2}),
\end{equation}
uniformly for $b\geq0$.
\end{corollary}

\begin{proof}
On the stratum $\kappa=k$, use the core of the entire crossed set,
not a selected proper subfamily.  Its complexity is at least $2k-1$
by Lemma~\ref{lem:core-gap}.  Choose the fixed cutoff
\begin{equation}\label{h:higher-cutoff}
 K_{m,k}=2m+3k+6,\qquad 4\leq k\leq m.
\end{equation}
Proposition~\ref{prop:higher-fixed-bands} controls
$2k-1\leq C<K_{m,k}$.  The sum of these fixed-band rates is
$O_{m,S}(g^{-7/4}H^{79/4})$, since
\[
 g^{-k+9/4}H^{6k-17/4}
 =g^{-7/4}H^{79/4}(H^6/g)^{k-4}.
\]
The last row of \eqref{h:tail-cutoffs} controls the remaining
complexities and, after summing $4\leq k\leq m$, gives
$O_{m,S}(g^{-1}H^{(1-m)/2})$.
The ratio of the fixed-band total to this last bound is
$g^{-3/4}H^{79/4+(m-1)/2}\to0$.  Only finitely many crossing
counts and constants are summed.  This proves the corollary.
\end{proof}

\subsection{Combining the witness estimates}

Return to $b=10\log L$. The local loss $Q_m$ in
\eqref{w:witness-set}--\eqref{w:local-loss} counts each disjoint tuple
once if it has a witness among the retained curves of length at most
$L$. We first bound the tuples with disconnected joint complement.
For the others, we divide each tuple's contribution among its
witnesses and sum the estimates for the individual classes.

\begin{theorem}[Local weighted deletion at every fixed order]
\label{thm:weighted-loss}
For every fixed $m\geq1$ and $S\in\mathbb R$,
\begin{equation}\label{h:weighted-loss}
 \mathbb E_g Q_m\leq A_{m,S}(H^{-2}+g^{-1})=o(1).
\end{equation}
The assertion includes disconnected joint marked complements.
Restricting the marked slots to arbitrary Borel windows in
$(-\infty,S]$ only decreases this nonnegative count.
\end{theorem}

\begin{proof}
For $m=1$, the number of nonisolated vertices is at most twice the
number of conflict edges.  Thus $Q_1\leq2B_{g,S}$ and
Proposition~\ref{prop:conflicts} proves the assertion.
Fix $m\geq2$.  First separate the disconnected marked branch.
Proposition~\ref{prop:disconnected} gives
\begin{equation}\label{h:disconnected-deletion}
 Q_m^{\rm disc}\leq D_m^{\rm disc},\qquad
 \mathbb E_gD_m^{\rm disc}=O_{m,S}(g^{-1}).
\end{equation}
All remaining sums are over connected marked complements.

Partition the original witness set into equal-support witnesses
(class $E$), pure-touching witnesses (class $T$), and witnesses
crossing exactly $k$ marked
supports, $1\leq k\leq m$.  Equal support at one mark precludes
contact with any other mark because the marked supports are
closed-disjoint.  Mixed touching and crossing signatures belong to
the crossing class.  These alternatives are disjoint and exhaustive.

Let $X_k$ denote the class of witnesses crossing exactly $k$ marked
supports. Using the denominator $d(\tau)$ from \eqref{w:normalized-partition}
for every class gives
\begin{equation}\label{h:exact-crossed-set}
 Q_m^{\rm conn}
 =Q_{E,m}^{\rm conn}+Q_{T,m}^{\rm conn}
  +\sum_{k=1}^m Q_{X_k,m}^{\rm conn}.
\end{equation}
By \eqref{w:raw-domination}, a bound on the raw incidence bounds each expected share in \eqref{h:exact-crossed-set}. For fixed $m$ only
rows with at most $m$ crossed marks contribute.
\begin{center}
\begin{minipage}{\linewidth}\centering
{\small\itshape Expected deletion by witness class for a connected marked complement\par}\smallskip
\small
\begin{tabular}{p{0.53\linewidth}p{0.39\linewidth}}
\hline
Witness class and input & Expected normalized share, at most\\
\hline
Equal support; Lemma~\ref{w:equal-support}
 & $O_{m,S}(H^{-6}+g^{-1})$\\
Pure touching; Proposition~\ref{prop:touching}
 & $O_{m,S}(\log H/H^{10})$\\
Exactly one crossed mark;
Proposition~\ref{prop:single-crossing}
 & $O_{m,S}(H^{-2})$\\
Exactly two crossed marks;
Proposition~\ref{prop:two-crossing}
 & $O_{m,S}(H^{-4})$\\
Exactly three crossed marks;
Proposition~\ref{prop:three-crossing}
 & $O_{m,S}(g^{-3/4}H^{55/4})$\\
At least four crossed marks;
Corollary~\ref{cor:higher-crossings}
 & $O_{m,S}(g^{-1}H^{(1-m)/2})$\\
\hline
\end{tabular}
\end{minipage}
\end{center}
The crossing estimates include every core complexity and all
noncrossed placements.  Each displayed term is
$O_{m,S}(H^{-2}+g^{-1})$, and there are finitely many classes at fixed
$m$.  Summing \eqref{h:exact-crossed-set} and adding
\eqref{h:disconnected-deletion} proves \eqref{h:weighted-loss}.
\end{proof}

\begin{corollary}[Bounded tuple and slot weights]\label{cor:bounded-weight-loss}
For a measurable real or complex function $F$ on ordered retained
disjoint $m$-tuples, uniformly bounded over the relevant surfaces
and tuples, define
\[
 D_m(F)=\sum_{\tau\in\mathcal T_m}F(\tau),\qquad
 D_m^{\rm iso}(F)
 =\sum_{\substack{\tau\in\mathcal T_m\\
                  \text{all marked vertices isolated}}}F(\tau).
\]
Then, at every fixed $m\geq1$,
\begin{equation}\label{h:bounded-weight-loss}
 \mathbb E_g\bigl|D_m(F)-D_m^{\rm iso}(F)\bigr|
 \leq A_{m,S}\|F\|_\infty(H^{-2}+g^{-1}).
\end{equation}
For product slot weights $F(\tau)=\prod_{i=1}^m
f_i(\ell(\alpha_i)-H+\log H)$, the norm on the right may be replaced
by $\prod_i\|f_i\|_\infty$.  Arbitrary Borel slot windows, including
repeated windows, are covered by taking indicator weights.
\end{corollary}

\begin{proof}
The selected marks have no mutual conflicts, so all marked vertices
are isolated exactly when the external witness count is zero.
Consequently
\begin{equation}\label{h:exact-weighted-difference}
 D_m(F)-D_m^{\rm iso}(F)
 =\sum_{\tau\in\mathcal T_m}F(\tau)\mathbf1_{\{d(\tau)>0\}}.
\end{equation}
The absolute value of this sum is at most $\|F\|_\infty Q_m$
pointwise.  Apply Theorem~\ref{thm:weighted-loss}.
\end{proof}

Section~\ref{sec:probability} combines Theorems~\ref{thm:disjoint}
and~\ref{thm:weighted-loss} to obtain the isolated factorial-moment
estimate \eqref{p:isolated-factorial-limit} and prove Poisson convergence.
Isolation is tested against every retained curve below the centered
ceiling $S$, not only against curves in the observation intervals.
All bounds above are for fixed marked order; their constants need
not remain bounded when that order grows with $g$.

% END sections/06_higher_cores.tex

% BEGIN sections/07_probability.tex
\section{Poisson convergence and the first arrivals}
\label{sec:probability}

To prove Theorems~\ref{thm:main} and~\ref{thm:arrivals}, we first count
curves in bounded intervals of centered length. For isolated curves,
a product of falling factorials counts exactly the disjoint tuples that
survive deletion. Taking expectations gives the mixed factorial moments
and hence the limiting Poisson law for the isolated counts. The coupling
transfers this limit to the full counts. To identify the first arrivals,
we also control curves below the left endpoint of the observation window.

We retain $H=\log g$ and $a_g=H-\log H$ from \eqref{i:center}.
Count-vector dimensions, moment orders, and observation windows are
fixed before $g\to\infty$. For probability measures on a common
measurable space, we use
\begin{equation}
 d_{\mathrm{TV}}(P,Q)=\sup_A|P(A)-Q(A)|.
 \label{p:tv-convention}
\end{equation}

\subsection{A factorial-moment criterion on a fixed lattice}

We use the following factorial-moment criterion for count vectors, allowing
zero Poisson parameters. Write $(n)_0=1$ and
$(n)_r=n(n-1)\cdots(n-r+1)$ for $r\ge1$.

\begin{lemma}
\label{lem:factorial-poisson}
Let $X_n=(X_{n,1},\ldots,X_{n,d})$ take values in $\mathbb N_0^d$, with
$d$ fixed, and let $\lambda_1,\ldots,\lambda_d\in[0,\infty)$ be fixed
parameters. Suppose that all mixed factorial moments are finite and,
for every $\mathbf r=(r_1,\ldots,r_d)\in\mathbb N_0^d$,
\begin{equation}
 \mathbb E\prod_{j=1}^d(X_{n,j})_{r_j}
 \longrightarrow\prod_{j=1}^d\lambda_j^{r_j},
 \qquad 0\le\lambda_j<\infty.
 \label{p:factorial-hypothesis}
\end{equation}
Then
\begin{equation}
 d_{\mathrm{TV}}\!\left(\mathcal L(X_n),
       \bigotimes_{j=1}^d\operatorname{Pois}(\lambda_j)\right)
 \longrightarrow0.
 \label{p:factorial-conclusion}
\end{equation}
Here $0^0=1$, and $\operatorname{Pois}(0)$ denotes the point mass at zero.
\end{lemma}

\begin{proof}
Let $J=\{j:\lambda_j>0\}$. If $J$ is nonempty, take $r_j=0$ for
$j\notin J$ in \eqref{p:factorial-hypothesis}. By the factorial-moment
criterion of Mirzakhani and Petri~\cite[Theorem~2.9]{MP19}, the joint
atom probabilities for $(X_{n,j})_{j\in J}$ converge to those of
independent Poisson variables with means $(\lambda_j)_{j\in J}$.
For the other coordinates, Markov's inequality and the first moments give
\[
 \mathbb P\{\text{some }j\notin J\text{ has }X_{n,j}>0\}
 \le \sum_{j\notin J}\mathbb EX_{n,j}\longrightarrow0.
\]
Thus setting these coordinates to zero changes each atom probability by
$o(1)$. This gives atom convergence to the stated product law, including
when $J$ is empty and the limit is concentrated at the zero vector.

Finally, let $p_n$ and $p$ be the prelimit and limiting mass functions.
For $\varepsilon>0$, choose a finite set $F$ with $p(F)>1-\varepsilon$.
Then
\begin{equation}
 \sum_{\mathbf a}|p_n(\mathbf a)-p(\mathbf a)|
 \le \sum_{\mathbf a\in F}|p_n(\mathbf a)-p(\mathbf a)|
       +p_n(F^c)+p(F^c).
 \label{p:lattice-tv}
\end{equation}
The finite sum tends to zero and $p_n(F^c)\to p(F^c)$, so the limsup is
at most $2\varepsilon$. Under \eqref{p:tv-convention}, total variation
is half this $\ell^1$ distance. Letting $\varepsilon\downarrow0$
proves \eqref{p:factorial-conclusion}.
\end{proof}

\subsection{Isolated retained counts}

Fix a centered ceiling $S\in\mathbb R$. For sufficiently large $g$,
the corresponding length and cuff cutoffs are $L=a_g+S$ and
$b=10\log L$. As in Section~\ref{sec:foundations}, write
$\Xi_{g,S}^{\rm tr}$ and $\widehat\Xi_{g,S}$ for the centered counting
measures of retained and isolated retained curves. Isolation is tested
in the entire retained conflict graph below $S$, including neighbors
outside the selected observation windows.

Let $I_j=(u_j,v_j]\subset(-\infty,S]$, $1\le j\le d$, be a fixed
pairwise disjoint family of bounded half-open intervals. Set
\[
 \widehat N_{g,j}=\widehat\Xi_{g,S}(I_j),
 \qquad \lambda_j=\lambda(I_j)=\int_{I_j}\frac{e^t}{8\pi^2}\,dt.
\]
Empty intervals are permitted. For a multiindex $\mathbf r$ of total
order $m\ge1$, let $D_{g,\mathbf r}$ count ordered retained selections
with $r_j$ slots in $I_j$ and pairwise disjoint closed supports. Let
$Q_{g,\mathbf r}$ count those selections containing at least one
nonisolated selected curve. These definitions include all joint
complementary topologies.

Distinct isolated vertices have disjoint closed
supports, and a disjoint selection contributes precisely when none
of its selected curves has a retained neighbor. Thus, on each surface,
\begin{equation}
 \prod_{j=1}^d(\widehat N_{g,j})_{r_j}
       =D_{g,\mathbf r}-Q_{g,\mathbf r}.
 \label{p:isolated-factorial-identity}
\end{equation}
The slots within each interval are already ordered by the falling
factorial. Since the intervals are disjoint, there is no additional
multinomial or factorial coefficient. At order zero the identity is
$1=1$.

Theorem~\ref{thm:disjoint} and Theorem~\ref{thm:weighted-loss} give, for
each fixed $\mathbf r$,
\begin{equation}
 \mathbb E_g\prod_{j=1}^d(\widehat N_{g,j})_{r_j}
   =\prod_{j=1}^d\lambda_j^{r_j}
       +O_{m,S,\mathbf I}\!\left(\frac{\log H}{H}\right).
 \label{p:isolated-factorial-limit}
\end{equation}
Indeed, every tuple lost from the chosen intervals also contributes
to the unrestricted loss in Theorem~\ref{thm:weighted-loss}.
All these moments are finite. One may also see this directly from the
metric-uniform bound
$N_X^{\rm prim}(L)\le(g-1)e^{L+7}$ for primitive unoriented geodesics
on a closed genus-$g$ surface, stated in~\cite[Theorem~2.8]{HSWX}.

Applying Lemma~\ref{lem:factorial-poisson} to
\eqref{p:isolated-factorial-limit}, we obtain
\begin{equation}
 d_{\mathrm{TV}}\!\left(
  \mathcal L_g(\widehat N_{g,1},\ldots,\widehat N_{g,d}),
  \bigotimes_{j=1}^d\operatorname{Pois}(\lambda_j)\right)
 \longrightarrow0.
 \label{p:isolated-tv}
\end{equation}
The error in \eqref{p:isolated-factorial-limit} is for each fixed
moment order. Since its constant depends on $m$, the factorial-moment
argument gives a qualitative limit in \eqref{p:isolated-tv}.

\subsection{Full counts and vague convergence}

We transfer the isolated-count limits to the full process using
Corollary~\ref{cor:f-isolation-coupling}, which bounds the mismatch
probability by
\begin{equation}
 \varepsilon_g(S)
 :=\mathbb P_g\!\left\{
    \Xi_g^{(8)}|_{(-\infty,S]}\ne\widehat\Xi_{g,S}\right\}
   =O_S\!\left(\frac{\log H}{H}\right).
 \label{p:cutoff-coupling}
\end{equation}

For a fixed disjoint bounded interval family $I_1,\ldots,I_d$, choose one
fixed $S\ge\max_j\sup I_j$. On the same sampled surface, the full and
isolated count vectors differ only on the event in
\eqref{p:cutoff-coupling}. The coupling and triangle inequalities give
\begin{equation}
 \begin{split}
 &d_{\mathrm{TV}}\!\left(
    \mathcal L_g\bigl((\Xi_g^{(8)}(I_j))_{j=1}^d\bigr),
    \bigotimes_{j=1}^d\operatorname{Pois}(\lambda(I_j))\right)\\
 &\quad\le\varepsilon_g(S)+
 d_{\mathrm{TV}}\!\left(
    \mathcal L_g\bigl((\widehat\Xi_{g,S}(I_j))_{j=1}^d\bigr),
    \bigotimes_{j=1}^d\operatorname{Pois}(\lambda(I_j))\right)
 \longrightarrow0.
 \end{split}
 \label{p:full-bounded-tv}
\end{equation}
For a finite family of possibly overlapping bounded half-open intervals,
refine their endpoints into disjoint half-open cells. Each original count
is a deterministic sum of cell counts, and total variation contracts
under deterministic maps. Thus every such finite count vector converges
in total variation to the corresponding vector of $\Pi$ counts. The
overlapping coordinates have the dependence induced by their common
cells. This argument uses a single ceiling for each finite family and
does not require consistency between isolated processes constructed at
different ceilings.

To pass to the vague topology and then to half-line counts, we use
the first-moment estimate~\eqref{i:external-mean}.
At every fixed $s\in\mathbb R$, substitution of
$L=a_g+s$ gives
\begin{equation}
 \mathbb E_g N_g(s)
 =\Lambda(s)\left(1-\frac{\log H+c_0-s}{H}\right)
   +O_s\!\left(g^{-1/2}H^{3/2}+g^{-1}H^3\right)
 \longrightarrow\Lambda(s),
 \label{p:critical-first-moment}
\end{equation}
where $N_g(s)=\Xi_g^{(8)}(( -\infty,s])$ and
$\Lambda(s)=e^s/(8\pi^2)$. The lower length condition in
\eqref{i:external-mean} holds eventually for each fixed $s$.

The interval arguments above use the convention $(u,v]$.  For a fixed
endpoint $t$ and $\delta>0$, Markov's inequality gives
\[
 \mathbb P_g\{\Xi_g^{(8)}(\{t\})>0\}
 \le \mathbb E_g N_g(t)-\mathbb E_g N_g(t-\delta).
\]
By \eqref{p:critical-first-moment}, the limit superior as $g\to\infty$
is at most $\Lambda(t)-\Lambda(t-\delta)$. This bound tends to zero
as $\delta\downarrow0$.  Changing the inclusion of finitely many fixed
endpoints therefore has vanishing coupling cost.  This justifies the
bounded-interval count conclusions for any endpoint convention.

If $K\subset\mathbb R$ is compact, choose $s\ge\sup K$.
Equation~\eqref{p:critical-first-moment} bounds its expected count
for all sufficiently large $g$. For the finitely many remaining genera,
the metric-uniform primitive-geodesic bound used above gives a finite
expectation; the count is zero if the length cutoff is nonpositive.
Hence
\begin{equation}
 C_K:=\sup_{g\ge3}\mathbb E_g\Xi_g^{(8)}(K)<\infty,
 \qquad
 \sup_{g\ge3}\mathbb P_g\{\Xi_g^{(8)}(K)>A\}\le C_K/A.
 \label{p:compact-tightness}
\end{equation}
The same count-tightness conclusion holds for every bounded Borel set
by containment in a compact interval.

We apply \cite[Theorem~1.1 and Lemma~1.2]{Kallenberg}. On a locally
compact second countable Hausdorff space, finite joint count convergence
on a dissecting-covering semiring of relatively compact Borel sets
whose boundaries have zero mass almost surely under the target law
implies convergence in distribution of random Radon measures for the
vague topology.
The count bounds in \eqref{p:compact-tightness} also verify the
tightness criterion. On $\mathbb R$, take
\begin{equation}
 \mathcal J=\{\varnothing\}\cup
       \{(a,b]:a,b\in\mathbb Q,\ a<b\}.
 \label{p:dc-semiring}
\end{equation}
Intersections are members of the semiring, and the difference of two
members is a finite disjoint union of members.
Every bounded set can be covered by finitely many of its intervals of
arbitrarily small diameter, by enclosing it in a larger rational interval
and partitioning that interval at rational points. For every rational
$q$, $\mathbb E\Pi(\{q\})=\lambda(\{q\})=0$; hence all rational
endpoints simultaneously have zero mass almost surely. Finally, finite disjoint refinement together with
\eqref{p:full-bounded-tv} supplies every required joint count limit. The criterion proves
\begin{equation}
 \Xi_g^{(8)}\Longrightarrow\Pi
 \quad\hbox{in }\mathcal M_p(\mathbb R)
 \quad\hbox{for the vague topology}.
 \label{p:vague-limit}
\end{equation}
Here the subspace of locally finite counting measures, allowing
multiplicities, is closed in the space of Radon measures for the vague
topology, as recalled in~\cite[Section~1]{Kallenberg}. Prelimit endpoint
atoms need not be excluded. The target is simple because its intensity
is diffuse.

For every nonnegative $f\in C_c(\mathbb R)$, the map
$\xi\mapsto e^{-\xi(f)}$ is bounded and continuous in the vague
topology. Thus \eqref{p:vague-limit} and the Poisson Laplace functional give
\begin{equation}
 \mathbb E_g\exp\!\left(-\int f\,d\Xi_g^{(8)}\right)
 \longrightarrow
 \exp\!\left[-\int_{\mathbb R}(1-e^{-f(t)})
                         \frac{e^t}{8\pi^2}\,dt\right].
 \label{p:laplace-limit}
\end{equation}

\subsection{Half-line counts and the order of truncation limits}

The half-line $(-\infty,s]$ is not relatively compact, so we control
its left tail before removing the truncation. For each fixed $R>0$,
Markov's inequality and
\eqref{p:critical-first-moment} at $s=-R$ imply
\begin{equation}
 \limsup_{g\to\infty}\mathbb P_g\{N_g(-R)>0\}
 \le\Lambda(-R)=\frac{e^{-R}}{8\pi^2}.
 \label{p:left-tail}
\end{equation}
Thus the limsup vanishes as $R\to\infty$. The estimate takes the
genus limit at fixed $R$; it asserts no uniformity for $R=R(g)$.

\begin{proposition}
\label{p:halfline-counts}
For any fixed real thresholds $s_1,\ldots,s_d$,
\begin{equation}
 d_{\mathrm{TV}}\!\left(
     \mathcal L_g((N_g(s_i))_{i=1}^d),
     \mathcal L((\Pi(( -\infty,s_i]))_{i=1}^d)\right)
 \longrightarrow0.
 \label{p:halfline-tv}
\end{equation}
\end{proposition}

\begin{proof}
Choose $R>0$ with $-R<\min_i s_i$, and let
\[
 \begin{aligned}
 Y_g&=(N_g(s_i))_{i=1}^d,
 &Y_g^R&=(\Xi_g^{(8)}((-R,s_i]))_{i=1}^d,\\
 Y&=(\Pi(( -\infty,s_i]))_{i=1}^d,
 &Y^R&=(\Pi((-R,s_i]))_{i=1}^d.
 \end{aligned}
\]
The same-surface identity $Y_g=Y_g^R+N_g(-R)\mathbf1$ shows that
the prelimit truncation fails only when $N_g(-R)>0$. The corresponding
target truncation fails with probability
$1-e^{-\Lambda(-R)}$. Hence
\begin{equation}
 \begin{split}
 d_{\mathrm{TV}}(\mathcal L_g(Y_g),\mathcal L(Y))
 &\le\mathbb P_g\{N_g(-R)>0\}
    +d_{\mathrm{TV}}(\mathcal L_g(Y_g^R),\mathcal L(Y^R))\\
 &\hspace{2em}+1-e^{-\Lambda(-R)}.
 \end{split}
 \label{p:halfline-triangle}
\end{equation}
For fixed $R$, the middle term tends to zero by
\eqref{p:full-bounded-tv} after refinement into disjoint intervals.
Together with \eqref{p:left-tail}, this gives
\begin{equation}
 \limsup_{g\to\infty}
 d_{\mathrm{TV}}(\mathcal L_g(Y_g),\mathcal L(Y))
 \le\Lambda(-R)+1-e^{-\Lambda(-R)}.
 \label{p:halfline-truncation-bound}
\end{equation}
Letting $R\to\infty$ proves the result. Repeated or unordered
thresholds cause no difficulty: refine at the distinct endpoints and
then take deterministic sums and repeated coordinates.
\end{proof}

The limiting half-line coordinates are cumulative counts, with independent
Poisson increments on disjoint intervals. They are generally dependent.
The same truncation and finite refinement also handle any fixed count
vector containing both bounded intervals and these half-lines.

\subsection{Existence and identification of the first arrivals}

We first ensure that the required prelimit arrivals exist on every
surface for sufficiently large genus.

\begin{lemma}
\label{p:arrival-existence}
Every closed genus-$g$ hyperbolic surface with $g\ge3$ has at least
$g-2$ distinct primitive unoriented figure-eight geodesics. In particular,
the first $k$ figure-eight lengths, ordered with multiplicities, exist
on every such surface once $g\ge k+2$.
\end{lemma}

\begin{proof}
Topologically decompose the surface into $g$ one-holed tori attached to
a sphere with $g$ boundary components. Straighten their separating cuffs
to disjoint geodesics and pants-decompose the central sphere. This gives
$g-2$ embedded geodesic pants, each with three distinct ambient cuffs.
Choose an intrinsic primitive figure-eight geodesic in each interior.
The inclusion is isometric on the interior, so each chosen curve remains
a closed geodesic with exactly one transverse self-intersection in the
ambient surface. Curves in different interiors are distinct. It also
remains primitive: any shorter ambient periodic traversal of that same
geodesic would lie in the same interior and contradict intrinsic
primitivity. Finally, on a fixed compact hyperbolic surface only
finitely many primitive closed geodesics have length below any fixed
bound. The existence of at least $k$ curves therefore gives the first
$k$ lengths, with multiplicities retained.
\end{proof}

We next justify the ordered list of limiting points. For each integer
$n$, the count $\Pi(( -\infty,n])$ is finite almost surely.
Intersecting these countably many events and using monotonicity gives
finiteness of all half-line counts simultaneously. Moreover,
\[
 \mathbb P\{\Pi(( -\infty,-R])>0\}
    =1-e^{-\Lambda(-R)}\longrightarrow0.
\]
For integer $R$ these events decrease, so almost surely some left
half-line contains no points. The process nevertheless has infinitely
many points to the right. Indeed,
for every fixed positive integer $k$,
\[
 \mathbb P\{\Pi(\mathbb R)<k\}
 \le\mathbb P\{\Pi(( -\infty,S])<k\}
 \longrightarrow0\qquad(S\to\infty).
\]
Local finiteness and the diffuse intensity therefore give an infinite
well-defined ordered list $S_1<S_2<\cdots$ of its points.

The increasing map $y=\Lambda(s)=e^s/(8\pi^2)$ sends $\mathbb R$
onto $(0,\infty)$ and satisfies $dy=\lambda(ds)$. Its image process
has independent Poisson counts of Lebesgue means on disjoint sets,
so it is a rate-one homogeneous Poisson process. For a Poisson process
with intensity $\rho$, fix a Borel set $B$ with $0<\rho(B)<\infty$.
Conditional on having $n$ points in $B$, its restriction to $B$ has the
law of the counting measure of $n$ independent samples from
$\rho|_B/\rho(B)$. Indeed, on every finite measurable partition of $B$,
the independent Poisson counts conditioned on their sum have the same
multinomial law as the sample counts. These counts determine the law
of the restriction.

In $(0,T]$ with $T>0$, the ordered positions given $n$ points are thus
the order statistics of $n$ independent uniform samples. When $n\ge k$
and $0<y_1<\cdots<y_k<T$, the remaining $n-k$ ordered
positions occupy a simplex of volume $(T-y_k)^{n-k}/(n-k)!$.
Integrating the joint density $n!/T^n$ over this simplex and
averaging over the Poisson count gives
\begin{equation}
 \sum_{n=k}^{\infty}
 \frac{e^{-T}T^n}{n!}
 \frac{n!}{(n-k)!T^n}(T-y_k)^{n-k}
 =e^{-y_k}.
 \label{p:homogeneous-arrival-density}
\end{equation}
This holds for every $T>y_k$. Changing variables to the increments
$e_1=y_1$, $e_j=y_j-y_{j-1}$ has Jacobian one and transforms the
density into $\prod_{j=1}^k e^{-e_j}\mathbf1_{\{e_j>0\}}$.
Consequently, if $E_1,E_2,\ldots$ are independent rate-one
exponential variables and $\Gamma_j=E_1+\cdots+E_j$, then
\begin{equation}
 (S_1,\ldots,S_k)
 \overset{d}=\bigl(\log(8\pi^2\Gamma_1),\ldots,
                  \log(8\pi^2\Gamma_k)\bigr).
 \label{p:gamma-arrivals}
\end{equation}
Changing variables back in \eqref{p:homogeneous-arrival-density}
gives the joint density
\begin{equation}
 p_k(s_1,\ldots,s_k)
 =\mathbf1_{\{s_1<\cdots<s_k\}}
     e^{-\Lambda(s_k)}\prod_{j=1}^k\frac{e^{s_j}}{8\pi^2}.
 \label{p:inhomogeneous-arrival-density}
\end{equation}
There is no extra $k!$ factor for this ordered-arrival density.
In particular,
\begin{equation}
 \mathbb P\{S_j>s\}
 =\mathbb P\{\Pi(( -\infty,s])<j\}
 =e^{-\Lambda(s)}\sum_{r=0}^{j-1}\frac{\Lambda(s)^r}{r!}.
 \label{p:arrival-survival}
\end{equation}
Every limiting arrival marginal is continuous.

Let $T_{g,j}=\ell_{g,j}^{(8)}-a_g$, defined for sufficiently large
$g$ by Lemma~\ref{p:arrival-existence}. For any threshold vector
$(s_1,\ldots,s_k)$, including repeated or unordered thresholds,
the exact event identity is
\begin{equation}
 \{T_{g,j}\le s_j\text{ for every }j\le k\}
 =\{N_g(s_j)\ge j\text{ for every }j\le k\}.
 \label{p:order-statistic-event}
\end{equation}
It remains valid at tied prelimit lengths because counts retain
multiplicities. Proposition~\ref{p:halfline-counts} shows that these
probabilities converge to the corresponding probabilities for
$(S_1,\ldots,S_k)$. The boundary of any lower orthant for the limit
vector lies in a finite union of events $\{S_j=s_j\}$, all null by
\eqref{p:arrival-survival}. The joint limiting distribution function
is therefore continuous, and convergence of its lower-orthant
probabilities proves
\begin{equation}
 (T_{g,1},\ldots,T_{g,k})\Longrightarrow(S_1,\ldots,S_k)
 \quad\hbox{in }\mathbb R^k.
 \label{p:arrival-vector-limit}
\end{equation}
Together with \eqref{p:gamma-arrivals}, this is the arrival-vector
assertion of Theorem~\ref{thm:arrivals}.

Escape at both ends can also be controlled directly. For fixed
$R>0$ and $S\in\mathbb R$,
\begin{equation}
 \begin{split}
 \limsup_{g\to\infty}\mathbb P_g\{T_{g,1}\le-R\}
 &\le\Lambda(-R),\\
 \lim_{g\to\infty}\mathbb P_g\{T_{g,k}>S\}
 &=e^{-\Lambda(S)}\sum_{r=0}^{k-1}\frac{\Lambda(S)^r}{r!}.
 \end{split}
 \label{p:arrival-tightness}
\end{equation}
The first bound tends to zero as $R\to\infty$, and the second
expression tends to zero as $S\to\infty$. Since the vector is ordered,
these bounds control every coordinate. Both truncation parameters are
sent to their limits after genus tends to infinity. No continuity of
the order-statistic map on the entire vague measure space is assumed.

Finally, the non-simple systole is realized by a figure-eight
by~\cite[Theorem~4.2.4]{Buser2010}; see also
\cite[Introduction and Section~2.4]{HSWX}. Thus
$M_g=\ell_{\rm sys}^{\rm ns}(X)-a_g=T_{g,1}$, and
\begin{equation}
 \{M_g>s\}=\{N_g(s)=0\}.
 \label{p:minimum-void}
\end{equation}
The one-dimensional case of \eqref{p:halfline-tv} gives
\begin{equation}
 \mathbb P_g\{M_g>s\}\longrightarrow e^{-\Lambda(s)}
   =\exp\!\left(-\frac{e^s}{8\pi^2}\right).
 \label{p:minimum-survival}
\end{equation}
If $G=\log(8\pi^2)-S_1$, then the continuity of the limiting
distribution implies
\begin{equation}
 \mathbb P\{G\le x\}
 =\mathbb P\{S_1\ge\log(8\pi^2)-x\}
 =\exp(-e^{-x}).
 \label{p:gumbel}
\end{equation}
Thus $\log(8\pi^2)-M_g$ converges to the standard Gumbel law.
Continuity of this limiting law permits either endpoint convention.

This completes the proofs of Theorems~\ref{thm:main} and~\ref{thm:arrivals}.
The higher arrivals describe the figure-eight spectrum at fixed order.
We next determine the cuff shapes accompanying these lengths.
Section~\ref{sec:tails} returns to the minimum and proves the uniform
integrability needed for convergence of its moments.
% END sections/07_probability.tex

% BEGIN sections/08_marked_intensity.tex
\section{Cuff shapes and marked first intensity}\label{sec:marked-intensity}

We now record the three cuff lengths of each supporting pair of pants.
For a single figure-eight sampled from the normalized expected count,
the third cuff contributes half the perimeter in the limit, and the two
winding cuffs divide the other half in a uniformly random proportion.
We derive this law from the one-curve integral. The joint shapes of the
first arrivals will also require the tuple estimates of
Section~\ref{sec:marked-factorials}.

\subsection{Expected counting measures}

We count each primitive unoriented figure-eight once. For each curve,
label its two winding cuffs in either order with probability $1/2$,
independently across curves, and use that order at every cutoff.
Write $x,y$ for their lengths and $z$ for the length of the third cuff.
These labels refer to intrinsic boundary circles, including when two
of them have the same ambient image. The expectation $\mathbb E_g^\circ$
includes the fair ordering. In a first-moment calculation, averaging
over the order simply replaces a cuff test $F(x,y,z)$ by
$[F(x,y,z)+F(y,x,z)]/2$.

Set
\begin{equation}\label{mk1:eq:scale}
 H=\log g,\qquad a_g=H-\log H,\qquad L_g=a_g+S,
 \qquad \Lambda(S)=\frac{e^S}{8\pi^2},
\end{equation}
where $S\in\mathbb R$ is fixed and $g$ is sufficiently large. For a curve
$\gamma$, put
\begin{equation}\label{mk1:eq:shape}
 Z_\gamma=\frac{(x_\gamma,y_\gamma,z_\gamma)}
 {x_\gamma+y_\gamma+z_\gamma}\in\Delta,
 \qquad
 \Delta=\{(u,v,w)\in[0,1]^3:u+v+w=1\}.
\end{equation}
Let $N_g(S)$ count all figure-eights of length at most $L_g$, retaining
multiplicities of equal lengths. The expected shape measure on $\Delta$
and the expected joint measure on $(-\infty,S]\times\Delta$ act on bounded
Borel tests $F$ and $\Psi$ by
\begin{align}
 \mathcal M_{g,S}(F)
 &=\mathbb E_g^\circ\sum_{\ell_\gamma\le L_g}F(Z_\gamma),
 &m_{g,S}&=\mathbb E_gN_g(S),\label{mk1:eq:markmeasure}\\
 \mathcal J_{g,S}(\Psi)
 &=\mathbb E_g^\circ\sum_{\ell_\gamma\le L_g}
 \Psi(\ell_\gamma-a_g,Z_\gamma).\label{mk1:eq:jointmeasure}
\end{align}
Both measures have total mass $m_{g,S}$, which is positive for all
sufficiently large $g$. Write
$\widehat{\mathcal M}_{g,S}=\mathcal M_{g,S}/m_{g,S}$.
To sample from this probability law, weight the surface law by
$N_g(S)/m_{g,S}$, choose uniformly among its figure-eights of length
at most $L_g$, and apply the fair cuff ordering.

A support is good if its closure is embedded, its three ambient cuffs
are distinct, and its complement is connected. Among curves of length at most $L_g$,
retain those with good support and all three cuff lengths at least
$10\log L_g$.

\begin{theorem}[Boundary shape and joint first intensity]\label{mk1:thm:m1}
For each fixed $S$, $m_{g,S}\to\Lambda(S)$. For every $F\in C(\Delta)$,
\begin{equation}\label{mk1:eq:shape-limit}
 \mathcal M_{g,S}(F)\longrightarrow
 \Lambda(S)\int_0^1 F\left(\frac t2,\frac{1-t}2,\frac12\right)dt.
\end{equation}
Consequently $\widehat{\mathcal M}_{g,S}$ converges weakly to the probability
law $\nu$ of $(T/2,(1-T)/2,1/2)$, $T$ uniform on $[0,1]$.
More precisely, for every bounded continuous $\Psi$ on
$(-\infty,S]\times\Delta$,
\begin{equation}\label{mk1:eq:joint-limit}
 \mathcal J_{g,S}(\Psi)\longrightarrow
 \int_{-\infty}^S\frac{e^s}{8\pi^2}
 \int_0^1\Psi\!\left(s,\left(\frac t2,\frac{1-t}2,\frac12\right)\right)dt\,ds.
\end{equation}
The expected count of retained curves also tends to $\Lambda(S)$, and
their shape and joint measures satisfy \eqref{mk1:eq:shape-limit}
and \eqref{mk1:eq:joint-limit}. Dividing each retained measure by the
expected retained count divides its limit by $\Lambda(S)$.
\end{theorem}

For the shape up to exchange of the winding cuffs, sort the first two
coordinates. The limiting shape is then $(V/2,(1-V)/2,1/2)$, where
$V=\min(T,1-T)$ has density $2$ on $[0,1/2]$. The third cuff retains
its distinguished role.

\subsection{Exact marked integral and analytic inputs}

For cuff lengths $x,y,z\ge0$, the exact figure-eight length is given by
\eqref{f:pants-length}.
Let $D_L=\{(x,y,z)\ge0:\ell(x,y,z)\le L\}$ and
$\Phi(x,y,z)=\sinh(x/2)\sinh(y/2)\sinh(z/2)$.
On $D_L$, the identity gives $x+y<L$ and $z<L$.

\begin{lemma}[Fair-order integration]\label{mk1:lem:integration}
For $g\ge3$, $L\ge0$, and any nonnegative Borel function $F$ of the
three cuff lengths in the specified order,
\begin{equation}\label{mk1:eq:exact-intensity}
 \mathbb E_g^\circ\sum_{\substack{\gamma\text{ good}\\\ell_\gamma\le L}}
 F(x_\gamma,y_\gamma,z_\gamma)
 =\frac1{2V_g}\int_{D_L}F(x,y,z)
 V_{g-2,3}(x,y,z)xyz\,dx\,dy\,dz.
\end{equation}
No permutation symmetry of $F$ is required. The identity extends to
signed or complex tests whose absolute integral is finite.
\end{lemma}

\begin{proof}
Fix one good pair of pants. Each choice of its third cuff determines
one figure-eight. Conditional on that choice, the two orders of the
winding cuffs receive equal weight. Conditional on the surface, the expected contribution of
the three figure-eights is half the sum of $F\mathbf1_{D_L}$ over the
six ordered triples of geometric cuffs. This counts the same six
labelings even when some cuff lengths agree.

The ordered-cuff integration formula used in
Lemma~\ref{lem:f-one-point-integral} has coefficient one for $g\ge3$,
and $V_{0,3}=1$. Applied to the ordered test just described, it gives
\eqref{mk1:eq:exact-intensity}. Monotone convergence covers nonnegative
Borel tests, and absolute integrability gives the signed and complex
versions.

\end{proof}

The ratio $\rho_g=V_{g-2,3}/V_g$ sets the scale of the good-support
intensity. By \eqref{f:one-pants-volume-ratio},

\begin{equation}\label{mk1:eq:rho}
 \rho_g=\frac{1+O(g^{-1})}{8\pi^2g}.
\end{equation}
Write $\psi(t)=\sinh(t/2)/(t/2)$ for $t>0$, with $\psi(0)=1$.
The boundary-profile estimate \eqref{a:boundary-profile} at genus
$g-2\ge1$ with three boundary components gives
\begin{equation}\label{mk1:eq:q}
 Q_g(x,y,z)=
 \frac{V_{g-2,3}(x,y,z)}
 {V_{g-2,3}\,\psi(x)\psi(y)\psi(z)},
 \quad 0<Q_g\le1,\quad
 1-Q_g\le\frac{c(3)}{g-2}(x^2+y^2+z^2).
\end{equation}
Since $xyz\,\psi(x)\psi(y)\psi(z)=8\Phi(x,y,z)$, the expected
good-support cuff measure $\eta_{g,L}$ in \eqref{mk1:eq:exact-intensity}
has density
\begin{equation}\label{mk1:eq:density}
 d\eta_{g,L}=4\rho_g Q_g\mathbf1_{D_L}\Phi\,dx\,dy\,dz,
 \qquad 1-\epsilon_{g,L}\le Q_g\le1,
 \qquad \epsilon_{g,L}=\frac{3c(3)L^2}{g-2}.
\end{equation}
The bound on $Q_g$ holds on $D_L$, where each cuff length is at most
$L$. At $L=L_g$ we have $\epsilon_{g,L}=o(1)$.

Let $E_{g,L}$ count figure-eights of length at most $L$ whose supports
are exceptional. For $g\ge3$ and $L\ge2\operatorname{arccosh}3$,
summing the three estimates in \eqref{f:exceptional-supports} gives
\begin{equation}\label{mk1:eq:bad}
 \mathbb E_gE_{g,L}\le C\left(
 \frac{L^2e^{L/2}}g+\frac{L^2e^L}{g^2}+\frac{L^2e^L}{g^3}\right).
\end{equation}
At $L=L_g$, this is
$O_S(g^{-1/2}H^{3/2}+H/g)=o(1)$. Thus exceptional supports make a
vanishing contribution to every bounded test in the first-intensity limit.

\subsection{An exact triangle law}

The cuff profile $\mathbf1_{D_L}\Phi\,dx\,dy\,dz$ occurs in the first-moment calculation of He, Shen, Wu and Xue~\cite[equations~(3.4)--(3.7)]{HSWX}. Anantharaman and Monk give the first-moment integral in the three cuff lengths and the corresponding figure-eight length formula in \cite[equations~(7.1)--(7.2)]{AMI}. They also express the cuff lengths in terms of the two loop lengths and the crossing angle, compute the Jacobian, and describe the admissible domain~\cite[Lemmas~7.6, 7.7 and~7.10]{AMI}. Their formulas imply that, under the normalized one-point cuff profile, the winding fraction $x/(x+y)$ converges to the uniform law on $[0,1]$ as $L\to\infty$.

Our coordinates separate the length constraint from the division
between the two winding cuffs. The identity \eqref{f:pants-length}
involves $x$ and $y$ only through $\cosh(x/2)\cosh(y/2)$. We use this
product and $\cosh(z/2)$ as coordinates, and record the fraction of the
logarithm of the product contributed by the first cuff.
For $L>2\operatorname{arccosh}3$, set
\begin{equation}\label{mk1:eq:K}
 C_L=\cosh(L/2),\qquad K=K_L=\frac{C_L-1}{2}>1.
\end{equation}
Introduce the coordinates
\begin{equation}\label{mk1:eq:coords}
 \begin{gathered}
 p=\cosh(x/2),\quad q=\cosh(y/2),\quad w=\cosh(z/2),\quad r=pq,\\
 a=\frac rK,\qquad b=\frac{w-1}{2K},\qquad
 \tau=\frac{\log p}{\log r}.
 \end{gathered}
\end{equation}
When $r=1$ both winding lengths vanish. This set has profile measure
zero, so we may set $\tau=0$ there. Write
\[
 \mathcal T=\{a,b\ge0:a+b\le1\},\qquad
 \mathcal T_K=\{(a,b)\in\mathcal T:a\ge K^{-1}\}.
\]

\begin{lemma}[Exact density and boundary control]\label{mk1:lem:triangle}
Let $I(L)=\int_{D_L}\Phi\,dx\,dy\,dz$. Let $P_K$ be the law of
$(a,b,\tau)$ under the probability measure with density
$\mathbf1_{D_L}\Phi/I(L)$ in $(x,y,z)$. Its density with respect to
$da\,db\,d\tau$ is
\begin{equation}\label{mk1:eq:profile-density}
 f_K(a,b,\tau)=\frac{\log(Ka)}{J_K}\,
 \mathbf1_{\mathcal T_K}(a,b)\mathbf1_{[0,1]}(\tau),
 \qquad
 J_K=\frac12\log K-\frac34+\frac1K-\frac1{4K^2}.
\end{equation}
Thus $\tau$ is uniform on $[0,1]$ and independent of $(a,b)$ for
every $K>1$. The unnormalized mass is
\begin{equation}\label{mk1:eq:mass}
 I(L)=16K^2J_K=8K^2\log K-12K^2+16K-4
 =\frac{(L-3-4\log2)e^L}{4}+O(Le^{L/2}).
\end{equation}
Let $\mathsf U$ be the product of the uniform probability measures on
$\mathcal T$ and $[0,1]$. Its joint density is
$2\mathbf1_{\mathcal T}(a,b)\mathbf1_{[0,1]}(\tau)$. Then
\begin{equation}\label{mk1:eq:triangle-tv}
 d_{\mathrm{TV}}(f_K\,da\,db\,d\tau,\mathsf U)
 \le\frac{3/2-2/K+1/(2K^2)}{\log K}
 \le\frac{3}{2\log K}.
\end{equation}
For $0<\delta,\varepsilon<1/2$, the $P_K$-probability that
$a<\delta$, $b<\delta$, or $\tau\notin[\varepsilon,1-\varepsilon]$
is at most
\begin{equation}\label{mk1:eq:layers}
 4\delta+2\varepsilon+\frac{3}{2\log K}.
\end{equation}
Also, $P_K(a+b<\delta)\le\delta^2+3/(2\log K)$.
\end{lemma}

\begin{proof}
In $(p,q,w)$ coordinates, $\Phi\,dx\,dy\,dz=8\,dp\,dq\,dw$ and
the length constraint is $2pq+w\le C_L$. For $r>1$, the inverse map
$(r,\tau)\mapsto(p,q)=(r^\tau,r^{1-\tau})$ has absolute Jacobian
$\log r$. Substituting $r=Ka$ and $w=1+2Kb$ gives
$dr\,dw=2K^2\,da\,db$. The inequalities $r\ge1$, $w\ge1$ and
$2r+w\le C_L$ give precisely $(a,b)\in\mathcal T_K$.
The transformed profile measure is therefore
\[
 16K^2\log(Ka)\mathbf1_{\mathcal T_K}(a,b)
 \mathbf1_{[0,1]}(\tau)\,da\,db\,d\tau.
\]
Integrating first in $\tau$ and $b$ leaves
$J_K=\int_{1/K}^1(1-a)\log(Ka)\,da$ after removing the factor
$16K^2$. Evaluation gives \eqref{mk1:eq:profile-density} and the exact
mass in \eqref{mk1:eq:mass}. The asymptotic follows on substituting
$K=e^{L/2}/4-1/2+e^{-L/2}/4$ and letting $L\to\infty$.

To compare $P_K$ with $\mathsf U$, accept a sample from $\mathsf U$
with probability $\log(Ka)/\log K$ on $\mathcal T_K$ and zero
elsewhere. The acceptance probability is $c_K=2J_K/\log K$.
Conditional on acceptance, the sample has density $f_K$. Thus
$\mathsf U=c_KP_K+(1-c_K)\mathsf R_K$, where $\mathsf R_K$ is the
conditional law given rejection. It follows that
$d_{\mathrm{TV}}(P_K,\mathsf U)\le1-c_K$, which gives
\eqref{mk1:eq:triangle-tv}.

Under $\mathsf U$, each of the events $a<\delta$ and $b<\delta$
has probability $2\delta-\delta^2$. The two end intervals for $\tau$
have combined probability $2\varepsilon$, and the smaller triangle
$a+b<\delta$ has probability $\delta^2$. Applying the union bound
under $\mathsf U$ and then \eqref{mk1:eq:triangle-tv} proves the
stated estimates under $P_K$.
\end{proof}

\subsection{Transfer to the expected curve measure}

Two errors separate the triangle calculation from the full expected cuff
measure. On good supports the volume factor $Q_g$ perturbs the profile.
The remaining curves contribute the exceptional-support mass in
\eqref{mk1:eq:bad}.

\begin{proposition}[Mass and variation comparison]\label{mk1:prop:transfer}
At $L=L_g$, the expected number of figure-eights with good support
tends to $\Lambda(S)$. For all sufficiently large $g$, the total
variation distance between the normalized good-support cuff measure
and the normalized profile measure is at most $\epsilon_{g,L_g}$.
The normalized full and good-support cuff measures have total variation
distance $o(1)$. In the coordinates \eqref{mk1:eq:coords}, both
therefore converge in total variation to $\mathsf U$.

\end{proposition}

\begin{proof}
Equation \eqref{mk1:eq:density} gives
\[
 (1-\epsilon_{g,L})4\rho_g I(L)
 \le \eta_{g,L}(1)\le4\rho_g I(L).
\]
By \eqref{mk1:eq:rho}, \eqref{mk1:eq:mass}, and $e^{L_g}=e^Sg/H$,
\begin{equation}\label{mk1:eq:good-mass}
 \eta_{g,L_g}(1)
 =\frac{(L_g-3-4\log2)e^{L_g}}{8\pi^2g}
 \left(1+O\left(\frac{1+L_g^2}{g}\right)\right)
 +O\left(\frac{L_ge^{L_g/2}}g\right)
 \longrightarrow\Lambda(S).
\end{equation}
Accept a sample from the normalized profile measure with probability
$Q_g$. By \eqref{mk1:eq:density}, the accepted sample has the
normalized good-support law. The rejection probability is at most
$\epsilon_{g,L_g}$, which bounds the total variation distance.

Write the full expected cuff measure as
$\eta^{\rm full}_{g,L}=\eta_{g,L}+\beta_{g,L}$, where $\beta_{g,L}$
counts exceptional occurrences in expectation. At $L=L_g$ its
mass is $o(1)$ by \eqref{mk1:eq:bad}, while the good-support mass
tends to $\Lambda(S)>0$. Under the normalized full law the probability
of selecting an exceptional occurrence is therefore $o(1)$.
Conditioning on a good occurrence gives the normalized good-support
law, so the two laws differ in total variation by at most that
probability. For each $g$, applying the coordinate map
\eqref{mk1:eq:coords} can only decrease both variation distances.
The triangle limit follows from \eqref{mk1:eq:triangle-tv}.

\end{proof}

\subsection{Returning to physical cuff lengths}

The variable $\tau$ gives the division between the winding cuffs in
logarithmic coordinates. To recover the fraction of the actual cuff
length $x+y$ contributed by $x$, we control the error in replacing
the inverse hyperbolic cosine by a logarithm. For $v\ge1$, write
\begin{equation}\label{mk1:eq:kappa}
 \operatorname{arccosh}v=\log v+\kappa(v),\qquad
 \kappa(v)=\log(1+\sqrt{1-v^{-2}}),\qquad
 0\le\kappa(v)\le\log2.
\end{equation}
In the coordinates of \eqref{mk1:eq:coords}, the inverse formulas are
\begin{align}
 x&=2\tau\log(Ka)+2\kappa((Ka)^\tau),\notag\\
 y&=2(1-\tau)\log(Ka)+2\kappa((Ka)^{1-\tau}),\label{mk1:eq:inverse}\\
 z&=2\log(1+2Kb)+2\kappa(1+2Kb).\notag
\end{align}
Fix $0<\delta<1/2$. On $a,b\ge\delta$, we have
$\log(Ka)=\log K+O_\delta(1)$ and
$\log(1+2Kb)=\log K+O_\delta(1)$. Since $L=2\log K+O(1)$,
it follows that, uniformly for $a,b\ge\delta$ and $0\le\tau\le1$,
\begin{equation}\label{mk1:eq:physical-shape}
 \frac{(x,y,z)}L\longrightarrow(\tau,1-\tau,1),\qquad
 \frac{(x,y,z)}{x+y+z}
 \longrightarrow\left(\frac\tau2,\frac{1-\tau}2,\frac12\right).
\end{equation}
The winding ratio is also controlled up to both endpoints. With $u=x+y$
and $r>1$, \eqref{mk1:eq:inverse} gives
\begin{equation}\label{mk1:eq:physical-t}
 \left|\frac{x}{u}-\tau\right|
 =\left|\frac{(1-\tau)\kappa(p)-\tau\kappa(q)}
 {\log r+\kappa(p)+\kappa(q)}\right|
 \le\frac{\log2}{\log r}.
\end{equation}

The figure-eight length depends only on $a+b$. In fact,
\begin{equation}\label{mk1:eq:radial-length}
 \ell=2\operatorname{arccosh}(1+2K(a+b)),\qquad
 L=2\operatorname{arccosh}(1+2K).
\end{equation}
Since $\kappa(v)\to\log2$ as $v\to\infty$, subtraction gives
$\ell-L\to2\log(a+b)$ uniformly on $a+b\ge\delta$.

\begin{proof}[Proof of Theorem~\ref{mk1:thm:m1}]
Set $L=L_g$ and let $\mu_g$ be the law of $(a,b,\tau)$ under the
normalized full expected cuff measure. Lemma~\ref{mk1:lem:triangle}
and Proposition~\ref{mk1:prop:transfer} give
$d_{\mathrm{TV}}(\mu_g,\mathsf U)\to0$.
On $a,b\ge\delta$, \eqref{mk1:eq:physical-shape} and
\eqref{mk1:eq:radial-length} give uniform convergence of the centered
length and shape to
\[
 \left(S+2\log(a+b),
       \left(\frac{\tau}{2},\frac{1-\tau}{2},\frac12\right)\right).
\]
For large $g$, the prelimit and limiting pairs on this region lie in a
common compact subset of $(-\infty,S]\times\Delta$. Hence composing a bounded continuous test $\Psi$ with the prelimit
map gives uniform convergence on this region to its composition with
the limiting map. Total-variation convergence allows us to
integrate the limiting test against $\mathsf U$. The omitted region
has $\mathsf U$-probability at most $4\delta$ and $\mu_g$-probability
at most $4\delta+o(1)$. Taking $g\to\infty$ and then
$\delta\downarrow0$ therefore gives joint weak convergence of the
normalized measures.

To identify the limit, put $R=a+b$ and $W=a/(a+b)$ under $\mathsf U$.
The inverse map $(r,w)\mapsto(rw,r(1-w))$ has absolute Jacobian $r$.
Hence $R$ has density $2r$ on $(0,1)$, and $W,\tau$ are independent
uniform variables on $[0,1]$, independent also of $R$. The centered
length $S+2\log R$ has density $e^{s-S}\mathbf1_{\{s\le S\}}$ and
is independent of the limiting shape. Multiplying by
$m_{g,S}\to\Lambda(S)$ gives \eqref{mk1:eq:joint-limit}.
Taking $\Psi(s,Z)=F(Z)$ gives \eqref{mk1:eq:shape-limit}.

Finally, with the cuff threshold $10\log L_g$, the full expected
cuff measure exceeds the retained measure by a positive measure of
mass $\mathbb E_gR_{g,S}=o(1)$ by \eqref{f:trimming-critical}.
The retained mass therefore also tends to $\Lambda(S)>0$. For large
$g$, the normalized full measure is a mixture in which the deleted
curves have weight $\mathbb E_gR_{g,S}/m_{g,S}$. Its total-variation
distance from the normalized retained measure is at most this weight,
which tends to zero. This bound holds after every measurable
coordinate map, including those that depend on $L_g$. It proves
the retained limits with their stated normalizations.
\end{proof}

\subsection{Triangle and finite-offset conclusions}

Define the exponential coordinates by
\begin{equation}\label{mk1:eq:physical-coordinates}
 A_L=e^{(x+y-L)/2},\qquad B_L=e^{(z-L)/2},\qquad t_L=\frac{x}{x+y}.
\end{equation}
They record the offsets of $x+y$ and $z$ from the length cutoff, together
with the winding ratio. For the normalized one-point cuff profile, the joint
limit in \eqref{mk1:eq:offset-limit} also follows from the cuff formulas
and Jacobian of Anantharaman and Monk~\cite[Lemmas~7.6--7.7]{AMI}.
The comparison in triangle coordinates transfers this law to the full
and retained expected cuff measures at $L=L_g$.

\begin{corollary}[Physical triangle and offsets]\label{mk1:cor:offset}
Under the normalized full expected cuff measure for $\ell\le L_g$,
$(A_{L_g},B_{L_g},t_{L_g})$ converges weakly to $\mathsf U$, the
product of the uniform laws on $\mathcal T$ and $[0,1]$. Moreover,
\begin{equation}\label{mk1:eq:offset-limit}
 \left(\ell-L_g,\ x+y-\ell,\ z-\ell,\ \frac{x}{x+y}\right)
 \Longrightarrow
 \bigl(2\log R,\ 2\log W,\ 2\log(1-W),\ T\bigr),
\end{equation}
where $R$ has density $2r$ on $(0,1)$ and $W,T$ are independent uniform
variables on $[0,1]$, independent of $R$. The normalized retained
measure has the same limits.
\end{corollary}

\begin{proof}
Take $L=L_g$ and fix $0<\delta,\varepsilon<1/2$. On the region
$a,b\ge\delta$, $\varepsilon\le\tau\le1-\varepsilon$, both $p$ and
$q$ are at least $(K\delta)^\varepsilon$ for large $g$, while
$w\ge1+2K\delta$. Thus $\kappa(p),\kappa(q),\kappa(w)\to\log2$
uniformly there. Since $Ke^{-L/2}\to1/4$,
\eqref{mk1:eq:inverse} and \eqref{mk1:eq:physical-t} give, uniformly
on this region,
\[
 A_L=Ke^{-L/2}a\,e^{\kappa(p)+\kappa(q)}\to a,\qquad
 B_L=e^{-L/2}(1+2Kb)e^{\kappa(w)}\to b,\qquad t_L\to\tau.
\]
Here $x+y-\ell=2\log A_L-(\ell-L)$ and
$z-\ell=2\log B_L-(\ell-L)$. Both $A_L$ and $B_L$ are bounded away
from zero for large $g$ on this region. The uniform radial limit in
\eqref{mk1:eq:radial-length} therefore gives the limits of the two
offsets as $2\log W$ and $2\log(1-W)$, respectively.

By \eqref{mk1:eq:layers} and the total-variation comparisons above,
the probability of the omitted region has limsup at most
$4\delta+2\varepsilon$ under either normalized expected measure.
For each bounded continuous test, apply uniform convergence on the
compact region above, then let $\delta,\varepsilon\downarrow0$
after the genus limit.
This gives both weak limits, with the joint law of $(R,W,T)$ computed
in the proof of Theorem~\ref{mk1:thm:m1}.
\end{proof}

\begin{remark}[Scope and the topology of convergence]\label{mk1:rem:scope}
The total-variation estimate \eqref{mk1:eq:triangle-tv} applies in
the exact triangle coordinates. For the physical shape, the limit law
$\nu$ is supported on $\{(u,v,w)\in\Delta:w=1/2\}$. Its preimage in
cuff space is the plane $z=x+y$, which has zero measure under the
good-support density \eqref{mk1:eq:density}. Thus the normalized
good-support shape law, whenever defined, has total-variation distance
one from $\nu$. The corresponding distance for the full shape law
tends to one because the normalized mass of exceptional supports tends
to zero.

The factorization in \eqref{mk1:eq:joint-limit} describes the length
and shape of a single occurrence sampled from the expected counting
measure. The factorial moments in Section~\ref{sec:marked-factorials}
and the process limit in Section~\ref{sec:marked-process} identify the joint limits
for distinct curves and the shape law of the first arrival.
\end{remark}

% END sections/08_marked_intensity.tex

% BEGIN sections/09_marked_factorials.tex
\section{Shape-marked factorial measures}\label{sec:marked-factorials}

The first intensity in Section~\ref{sec:marked-intensity} describes
one curve sampled from the expected count. To obtain the joint law of
counts in several marked cells, we compare the expected factorial measure
of isolated curves with a product measure. At finite genus the factors
are the exact retained one-curve intensities. Averaging the winding
orders lets us use the disjoint-support integral of
Section~\ref{sec:disjoint} for this comparison.
Section~\ref{sec:higher} controls the tuples lost through isolation.

\subsection{Processes and mixed factorial moments}

On a Weil--Petersson random surface \(X\) of genus \(g\), fix the
centered length ceiling \(S\in\mathbb R\). For sufficiently large \(g\), set
\begin{equation}\label{mk2:eq:scales}
 H=\log g,\qquad a_g=H-\log H,\qquad L=L_g=a_g+S,\qquad b=10\log L.
\end{equation}
As in Section~\ref{sec:foundations}, retain the figure-eights of length
at most \(L\) with good support and all three cuffs of length at least
\(b\). A retained curve is isolated if its closed support is disjoint
from those of all other retained curves below this ceiling.

Use the fair order of the two winding cuffs assigned to each curve
in Section~\ref{sec:marked-intensity}. Given \(X\), different curves
receive independent orders. A curve keeps its order whenever it appears
at another ceiling or in another tuple. With \(x_\gamma,y_\gamma\) the
ordered winding lengths and \(z_\gamma\) the third cuff length, put
\[
Z_\gamma=\frac{(x_\gamma,y_\gamma,z_\gamma)}
 {x_\gamma+y_\gamma+z_\gamma},\qquad
 Y_\gamma=(\ell_\gamma-a_g,Z_\gamma).
\]
If the winding lengths are equal, either order gives the same mark
for that curve. Distinct curves with identical marks each contribute
one atom. On \(E_S=(-\infty,S]\times\Delta\), define
\begin{equation}\label{mk2:eq:processes}
 \mathsf N^{\mathrm{ret}}_{g,S}=\sum_{\gamma\ \mathrm{retained}}\delta_{Y_\gamma},
 \qquad
 \mathsf N^{\mathrm{iso}}_{g,S}=\sum_{\gamma\ \mathrm{isolated}}\delta_{Y_\gamma},
 \qquad \mu_{g,S}=\mathbb E_g^\circ\mathsf N^{\mathrm{ret}}_{g,S},
\end{equation}
where \(\mathbb E_g^\circ\) averages over both the surface and its
marks. Each sum is finite for a fixed surface. The measurable marking
of the curve classes is described in Section~\ref{mk3:subsec:marking}.

Let \(\nu\) be the law of \((T/2,(1-T)/2,1/2)\) for \(T\) uniform on
\([0,1]\). The limiting retained intensity is
\begin{equation}\label{mk2:eq:limit-measure}
 \mu_S(ds,dz)=\frac{e^s}{8\pi^2}\mathbf1_{\{s\le S\}}\,ds\,\nu(dz).
\end{equation}
This is a finite measure of mass \(e^S/(8\pi^2)\).
Write \((n)_r=n(n-1)\cdots(n-r+1)\) and \((n)_0=1\).

\begin{theorem}[Marked mixed factorial limits]\label{mk2:thm:mixed}
Fix \(S\in\mathbb R\) and pairwise disjoint cells
\(C_j=(u_j,v_j]\times A_j\subset E_S\), \(1\le j\le d\), where
\(-\infty<u_j<v_j\le S\) and each \(A_j\subset\Delta\) is Borel with
\(\nu(\partial_\Delta A_j)=0\). For every fixed multiindex
\(\mathbf r\in\mathbb N_0^d\),
\begin{equation}\label{mk2:eq:mixed-limit}
 \mathbb E_g^\circ\prod_{j=1}^d
       \bigl(\mathsf N^{\mathrm{iso}}_{g,S}(C_j)\bigr)_{r_j}
 \longrightarrow
 \prod_{j=1}^d
 \left[
  \frac{e^{v_j}-e^{u_j}}{8\pi^2}\,
  \int_0^1\mathbf1_{A_j}
       \left(\frac t2,\frac{1-t}2,\frac12\right)dt
 \right]^{r_j}.
\end{equation}
Each such mixed moment is finite for all sufficiently large \(g\).
We use \(0^0=1\), so the formula includes cells of zero limiting mass
and gives \(1=1\) at total order zero.
\end{theorem}

Disjointness of the product cells ensures that an individual curve can
occupy at most one cell. Their projections onto length or shape may
overlap. The boundary condition on \(A_j\) uses the relative topology
of \(\Delta\). For example, the segment
\(\{(u,v,w)\in\Delta:w=1/2\}\) is its own relative boundary and has
\(\nu\)-mass one.

\subsection{The retained intensity and the geometric estimates}

We express the geometric inputs in cuff coordinates. As in
Section~\ref{sec:disjoint}, bordered volumes \(V_{h,n}(\mathbf x)\)
have labeled boundaries, and \(V_g=V_{g,0}\).

For a cuff triple \(\mathbf x=(x,y,z)>0\), define
\[
 \ell(\mathbf x)=2\operatorname{arccosh}
 \bigl(\cosh(z/2)+2\cosh(x/2)\cosh(y/2)\bigr),\qquad
 q_g(\mathbf x)=\left(\ell(\mathbf x)-a_g,
                 \frac{\mathbf x}{x+y+z}\right).
\]
On the retained domain \(\mathcal D_{L,b}\) from
\eqref{a:length-domain}, every cuff length is at least \(b\) and
\(\ell(\mathbf x)-a_g\le S\). Thus \(q_g\) maps this domain into \(E_S\).
The exact retained cuff measure is
\begin{equation}\label{mk2:eq:eta}
 d\eta_g(\mathbf x)=
 \frac{\mathbf1_{\mathcal D_{L,b}}}{2V_g}
 V_{g-2,3}(x,y,z)\,xyz\,dx\,dy\,dz.
\end{equation}
Lemma~\ref{mk1:lem:integration} and Theorem~\ref{mk1:thm:m1} give
\begin{equation}\label{mk2:eq:first-intensity}
 \mu_{g,S}=(q_g)_*\eta_g,\qquad
 \mu_{g,S}\Longrightarrow\mu_S,\qquad
 \mathfrak m_{g,S}:=\mu_{g,S}(E_S)\longrightarrow\frac{e^S}{8\pi^2}.
\end{equation}
Here convergence is weak convergence of finite measures on \(E_S\).
Exchanging the winding coordinates \(x,y\) preserves the retained
domain \(\mathcal D_{L,b}\) and its cuff measure \(\eta_g\).

Now fix \(m\ge1\). Let \(\mathcal T_m(X)\) be the ordered tuples of
\(m\) retained curves with pairwise disjoint closed supports. Write
\(\mathbf x=(\mathbf x_1,\ldots,\mathbf x_m)\) for their cuff triples.
If a Borel function \(G(\mathbf x)\) is invariant under swapping the
two winding coordinates in each slot separately, \(G(\tau)\) is well
defined on the geometric tuple \(\tau\) before the random orders are
chosen. Proposition~\ref{prop:exact-integration} and the cuff-kernel
extension following \eqref{a:one-pants} give, for bounded such \(G\),
\begin{equation}\label{mk2:eq:exact-disjoint}
 \mathbb E_g\sum_{\tau\in\mathcal T_m}G(\tau)
 =\frac1{2^mV_g}\int_{\mathcal D_{L,b}^{\,m}}
 G(\mathbf x)\,K_{g,m}(\mathbf x)\prod_{i,a}x_{i,a}\,d\mathbf x.
\end{equation}
The kernel \(K_{g,m}\) in \eqref{a:volume-kernel} sums over the
admissible joint complements of Definition~\ref{a:def-admissible}.
The integration formula holds for \(g\ge3\). Its integral is absolutely
finite for bounded signed or complex \(G\), since the cuff lengths are
bounded by \(L\) and the kernel is a finite sum of volume polynomials.

For \(g\ge2m+1\), the term with connected joint complement is the
bordered volume \(V_{g-2m,3m}(\mathbf x)\), evaluated at all \(3m\)
cuff lengths. Denote its contribution to \eqref{mk2:eq:exact-disjoint}
by \(\mathcal C_{g,m}(G)\). Proposition~\ref{prop:factorization}
compares its density with the product of the one-slot densities
in \eqref{mk2:eq:eta}. Integrating \eqref{a:pointwise-density} gives
\begin{equation}\label{mk2:eq:factorization}
 \left|\mathcal C_{g,m}(G)-\eta_g^{\otimes m}(G)\right|
 \le \epsilon_{g,m}\eta_g^{\otimes m}(|G|),
 \qquad \epsilon_{g,m}=C_{m,S}\frac{1+L^2}{g},
\end{equation}
for fixed \(m\) and all sufficiently large \(g\). The proposition
applies whenever \(0\le b\le L\le A\log g\) for fixed \(A>0\).
The parameters in \eqref{mk2:eq:scales} satisfy this condition
eventually with \(A=1\).

For \(m\ge2\), let \(D_{g,m}^{\mathrm{disc}}(1)\) count tuples whose
joint complement is disconnected. Applying
Proposition~\ref{prop:disconnected} to the unit kernel gives
\begin{equation}\label{mk2:eq:disc}
 \mathbb E_g D_{g,m}^{\mathrm{disc}}(1)
 \le \frac{C_m}{g}\left(\frac{Le^L}{g}\right)^m
 =O_{m,S}(g^{-1}).
\end{equation}
At \(m=1\) the disconnected count is zero by the good-support
condition. For a bounded weight \(G\), the absolute disconnected
contribution is bounded by \(\|G\|_\infty\) times the expected
unit count.

Finally, the curves within a tuple \(\tau\in\mathcal T_m(X)\) have
no mutual conflicts. Let \(\mathcal W(\tau)\) be the set of retained
curves outside \(\tau\) whose closed supports meet at least one
selected support. Write
\begin{equation}\label{mk2:eq:Q}
 d(\tau)=|\mathcal W(\tau)|,\qquad
 Q_m=\sum_{\tau\in\mathcal T_m}\mathbf1_{\{d(\tau)>0\}}.
\end{equation}
Theorem~\ref{thm:weighted-loss} proves
\begin{equation}\label{mk2:eq:loss-input}
 \mathbb E_gQ_m\le A_{m,S}(H^{-2}+g^{-1}).
\end{equation}
A witness may lie outside the observed cells. Each affected tuple is
counted once, regardless of its joint complement or its number of
witnesses. Membership in the tuple and witness sets depends only on the
surface, so both sets stay fixed when we average the winding orders.

\subsection{Conditional marking and the quantitative tuple comparison}

Fix an integer \(m\ge1\). For
\(\tau=(\gamma_1,\ldots,\gamma_m)\in\mathcal T_m(X)\), write
\(Y_\tau=(Y_{\gamma_1},\ldots,Y_{\gamma_m})\).
To compute the conditional expectation of \(F(Y_\tau)\) given \(X\),
start with either order of the two winding cuffs in each triple
\(\mathbf x_i\), and let \(\sigma(x,y,z)=(y,x,z)\). For a bounded
Borel function \(F:E_S^m\to\mathbb C\), define on \(\mathcal D_{L,b}^{\,m}\)
\begin{equation}\label{mk2:eq:sym}
 G_{g,F}(\mathbf x_1,\ldots,\mathbf x_m)
 =2^{-m}\sum_{\varepsilon\in\{0,1\}^m}
 F\bigl(q_g(\sigma^{\varepsilon_1}\mathbf x_1),\ldots,
        q_g(\sigma^{\varepsilon_m}\mathbf x_m)\bigr).
\end{equation}
Given \(X\), the \(m\) distinct curves have independent fair winding
orders. The average in \eqref{mk2:eq:sym} is therefore the conditional
expectation of \(F(Y_\tau)\). It is invariant under a swap in any
slot, so its value \(G_{g,F}(\tau)\) is independent of the initial
orders. Also, \(\|G_{g,F}\|_\infty\le\|F\|_\infty\). Applying
\eqref{mk2:eq:factorization} and \eqref{mk2:eq:disc} to the unit
kernel gives \(\mathbb E_g|\mathcal T_m(X)|<\infty\) for large \(g\).
This justifies taking expectations of the bounded real or complex
weights term by term.

\begin{lemma}[Product integral and multiplicities]\label{mk2:lem:average}
For every such \(F\),
\begin{equation}\label{mk2:eq:sym-integral}
 \eta_g^{\otimes m}(G_{g,F})=\mu_{g,S}^{\otimes m}(F),
 \qquad
 \eta_g^{\otimes m}(|G_{g,F}|)
       \le\mu_{g,S}^{\otimes m}(|F|).
\end{equation}
The expected marked sum
\(\mathbb E_g^\circ\sum_{\tau\in\mathcal T_m(X)}F(Y_\tau)\)
is given by the integral in \eqref{mk2:eq:exact-disjoint} with
\(G=G_{g,F}\).
\end{lemma}
\begin{proof}
Each winding swap preserves \(\eta_g\). Hence every one of the
\(2^m\) terms in \eqref{mk2:eq:sym} has integral
\(\mu_{g,S}^{\otimes m}(F)\) by the pushforward identity
\eqref{mk2:eq:first-intensity}. Averaging gives the first equality.
The pointwise inequality \(|G_{g,F}|\le G_{g,|F|}\) gives the bound.

Conditioning the marked sum on \(X\) gives the sum of \(G_{g,F}\)
over \(\mathcal T_m(X)\), to which \eqref{mk2:eq:exact-disjoint}
applies. The two factors \(2^{-m}\) have different roles.
In \eqref{mk2:eq:sym} each joint choice of winding orders has probability
\(2^{-m}\). In \eqref{mk2:eq:exact-disjoint} the factor divides out
the \(2^m\) cuff labelings of each ordered curve tuple. Geometric
labelings retain this multiplicity even when cuff lengths coincide.
Taking \(F=1\) recovers the unmarked integral.
\end{proof}

For the isolated marked process, define its expected factorial measure by
\begin{equation}\label{mk2:eq:alpha}
 \alpha^{\mathrm{iso}}_{g,m}(F)=
 \mathbb E_g^\circ
 \sum_{\substack{(\gamma_1,\ldots,\gamma_m)\\
                  \text{distinct isolated retained curves}}}
 F(Y_{\gamma_1},\ldots,Y_{\gamma_m}).
\end{equation}
This sum runs over ordered tuples of distinct curves. Two distinct
curves may have the same marked value. Since any two isolated curves
have disjoint closed supports, every tuple in \eqref{mk2:eq:alpha}
belongs to \(\mathcal T_m(X)\).

\begin{proposition}[Uniform finite-genus factorial comparison]
\label{mk2:prop:comparison}
For each fixed integer \(m\ge1\) and \(S\in\mathbb R\), all sufficiently
large \(g\), and every bounded Borel \(F:E_S^m\to\mathbb C\),
\begin{align}
 \left|\alpha^{\mathrm{iso}}_{g,m}(F)-\mu_{g,S}^{\otimes m}(F)\right|
 &\le C_{m,S}\frac{1+L^2}{g}\,
                      \mu_{g,S}^{\otimes m}(|F|)
 \notag\\[-1mm]
 &\quad+C_{m,S}(H^{-2}+g^{-1})\|F\|_\infty.
 \label{mk2:eq:master}
\end{align}
For a finite signed measure \(\zeta\) on \(E_S^m\), let
\(\|\zeta\|_{\mathrm{var}}=\sup_{\|F\|_\infty\le1}|\zeta(F)|\),
with the supremum over Borel tests. Then
\begin{equation}\label{mk2:eq:finite-tv}
 \|\alpha^{\mathrm{iso}}_{g,m}-\mu_{g,S}^{\otimes m}\|_{\mathrm{var}}
 =O_{m,S}\left(H^{-2}+\frac{1+L^2}{g}\right).
\end{equation}
\end{proposition}
\begin{proof}
First sum over all \(\mathcal T_m(X)\). By \eqref{mk2:eq:factorization}
applied to \(G_{g,F}\) and the identities in \eqref{mk2:eq:sym-integral},
the expected contribution from tuples with connected joint complement
differs from \(\mu_{g,S}^{\otimes m}(F)\) by at most the first term in
\eqref{mk2:eq:master}. The absolute contribution
from tuples with disconnected joint complement is at most
\(\|F\|_\infty\mathbb E_gD^{\mathrm{disc}}_{g,m}(1)
=O_{m,S}(g^{-1}\|F\|_\infty)\) by \eqref{mk2:eq:disc}.

It remains to remove tuples with a nonisolated curve. Conditional
averaging gives the removed weight
\begin{equation}\label{mk2:eq:loss-identity}
 D_m(G_{g,F})-D_m^{\rm iso}(G_{g,F})
 =\sum_{\tau\in\mathcal T_m}G_{g,F}(\tau)
                  \mathbf1_{\{d(\tau)>0\}}.
\end{equation}
The measurable tuple weight $G_{g,F}$ satisfies
$\|G_{g,F}\|_\infty\le\|F\|_\infty$. Corollary~\ref{cor:bounded-weight-loss}
therefore bounds the expected absolute value of
\eqref{mk2:eq:loss-identity} by
$C_{m,S}(H^{-2}+g^{-1})\|F\|_\infty$. This gives the second error
bound in \eqref{mk2:eq:master}, after absorbing the disconnected
contribution into the constant. Finally,
\(\mu_{g,S}^{\otimes m}(|F|)\le\mathfrak m_{g,S}^m\|F\|_\infty\)
and \(\mathfrak m_{g,S}=O_S(1)\) by \eqref{mk2:eq:first-intensity}.
Taking the supremum over \(\|F\|_\infty\le1\) gives
\eqref{mk2:eq:finite-tv}.
\end{proof}

\subsection{From tuple integrals to the mixed factorial limits}

\begin{proof}[Proof of Theorem~\ref{mk2:thm:mixed}]
For \(m=\sum_jr_j\ge1\), assign \(r_j\) labeled slots to \(C_j\) in
successive blocks. Let \(j(i)\) be the index of the cell assigned to
slot \(i\), and set
\[
 F_{\mathbf r}(y_1,\ldots,y_m)
       =\prod_{i=1}^m\mathbf1_{C_{j(i)}}(y_i).
\]
Fix a surface and a realization of its marks. The falling factorial
for \(C_j\) counts ordered lists of \(r_j\) distinct isolated curves
in that cell. Because the cells are pairwise disjoint, the lists for
different cells contain distinct curves as well. Concatenating the lists
in slot order gives a bijection with the tuples counted on the right of
\begin{equation}\label{mk2:eq:pointwise-factorial}
 \prod_j\bigl(\mathsf N^{\mathrm{iso}}_{g,S}(C_j)\bigr)_{r_j}
 =
 \sum_{\substack{(\gamma_1,\ldots,\gamma_m)\\
                 \text{distinct isolated retained curves}}}
 F_{\mathbf r}(Y_{\gamma_1},\ldots,Y_{\gamma_m}).
\end{equation}

It matters that we first selected distinct curves. Suppose the two
winding orders for one curve put it in two different cells. The product
of its two cell indicators is then always zero, although each indicator
has conditional mean \(1/2\). In \eqref{mk2:eq:pointwise-factorial}
each slot corresponds to a different curve, so conditional independence
does give a product after averaging.
\[
 G_{g,F_{\mathbf r}}(\mathbf x_1,\ldots,\mathbf x_m)
   =\prod_i h_{g,j(i)}(\mathbf x_i),\qquad
 h_{g,j}(\mathbf x)
   =\frac{\mathbf1_{C_j}(q_g(\mathbf x))
           +\mathbf1_{C_j}(q_g(\sigma\mathbf x))}{2}.
\]
Here \(h_{g,j}\) is the conditional probability that a curve with cuff
triple \(\mathbf x\) falls in \(C_j\). Taking expectations in
\eqref{mk2:eq:pointwise-factorial} and applying
Proposition~\ref{mk2:prop:comparison} to \(F_{\mathbf r}\) gives
\begin{equation}\label{mk2:eq:finite-product}
 \mathbb E_g^\circ\prod_j\bigl(\mathsf N^{\mathrm{iso}}_{g,S}(C_j)\bigr)_{r_j}
 =\prod_j\mu_{g,S}(C_j)^{r_j}
   +O_{m,S}\left(H^{-2}+\frac{1+L^2}{g}\right).
\end{equation}
The error is uniform over pairwise disjoint Borel cells at fixed \(m,S\),
since \(0\le F_{\mathbf r}\le1\) and the total masses
\(\mathfrak m_{g,S}\) are bounded.

It remains to pass to the limit in the cell masses. Taking the boundary
of \(C_j\) in \(E_S\), and the closure and boundary of \(A_j\) in
\(\Delta\), gives
\[
 \partial C_j\subset
 \bigl(\{u_j,v_j\}\times\overline{A_j}\bigr)
 \ \cup\ \bigl([u_j,v_j]\times\partial_\Delta A_j\bigr).
\]
The length intensity assigns zero mass to the endpoints, and
\(\nu(\partial_\Delta A_j)=0\) by hypothesis. Thus \(C_j\) is a
continuity set for \(\mu_S\). In \eqref{mk2:eq:first-intensity}, the
total masses converge to \(\mu_S(E_S)>0\). Applying Portmanteau to
the normalized measures and restoring these masses gives
\begin{equation}\label{mk2:eq:cell-limit}
 \mu_{g,S}(C_j)\longrightarrow
 \lambda_j:=\mu_S(C_j)
 =\frac{e^{v_j}-e^{u_j}}{8\pi^2}\nu(A_j).
\end{equation}
Substituting these limits in \eqref{mk2:eq:finite-product} proves
\eqref{mk2:eq:mixed-limit}, including cells of zero limiting mass with
the convention \(0^0=1\). For \(m=0\) both sides equal one.
\end{proof}

\subsection{Comparison with the limiting cell masses}

Equation~\eqref{mk2:eq:finite-product} compares factorial moments with
the product of the exact retained cell masses. To compare instead with
the limiting masses, we must also control their individual errors. Put
\(\delta_g=\max_j|\mu_{g,S}(C_j)-\lambda_j|\). Choose \(K\ge1\) that bounds
both \(\mathfrak m_{g,S}\) and \(\mu_S(E_S)\) for all sufficiently large
\(g\). Then every cell mass under either measure is at most \(K\). Replacing the \(m\ge1\) factors one at a time
in \eqref{mk2:eq:finite-product} gives
\begin{equation}\label{mk2:eq:cell-error}
 \left|
 \mathbb E_g^\circ\prod_j\bigl(\mathsf N^{\mathrm{iso}}_{g,S}(C_j)\bigr)_{r_j}
       -\prod_j\lambda_j^{r_j}\right|
 \le C_{m,S}\left(H^{-2}+\frac{1+L^2}{g}\right)
       +mK^{m-1}\delta_g.
\end{equation}
For each fixed family of continuity cells, \eqref{mk2:eq:cell-limit}
gives \(\delta_g\to0\). A rate for this term would require a
quantitative estimate for the chosen cells, including their shape
boundaries.

Section~\ref{sec:marked-process} applies the factorial-moment criterion to
these isolated counts and uses the coupling to pass to all figure-eights.
% END sections/09_marked_factorials.tex

% BEGIN sections/10_marked_process.tex
\section{The marked process and marked arrivals}\label{sec:marked-process}

The mixed factorial moments of Section~\ref{sec:marked-factorials}
give the limiting counts of isolated curves in marked cells. We use
the same winding order for a curve before and after isolation. The
unmarked coupling of Section~\ref{sec:probability} then transfers these
count limits to the full marked process. To read off the first arrivals,
we sort the marks by length in a bounded window and then remove the two
endpoints. Each genus limit is taken at fixed cutoffs, observation
cells and moment orders.

\subsection{Persistent marks and the coupling}\label{mk3:subsec:marking}

On a fixed surface, the set of actual figure-eights is countable, with
only finitely many curves below each length bound. To construct their
marks measurably, enumerate the primitive unoriented free-homotopy
classes with self-intersection number one on a reference surface, and
choose a reference order of each pair of winding cuffs. On Teichm\"uller
space, attach independent fair bits to these classes and form the marked
counting measure from their geodesic lengths and intrinsic cuff lengths.
These length functions are measurable, so the conditional marking law
is a measurable probability kernel. A change of marking permutes the
classes and may exchange their reference winding orders. The product
fair law is invariant under both operations, so this kernel descends to
moduli space. Each curve keeps its bit at every cutoff. Intrinsic
boundary occurrences remain distinct even when their ambient images
coincide, as in Lemma~\ref{lem:f-support-types}.

Recall $H=\log g$ and $a_g=H-\log H$. For fixed $S$ and sufficiently
large $g$, put $L=a_g+S$ and $b=10\log L$. The retained and isolated marked measures
$\mathsf N^{\mathrm{ret}}_{g,S}$ and
$\mathsf N^{\mathrm{iso}}_{g,S}$ inherit the marks already assigned to
the full curve set, with the conditional law used in
Sections~\ref{sec:marked-intensity} and~\ref{sec:marked-factorials}.
Isolation is tested against the entire retained graph below the centered
ceiling $S$, including curves outside the observed cells. We use the
isolated process at each fixed ceiling separately. Pointwise on the
enlarged probability space,
\begin{equation}\label{mk3:eq:nested}
 \mathsf N^{\mathrm{iso}}_{g,S}
 \le\mathsf N^{\mathrm{ret}}_{g,S}
 \le\mathsf N_g|_{(-\infty,S]\times\Delta}.
\end{equation}
Persistent marks make the full-to-isolated mismatch event identical
to the unmarked event in Corollary~\ref{cor:f-isolation-coupling}:
retained atoms keep their marks, and every deletion reduces the total
mass. Hence
\begin{equation}\label{mk3:eq:coupling}
 \begin{split}
 \mathbb P_g^\circ\{\mathsf N_g|_{(-\infty,S]\times\Delta}
                  \ne\mathsf N^{\mathrm{iso}}_{g,S}\}
 &=\varepsilon_g(S)\\
 &\le\mathbb E_gR_{g,S}+\mathbb E_gB_{g,S}
       =O_S(\log H/H).
 \end{split}
\end{equation}
An edge can remove both endpoints, so the number of atoms lost is at
most $R_{g,S}+2B_{g,S}$, as in Corollary~\ref{cor:f-isolation-coupling}.

The projection identity following~\eqref{mk:full} shows that the
unmarked first-moment and tail estimates apply unchanged. In particular,
for each compact $K\subset E$ and each fixed $R>0$,
\begin{equation}\label{mk3:eq:compact-tail}
 \begin{split}
 \sup_{g\ge3}\mathbb E_g^\circ\mathsf N_g(K)&<\infty,\\
 \limsup_{g\to\infty}\mathbb P_g^\circ\{
       \mathsf N_g((-\infty,-R]\times\Delta)>0\}
       &\le\Lambda(-R).
 \end{split}
\end{equation}
The first line follows by projection from~\eqref{p:compact-tightness},
or by comparison with a fixed ceiling in~\eqref{p:critical-first-moment};
the second is~\eqref{p:left-tail}. The parameter $R$ is sent to infinity
only after the genus limit.

\subsection{Fixed count vectors}

Fix a finite family of pairwise disjoint bounded cells
$C_j=(u_j,v_j]\times A_j\subset(-\infty,S]\times\Delta$, where
each $A_j$ is Borel and $\nu(\partial_\Delta A_j)=0$.
Theorem~\ref{mk2:thm:mixed} gives, for every fixed
$\mathbf r\in\mathbb N_0^d$,
\begin{equation}\label{mk3:eq:moment-input}
 \mathbb E_g^\circ\prod_j
       (\mathsf N^{\mathrm{iso}}_{g,S}(C_j))_{r_j}
       \longrightarrow\prod_j\mu(C_j)^{r_j}.
\end{equation}
The factorial-moment criterion, Lemma~\ref{lem:factorial-poisson},
therefore gives convergence in total variation of the isolated count
vector to independent Poisson variables of means $\mu(C_j)$, including
zero means. The needed moment finiteness also follows from the
metric-uniform primitive-count bound recalled after
\eqref{p:isolated-factorial-limit}; adding marks cannot increase that
count. By the coupling inequality, full and isolated count-vector laws
have distance at most $\varepsilon_g(S)$. Together with
\eqref{mk3:eq:coupling}, this proves~\eqref{mk:count-tv} for disjoint
bounded cells.

For a finite overlapping family, partition the length intervals at all
their endpoints and partition $\Delta$ into the Boolean atoms generated
by $A_1,\ldots,A_d$. The relative boundary of each Boolean atom is
contained in $\bigcup_j\partial_\Delta A_j$ and is therefore
$\nu$-null. Taking products with the length subintervals gives finitely many
disjoint continuity cells. Each original count is a deterministic sum
of refined counts, for both the prelimit and target processes.
Total variation contracts under this summation map, which yields the
bounded overlapping case with its prescribed dependence.

For a fixed family that includes left half-lines, let $Y_g,Y$ be the
count vectors for the full and limiting marked processes. Truncate every
cell to $(-R,\infty)\times\Delta$. A change in any coordinate requires
a point in the length projection at or below $-R$. The truncated count
vectors converge in total variation for each fixed $R$, so the argument
of Proposition~\ref{p:halfline-counts}, using~\eqref{p:left-tail}, gives
\begin{equation}\label{mk3:eq:halfline}
 \limsup_{g\to\infty}
 d_{\mathrm{TV}}(\mathcal L_g^\circ(Y_g),\mathcal L(Y))
       \le\Lambda(-R)+1-e^{-\Lambda(-R)}.
\end{equation}
Letting $R\to\infty$ proves~\eqref{mk:count-tv}.

\subsection{A continuity semiring and vague convergence}

The space $E=\mathbb R\times\Delta$ is locally compact, second
countable and metrizable. To apply the same convergence criterion as for unmarked lengths, we
choose cells whose boundaries have zero limiting mass. The third
coordinate under $\nu$ equals $1/2$ almost surely, so we avoid that
level when choosing the shape endpoints. Set
$D=\mathbb Q\setminus\{1/2\}$ and define
\[
 \mathcal A=\{\varnothing\}\cup
 \left\{\Delta\cap\prod_{i=1}^3(a_i,b_i]:
       a_i,b_i\in D,\ a_i<b_i\right\}.
\]
This is a countable semiring in $\Delta$. Every relative boundary is
contained in its coordinate endpoint levels. Under $\nu$, each of the first
two coordinates has the uniform law on $[0,1/2]$, and the third equals
$1/2$ almost surely. Thus $\nu(\partial_\Delta A)=0$ for every
$A\in\mathcal A$.

The product semiring
\begin{equation}\label{mk3:eq:semiring}
 \mathcal J_{\mathrm{mark}}=\{\varnothing\}\cup
 \{(a,b]\times A:a,b\in\mathbb Q,\ a<b,\ A\in\mathcal A\}
\end{equation}
is dissecting and covering in the sense of
\cite[Section~1]{Kallenberg}. To see this, enclose a relatively compact
Borel set in a bounded rational length interval times a coordinate box
containing $\Delta$. Subdivide this product into a finite grid of
arbitrarily small mesh using the indicated dense endpoint sets.
Intersecting the shape boxes with $\Delta$ gives a finite cover by
semiring members of arbitrarily small diameter. Intersections remain in
the semiring, and differences are finite disjoint unions of its members.
Each member has $\mu$-null boundary. Countability then makes all these
boundaries simultaneously $\Pi^{\mathrm{mark}}$-null almost surely.

The bounded count-vector result provides every finite joint count limit
on this semiring. The random-measure convergence criterion
\cite[Theorem~1.1]{Kallenberg}, already used in
Section~\ref{sec:probability}, therefore proves~\eqref{mk:vague}.
The first bound in~\eqref{mk3:eq:compact-tail} also gives local count
tightness as in~\cite[Lemma~1.2]{Kallenberg}. Counting measures with
multiplicities form a closed subspace of the Radon measures, so the
conclusion holds in $\mathcal M_p(E)$.

For nonnegative $f\in C_c(E)$, the map
$\xi\mapsto e^{-\xi(f)}$ is bounded and continuous in the vague
topology. Its limiting expectation is the Poisson Laplace functional
in~\eqref{mk:laplace}. For compactly supported nonnegative step functions, the functional
follows from independent Poisson counts. Uniform approximation within one compact
neighborhood of the support, whose $\mu$-mass is finite, gives the
formula for $f$. This completes the proof of Theorem~\ref{mk:main}.

\subsection{Sorting in a fixed window and the two tails}

Half-line counts locate each arrival length. To recover its accompanying
mark, we follow the entire atom as its position in the length ordering
stabilizes. We first do this on a bounded length window.

\begin{lemma}[Sorting in a fixed window]\label{mk3:lem:sorting}
Fix an integer $k\ge1$ and finite $a<b$. For
$\xi\in\mathcal M_p(E)$, list the first $k$ atoms in
$(a,b]\times\Delta$ in increasing length order, or return $\dagger$
if the window contains fewer than $k$ atoms. Give the target
$E^k\sqcup\{\dagger\}$ the topology in which $\dagger$ is isolated.
This map is continuous at every $\xi$ having no atom on
$\{a,b\}\times\Delta$ and whose atoms in the window, counted with
multiplicity, have pairwise distinct length coordinates. Any fixed
measurable rule may be used to order length ties.
\end{lemma}
\begin{proof}
Let $\xi_n\to\xi$ vaguely, with $\xi$ satisfying the stated hypotheses.
Since the window boundary has zero $\xi$-mass, the restricted measures
converge weakly as finite measures on the compact closure of the window. Their total masses converge to the integer
$N=\xi((a,b]\times\Delta)$ and hence equal $N$ for all sufficiently
large indices. If $N<k$, the value is eventually $\dagger$.
Otherwise choose disjoint small neighborhoods of the $N$ limiting atoms, with $\xi$-null boundaries and pairwise
disjoint length projections. For all large indices each neighborhood
contains exactly one prelimit atom and the remaining part of the window
contains none. Shrinking the neighborhoods gives convergence of their
positions and marks. Their length order is fixed, proving convergence
of the first $k$ in that order. This also proves that prelimit ties do
not affect continuity at a distinct-length limit.
\end{proof}

The length projection of $\Pi^{\mathrm{mark}}$ has the law of $\Pi$.
Section~\ref{sec:probability} gives an infinite increasing enumeration,
with finitely many points below each fixed ceiling. Conditional on the
count in a bounded window, the length coordinates have a continuous
joint density. Applying this observation to the countable family
$(-m,m]\times\Delta$, $m\in\mathbb N$, shows that all length coordinates are
distinct almost surely. For each fixed window, its two endpoints also
carry no atoms almost surely. Lemma~\ref{p:arrival-existence} supplies
the first $k$ actual figure-eights whenever $g\ge k+2$.

Fix $k$, then $R>0$ and $S>-R$. By~\eqref{mk:vague} and
Lemma~\ref{mk3:lem:sorting}, the first $k$ marked atoms in
$(-R,S]\times\Delta$, with the cemetery convention, converge in law.
The global and truncated first-$k$ vectors agree whenever there is no
point at or below $-R$ and at least $k$ points at or below $S$.
By~\eqref{p:left-tail} and the half-line count limit
\eqref{mk:count-tv},
\begin{equation}\label{mk3:eq:two-tails}
 \begin{split}
 \limsup_{g\to\infty}\mathbb P_g^\circ\{
       \text{global and truncated vectors differ}\}
   &\le\Lambda(-R)+p_k(S),\\
 p_k(S)&=e^{-\Lambda(S)}
               \sum_{r=0}^{k-1}\frac{\Lambda(S)^r}{r!}.
 \end{split}
\end{equation}
The corresponding target discrepancy is at most
$1-e^{-\Lambda(-R)}+p_k(S)$. Extend any bounded continuous test on
$E^k$ by zero at the isolated cemetery point. The difference between
its global and truncated expectations is bounded by twice its supremum
norm times the discrepancy probability. First take $g\to\infty$ at
fixed $R,S$, then let $R\to\infty$ and $S\to\infty$. Both errors
vanish, proving convergence of the global marked arrival vector to the
ordered atoms of $\Pi^{\mathrm{mark}}$.

\subsection{Independence and the marked arrival law}

The increasing map $y=\Lambda(s)$ sends $\mu$ to $dy\,\nu(dz)$ on
$(0,\infty)\times\Delta$. For $T>0$, the transformed process has a
Poisson total of mean $T$ in $(0,T]\times\Delta$. By the conditional sampling
argument preceding~\eqref{p:homogeneous-arrival-density}, given a total
of $n$ its points have the law of an unordered sample of $n$ independent
pairs with distribution $T^{-1}dy\,\nu(dz)$.

Conditional on all sampled lengths, the marks have product law
$\nu^{\otimes n}$. Sorting by length permutes the marks and preserves
this conditional law. The length density therefore gives the joint first-$k$
density
\begin{equation}\label{mk3:eq:arrival-density}
 e^{-y_k}\,dy_1\cdots dy_k\,\nu^{\otimes k}(d\mathbf z),
 \qquad 0<y_1<\cdots<y_k.
\end{equation}
Thus the first $k$ transformed arrival lengths have the law of
$(\Gamma_1,\ldots,\Gamma_k)$, independently of their marks, whose law
is $\nu^{\otimes k}$. Every finite selection of coordinates from the
two sequences lies in a common initial segment. The factorization for
all $k$ therefore makes the entire mark sequence independent of the
entire arrival sequence. The preceding subsection now proves
\eqref{mk:arrival-vector}.

The first limiting length $\log(8\pi^2E_1)$ has a continuous law.
For $\nu(\partial_\Delta A)=0$, the joint continuity-set assertion
therefore gives~\eqref{mk:shortest}. Equivalently,
$\log(8\pi^2)-T_{g,1}$ converges to the standard Gumbel law, with its
limiting shape independent of that variable. The map $Q$
in~\eqref{mk:quotient} is continuous, and
$\mathbb P\{\min(U,1-U)\le v\}=2v$ for $0\le v\le1/2$.
The continuous mapping theorem proves the quotient assertion, completing
Theorem~\ref{mk:arrivals}.

\subsection{Failure of process-law total variation}

\begin{proposition}\label{mk3:prop:tv-fails}
For every $g\ge3$, the full marked processes on $E$ satisfy
\begin{equation}\label{mk3:eq:tv-fails}
 d_{\mathrm{TV}}\!\left(\mathcal L_g^\circ(\mathsf N_g),
                     \mathcal L(\Pi^{\mathrm{mark}})\right)=1.
\end{equation}
\end{proposition}
\begin{proof}
Let $A=\{z\in\Delta:z_3<1/2\}$. The event that some atom belongs to
$\mathbb R\times A$ is measurable, since it is the union of
$\{\xi((-m,m]\times A)>0\}$ over positive integers $m$.
By the construction in Lemma~\ref{p:arrival-existence}, every such
surface contains an embedded geodesic pair of pants with three distinct
cuffs. Choose a shortest cuff as the third cuff and take the figure-eight
winding around the other two. Its third cuff accounts for at most a
third of the perimeter, regardless of the winding order. Thus
$\mathsf N_g$ has an atom in $\mathbb R\times A$ on every surface.
Since $\nu(A)=0$, the limiting Poisson process almost surely has no
such atom. The two event probabilities differ by one, which gives
\eqref{mk3:eq:tv-fails} by~\eqref{p:tv-convention}.
\end{proof}

The set $\{z\in\Delta:z_3<1/2\}$ has relative boundary $z_3=1/2$,
which has $\nu$-mass one. The separating event therefore uses a shape
set outside the continuity class in~\eqref{mk:count-tv}.
Count-vector convergence in total variation holds for each fixed
continuity-cell family, and marked arrivals converge weakly at each
fixed order. Section~\ref{sec:tails} treats the fixed finite moments
of the centered minimum using separate tail estimates.
% END sections/10_marked_process.tex

% BEGIN sections/11_tail_integrability.tex
\section{Exponential tails, moments and the expected minimum}\label{sec:tails}

The limiting law of the shortest figure-eight describes fluctuations
in a fixed window about $a_g$. To pass to moments, we must also control
the contribution from surfaces on which the minimum lies outside that
window. Write
\begin{equation}\label{tail:notation}
 M_g=\ell_{\rm sys}^{\rm ns}-a_g,\qquad
 Y_g=(-M_g)_+,\qquad Z_g=(M_g)_+,\qquad q=\frac1{8\pi^2},
\end{equation}
where $x_+=\max\{x,0\}$, $H=\log g$, and $a_g=H-\log H$ as in
\eqref{i:center}. A small minimum means that at least one short
figure-eight occurs, so the first moment of the full count bounds the
negative tail. A large minimum means that there are no short
figure-eights. We bound the probability of this event through the
variance of a count restricted by a lower cuff cutoff $b$. Keeping $b$
free lets us choose $b=L/5$ when the length window grows with the genus.

By Theorem~\ref{thm:arrivals},
\begin{equation}\label{tail:weak-limit}
 M_g\Longrightarrow M=\log(8\pi^2E),\qquad
 E\sim\operatorname{Exp}(1),\qquad
 \mathbb P(M>s)=e^{-qe^s}.
\end{equation}
Fix a universal constant $D\geq1$ for which
$\ell_{\rm sys}^{\rm ns}\le D\log g$ for every $g\ge3$, as supplied
by \cite[Lemma~2.7]{HSWX}. Since $a_g>0$, this gives $Z_g\le DH$.

\subsection{Uniform negative tails}

Write $N_{\rm f8}(X,L)=0$ for $L<0$. The non-simple systole is the
shortest figure-eight length, as used in Theorem~\ref{thm:arrivals}.
Hence, for $t\geq0$,
\begin{equation}\label{tail:count-event}
 \{M_g\leq-t\}=\{N_{\rm f8}(X,a_g-t)\geq1\}.
\end{equation}
Recall the lower length bound $L_*=2\operatorname{arccosh}3$ from
\eqref{f:pants-length}. Choose an absolute constant $C_0$ for which
the first-moment estimate~\eqref{i:external-mean}, from
\cite[Proposition~3.1]{HSWX}, takes the form
\begin{equation}\label{tail:first-moment}
 \left|\mathbb E_gN_{\rm f8}(X,L)
       -q\frac{(L-c_0)e^L}{g}\right|
 \leq C_0\left(\frac{L^2e^{L/2}}g+\frac{L^4e^L}{g^2}\right),
 \qquad c_0=3+4\log2.
\end{equation}
This estimate holds for all $g\geq3$ and $L\geq L_*$, and the count
includes every primitive unoriented figure-eight of length at most $L$.

\begin{proposition}[Uniform negative-tail bounds]\label{tail:negative}
Define
\begin{equation}\label{tail:constants}
 \begin{split}
 K_g&=q+C_0\left(\frac{16}{e^2\log g}
                       +\frac{(\log g)^3}{g}\right),\\
 K_*&=q+C_0\left(\frac{16}{e^2\log3}+\frac{27}{e^3}\right).
 \end{split}
\end{equation}
Then $K_g\leq K_*$, $K_g\to q$, and, for every $g\geq3$ and
$t\geq0$,
\begin{equation}\label{tail:negative-probability}
 \mathbb P_g(M_g<-t)=\mathbb P_g(Y_g>t)\leq K_ge^{-t}.
\end{equation}
For every $R\geq0$,
\begin{align}
 \mathbb E_g[(-M_g-R)_+]&\leq K_ge^{-R},
       \label{tail:negative-stop-bound}\\
 \mathbb E_g[Y_g\mathbf1_{\{Y_g>R\}}]
       &\leq K_g(R+1)e^{-R}.\label{tail:negative-ui}
\end{align}
In particular $\{Y_g:g\geq3\}$ is uniformly integrable, and
\begin{equation}\label{tail:negative-uniform-tail}
 \lim_{R\to\infty}\sup_{g\geq3}\mathbb E_g[(-M_g-R)_+]=0.
\end{equation}
\end{proposition}

\begin{proof}
Fix $g\geq3$ and $t\geq0$, and put $L=a_g-t$. If $L<L_*$, the
event $\{M_g<-t\}$ is empty. For $L\geq L_*$, the inclusive count
in~\eqref{tail:count-event} bounds this strict tail event, including at
$L=L_*$. Since $H\geq\log3>1$, we have $0<L\leq a_g\leq H$ and
\begin{equation}\label{tail:negative-normalization}
 \frac{e^L}{g}=\frac{e^{-t}}H.
\end{equation}
Markov's inequality, \eqref{tail:count-event}, and
\eqref{tail:first-moment} give
\begin{align}
 \mathbb P_g(M_g<-t)
 &\leq\mathbb E_gN_{\rm f8}(X,L)\notag\\
 &\leq e^{-t}\left[
    q\frac LH+C_0\left(\frac{L^2e^{-L/2}}H+
                              \frac{L^4}{gH}\right)\right].
       \label{tail:negative-normalized-bound}
\end{align}
Dropping the nonpositive contribution $-qc_0e^{-t}/H$ gives the
second inequality even when $L-c_0<0$. The elementary maxima
\begin{equation}\label{tail:maxima}
 \sup_{x\geq0}x^2e^{-x/2}=16e^{-2},\qquad
 \sup_{h\geq0}h^3e^{-h}=27e^{-3}
\end{equation}
are attained at $x=4$ and $h=3$.  Since $L/H\leq1$ and
$L^4/(gH)\leq H^3/g$, these maxima prove
\eqref{tail:negative-probability} and $K_g\leq K_*$.
The convergence $K_g\to q$ follows from $H\to\infty$ and
$H^3e^{-H}\to0$.

For any nonnegative random variable $Y$ and $R\geq0$, Tonelli's
theorem gives
\begin{align*}
 \mathbb E[(Y-R)_+]&=\int_R^\infty\mathbb P(Y>t)\,dt,\\
 \mathbb E[Y\mathbf1_{\{Y>R\}}]
 &=R\mathbb P(Y>R)+\int_R^\infty\mathbb P(Y>t)\,dt.
\end{align*}
For $R\geq0$, $(Y_g-R)_+=(-M_g-R)_+$.  Integrating
\eqref{tail:negative-probability} proves
\eqref{tail:negative-stop-bound}--\eqref{tail:negative-ui}.
Their uniform versions with $K_*$ give uniform integrability and
\eqref{tail:negative-uniform-tail}.
\end{proof}

\subsection{An auxiliary count for the positive tail}

For $L>1$ and $0\le b\le L$, let $N=N_g^{\rm aux}(L,b)$ count the
figure-eights of length at most $L$ on good supports whose three cuffs
all have length at least $b$. Each such support has embedded closure,
three distinct ambient cuffs, and connected complement. If $M_g>t$ and
$L=a_g+t$, then $N=0$ for every choice of $b$. Increasing $b$ helps
control intersections between the counted supports. The lower bound
on $\mathbb E_gN$ below shows how large a cutoff we can use in the
variance estimate.

For fixed real parameters $A\geq1$ and $B\geq2A+8$, consider
\begin{equation}\label{tail:window}
 T_g=A\log H,\qquad L=a_g+t,\quad0\leq t\leq T_g,
 \qquad b=B\log L.
\end{equation}
We will also use estimates valid throughout the larger range
\begin{equation}\label{tail:window-range}
 t\geq0,\qquad H/2\leq L=a_g+t\leq3H/2,
 \qquad 10\log L\leq b<L/4.
\end{equation}
For sufficiently large $g$ depending on $A,B$, \eqref{tail:window}
lies in \eqref{tail:window-range}, since $L\le H+(A-1)\log H$
and $B\log L<L/4$.
For $L=a_g+t$, write
\begin{equation}\label{tail:mean-scale}
 \mu=\mathbb E_gN,\qquad \mathcal A=\frac{Le^L}{g}=\frac LH e^t.
\end{equation}
Here $\mathcal A$ denotes the mean scale, while $A$ controls the window
width in~\eqref{tail:window}.

\begin{lemma}[Uniform mean bounds]\label{tail:mean}
There are absolute constants $c,C>0$ such that, for all sufficiently
large $g$, uniformly in~\eqref{tail:window-range},
\begin{equation}\label{tail:mean-bounds}
 c\mathcal A\leq\mu\leq C\mathcal A.
\end{equation}
\end{lemma}

\begin{proof}
For winding cuffs $x,y$ and third cuff $z$, rewrite
\eqref{f:pants-length} as
\begin{equation}\label{tail:length}
 \cosh\frac{\ell(x,y,z)}2
 =\cosh\frac z2+\cosh\frac{x+y}2+\cosh\frac{x-y}2.
\end{equation}
Lemma~\ref{lem:f-one-point-integral}, applied to the symmetric
indicator of the cuff lower bounds, gives the exact integral
\begin{equation}\label{tail:one-point}
 \mu=\frac1{2V_g}\int_{\substack{\ell(x,y,z)\leq L\\x,y,z\geq b}}
 V_{g-2,3}(x,y,z)xyz\,dx\,dy\,dz.
\end{equation}
The boundary profile~\eqref{a:boundary-profile} and volume
ratio~\eqref{f:one-pants-volume-ratio}, based on \cite{NWX,MZ},
compare the cuff density $V_{g-2,3}(x,y,z)xyz/(2V_g)$ with
$\Phi/g$ on $[0,L]^3$, where
$\Phi=\sinh(x/2)\sinh(y/2)\sinh(z/2)$.
The profile's relative deficit is $O(L^2/g)=o(1)$ uniformly for
$L\leq2H$, so this density lies between $c\Phi/g$ and $C\Phi/g$.
Equation~\eqref{tail:length} implies $x+y\leq L$ and $z\leq L$,
which gives
\[
 \mu\leq\frac Cg
 \left(\int_0^L\int_0^{L-x}e^{(x+y)/2}\,dy\,dx\right)
 \left(\int_0^Le^{z/2}\,dz\right)
 \leq C\frac{Le^L}{g}.
\]

For the lower bound put $d=2\log6$ and $r=x+y$.  Restrict to
\begin{equation}\label{tail:true-box}
 L-d-1\leq r,z\leq L-d,\qquad L/4\leq x\leq L/2.
\end{equation}
For all sufficiently large $L$, we have
$y=r-x\geq L/2-d-1\geq L/4$ and $z\geq L-d-1\geq L/4$.
Together with $x\geq L/4>b$, these inequalities verify all three
cuff lower bounds. Since $|x-y|\leq r$, the box also satisfies
\[
 \cosh(\ell/2)\leq3\cosh((L-d)/2)
 \leq3e^{(L-d)/2}=\tfrac12e^{L/2}\leq\cosh(L/2).
\]
Thus the true figure-eight length is at most $L$ throughout the box.
The change from $(x,y,z)$ to $(x,r,z)$ has Jacobian one. The box has
volume $L/4$, and $r+z\geq2L-2d-2$ gives $\Phi\geq ce^L$ there.
Integrating the lower density bound proves $\mu\geq cLe^L/g$.
All estimates are uniform in~\eqref{tail:window-range}, since
$H/2\le L\le3H/2$ and $b<L/4$.
\end{proof}

Let $D_2=D_{g,2}^{L,b}(1;X)$ count ordered pairs of these curves
whose closed supports are disjoint. Since $L\leq2\log g$,
Proposition~\ref{prop:factorization} applies with $m=2$ and compares
the connected-complement contribution to the square of the exact mass
in~\eqref{tail:one-point}. Proposition~\ref{prop:disconnected}
bounds the remaining contribution by $C\mathcal A^2/g$.
Using Lemma~\ref{tail:mean}, we obtain
\begin{equation}\label{tail:disjoint}
 \mathbb E_gD_2\leq
 \mu^2\left(1+C\frac{1+L^2}{g}\right).
\end{equation}
The leading term is the square of the exact mean at the same genus.
It cancels with $\mu^2$ when we estimate the variance below.

\subsection{Crossing pairs in an embedded four-holed sphere}

All pair counts $Q$ below count ordered pairs of distinct actual
curves from the auxiliary family. We first consider crossing pairs,
whose distinct supports have overlapping interiors. Each actual ordered
pair gives four ordered source pairs, from the two winding orders in
each slot. Thus Lemma~\ref{lem:f-eight-transfer} applies with the
new cuff threshold. The length identity shows that every distinguished
cuff is at most $L$, the sum of the two winding cuffs is at most $L$,
and each support perimeter is at most $2L$.

Let $Q_{04^0}(L,b)$ count crossing pairs whose intrinsic filled core
is an embedded four-holed sphere with four distinct ambient boundary
geodesics. Let $Q_{04^1}(L,b)$ count crossing pairs with intrinsic core
$S_{0,4}$ and one identified pair of boundary copies. We estimate these
two classes separately.

\begin{lemma}[Planar refinement]\label{tail:planar}
For $L$ sufficiently large, $10\log L\leq b\leq L$, and all
sufficiently large $g$, with constants independent of $L,b,g$,
\begin{equation}\label{tail:planar-bound}
 \mathbb E_gQ_{04^0}(L,b)
 \leq\frac C{g^2}\left(L^5e^{2L-b/2}+L^6e^{3L/2}\right).
\end{equation}
At $b=10\log L$ this is $O(e^{2L}/g^2)$. At the auxiliary threshold $b=B\log L$ and on
\eqref{tail:window},
\begin{equation}\label{tail:planar-relative}
 \frac{\mathbb E_gQ_{04^0}}{\mu^2}
 \leq C\left(L^{3-B/2}+L^4e^{-L/2}\right).
\end{equation}
\end{lemma}

\begin{proof}
We use the cuff domains and integrals of
\cite[Lemmas~4.13--4.15 and equations~(4.57)--(4.61)]{HSWX},
keeping track of the lower cuff bound $b$ in each domain.

Label the four outer cuff lengths $x_1,x_2,x_3,x_4$ and the internal
cuff length of the first support $y$. Each pant contains exactly two
outer cuffs of $S_{0,4}$. Thus the two supports together contain four,
three, or two distinct outer cuffs. After relabeling, their
representative domains lie in the following sets, where all coordinates
are nonnegative unless a stronger lower bound is displayed.
\begin{align}
 \mathcal D_1:\quad&x_i,y\leq L,\quad
 x_1+x_2\leq2L-b,\quad x_3+x_4\leq2L-b;
       \label{tail:domain-one}\\
 \mathcal D_2:\quad&b\leq x_1\leq L,\quad x_2,x_3,y\leq L,\quad
 2x_1+x_2+x_3+x_4\leq4L;
       \label{tail:domain-two}\\
 \mathcal D_3:\quad&b\leq x_1,x_2\leq L,\quad y\leq L,\quad
 2x_1+2x_2+x_3+x_4\leq4L.
       \label{tail:domain-three}
\end{align}
In the first case each support has one internal cuff of length at
least $b$, so its two outer cuffs have sum at most $2L-b$.  In the
second case the shared outer cuff has length $x_1\geq b$.  If
$\eta,\xi$ are the two relevant internal cuffs, shortening the
complementary cylinder boundary gives
$x_4\leq\ell(\eta)+\ell(\xi)$; adding the two support-perimeter
inequalities gives \eqref{tail:domain-two}.  In the third case the
shared outer cuffs have $x_1,x_2\geq b$, and shortening gives
$x_3+x_4\leq\ell(\eta)+\ell(\xi)$, yielding
\eqref{tail:domain-three}.

Once a representative is fixed, count the missing internal curve by
the four-holed-sphere McShane estimate.  The relevant outer pairs
$(x_j,x_k)$ are $(x_3,x_4)$, $(x_2,x_4)$, and $(x_3,x_4)$,
respectively.  In the second and third cases, these are the outer cuffs of the pant
complementary to the missing second pant. The chosen pairs agree with
the integration kernels in the cited equations. By
\cite[Theorem~2.13 and Lemma~2.12]{HSWX}, the multiplicity is at most
\begin{equation}\label{tail:mcshane}
 C(1+x_j)\left(1+e^{(L-x_j-x_k)/2}\right).
\end{equation}
There are only a fixed number of relabelings.

Cutting the four outer cuffs and $y$ leaves two pants and the
complement. The sum of zero-boundary volumes over all complementary
boundary partitions, divided by $V_g$, is at most $C/g^2$ by
\cite[equation~(4.57)]{HSWX}. For sufficiently large $g$, the only
possible genus-zero complementary component is $S_{0,3}$ in the
partition $3+1$. For this component, $V_{0,3}=1$ and
$x\leq2\sinh(x/2)$ give the required cuff-weighted bound directly.
Apply the upper boundary profiles to the positive-genus components,
as in \cite[equations~(4.59)--(4.61)]{HSWX}. Each integrated
representative is then bounded by
\begin{equation}\label{tail:planar-kernel}
 \frac C{g^2}\int_{\mathcal D_i}
 (1+x_j)\left(1+e^{(L-x_j-x_k)/2}\right)
 e^{(x_1+x_2+x_3+x_4)/2}
 y\,dx_1\,dx_2\,dx_3\,dx_4\,dy.
\end{equation}
The exterior profiles involve only the outer cuffs. The internal cuff
contributes $y\,dy$. Since $y\leq L$ and all outer coordinates are
at most $4L$, integrating $y$ and bounding $1+x_j$ gives a factor
at most $CL^3$.

Put $s=x_1+x_2+x_3+x_4$. The columns below bound the contributions
of the two terms in the McShane estimate~\eqref{tail:mcshane}. In
each row the integral is over the projection of the indicated domain
onto the four outer cuff lengths.
\begin{equation}\label{tail:planar-integrals}
 \begin{array}{c|cc}
 &\displaystyle\int e^{s/2}
 &\displaystyle\int e^{(s+L-x_j-x_k)/2}\\ \hline
 \mathcal D_1&CL^2e^{2L-b}&CL^3e^{3L/2-b/2}\\
 \mathcal D_2&CL^2e^{2L-b/2}&CL^3e^{3L/2}\\
 \mathcal D_3&CL^3e^{2L-b}&CL^3e^{3L/2}.
 \end{array}
\end{equation}
For the first row set
$u=x_1+x_2$, $v=x_3+x_4$, and enlarge to
$0\leq u,v\leq2L-b$.  Their slices have lengths at most $u,v$.
Thus the first integral is at most
$(\int_0^{2L-b}u e^{u/2}\,du)^2\leq CL^2e^{2L-b}$.
In the second integral the extra exponential cancels $v$ from the
exponent, leaving $e^{L/2}e^{u/2}$; the $v$ slice then gives a
factor $O(L^2)$, proving its displayed bound.

On $\mathcal D_2$, integrating $e^{s/2}$ in $x_4$ up to
$4L-2x_1-x_2-x_3$ gives at most $2e^{2L-x_1/2}$. Since
$x_2,x_3\le L$, the first integral is bounded by
$2L^2e^{2L}\int_b^L e^{-x_1/2}\,dx_1\le 4L^2e^{2L-b/2}$.
On $\mathcal D_3$, we have $s\leq4L-2b$. A simplex slice at
fixed $s$ has volume $O(s^3)$, giving the third-row bound.
For the second column in these two rows, the exponents become
$L/2+(x_1+x_3)/2$ and $L/2+(x_1+x_2)/2$, respectively.  Each of
these two coordinates is at most $L$, and their exponential integral
is at most $Ce^L$.  The remaining two coordinates have domain area
$O(L^2)$.  This yields even $CL^2e^{3L/2}$, which implies the
displayed estimates for $L>1$.

Multiplying \eqref{tail:planar-integrals} by $CL^3/g^2$ and summing
gives the upper bound
\[
 \frac C{g^2}\left[
 L^5e^{2L-b}+L^6e^{3L/2-b/2}
 +L^5e^{2L-b/2}+L^6e^{3L/2}+L^6e^{2L-b}\right].
\]
Since $b\ge 10\log L$ implies $Le^{-b/2}\le1$, the last term is at most
$L^5e^{2L-b/2}$. This proves
\eqref{tail:planar-bound}.  Divide by
$\mu^2\geq cL^2e^{2L}/g^2$ from Lemma~\ref{tail:mean} and use
$e^{-b/2}=L^{-B/2}$ to obtain \eqref{tail:planar-relative}.
\end{proof}

\subsection{Completion and the remaining small cores}

Let $Q_{12}(L,b)$ count crossing pairs with intrinsic filled core
$S_{1,2}$. Gluing the identified boundary copies of a core counted by
$Q_{04^1}$ gives an embedded $S_{1,2}$. We can therefore estimate
$Q_{12}$ and $Q_{04^1}$ together by counting pairs of supporting
pants in this common completion.

\begin{lemma}[Completion and witnesses]\label{tail:witness}
For either of these two pair types, let $\widehat Y\cong S_{1,2}$ be
the embedded completed image, with true outer cuffs $\gamma_1,\gamma_2$
of lengths $x_1,x_2$.
Each support contains exactly one true outer cuff. After relabeling cuffs,
every pair is covered by one of the following witness families.
\begin{enumerate}[label=(\roman*),leftmargin=*]
\item The pants contain different true outer cuffs. The first pant contains $\gamma_1$ and internal
cuffs of lengths $y_1,y_2$. Cutting them leaves two pants. The witness lengths all lie in $[0,L]$.
\item Both pants contain the same true outer cuff $\gamma_1$, and
neither contains $\gamma_2$. The first pant and its complementary
pant again have internal cuffs $y_1,y_2$. The witness domain is
\begin{equation}\label{tail:witness-domain}
 b\le x_1\le L,\qquad 0\le x_2\le4L-2x_1,\qquad
 0\le y_1,y_2\le L.
\end{equation}
\end{enumerate}
In both cases the McShane count in $\widehat Y$ bounds the number
of possible second pants by $C K_D(x_2)$, where
\begin{equation}\label{tail:witness-kernels}
 K_D(x)=(1+x)(1+e^{L-x/2}).
\end{equation}
\end{lemma}
\begin{proof}
An embedded geodesic pant in $S_{1,2}$ with distinct cuffs must meet
its outer boundary: otherwise it would have three distinct essential
nonperipheral disjoint curves in a surface of pants complexity two.
If a support contained both outer cuffs, its third cuff would cut off
a one-holed torus in $\widehat Y$. Removing the support would separate
this torus from the nonempty exterior of $\widehat Y$ in $X$, contrary
to the connected-complement hypothesis. The two cases follow according
to whether the outer cuffs of the supports differ or agree.

For the identified case, cut the seam to recover the intrinsic $S_{0,4}$,
as in the description of $\mathcal C_{0,4}^1$ in \cite[\S4.3.3, paragraph
before Proposition~4.12]{HSWX}. Each pant has exactly one interior cuff and two
outer cuffs. There is at most one interior cuff because distinct
essential interior curves of a four-holed sphere cannot be disjoint.
There is at least one because a pant with all three cuffs on the outer
boundary would be a connected component of the four-holed sphere.
The auxiliary family requires distinct ambient cuffs, so a pant cannot
contain both copies of the identified seam. Gluing these copies
therefore leaves embedded pants in $\widehat Y$, and the classification
applies. A second identification would produce a closed genus-two
component, which cannot lie properly in a connected closed surface of
genus $g\ge3$.

In (ii), write the other cuffs of the two supports as
$\alpha_i,\beta_i$ for $i=1,2$. The true outer cuff $\gamma_2$ is a
boundary of the filled intrinsic core and is not a support cuff.
The boundary-shortening construction supplies a piecewise geodesic
representative made of distinct boundary arcs of those four cuffs.
Consequently
\[
 x_2\le\sum_{i=1}^2(\ell(\alpha_i)+\ell(\beta_i))
       \le4L-2x_1.
\]
The same construction takes place before gluing in the identified case:
$\gamma_2$ is not a glued seam, so gluing the other boundary copies
does not change this inequality. This proves~\eqref{tail:witness-domain}.
It is the domain argument of \cite[Lemma~4.10]{HSWX} with $b$ retained.

Apply the McShane--Mirzakhani bound
\cite[equations~(2.5)--(2.6) and Lemma~2.12]{HSWX} at $\gamma_2$.
In (i) it counts the second pant directly. In (ii) it counts its
complementary pant in $\widehat Y$, which determines the second pant
and has the same two internal cuffs. In both cases these cuffs have
total length at most $2L$, giving~\eqref{tail:witness-kernels}.

All counting takes place in the completed $S_{1,2}$. An actual pair
uniquely determines its supports and intrinsic filled core. Forgetting
the identified seam records only its two supporting pants in the
completion. Bounding all possible second pants therefore gives an
upper bound without a seam multiplicity. Each pant carries at most three actual figure-eights.
The cuff orders and relabelings contribute only an
absolute factor, which also absorbs the source-to-actual counting
convention. The witness bounds therefore apply to both pair types.
\end{proof}

\begin{lemma}[Refined small-core mass]\label{tail:small}
For sufficiently large $g$, $L>1$ and $0\le b\le L$,
\begin{equation}\label{tail:other-small}
 \mathbb E_g(Q_{12}+Q_{04^1})\le\frac C{g^2}
 \left[L^5e^{2L-b/2}+L^6e^{3L/2}+L^5e^L\right].
\end{equation}
For the auxiliary parameters~\eqref{tail:window}, this gives
\begin{equation}\label{tail:small-relative}
 \frac{\mathbb E_g(Q_{12}+Q_{04^1})}{\mu^2}
 \le C\left[L^{3-B/2}+L^4e^{-L/2}+L^3e^{-L}\right].
\end{equation}
\end{lemma}
\begin{proof}
The exterior of $\widehat Y$ is either $S_{g-2,2}$ or
$S_{g_1,1}\sqcup S_{g_2,1}$, where $g_1,g_2\ge1$ and
$g_1+g_2=g-1$. Summing over all these exteriors gives the zero-boundary volume bound
\begin{equation}\label{tail:small-exterior}
 B_g^{\rm ext}:=
 \frac{V_{g-2,2}+\sum_{g_1+g_2=g-1}V_{g_1,1}V_{g_2,1}}{V_g}
 \le Cg^{-2}.
\end{equation}
This is \cite[equation~(4.41)]{HSWX}, based on its Theorems~2.4--2.5.
The boundary upper profiles cancel the two exterior cuff factors
against their denominators and leave at most
$C B_g^{\rm ext} e^{(x_1+x_2)/2}$.
Internal cuffs contribute their Fenchel--Nielsen factors. Automorphism
and side-label factors in the unfolding are bounded absolute constants
for these fixed witness types.

In both cases the two local components are pants, so their volumes
are one. The $y_1y_2$ integral is at most $CL^4$. Since $x_2\le4L$,
the polynomial factor in $K_D(x_2)$ costs at most $CL$. In (i) the
two outer integrals over $[0,L]^2$ are bounded by
$Ce^L+CL e^{3L/2}$, giving
$C(L^5e^L+L^6e^{3L/2})$.
For (ii) the remaining integrals are
\begin{align}
 \int_b^L\!\int_0^{4L-2x_1}e^{(x_1+x_2)/2}\,dx_2\,dx_1
 &\le 4e^{2L-b/2},\label{tail:small-integral-one}\\
 \int_b^L\!\int_0^{4L-2x_1}e^{L+x_1/2}\,dx_2\,dx_1
 &\le 8L e^{3L/2}.\label{tail:small-integral-two}
\end{align}
For the first inequality integrate $x_2$ and then use
$\int_b^L e^{-x_1/2}dx_1\le2e^{-b/2}$.
For the second use $4L-2x_1\le4L$ and
$\int_b^L e^{x_1/2}dx_1\le2e^{L/2}$.
Multiplying by $CL^5B_g^{\rm ext}$ proves~\eqref{tail:other-small}.
Finally, \eqref{tail:mean-bounds} and $e^{-b/2}=L^{-B/2}$ give
\eqref{tail:small-relative}.
\end{proof}

Lemma~\ref{s:identified-pair} uses the actual-pair estimate in
Lemmas~\ref{tail:witness} and~\ref{tail:small} to bound the
ordered cuff triples of \cite[Proposition~4.12]{HSWX}.

\subsection{Touching supports with a variable cuff cutoff}

Let $Q_{\rm touch}(L,b)$ count ordered auxiliary pairs whose distinct
support pants have disjoint interiors and share a boundary cuff.

\begin{lemma}[Symbolic touching bound]\label{tail:touch}
There are absolute constants $C$ and $g_*$ such that for $g\ge g_*$,
$L>1$ and $0\le b\le L$,
\begin{equation}\label{tail:touch-bound}
 \mathbb E_gQ_{\rm touch}(L,b)
 \le \frac{CL^4}{g^2}\left(e^{2L-b}+e^L\right).
\end{equation}
Every connected or disconnected complementary realization is included.
\end{lemma}
\begin{proof}
In an ordered touching pair, take the first curve as the mark and the
second as the witness. The marked support has connected complement by
goodness, so $Q_{\rm touch}(L,b)=I_{T,1}^{\rm conn}(L,b)$.
Apply \eqref{w:touching-bound} with $m=1$. Since
$\Theta(b)=(1+b)e^{-b}\le1$ and $1+b\le 2L$, we obtain
\[
 \mathbb E_gQ_{\rm touch}(L,b)
 \le Cg^{-2}L^2e^{2L}(1+b)e^{-b}
 \le 2Cg^{-2}L^3e^{2L-b}.
\]
As $L>1$, this implies \eqref{tail:touch-bound}.
\end{proof}

\subsection{Larger filled cores on the moving window}

For crossing supports, the geometric argument of
Lemma~\ref{lem:f-pair-core} gives
\begin{equation}\label{tail:core-geometry}
 p(Y)\leq4L,\qquad \operatorname{Area}(Y)\leq4L+4\pi.
\end{equation}
These bounds follow from the two support perimeters, independently of
\cite[Proposition~4.6]{HSWX}, and apply throughout
\eqref{tail:window-range}. In this range the area is $O(H)$, whereas
$\operatorname{Area}(X)=4\pi(g-1)$. Thus the core is proper and has
nonempty boundary for all sufficiently large $g$. The two actual curves fill
the core with sum of lengths at most $2L$.  Intrinsic boundary
circles remain distinct when their ambient geodesic images coincide.

For the estimates below, recall that for fixed $a>0$ and $T=aH$,
the family $\mathcal S_T(X)$ from \eqref{f:subsurface-family}
imposes $p(Y)\leq2T$ and $\operatorname{Area}(Y)\leq4T$.
The estimates \eqref{f:subsurface-tail} and
\eqref{f:subsurface-fixed-type}, from
\cite[Propositions~33--34]{WX}, are
\begin{align}
 \mathbb E_g\sum_{\substack{Y\in\mathcal S_T(X)\\|\chi(Y)|\geq k}}
 e^{-p(Y)/4}
 &\leq C_{a,k}\frac{Te^{7T/2}}{g^k},
       \label{tail:subsurface-tail}\\
 \mathbb E_g\sum_{\substack{Y\in\mathcal S_T(X)\\Y\cong S_{h,n}}}
 e^{-(1-\delta)p(Y)/2}
 &\leq C_{a,\delta,h,n}\frac{T^{66}e^{\delta T}}{g^m},
 \quad m=2h-2+n\leq16.
       \label{tail:subsurface-fixed}
\end{align}
Here $k\geq1$ is a fixed integer. In the second estimate, the bordered
type is fixed with $1\leq m\leq16$, and $\delta>0$ is fixed.
Both estimates include all complementary realizations. The filling-tuple estimate
\eqref{f:filling-tuple-count}, from
\cite[Definition~17 and Theorem~18]{WXPGT}, gives, for fixed $m\geq1$
and $0<\delta<1/2$,
\begin{equation}\label{tail:filling-pairs}
 N_2^{\rm fill}(Y,2L)\leq
 C_{m,\delta}(1+2L)e^{2L-(1-\delta)p(Y)/2}.
\end{equation}
Here $2L$ bounds the sum of the two curve lengths.

\begin{lemma}[Larger cores in the moving window]\label{tail:high}
Let $Q_{\geq3}$ count the crossing auxiliary ordered pairs with
intrinsic core complexity at least three. For all sufficiently large $g$,
uniformly on \eqref{tail:window-range},
\begin{equation}\label{tail:high-bound}
 \frac{\mathbb E_gQ_{\geq3}}{\mu^2}
 \leq C\left(H^{65}g^{-3/4}+e^tg^{-3/2}\right).
\end{equation}
\end{lemma}

\begin{proof}
Take $T=3H$. Since $L\leq3H/2$, for all sufficiently large $g$ we have
$p(Y)\leq6H=2T$ and
$\operatorname{Area}(Y)\leq6H+4\pi\leq12H=4T$, placing every
actual core in the single family $\mathcal S_{3H}(X)$ throughout the
moving window.

For $3\leq m\leq14$, combine \eqref{tail:filling-pairs} and
\eqref{tail:subsurface-fixed} with $\delta=1/12$.  There are only
finitely many types, all within the stated range $m\leq16$.
Since $e^{\delta T}=e^{H/4}$, their sum satisfies
\begin{equation}\label{tail:fixed-high}
 \mathbb E_gQ_{3\leq m\leq14}
 \leq C\frac{Le^{2L}H^{66}e^{H/4}}{g^3}.
\end{equation}
Dividing by $\mu^2\geq cL^2e^{2L}/g^2$ and using $L\asymp H$
gives $CH^{65}g^{-3/4}$.

For $m\geq15$, apply the doubling and retraction argument from
Lemma~\ref{lem:f-high-complexity}. The double has genus $m+1$;
the actual curves retain their lengths, primitivity and distinct
unoriented classes. The metric-uniform bound
\eqref{f:primitive-count}, from \cite[Theorem~2.8]{HSWX}, therefore
bounds the eligible curves by $me^{L+7}$.

For actual cores, \eqref{tail:core-geometry} and Gauss--Bonnet give
$2\pi m\leq4L+4\pi$, hence $m\leq CL$.  The number of ordered
pairs filling a fixed actual core is at most
\begin{equation}\label{tail:double-bound}
 CL^2e^{2L}\leq CL^2e^{3L}e^{-p(Y)/4}.
\end{equation}
The second inequality follows from $p(Y)\leq4L$.
Sum first over actual cores and then enlarge the nonnegative weighted
sum to $\mathcal S_{3H}(X)$. The bounds $m\leq CL$ and $p(Y)\leq4L$
are used before this enlargement. Applying
\eqref{tail:subsurface-tail} with fixed $k=15$ gives
\begin{equation}\label{tail:growing-high}
 \mathbb E_gQ_{m\geq15}
 \leq CL^2e^{3L}\frac{He^{21H/2}}{g^{15}}.
\end{equation}
After division by $\mu^2$, this is at most
$CHe^Lg^{-5/2}=Ce^tg^{-3/2}$, since $e^L=ge^t/H$.
This proves \eqref{tail:high-bound}. The filling-pair constants range
over the finite band $3\leq m\leq14$, and the primitive-count bound
handles all $m\geq15$ uniformly. Both subsurface estimates use the
fixed value $a=3$, so the bound is uniform in \eqref{tail:window-range}.
\end{proof}

\subsection{Variance and the positive tail of the minimum}

The disjoint pairs supply the $\mu^2$ term in the second moment.
After it cancels in the variance, each estimate for intersecting
supports contributes one of the error terms below. We keep the bound
in terms of $b$ so that it applies both to the logarithmic cutoff
in~\eqref{tail:window} and to the choice $b=L/5$ in the moment proof.

\begin{proposition}[Parameterized variance and void bounds]\label{tail:variance}
There are absolute constants $C$ and $g_*$ such that for $g\ge g_*$,
uniformly on \eqref{tail:window-range},
\begin{equation}\label{tail:variance-before}
\begin{split}
 \frac{\operatorname{Var}_gN}{\mu^2}\leq C\bigg[
 &e^{-t}+\frac{1+L^2}{g}
       +L^3e^{-b/2}+L^2e^{-b}\\
 &+L^4e^{-L/2}+L^3e^{-L}+L^2e^{-L}\\
 &+H^{65}g^{-3/4}+e^tg^{-3/2}\bigg].
\end{split}
\end{equation}
In particular, for every fixed $A\geq1$ and $B\geq2A+8$, there are finite constants
$C_{A,B}$ and $g_0(A,B)$ such that, uniformly on
\eqref{tail:window} for $g\geq g_0(A,B)$,
\begin{equation}\label{tail:variance-bound}
 \frac{\operatorname{Var}_gN}{\mu^2}
 \leq C_{A,B}(e^{-t}+H^{-A}).
\end{equation}
Consequently, simultaneously throughout this window,
\begin{equation}\label{tail:parameterized-positive}
 \mathbb P_g(M_g>t)\leq C_{A,B}(e^{-t}+H^{-A}).
\end{equation}
\end{proposition}

\begin{proof}
We first estimate the variance throughout \eqref{tail:window-range}.
Let $Q_{\rm same}$ count ordered pairs of distinct auxiliary curves
with the same support pant. Classifying the remaining pairs by how
their supports meet and by their intrinsic filled cores gives the exact
partition
\[
 N(N-1)=D_2+Q_{\rm same}+Q_{\rm touch}
             +Q_{12}+Q_{04^0}+Q_{04^1}+Q_{\geq3}.
\]
Each pant carries at most three actual figure-eights. If it contributes
$n\leq3$ auxiliary curves, then $n(n-1)\leq2n$. Summing over supports
gives $Q_{\rm same}\leq2N$, so the diagonal and same-support terms
together have expectation at most $3\mu$. Apply \eqref{tail:disjoint},
\eqref{tail:planar-bound}, \eqref{tail:other-small},
\eqref{tail:touch-bound} and Lemma~\ref{tail:high}. Subtract the exact
$\mu^2$ and use $\mu\asymp Le^L/g\asymp e^t$ from
Lemma~\ref{tail:mean} to obtain \eqref{tail:variance-before}.
Now specialize to \eqref{tail:window}. At $b=B\log L$, the two
cutoff terms are $L^{3-B/2}$ and $L^{2-B}$, each
$O_{A,B}(H^{-A})$ since $B\ge 2A+8$ and $L\asymp H$. For the other
terms use $e^t\leq H^A$, \eqref{tail:window-range}, and the fact
that a fixed negative exponential in $H$ dominates every fixed power
of $H$. After increasing the genus threshold, each remaining error
other than $e^{-t}$ is $O_{A,B}(H^{-A})$. This proves
\eqref{tail:variance-bound}.

If $M_g>t$, there is no non-simple geodesic of length at most
$L=a_g+t$, so the auxiliary count is zero. Chebyshev's inequality gives
\[
 \mathbb P_g(M_g>t)\leq\mathbb P_g(N=0)
 \leq\operatorname{Var}_g(N)/\mu^2,
\]
proving \eqref{tail:parameterized-positive} throughout \eqref{tail:window}.
\end{proof}

\subsection{Exponential moments and Wasserstein limits}

\begin{proof}[Proof of Theorem~\ref{thm:all-moments}]
Take $0\le t\le H/20$, $L=a_g+t$ and $b=L/5$. For all sufficiently
large $g$, these parameters lie in \eqref{tail:window-range}, since
$H-\log H\le L\le21H/20$ and $10\log L\le L/5$. In
\eqref{tail:variance-before} the largest cutoff term satisfies
\[
 L^3e^{-L/10}\le C H^{3+1/10}g^{-1/10}=O(g^{-1/20}).
\]
The terms from larger cores are $H^{65}g^{-3/4}$ and
$e^tg^{-3/2}\le g^{-29/20}$. These and all remaining terms other
than $e^{-t}$ are $O(g^{-1/20})$, uniformly over the interval.
Since $g^{-1/20}\le e^{-t}$ there, the variance bound and the
inclusion $\{M_g>t\}\subseteq\{N=0\}$ yield
\begin{equation}\label{tail:linear-void}
 \mathbb P_g(M_g>t)\le \operatorname{Var}_g(N)/\mu^2
 \le Ce^{-t},\qquad 0\le t\le H/20.
\end{equation}

Set $\kappa=1/(20D)$. For $H/20<t\le DH$, monotonicity gives
$\mathbb P_g(M_g>t)\le Ce^{-H/20}\le Ce^{-\kappa t}$.
For $t>DH$ the probability vanishes. Enlarging $C$ for the finitely
many smaller genera, using $Z_g\le DH$, gives
\begin{equation}\label{tail:global-positive}
 \mathbb P_g(M_g>t)\le Ce^{-\kappa t}
 \qquad(g\ge3,\ t\ge0).
\end{equation}
Together with \eqref{tail:negative-probability}, this implies
$\mathbb P_g(|M_g|>t)\le (K_*+C)e^{-\kappa t}$, since $\kappa<1$.
Integrating against $\theta e^{\theta t}\,dt$ gives
\begin{equation}\label{tail:exp-integrability}
 \sup_{g\ge3}\mathbb E_g e^{\theta|M_g|}
 \le 1+\frac{(K_*+C)\theta}{\kappa-\theta},
 \qquad 0<\theta<\kappa.
\end{equation}

For $-1<s<\kappa$, choose $p>1$ with $-1<ps<\kappa$. Integrating the
two one-sided bounds \eqref{tail:negative-probability} and
\eqref{tail:global-positive} separately gives
$\sup_{g\ge3}\mathbb E_g e^{psM_g}<\infty$. Hence $e^{sM_g}$ is
uniformly integrable. By \eqref{tail:weak-limit},
\begin{equation}\label{tail:laplace-limit}
 \lim_{g\to\infty}\mathbb E_g e^{sM_g}
 = (8\pi^2)^s\int_0^\infty u^s e^{-u}\,du
 = (8\pi^2)^s\Gamma(1+s).
\end{equation}
For $s\le-1$ the integral diverges at zero. Applying weak convergence
to the bounded continuous function $x\mapsto\min(e^{sx},K)$ and then
letting $K\to\infty$ gives $\liminf_{g\to\infty}\mathbb E_g e^{sM_g}=+\infty$.

To find the rate of divergence at $s=-1$, we need the probability of
a minimum well below the fixed window about $a_g$. Write
$F_g(L)=\mathbb P_g(\ell_{\rm sys}^{\rm ns}\le L)$. For a fixed $0<u<1$, set
$L=uH$, $b=10\log L$ and $\mathcal A=Le^L/g$. The first moment
in~\eqref{tail:first-moment} is asymptotic to $q\mathcal A$. The
trimming loss in~\eqref{f:trimming-uniform} is $o(\mathcal A)$ since
$b=o(L)$ and $L=u\log g$. Thus the count $N$ of retained curves at
these cutoffs has mean $\mu\sim q\mathcal A\to0$.

Here the mean tends to zero. To use it as a first-order estimate for
the probability of seeing a curve, we show that the contribution from
ordered pairs, $\mathbb E_g[N(N-1)]$, is $o(\mathcal A)$. Since $L\le H$,
Propositions~\ref{prop:factorization} and~\ref{prop:disconnected}
bound the disjoint-support contribution by $O(\mathcal A^2)$.
Pairs on the same support contribute $O(e^{L-b/2}/g)=o(\mathcal A)$
by~\eqref{f:same-pants-bound}. Substituting $e^{-b/2}=L^{-5}$ in the
absolute bounds~\eqref{tail:planar-bound}, \eqref{tail:other-small}
and~\eqref{tail:touch-bound} gives $O(\mathcal A^2)$ for crossing
pairs of complexity two and for touching pairs. For larger cores, the
perimeter and area bounds~\eqref{f:pair-core-geometry} place each
core in $\mathcal S_{2H}(X)$ for all sufficiently large $g$, since
$L\le H$. Apply~\eqref{f:fixed-complexity-pairs} with $\delta=1/8$
and the first inequality in~\eqref{f:growing-complexity-pairs},
counting both orders of each pair. After division by $\mathcal A$, these bounds are
$O(H^{66}g^{u-7/4})$ for complexities three through ten and
$O(H^2g^{2u-3})$ for complexities at least eleven. Both tend to zero.
It follows from
\[
 \mu-\tfrac12\mathbb E_g[N(N-1)]
 \le \mathbb P_g(N\ge1)\le F_g(L)\le\mathbb E_gN_{\rm f8}(X,L)
\]
that, for each fixed $0<u<1$,
\[
 F_g(uH)\sim \frac{quH e^{uH}}{g}\qquad(g\to\infty).
\]

To integrate these probabilities, write $e^{-\ell_{\rm sys}^{\rm ns}}$
as the integral of $e^{-L}$ over $L\ge\ell_{\rm sys}^{\rm ns}$.
Tonelli's theorem and $e^{a_g}=g/H$ give
\[
 \frac{\mathbb E_g e^{-M_g}}{H}
 =\frac g{H^2}\int_0^\infty e^{-L}F_g(L)\,dL.
\]
The contribution from $L\ge H$ is at most $H^{-2}$. For $0\le L\le H$,
\eqref{tail:first-moment} gives $F_g(L)\le CLe^L/g$ uniformly.
Indeed, the count vanishes below $L_*$, and for $L_*\le L\le H$ both
$Le^{-L/2}$ and $L^3/g$ are bounded. Substituting $L=uH$ in the
remaining integral gives an integrand $(g/H)e^{-uH}F_g(uH)$ bounded
by $Cu$ and converging to $qu$ for $0<u<1$. Dominated convergence
gives $\mathbb E_g e^{-M_g}/H\to \int_0^1qu\,du=1/(16\pi^2)$.

For every $r>0$, integrating \eqref{tail:negative-probability} and
\eqref{tail:global-positive} against
$r t^{r-1}\,dt$ gives
\begin{equation}\label{tail:all-bounded}
 \sup_{g\ge3}\mathbb E_g|M_g|^r
 \le (K_*+C\kappa^{-r})\Gamma(1+r)=:C_r.
\end{equation}

Now fix $p>0$ and $R>0$, and apply~\eqref{tail:all-bounded} at $r=p+1$. The inequality
$|x|^p\mathbf1_{\{|x|>R\}}\le |x|^{p+1}/R$ gives
\begin{equation}\label{tail:power-ui}
 \sup_{g\geq3}\mathbb E_g\left[|M_g|^p
                  \mathbf1_{\{|M_g|>R\}}\right]\leq\frac{C_{p+1}}R.
\end{equation}
This gives uniform integrability of $|M_g|^p$ over all genera for
every fixed $p>0$. The limit variable also has every finite absolute
moment, since
\[
 \int_0^\infty |\log(u/q)|^p e^{-u}\,du<\infty.
\]
Near zero, the substitution $u=e^{-v}$ gives integrability. At
infinity, the exponential factor is integrable against every power
of $\log u$.
Approximate $|x|^p$ by bounded continuous truncations. Weak
convergence~\eqref{tail:weak-limit} and~\eqref{tail:power-ui} then give
\begin{equation}\label{tail:absolute-moment-limit}
 \mathbb E_g|M_g|^p\longrightarrow\mathbb E|M|^p
 =\int_0^\infty|\log(u/q)|^p e^{-u}\,du.
\end{equation}
For a continuous test $f$ with the stated polynomial growth, use
continuous cutoffs equal to one on $[-R,R]$ and zero outside
$[-R-1,R+1]$. The truncated expectations converge by weak convergence.
The omitted parts tend to zero uniformly in $g$ by
\eqref{tail:power-ui} and the corresponding tail probability bound.
This proves the assertion for $f$, including all fixed positive integer
moments of $M_g$ and all fixed positive real moments of either part.

Differentiating the Euler gamma integral at one yields
\[
 \mathbb E\log E=-\gamma_{\!E},\qquad
 \mathbb E(\log E)^2=\gamma_{\!E}^2+\frac{\pi^2}{6}.
\]
Since $M=\log(8\pi^2)+\log E$, the first and second moment limits give
\begin{align}
 \mathbb E_gM_g^2&\longrightarrow
  [\log(8\pi^2)-\gamma_{\!E}]^2+\frac{\pi^2}{6},
       \label{tail:second-moment}\\
 \operatorname{Var}_g(M_g)
 =\operatorname{Var}_g(\ell_{\rm sys}^{\rm ns})
 &\longrightarrow\frac{\pi^2}{6}.
       \label{tail:limiting-variance}
\end{align}

Finally fix $p\geq1$ and set $r=p+1$. Write $Q_g,Q$ for generalized
quantiles of $M_g,M$. The limiting CDF
$F(x)=1-e^{-qe^x}$ is continuous and strictly increasing, so
$Q_g(u)\to Q(u)$ for every $0<u<1$, where
\[
 Q(u)=\log[-\log(1-u)/q].
\]
Let $U$ be uniform on $(0,1)$. Then $Q_g(U)$ and $Q(U)$ have the
required marginal laws, their difference tends to zero almost surely,
and
\[
 \mathbb E|Q_g(U)-Q(U)|^r\le2^{r-1}(C_r+\mathbb E|M|^r).
\]
The costs $|Q_g(U)-Q(U)|^p$ therefore have uniformly bounded $r/p>1$
moments and are uniformly integrable. Their expectations tend to zero.
Since this admissible coupling bounds $W_p^p$ from above,
we obtain
\begin{equation}\label{tail:wp-limit}
 W_p\bigl(\operatorname{Law}(M_g),\operatorname{Law}(M)\bigr)
 \longrightarrow0.
\end{equation}
\end{proof}

\subsection{Expected excess above a threshold}

\begin{corollary}[Uniform positive-tail bounds]\label{tail:positive}
There are constants $C_1,C_2$ and $g_0$ such that, for every
$g\geq g_0$ and every $R\geq0$,
\begin{equation}\label{tail:positive-stop-bound}
 \mathbb E_g[(M_g-R)_+]\leq C_1e^{-R}+C_2H g^{-1/20}.
\end{equation}
The family $\{Z_g:g\geq3\}$ is uniformly integrable.  In particular,
\begin{equation}\label{tail:positive-uniform-tail}
 \lim_{R\to\infty}\sup_{g\geq3}\mathbb E_g[(M_g-R)_+]=0.
\end{equation}
\end{corollary}

\begin{proof}
Put $T=H/20$. For $R<T$, integrating~\eqref{tail:linear-void} over
$[R,T]$ contributes at most
$Ce^{-R}$ to $\mathbb E_g[(M_g-R)_+]$. For every $R\ge0$, monotonicity
and $Z_g\le DH$ bound the contribution above $\max(R,T)$ by
\[
 DH\,\mathbb P_g(M_g>T)\le CDH e^{-T}=CDH g^{-1/20}.
\]
This proves \eqref{tail:positive-stop-bound} for every $R\ge0$.
Theorem~\ref{thm:all-moments} gives uniform integrability of $|M_g|$
and hence of $Z_g$. Since $(Z_g-R)_+\le Z_g\mathbf1_{\{Z_g>R\}}$,
\eqref{tail:positive-uniform-tail} follows.
\end{proof}

\subsection{Uniform stop-loss limits and the mean}

\begin{proof}[Proof of Corollary~\ref{thm:expected-minimum}]
Theorem~\ref{thm:all-moments} gives uniform integrability of $|M_g|$
and convergence in $W_1$ to $M$ in \eqref{tail:weak-limit}.
For every $R\ge0$, both $x\mapsto(x-R)_+$ and
$x\mapsto(-x-R)_+$ are $1$-Lipschitz. Applying either function to
any coupling of $M_g$ and $M$, then taking the infimum over couplings,
gives
\[
 \sup_{R\ge0}\left|
 \mathbb E_g[(\pm M_g-R)_+]-\mathbb E[(\pm M-R)_+]\right|
 \le W_1(\mathcal L_g(M_g),\mathcal L(M))\longrightarrow0.
\]
Using the law in \eqref{tail:weak-limit}, Tonelli's theorem and the
substitutions $u=qe^{-t}$ and $u=qe^t$ give
\begin{align*}
 \mathbb E[(-M-R)_+]
 &=\int_R^\infty(1-e^{-qe^{-t}})\,dt
   =\operatorname{Ein}(qe^{-R}),\\
 \mathbb E[(M-R)_+]
 &=\int_R^\infty e^{-qe^t}\,dt
   =\operatorname E_1(qe^R),
\end{align*}
with the functions defined in \eqref{i:exponential-integrals}.
This proves \eqref{i:negative-stop}--\eqref{i:positive-stop}.
Theorem~\ref{thm:all-moments}, applied to $f(x)=x$, also gives
\[
 \mathbb E_gM_g\longrightarrow\mathbb EM
 =\log(8\pi^2)+\int_0^\infty e^{-u}\log u\,du
 =\log(8\pi^2)-\gamma_{\!E}.
\]
Adding $a_g$ proves \eqref{i:expected-minimum}.
\end{proof}

For any fixed number of arrivals, the marked Poisson limit gives the
joint limiting law of the centered figure-eight lengths and their cuff
shapes. For the minimum, the uniform tail estimates turn weak convergence
in~\eqref{tail:weak-limit} into convergence of polynomial moments,
stop-loss functions, and exponential moments for $-1<s<\kappa$, as
well as convergence in every finite-order Wasserstein distance. These
convergence statements are qualitative. The exponent $\kappa$ comes
from the positive-tail bound. At the negative endpoint $s=-1$, the
estimate for short curves at lengths $u\log g$, $0<u<1$, determines
the leading term $\mathbb E_g e^{-M_g}\sim(\log g)/(16\pi^2)$.

% END sections/11_tail_integrability.tex

\appendix

% BEGIN sections/appendix_a_genus.tex
\section{Complementary topologies and volume convolution}
\label{app:genus}

Two pairs of pants can each have connected complement even when
their union separates the surface. We first list the possible
complements and their contributions to the kernel in
Proposition~\ref{prop:exact-integration}. We then bound the sum of
the complementary volumes for any fixed tuple order, including stable
genus-zero components, as needed in Proposition~\ref{prop:disconnected}.

\subsection{The complete two-pants kernel}

Fix $g\geq3$ and consider an ordered pair of good pants with disjoint
closures.  Write their cuff triples as $\mathbf a=(a_1,a_2,a_3)$ and
$\mathbf b=(b_1,b_2,b_3)$, and let $A$ and $B$ be the corresponding
disjoint sets of cuff labels.  Removing the first pants vertex from the
incidence graph leaves the second pants vertex and the complementary
vertices.  The remaining graph is connected exactly when every block
meets $B$.  Removing the second pants vertex similarly requires every
block to meet $A$.  Thus the connectivity requirements in
Definition~\ref{a:def-admissible} are equivalent to
\begin{equation}\label{a:two-pants-incidence}
 C\cap A\neq\varnothing,\qquad C\cap B\neq\varnothing
 \quad\text{for every block }C.
\end{equation}
These conditions also make the full graph connected.  Let $q$ be the
number of complementary blocks.  Each block meets both triples, so
$q\leq3$. In the table, a block of incidence $(i,j)$ contains $i$ cuffs
from $A$ and $j$ from $B$. The third column counts partitions of the
labeled cuffs into unordered blocks.
\begin{center}
\begin{minipage}{\linewidth}\centering
{\small\itshape Complementary topologies for two disjoint pairs of pants\par}\smallskip
\begin{tabular}{c c c l}
\hline
$q$ & Block incidences & Partitions & Stable genus conditions\\
\hline
$1$ & $(3,3)$ & $1$ & $h=g-4\geq0$\\
$2$ & $(1,1)+(2,2)$ & $9$
    & $h_2\geq1,\ h_4\geq0,\ h_2+h_4=g-3$\\
$2$ & $(1,2)+(2,1)$ & $9$
    & $h_3,h'_3\geq0,\ h_3+h'_3=g-3$\\
$3$ & $(1,1)+(1,1)+(1,1)$ & $6$
    & $h_1,h_2,h_3\geq1,\ \sum_i h_i=g-2$\\
\hline
\end{tabular}
\end{minipage}
\end{center}
In the second row, $h_2$ and $h_4$ are the genera of the components with two
and four boundaries.  There are $3\cdot3$ choices for the size-two
block in that row.  In the third row, the block with incidence $(1,2)$
is determined by its unique $A$-cuff and the $B$-cuff it omits, again
giving $3\cdot3$ choices.  In the last row, a partition is a matching
between the triples, giving $3!=6$ choices.  The genus conditions follow
from \eqref{a:genus-identity} and stability.  In particular, a
two-boundary component has positive genus.

Write $\mathbf a_{\widehat i}$ for $\mathbf a$ with $a_i$ omitted,
and similarly for $\mathbf b_{\widehat j}$. Reading the four rows
in order gives
\begin{equation}\label{a:two-pants-kernel}
\begin{split}
 K_{g,2}(\mathbf a,\mathbf b)
 ={}&\mathbf1_{\{g\geq4\}}V_{g-4,6}(\mathbf a,\mathbf b)\\
 &+\sum_{i,j=1}^3\ \sum_{h=1}^{g-3}
      V_{h,2}(a_i,b_j)
      V_{g-3-h,4}(\mathbf a_{\widehat i},\mathbf b_{\widehat j})\\
 &+\sum_{i,j=1}^3\ \sum_{h=0}^{g-3}
      V_{h,3}(a_i,\mathbf b_{\widehat j})
      V_{g-3-h,3}(\mathbf a_{\widehat i},b_j)\\
 &+\sum_{\sigma\in\mathfrak S_3}
     \ \sum_{\substack{h_1+h_2+h_3=g-2\\h_1,h_2,h_3\geq1}}
             \prod_{i=1}^3V_{h_i,2}(a_i,b_{\sigma(i)}).
\end{split}
\end{equation}
At $g=3$ the first term is omitted.  All empty sums are zero.
Proposition~\ref{prop:exact-integration}, with $m=2$, now gives the
following identity in $[0,\infty]$ for $L,b\geq0$ and nonnegative Borel
$F$.
\begin{equation}\label{a:two-pants-expectation}
 \mathbb E_gD_{g,2}^{L,b}(F)
 =\frac1{4V_g}\int_{\mathcal D_{L,b}^2}
 F\bigl(\ell(\mathbf a),\ell(\mathbf b)\bigr)
 K_{g,2}(\mathbf a,\mathbf b)
 \prod_{i=1}^3a_ib_i\,d\mathbf a\,d\mathbf b.
\end{equation}
The same identity holds for bounded signed or complex Borel $F$.
The factor $1/4$ accounts for the two orders of the winding cuffs in
each selected pair of pants.  At genus three, only the third row of the table
is feasible, with both complementary genera zero.  Since $V_{0,3}=1$,
we obtain
\begin{equation}\label{a:genus-three-kernel}
 K_{3,2}(\mathbf a,\mathbf b)=9.
\end{equation}
All four cut pieces are pants in this case, but their labels and selected
sides distinguish their roles, so there is no additional symmetry factor.

The last three lines of \eqref{a:two-pants-kernel} account for
disconnected joint complements.  For example, when $g\geq4$, attach
one cuff of each selected pair of pants to $S_{1,2}$ and the remaining four
cuffs to $S_{g-4,4}$.  Removing either selected pair of pants leaves a connected
complement, whereas removing both disconnects it.  The cutting formula
in the proof of \cite[Proposition~4.4]{HSWX} accounts only for connected
joint complements, corresponding to the first line of
\eqref{a:two-pants-kernel}.  The other three lines are also admissible.
Their contribution is lower order in the regime used here, by
Proposition~\ref{prop:disconnected}.

The multiplicities $9,9,6$ count partitions. Within one type the
integrals in \eqref{a:two-pants-expectation} can still differ.
Exchanging the two winding cuffs preserves the length function,
whereas exchanging a winding cuff with the third cuff can change it.

\subsection{Uniform summation over stable genera}

For a fixed partition of the cuffs, the genera of the complementary
components can vary with $g$. The next estimate sums over these genera
and keeps track of the power of $g$ contributed by each additional
component.

\begin{lemma}[Stable-genus convolution]\label{lem:genus-convolution}
Fix integers $a,N\geq1$.  There are $C_{a,N}<\infty$ and $g_0(a,N)$
such that, for every $g\geq g_0(a,N)$, every integer $1\leq q\leq N$,
and positive integers $n_1,\ldots,n_q$ with $\sum_i n_i\leq N$,
\begin{equation}\label{a:genus-convolution}
 \frac1{V_g}
 \sum_{\substack{h_i\geq0,\ 2h_i-2+n_i\geq1\ (1\leq i\leq q)\\
                 \sum_i(2h_i-2+n_i)=2g-2-a}}
       \prod_{i=1}^qV_{h_i,n_i}
 \leq C_{a,N}g^{-a-q+1}.
\end{equation}
The sum runs over nonnegative integer genera and includes every stable
genus-zero allocation.  The disjoint-pants application has $a=m$ and
$N=3m$.
\end{lemma}

\begin{proof}
For each integer $1\leq n\leq N$, the fixed-$n$ asymptotic
\eqref{f:volume-asymptotic} from \cite[Theorem~1.2]{MZ} bounds the volume
for all sufficiently large $h$.  For a stable pair $(h,n)$, put
$d=2h-3+n\geq0$.  Since $d+1\asymp h$ for large $h$ with fixed $n$,
enlarging the constant over the finitely many remaining stable types gives
\begin{equation}\label{a:factorial-volume-envelope}
 V_{h,n}\leq C_N(4\pi^2)^d\frac{d!}{\sqrt{d+1}}
 \quad\text{for every stable }h\geq0,\quad1\leq n\leq N.
\end{equation}
This includes all stable genus-zero types in the stated range,
in particular $V_{0,3}=1$ when $N\geq3$.

Set
\begin{equation}\label{a:factorial-indices}
 D=2g-3,\qquad d_i=2h_i-3+n_i,\qquad
 F=\sum_i d_i=2g-2-a-q.
\end{equation}
For fixed boundary counts, genus allocations inject into the weak
compositions of $F$ into $q$ nonnegative parts.  At least one $d_i+1$
is at least $(F+q)/q$, and the remaining inverse square-root factors
are at most one.  Thus
\begin{equation}\label{a:square-root-retention}
 \prod_{i=1}^q(d_i+1)^{-1/2}
 \leq\sqrt{\frac q{F+q}}.
\end{equation}
For integers $F\geq0$ and $q\geq1$,
\begin{equation}\label{a:factorial-convolution}
 \sum_{d_1+\cdots+d_q=F}\prod_{i=1}^qd_i!
 \leq3^{q-1}F!.
\end{equation}
For two factors, division by $F!$ gives
$\sum_{j=0}^F\binom Fj^{-1}\leq3$: the bound is immediate for
$F=0,1$, and for $F\geq2$ the two endpoints contribute two while
each interior binomial coefficient is at least $F$.  The general
case follows by induction, grouping the first $q-1$ parts; $q=1$
is equality.
Combining these estimates, with $q\leq N$, yields
\begin{equation}\label{a:genus-sum-envelope}
 \sum_{\mathbf h}\prod_iV_{h_i,n_i}
 \leq C_N(4\pi^2)^F\frac{F!}{\sqrt{F+q}}.
\end{equation}
The closed-surface asymptotic gives
$V_g\geq c(4\pi^2)^DD!/\sqrt g$ for large $g$.  Here
\begin{equation}\label{a:factorial-gap}
 D-F=a+q-1,\qquad F+q=2g-2-a.
\end{equation}
Once $g\geq a+N+2$, we have $F\geq g$; hence the $a+q-1$
factors in $D!/F!$ are at least $g$, while $g/(F+q)$ is bounded.
Dividing \eqref{a:genus-sum-envelope} by the lower bound for $V_g$
proves \eqref{a:genus-convolution}, after enlarging the genus
threshold.  All constants are uniform over the finitely many allowed
boundary counts and values of $q$.
\end{proof}

For $a=m$ and $N=3m$, this agrees with the convolution estimate in
\cite[Lemma~24 and the proof of Proposition~25]{NWX}.  For integers
$r\geq1$, their reference
volumes are
\begin{equation}\label{a:nwx-reference-volume}
 W_r=\begin{cases}
       V_{r/2+1,0},&r\text{ even},\\
       V_{(r+1)/2,1},&r\text{ odd}.
     \end{cases}
\end{equation}
For $r=2g-2-m$, if $m=2k$ then $W_r=V_{g-k,0}$, while if $m=2k+1$
then $W_r=V_{g-k-1,1}$.  In either case \eqref{f:volume-asymptotic} gives
$W_r/V_g=(8\pi^2g)^{-m}(1+O_m(g^{-1}))$.  Their convolution estimate
supplies the further factor $r^{-(q-1)}$, giving the same power of $g$
as \eqref{a:genus-convolution}.  The direct argument above makes both
the dependence on fixed $m$ and the inclusion of stable genus-zero
terms explicit.

\subsection{A boundary-volume bound including genus zero}

\begin{lemma}\label{a:global-boundary-bound}
For each fixed integer $m\geq1$ there is $B_m\geq1$ such that every stable
type with $h\geq0$ and $2\leq n\leq3m$ satisfies
\begin{equation}\label{a:global-boundary-domination}
 V_{h,n}(\mathbf x)\leq B_mV_{h,n}
             \prod_{j=1}^n\frac{\sinh(x_j/2)}{x_j/2}
 \qquad(\mathbf x\in[0,\infty)^n).
\end{equation}
Each quotient is defined to be one at zero.  The constant is uniform
in the genus $h$ and in all boundary lengths.
\end{lemma}

\begin{proof}
For $h\geq1$, the assertion with constant one is the upper bound in
\cite[Lemma~22]{NWX}.  It remains to treat genus zero, for which
stability requires $n\geq3$.  For each such $n$, the volume-polynomial
theorem \cite[Theorem~1.1]{Mirzakhani} gives
\[
 V_{0,n}(\mathbf x)=\sum_\alpha c_\alpha
                       \prod_jx_j^{2\alpha_j},\qquad c_\alpha\geq0,
\]
where the sum runs over finitely many multiindices
$\alpha\in\mathbb Z_{\geq0}^n$ and $V_{0,n}>0$.
Put $p(x)=\sinh(x/2)/(x/2)$ for $x>0$ and $p(0)=1$.
For each integer $d\geq0$, the coefficient of $x^{2d}$ in its
power series is $1/(2^{2d}(2d+1)!)$.  Positivity of the coefficients
therefore gives
\begin{equation}\label{a:monomial-profile}
 x^{2d}\leq2^{2d}(2d+1)!\,p(x)
 \quad(x\geq0,\ d\geq0).
\end{equation}
Applying this inequality to each monomial gives
\[
 V_{0,n}(\mathbf x)
 \leq\left(\sum_\alpha c_\alpha
                   \prod_j2^{2\alpha_j}(2\alpha_j+1)!\right)
                 \prod_jp(x_j).
\]
The coefficient sum is finite and depends only on $n$.  Dividing it
by $V_{0,n}$ and taking the maximum over $3\leq n\leq3m$, together
with one for positive genus, gives $B_m$.  Stability excludes $(0,2)$,
so this covers every type in the statement.  The resulting bound is
uniform for all boundary lengths and introduces no power of a length
cutoff.
\end{proof}

% END sections/appendix_a_genus.tex

% BEGIN sections/appendix_b_attachments.tex
\section{Cuff attachments and fixed-core integrals}\label{app:attachments}

When a marked pair of pants shares a cuff with the core, that length
occurs in both pieces and is integrated only once. We first bound
the local mass of the core. We then hold the shared lengths fixed
while integrating the other cuffs of the marked pants. This gives
the attachment estimates used in Section~\ref{sec:witnesses}.
Boundaries of bordered moduli spaces are labeled throughout. We write
\[
 \psi(t)=\frac{\sinh(t/2)}{t/2},\qquad
 \rho(t)=\frac{t}{e^t-1},\qquad
 \psi(0)=\rho(0)=1,\qquad
 \mathcal A_g(L)=\frac{Le^L}{g}.
\]
Each attachment cuff and each internal identification seam contributes
one length/twist Jacobian.  An attachment cuff joins the core to a
selected marked pair of pants.  An identification seam lies inside the
completed core and is determined by the generating pair.

\subsection{A kernel bound for a fixed number of filling curves}

Consider an ordered pair of primitive unoriented figure-eights with
distinct embedded pants supports whose interiors intersect.  As in
Lemma~\ref{lem:f-pair-core}, write its intrinsic geodesic filled surface
and completed image as
\[
 Y=S_{h,n},\qquad c=2h-2+n,\qquad
 Z=S_{h+j,\nu},\qquad \nu=n-2j\geq1.
\]
Here \(h,j\geq0\) and \(n\geq1\) are integers.  Identifying \(j\)
disjoint pairs of intrinsic boundary copies produces the seams of
\(Z\).  This seam system is determined by the generating pair, and
cutting it recovers connected \(Y\).  Each curve has length at most
\(L\), and every cuff of either support has length at least \(b\).
The local count allows all ambient complementary topologies.

For a metric \(Z\in\mathcal M_{h+j,\nu}(\mathbf t)\), let
\(N_{L,b}^{h,n,j}(Z)\) count these ordered pairs whose intrinsic
filled surface has type \((h,n)\) and whose completed image is all
of \(Z\).  Each support is embedded and has three distinct actual
cuffs.  Define
\begin{equation}\label{att:fixed-kernel-def}
 \begin{split}
 K_{L,b}^{h,n,j}(\mathbf t)
   &=\int_{\mathcal M_{h+j,\nu}(\mathbf t)}
           N_{L,b}^{h,n,j}(Z)\,dV_{\rm WP}(Z),\\
 \mathcal J_{L,b}^{h,n,j}
   &=\int_{\mathbb R_+^\nu}
           K_{L,b}^{h,n,j}(\mathbf t)
           \prod_{i=1}^\nu t_i\psi(t_i)\,d\mathbf t .
 \end{split}
\end{equation}
These nonnegative integrals initially take values in \([0,\infty]\).
All internal curves and twists are included in the local count.

More generally, fix an integer \(q\geq2\).  On
\(Z\in\mathcal M_{h+j,\nu}(\mathbf t)\), let
\(N_{L,b}^{(q),h,n,j}(Z)\) count an admissible family of ordered
\(q\)-tuples of distinct primitive unoriented closed geodesics, each
tuple counted once, and assume this count is measurable.  Suppose that
every counted tuple determines a
system of \(j\) disjoint simple geodesic seams, disjoint from its
curves, whose simultaneous cut is a connected surface
\(Y\cong S_{h,n}\).  On \(Y\), the tuple is nontrivial,
nonperipheral and filling, its total length is at most \(qL\), and
\(p(Y)\leq2qL\).  The perimeter includes both copies of every seam.
Define \(K_{L,b}^{(q),h,n,j}\) and
\(\mathcal J_{L,b}^{(q),h,n,j}\) by
\eqref{att:fixed-kernel-def}, with this count in place of
\(N_{L,b}^{h,n,j}\).  Additional restrictions, including cuff lower
bounds depending on \(b\), may be imposed on the counted family.
The pair kernel above is the case \(q=2\).

\begin{lemma}[Fixed-tuple core kernel]
Fix\label{lem:fixed-core-kernel} integers \(c\geq3\), \(q\geq2\), an allowable type
\((h,n,j)\) with \(2h-2+n=c\) and \(n-2j\geq1\), and
\(0<\delta<1/2\).  Under the preceding hypotheses, for
\(L\geq1\) and \(b\geq0\),
\begin{equation}\label{att:tuple-kernel-bound}
 \mathcal J_{L,b}^{(q),h,n,j}
 \leq C_{c,q,\delta}(1+L)^{3c+q-1-2j}
                    e^{q(1+\delta)L}.
\end{equation}
The constant is uniform in \(L,b\), in the admissible subfamily,
and in the finitely many allowable types of complexity \(c\).
In particular, the pair kernel above satisfies, for \(c\geq3\),
\begin{equation}\label{att:fixed-core-bound}
 \mathcal J_{L,b}^{h,n,j}
 \leq C_{c,\delta}(1+L)^{3c+1-2j}
                    e^{(2+2\delta)L}.
\end{equation}
\end{lemma}

\begin{proof}
By \cite[Definition~17 and Theorem~18]{WXPGT}, the number of
ordered nontrivial, nonperipheral, primitive unoriented filling
\(q\)-tuples on \(Y\), with total length at most \(qL\), is bounded by
\begin{equation}\label{att:filling-tuple-input}
 N_q^{\rm fill}(Y,qL)
 \leq C_{c,q,\delta}(1+L)^{q-1}
        e^{qL-(1-\delta)p(Y)/2}.
\end{equation}
The constant is uniform in the metric and boundary lengths.
Appendix~A of the cited paper sketches the boundary-shortening
proof and refers to~\cite{WX} for details.

Cutting the determined seams leaves the curves and their lengths
unchanged.  They remain primitive, since a root on the cut surface
would be a root in \(Z\), and nonperipherality and filling are
hypotheses.  Enlarge the count to all seam multicurves whose cut is
connected of type \((h,n)\), and then to all filling tuples on that
cut surface with the stated total-length and perimeter bounds.
There are finitely many seam multicurve orbits at fixed type.
Labeling seams and their sides introduces only a fixed finite factor.
This is a statement about seam systems, not about the number of
filling-tuple orbits.

Put \(T=\sum_{i=1}^{\nu}t_i\) and
\(U=\sum_{a=1}^j s_a\), where \(s_a\) are the seam lengths.
Multicurve integration on the bordered core and
\eqref{att:filling-tuple-input} give
\begin{equation}\label{att:seam-integral}
 \begin{split}
 K_{L,b}^{(q),h,n,j}(\mathbf t)
 \leq{}&C_{c,q,\delta}(1+L)^{q-1}e^{qL}
 \int_{\mathbb R_+^j}
   \mathbf1_{\{T+2U\leq2qL\}}
   e^{-(1-\delta)(T+2U)/2}\\
 &\hspace{12mm}\cdot
 V_{h,n}(t_1,\ldots,t_\nu,s_1,s_1,\ldots,s_j,s_j)
 \prod_{a=1}^j s_a\,d\mathbf s .
 \end{split}
\end{equation}
For \(j=0\), the integral denotes the remaining integrand.
Each seam appears twice among the volume's boundary arguments but
contributes only one twist integral \(s_a\,ds_a\).  Its two sides
belong to the core, so it carries no exterior profile \(\psi(s_a)\).
Tonelli's theorem justifies the nonnegative unfolding before
finiteness is established.

The fixed-type volume polynomial has length degree
\(D=6h-6+2n\), and every boundary argument in
\eqref{att:seam-integral} is at most \(2qL\).  Hence
\(V_{h,n}(\mathbf t,\mathbf s,\mathbf s)
\leq C_{c,q}(1+L)^D\).
Multiplication by the outer weights
\(t_i\psi(t_i)\leq e^{t_i/2}\) leaves
\[
 e^{qL+\delta T/2-(1-\delta)U}
 \leq e^{q(1+\delta)L}e^{-(1-\delta)U}.
\]
After extracting this bound, enlarge the integration domain to
\(\{T\leq2qL\}\times\mathbb R_+^j\).  Its outer simplex has
volume \((2qL)^\nu/\nu!\), and the seam integral is
\[
 \int_{\mathbb R_+^j}\prod_{a=1}^j
         s_a e^{-(1-\delta)s_a}\,d\mathbf s
 = (1-\delta)^{-2j}.
\]
The resulting polynomial power is
\[
 q-1+D+\nu=q-1+6h-6+3n-2j=3c+q-1-2j.
\]
This proves \eqref{att:tuple-kernel-bound}.  Dropping additional
counting restrictions makes the bound uniform in \(b\) and in the
admissible subfamily.

For the pair kernel in \eqref{att:fixed-kernel-def}, the intrinsic
core construction in Lemma~\ref{lem:f-pair-core} supplies the filling
property and determined seams.  Each figure-eight is nonperipheral,
and a root after cutting would give an ambient root.  The pants length
identity bounds each support perimeter by \(2L\); straightening the
union therefore gives
\begin{equation}\label{att:intrinsic-perimeter}
 p(Y)\leq4L.
\end{equation}
The total curve length is at most \(2L\), so \(q=2\) in
\eqref{att:tuple-kernel-bound} proves \eqref{att:fixed-core-bound}.
\end{proof}

For each intrinsic type, $j$ counts the identified boundary pairs and $\nu$ counts the exposed cuffs. The last column gives the power of $1+L$ in \eqref{att:fixed-core-bound}.
\begin{center}
\begin{minipage}{\linewidth}\centering
{\small\itshape Boundary identifications for cores of complexity three\par}\smallskip
\begin{tabular}{c c c c c}
\hline
Intrinsic \(Y\)&\(j\)&Completed \(Z\)&\(\nu\)&Power in
\eqref{att:fixed-core-bound}\\
\hline
\(S_{0,5}\)&0&\(S_{0,5}\)&5&10\\
\(S_{0,5}\)&1&\(S_{1,3}\)&3&8\\
\(S_{0,5}\)&2&\(S_{2,1}\)&1&6\\
\(S_{1,3}\)&0&\(S_{1,3}\)&3&10\\
\(S_{1,3}\)&1&\(S_{2,1}\)&1&8\\
\(S_{2,1}\)&0&\(S_{2,1}\)&1&10\\
\hline
\end{tabular}
\end{minipage}
\end{center}
The type \(S_{2,1}\) is compatible with connected individual pants
complements.  A pair of pants embedded in \(S_{2,1}\) can have connected
complement \(S_{0,4}\) within that surface.

\subsection{Prescribing one, two, or three cuffs of a mark}

We first estimate the contribution of a marked pair of pants when some of
its cuff lengths are fixed.  Let \(\mathcal D_{L,b}\) be the true
pants domain \eqref{w:true-pants-domain}.  Prescribe
\(s\in\{0,1,2,3\}\) cuff lengths and their winding or third-cuff
roles.  Let \(B_s\) be the integral of the domain indicator with
measure \(u\psi(u)\,du\) for each of the \(3-s\) free cuffs.
Thus \(B_3\) is the indicator itself, and \(B_s=0\) whenever a
prescribed cuff is below \(b\).  We suppress the parameters \(L,b\)
and the role assignment from the notation.  For \(L\geq1\), every
assignment satisfies
\begin{equation}\label{att:prescribed-cuffs}
 \begin{aligned}
 B_0&\leq4Le^L,&
 B_1(t)&\leq4Le^{L-t/2},\\
 B_2(t_1,t_2)&\leq2e^{L-(t_1+t_2)/2},&
 B_3(t_1,t_2,t_3)&\leq e^{L-(t_1+t_2+t_3)/2}.
 \end{aligned}
\end{equation}
The first three bounds follow from the integrations for each cuff
role in \eqref{w:touching-slices}, using
\(u\psi(u)\leq e^{u/2}\).  For the last bound, let
\(T=t_1+t_2+t_3\).  The true domain implies that the two winding
cuffs sum to at most \(L\) and the third cuff is at most \(L\).
Thus \(T\leq2L\) on feasible data, where \(1\leq e^{L-T/2}\).
On infeasible data the indicator vanishes.

The identity
\begin{equation}\label{att:profile-ratio}
 \frac{e^{-t/2}}{\psi(t)}=\rho(t)
\end{equation}
gives the saving at each attachment cuff.  The function \(\rho\)
decreases from one.  For \(t>0\), the numerator of its derivative
is \(e^t(1-t)-1\).  This expression is negative because it vanishes
at zero and has derivative \(-te^t<0\).

Now attach a marked pair of pants along \(s\in\{1,2,3\}\) outer core
cuffs.  Each shared cuff contributes one length/twist Jacobian and
receives no complementary boundary profile.  Relative to the core
measure \(\prod_i t_i\psi(t_i)\,dt_i\), the multiplier is therefore
\(B_s/\prod_{\rm shared}\psi(t_i)\).  On the attachment domain,
each shared cuff has length at least \(b\), so
\eqref{att:prescribed-cuffs}--\eqref{att:profile-ratio} give
\begin{equation}\label{att:single-pants-multiplier}
 \frac{B_s}{\prod_{\rm shared}\psi(t_i)}
 \leq C_s L^{\mathbf1_{\{s=1\}}}e^L
          \prod_{\rm shared}\rho(t_i)
 \leq C_s L^{\mathbf1_{\{s=1\}}}e^L\rho(b)^s .
\end{equation}
The threshold \(b\) is imposed on the shared core cuffs because they
belong to the attached mark.  This argument imposes no lower bound on
the unshared outer core cuffs.

\subsection{The joint cutting estimate with arbitrary attachments}

The following kernel formulation applies to the low-complexity
masses of Lemma~\ref{w:low-local-mass}, the fixed pair masses above,
and the decorated cores used later.  Fix a completed embedded core
of type \(S_{h_0,\nu}\), with \(\nu\geq1\) and
\(c=2h_0-2+\nu\geq1\).  Each decoration determines its selected
core.  On the bordered moduli space, take a nonnegative measurable
function that pointwise dominates the number of admissible local
decorations, including all their internal choices.  Assume this
function is independent of the outer gluing twists.  Let
\(K(\mathbf t)\) be its bordered moduli integral, and put
\begin{equation}\label{att:abstract-core-mass}
 \mathcal J(K)=\int_{\mathbb R_+^\nu}
                  K(\mathbf t)\prod_{i=1}^{\nu}t_i\psi(t_i)
                  \,d\mathbf t .
\end{equation}

Let \(p_s\in\mathbb Z_{\geq0}\) count the marked pants attached to
exactly \(s\) distinct outer core cuffs, for \(s=1,2,3\).  These
pants occupy ordered marked slots, and their interiors avoid the
core.  Put
\begin{equation}\label{att:contact-parameters}
 p=p_1+p_2+p_3,\qquad
 k=p_1+2p_2+3p_3,\qquad d=p_2+p_3.
\end{equation}
Let \(r\in\mathbb Z_{\geq0}\) further ordered marked pants be
exterior, meaning that their closed supports are disjoint from the
core.  All marked closed supports are pairwise disjoint.  In
applications with an anchor in the core, they are also
closed-disjoint from that anchor.  Each marked curve is a primitive
unoriented figure-eight of length at most \(L\), with three
distinct actual support cuffs of length at least \(b\).
Additional goodness conditions may be imposed.  Write
\(\mathbf p=(p_1,p_2,p_3)\) and denote the count of these decorated
core and mark configurations by \(I_{K;\mathbf p,r}\).  Add a
subscript \(q\) when the complement of the selected union has
exactly \(q\) components.

\begin{lemma}[Cuff attachment bound]\label{lem:attachment-bounds}
For fixed core type and nonnegative integers \(p,r\), there are
constants \(C_{c,p,r}\) and \(g_0(c,p,r)\) such that, for every
\(g\geq g_0(c,p,r)\) and integer \(q\geq1\),
\begin{equation}\label{att:joint-bound}
 \mathbb E_g I_{K;\mathbf p,r;q}
 \leq C_{c,p,r}\,
 g^{-(q-1)}\mathcal A_g(L)^{p+r}g^{-c}
 L^{-d}\rho(b)^k\,\mathcal J(K)
\end{equation}
for \(L\geq1,b\geq0\).  Summing over all complementary topologies
gives
\begin{equation}\label{att:all-complements-bound}
 \mathbb E_g I_{K;\mathbf p,r}
 \leq C_{c,p,r}\,
 \mathcal A_g(L)^{p+r}g^{-c}
 L^{-d}\rho(b)^k\,\mathcal J(K).
\end{equation}
Impossible contact data contribute zero.  The constants are uniform
in the contact data, \(L,b\), and the choice of local majorant
defining \(K\).
\end{lemma}

\begin{proof}
We first enumerate the possible outer attachments.  Label the
\(\nu\) outer core cuffs.  Each attached pair of pants specifies a subset
of size one, two or three, and distinct attached pants use disjoint
subsets.  Assign the subsets to ordered marked slots and specify
the cuff roles.  At fixed \(c,p,r\), there are finitely many such
data, with
\[
 p\leq k\leq\nu\leq c+2.
\]
These data cover every actual attachment, and unrealizable choices
contribute zero.  Internal seams are unavailable for attachment
because the completed core occupies both sides of each seam.

Let \(R\) be the union of the core and its attached pants.
It is connected and has
\begin{equation}\label{att:enlarged-core-topology}
 \chi(R)=-c-p,\qquad
 n_R=\nu+3p-2k,\qquad
 h_R=h_0+k-p.
\end{equation}
If \(n_R=0\), the embedded closed subsurface \(R\) equals the
connected ambient surface \(X\).  This would force
\(2g-2=c+p\), which is excluded by taking \(g\) sufficiently large
at fixed \(c,p\).
Every exposed cuff is a distinct actual geodesic.  Equality of a
free attached cuff with an outer core cuff would be another
recorded attachment.  Equality with a cuff of another selected
mark violates marked closed-disjointness, and equality within one
marked pair of pants violates its three-distinct-cuff condition.
The completed core already has distinct outer cuffs.

With the \(r\) exterior pants, the number of exposed cuff geodesics
and the number of distinct cuts are respectively
\begin{equation}\label{att:cut-numbers}
 N=\nu+3(p+r)-2k,\qquad N+k .
\end{equation}
The \(N+k\) cuts include the attachment cuffs.  Internal
identification seams have already been integrated in \(K\).
Every complementary component has nonempty boundary, since a
boundaryless component would be a separate component of \(X\).
The exposed boundaries are distinct geodesics, which excludes disks,
annuli and collapsed annuli.  Thus every complementary component is
stable.

Partition the \(N\) exposed labels into \(q\) nonempty blocks
\(A\), with boundary counts \(n_A\) and nonnegative integer genera
\(h_A\).  Euler characteristic gives
\begin{equation}\label{att:complementary-equations}
 \begin{gathered}
  \sum_A n_A=\nu+3(p+r)-2k,\qquad
  2h_A-2+n_A>0,\\
  \sum_A h_A=g+q+h_0+k-2(p+r)-c-2 .
 \end{gathered}
\end{equation}
For actual data the incidence graph is connected,
and the classification of labeled oriented bordered surfaces gives
the gluing orbit.  For an upper bound, sum over all stable genus
allocations and partitions satisfying
\eqref{att:complementary-equations}.  Dropping the remaining graph,
anchor-disjointness and individual-complement conditions enlarges
the sum.  There are finitely many partitions, with
\(1\leq q\leq N\leq c+2+3(p+r)\), including blocks with one
boundary.

The total complementary complexity is \(2g-2-(c+p+r)\).
Lemma~\ref{lem:genus-convolution},
with \(a=c+p+r\) and \(N=c+2+3(p+r)\), gives
\begin{equation}\label{att:genus-ratio}
 \frac1{V_g}\sum_{\mathbf h\ {\rm satisfying}\
                         \eqref{att:complementary-equations}}
      \prod_A V_{h_A,n_A}
 \leq C_{c,p,r}g^{-c-p-r-q+1}.
\end{equation}
For the same finite boundary range,
\begin{equation}\label{att:complementary-profiles}
 V_{h_A,n_A}(\mathbf u_A)
 \leq C_{c,p,r}V_{h_A,n_A}
                  \prod_{\ell\in A}\psi(u_\ell).
\end{equation}
For \(n_A\geq2\), this is Lemma~\ref{a:global-boundary-bound} with
an enlarged fixed boundary parameter.  A stable one-boundary
component has positive genus, so \cite[Lemma~22]{NWX} supplies that
case.  Thus both estimates include every stable genus-zero allocation
and every permitted one-boundary component.

It remains to integrate the boundary lengths.  Cut all exposed and
attachment cuffs.  Disintegration of the Weil--Petersson measure
gives the core moduli factor, one complementary moduli factor per
block, and one length/twist measure per distinct cut.  The local
majorant is independent of the outer twists.  Each attached or
exterior pair of pants has volume \(V_{0,3}=1\) and retains its true
curve-length indicator.
All labeling and orbit factors are bounded at fixed \(c,p,r\).
For a one-holed-torus component, the convention in
\cite[Remark~2.7]{AMI} halves the ordinary volume factor, so the
ordinary volume gives an upper bound.  All remaining finite factors
are absorbed into \(C_{c,p,r}\).

After applying \eqref{att:complementary-profiles}, each exposed
cuff receives one profile and each attachment cuff receives none.
With the core lengths fixed, the free variables of distinct marked
pants are disjoint, so their integrals separate.
Each exterior pair of pants contributes
\begin{equation}\label{att:exterior-pants-integral}
 \mathcal P_{L,b}:=
 \int_{\mathcal D_{L,b}}\prod_{i=1}^3u_i\psi(u_i)\,d\mathbf u
 \leq4Le^L.
\end{equation}
For the attached pants, multiply
\eqref{att:single-pants-multiplier} over the disjoint shared-cuff
subsets.  The combined multiplier relative to the core measure is
\begin{equation}\label{att:combined-multiplier}
 C(Le^L)^pL^{-d}\rho(b)^k,\qquad d=p_2+p_3.
\end{equation}
The unique shared-cuff Jacobian \(t_i\,dt_i\) is retained in
\(\mathcal J(K)\).  Extracting the multiplier in
\eqref{att:combined-multiplier} and dropping the remaining
compatibility restrictions leaves precisely this integrated core
mass.  Consequently the length part of the upper integral is at most
\[
 C(Le^L)^{p+r}L^{-d}\rho(b)^k\mathcal J(K).
\]
Combining this with \eqref{att:genus-ratio} proves
\eqref{att:joint-bound}.  Summing over the finitely many partitions
and values of \(q\) gives \eqref{att:all-complements-bound}.
All sums and integrands are nonnegative, so Tonelli justifies the
argument even when the core mass is initially allowed to be infinite.
\end{proof}

The exponent \(d=p_2+p_3\) counts attached pants with at least
two prescribed cuffs. Each contributes at most a constant times
\(e^L\), saving one factor of \(L\) relative to \(Le^L\).
Prescribing the third cuff gives the same saving as prescribing two.

\subsection{The complete complexity-two boundary catalog}

We now classify the complexity-two boundary contacts used in
Section~\ref{sec:witnesses}.  Consider an anchor/witness pair of intrinsic complexity two, with
completed image \(Z\), and write \(P\) for the anchor pants.  Both
generating supports are good.  All other marked closed supports
are pairwise disjoint and disjoint from \(P\), and their interiors avoid \(Z\).

\begin{lemma}[Available outer cuffs at complexity two]
\label{att:low-catalog}
If \(Z=S_{1,2}\), the anchor occupies exactly one outer cuff.
If \(Z=S_{0,4}\), it occupies exactly two outer cuffs.
The possible nonempty boundary-contact patterns are therefore
as follows.
\begin{enumerate}
\item one attached pair of pants using one cuff for intrinsic \(S_{1,2}\);
\item one pair of pants using one cuff, one pair of pants using two cuffs, or two
      pants using one cuff each for intrinsic \(S_{0,4}\) with
      completed image \(S_{0,4}\);
\item one attached pair of pants using one outer cuff for intrinsic
      \(S_{0,4}\) with completed image \(S_{1,2}\).
\end{enumerate}
\end{lemma}

\begin{proof}
Let \(a\) be the number of outer cuffs occupied by the anchor.
The metric completion of \(Z\setminus\overline P\) has Euler
characteristic \(-1\) and boundary number \(5-2a\) for
\(Z=S_{1,2}\), or \(7-2a\) for \(Z=S_{0,4}\).
After removing shared outer cuffs, its boundary curves are distinct.
These geodesic boundaries exclude disks and annuli, so every
component is stable.  Their total negative Euler characteristic
is one, and hence there is exactly one component.

In the planar case the component is a pair of pants, so \(7-2a=3\) and
\(a=2\).  For \(Z=S_{1,2}\), it is either a pair of pants or a one-holed
torus, giving \(a=1\) or \(2\).  The latter case leaves no
unoccupied outer cuff connecting the torus to the nonempty
exterior of \(Z\).  The anchor would then have disconnected
ambient complement, contrary to goodness.  Thus \(a=1\).
The argument uses the completed image \(Z\), so it also covers
the identified \(S_{0,4}\) subtype.

Geodesic-boundary subsurfaces with disjoint interiors can meet
only along common cuffs.  Closed-disjointness prevents another
mark from using an anchor-occupied cuff and prevents distinct
marks from using the same cuff.  Thus the \(S_{1,2}\) completion
has one available cuff, and the \(S_{0,4}\) completion has two.
This gives the listed patterns.  Internal identification seams
are unavailable because the completed core occupies both sides.

Goodness imposes no further topological exclusion of these
patterns, including the double-cuff attachment.  Attach all
remaining exposed cuffs to one connected stable exterior.
Each attached pair of pants has a free cuff, so every outer connection
required by extendibility still reaches this exterior, either
directly or through an attached pair of pants.  Both generating supports
therefore retain connected complements.  Removing an attached
pair of pants also leaves a connected complement, since the core still
joins the exterior along the anchor-occupied outer cuffs.
The length restrictions are imposed separately.
\end{proof}

For \(L>1\), use the notation \(\zeta\) and \(K_L^\zeta\) of
\eqref{w:low-kernels}, and put \(b_L=10\log L\).
Goodness makes both supports of every actual pair extendible.
The local mass in Lemma~\ref{w:low-local-mass} therefore bounds
these configurations for every arrangement of attachments and
exterior complementary components.  In type \(04^1\), the kernel
\(K_L^{04^1}\) retains the determined identification seam and all
its lengths and twists.

For intrinsic \(S_{1,2}\), let \(I_{12,\rm att;r}(L)\) count
configurations with one single-cuff attachment and \(r\) exterior
marks.  Lemma~\ref{lem:attachment-bounds}, with
\((c,p_1,p_2,p_3)=(2,1,0,0)\), gives
\begin{equation}\label{att:s12-one-cuff}
 \mathbb E_g I_{12,\rm att;r}(L)
 \leq C_r\mathcal A_g(L)^{r+1}
                \frac{e^{2L}}{g^2}\rho(b_L).
\end{equation}
If the outer core lengths are \(a,t\) and the cuff of length
\(t\) is shared, the profile division takes the explicit form
\[
 \begin{split}
 \int K_L^{12}(a,t)\,at\psi(a) B_1(t)\,da\,dt
 &=\int_{t\geq b_L}
    K_L^{12}(a,t)\,at\psi(a)\psi(t)
           \frac{B_1(t)}{\psi(t)}\,da\,dt\\
 &\leq4Le^L\rho(b_L)\mathcal J_L^{12}.
 \end{split}
\]
The shared cuff contributes its single Jacobian factor \(t\).
The enlarged core is \(S_{1,3}\).  With \(r\) exterior pants,
there are \(3+3r\) exposed boundaries, and
\eqref{att:genus-ratio} gives the genus exponent \(r+q+2\).

For intrinsic \(S_{0,4}\), let \(\epsilon=0\) for completed image
\(S_{0,4}\) and \(\epsilon=1\) for completed image \(S_{1,2}\).
Let \(p\) count attached pants and \(k\) count shared cuffs, and
write \(I^\epsilon_{p,k;r}(L)\) for the corresponding count with
\(r\) exterior marks. When \(k=2\), the two shared cuffs can belong to
one attached pair of pants or to two different pairs of pants.
\begin{center}
\begin{minipage}{\linewidth}\centering
{\small\itshape Attachments to a four-holed-sphere core\par}\smallskip
\begin{tabular}{c c c}
\hline
\((\epsilon,p,k)\)&Attachment pattern&Enlarged surface \(R\)\\
\hline
\((0,1,1)\)&one single-cuff attachment&\(S_{0,5}\)\\
\((0,1,2)\)&one double-cuff attachment&\(S_{1,3}\)\\
\((0,2,2)\)&two single-cuff pants&\(S_{0,6}\)\\
\((1,1,1)\)&one single-cuff attachment&\(S_{1,3}\)\\
\hline
\end{tabular}
\end{minipage}
\end{center}
Here
\[
 \chi(R)=-2-p,\qquad
 n_R=4-2\epsilon+3p-2k,\qquad h_R=\epsilon+k-p.
\]
Every attachment in this catalog uses at most two cuffs, so
\(d=k-p\). For either completion the local-mass and attachment bounds give
\begin{equation}\label{att:s04-bound}
 \mathbb E_g I^\epsilon_{p,k;r}(L)
 \leq C_r\mathcal A_g(L)^{p+r}
       \frac{e^{2L}}{g^2}
       L^{-(k-p)}\rho(b_L)^k .
\end{equation}
The bound includes all winding-cuff and third-cuff assignments,
all stable complementary genera, and both completion subtypes.
In the identified subtype, the seam has already been integrated
inside the local kernel.

Fix \(S\) and let \(g\to\infty\), with \(H=\log g\) and
\(L=H-\log H+S\).  Then
\begin{equation}\label{att:critical-basic}
 \mathcal A_g(L)=O_S(1),\qquad
 \frac{e^{2L}}{g^2}=e^{2S}H^{-2},\qquad
 \rho(b_L)=\frac{10\log L}{L^{10}-1}
          =O_S\!\left(\frac{\log H}{H^{10}}\right).
\end{equation}
Equation~\eqref{att:s12-one-cuff} is therefore
\(O_{r,S}(\log H/H^{12})\).  Applying \eqref{att:s04-bound}
to the four rows of the table, in order, gives
\begin{equation}\label{att:low-critical-rates}
 O_{r,S}\!\left(\frac{\log H}{H^{12}}\right),\qquad
 O_{r,S}\!\left(\frac{(\log H)^2}{H^{23}}\right),\qquad
 O_{r,S}\!\left(\frac{(\log H)^2}{H^{22}}\right),\qquad
 O_{r,S}\!\left(\frac{\log H}{H^{12}}\right).
\end{equation}
For fixed marked order \(m\), set \(r=m-1-p\) and sum over
feasible patterns and ordered slot assignments.  Only patterns
with \(r\geq0\) occur.  The resulting raw bound for the entire
complexity-two boundary-only class is
\(O_{m,S}(\log H/H^{12})\).

\subsection{Fixed higher complexity and assembly of all boundary placements}

Fix an intrinsic type \((h,n,j)\) with \(c=2h-2+n\geq3\) and
contact parameters as in \eqref{att:contact-parameters}.  Applying
Lemma~\ref{lem:attachment-bounds} to \(K_{L,b}^{h,n,j}\) and using
Lemma~\ref{lem:fixed-core-kernel} gives
\begin{equation}\label{att:fixed-attachments}
 \mathbb E_g I_{h,n,j;\mathbf p,r}(L,b)
 \leq C_{c,p,r,\delta}\mathcal A_g(L)^{p+r}
 \frac{(1+L)^{3c+1-2j}e^{(2+2\delta)L}}{g^c}
 L^{-d}\rho(b)^k .
\end{equation}
At fixed \(c,p,r\), there are finitely many choices of disjoint
subsets of the \(\nu=n-2j\) outer labels, of sizes one, two or three,
together with ordered slot assignments and winding roles.  These
choices cover every actual attachment pattern.  The core kernel
retains all internal curve configurations, whose number need not
be bounded by the size of this outer catalog.

Some entries are impossible.  For example, a pair of pants using all three
cuffs of a three-boundary core would close the enlarged core, which
is excluded for sufficiently large genus.  An internal identification
seam is also unavailable for attachment.  Further restrictions from
anchor-disjointness or individual connected complements may be
dropped in the nonnegative upper integral.  Thus the estimate covers
all boundary-only placements, including feasible three-cuff attachments.

The \(p+r\) additional marks occupy the slots other than the anchor,
so \(m=1+p+r\).  Fix \(m\) and \(S\), set \(H=\log g\) and
\(L=H-\log H+S\), and let \(g\to\infty\).  Taking \(\delta=1/8\)
in \eqref{att:fixed-attachments} gives
\begin{equation}\label{att:fixed-attachments-critical}
 O_{c,m,S}\!\left(
 g^{-c+9/4}H^{3c-5/4-2j-d}\rho(b)^k
 \right),\qquad m=1+p+r.
\end{equation}
A nonempty attachment pattern has \(k\geq1\).  At
\(b=b_L=10\log L\), use \(\rho(b_L)^k\leq\rho(b_L)\) and
\(H^{-2j-d}\leq1\) for large \(g\).  The complexity-\(c\) term
is then bounded by
\[
 C_{c,m,S}g^{-c+9/4}H^{3c-45/4}\log H.
\]
For a fixed integer \(M\geq3\), absorb all constants and type
multiplicities over \(3\leq c\leq M\) into one constant.  The
remaining factors form a geometric progression with ratio \(H^3/g\),
which tends to zero.  Summing therefore gives the order of the
\(c=3\) term,
\begin{equation}\label{att:finite-band-boundary}
 \mathbb E_g I^{\rm bd}_{3\leq c\leq M,m}
 =O_{M,m,S}\!\left(g^{-3/4}H^{-9/4}\log H\right).
\end{equation}
The left-hand count includes each tuple/witness pair once.
Summing over its possible anchor and pattern realizations gives
an upper bound.

Take \(K_m=2m+9\) and partition the boundary-only class by
intrinsic pair complexity.  Equations \eqref{att:low-critical-rates}
and \eqref{att:s12-one-cuff} cover complexity two, while
\eqref{att:finite-band-boundary} covers \(3\leq c<K_m\).
For \(c\geq K_m\), apply Proposition~\ref{prop:witness-pair-tail}
with \(r=m-1\) additional marks in arbitrary positions.  Its
threshold \(11+2r\) equals \(K_m\), so this tail contributes
\(O_{m,S}(g^{-1}H^{-(m-1)})\).  Hence
\begin{equation}\label{att:all-boundary-incidence}
 \mathbb E_g I_m^{\rm bd}
 \leq C_{m,S}\left(
       \frac{\log H}{H^{12}}
       +g^{-3/4}H^{-9/4}\log H
       +g^{-1}H^{-(m-1)}
       \right)
 =O_{m,S}\!\left(\frac{\log H}{H^{12}}\right).
\end{equation}
At \(m=1\) this class is empty.  The fixed-core constants are used
only for \(c<K_m\).  The last term comes from the separate tail
estimate, which covers the entire remaining complexity range.

For \(\tau\in\mathcal T_m(X)\), let \(d_{\rm bd}(\tau)\) count
its boundary-only witnesses.  The completed pair image contains the
witness support, while all other marked interiors avoid that image.
The crossed anchor is therefore unique, and
\(I_m^{\rm bd}=\sum_\tau d_{\rm bd}(\tau)\).
By \eqref{w:raw-domination}, retaining the original denominator
\(d(\tau)\) from \eqref{w:witness-set} and the convention \(0/0=0\),
we obtain
\[
 \mathbb E_g
 \sum_{\tau}\frac{d_{\rm bd}(\tau)}{d(\tau)}
 \leq\mathbb E_g I_m^{\rm bd}.
\]
Equation \eqref{att:all-boundary-incidence} includes disconnected
marked complements.  Marks whose interiors enter the completed
image are treated by the barrier lemma and the subsequent
multiple-crossing analysis.

% END sections/appendix_b_attachments.tex

% BEGIN sections/appendix_c_planar.tex
\section{Planar core calculations}\label{app:planar}

On the five-holed sphere, Lemma~\ref{lem:planar-mass} covers the
relevant triples by two families. In the first, an anchor and its
witness fill the entire core. In the second, the witness is supported
on a pair of pants with two outer cuffs. We bound the two masses
separately and then add the exterior marks in
\eqref{l:planar-exterior-incidence}.
Throughout, $L\ge1$, $b\ge0$, and $0<\delta<1/4$ is fixed.
We use the profile $\psi$ and the scale $\mathcal A_g$
from \eqref{l:profiles}.
All outer boundaries are labeled.  Constants absorb the finite choices
of boundary labels, the three actual figure-eight choices per supporting
pair of pants, and multicurve symmetry factors. The unfolding retains
all internal isotopy classes and twists. Nonnegativity permits the use
of Tonelli's theorem throughout.

\subsection{Counting input and the whole-core-pair integral}

For positive lengths $x,y,z$, the boundary summand in the
McShane--Mirzakhani identity is
\begin{equation}\label{pl:mcshane-summand}
 \mathcal R(x,y,z)
 =x-\log\frac{\cosh(y/2)+\cosh((x+z)/2)}
                    {\cosh(y/2)+\cosh((x-z)/2)}.
\end{equation}
Fix two outer boundaries of lengths $x,y$.  The identity and nonnegativity
of its summands imply that the sum of $\mathcal R(x,y,z)$ over all simple internal geodesics of
length $z$ bounding a pair of pants with these boundaries is at most $x$.
Because $\mathcal R(x,y,z)$ decreases in $z$, each term with $z\le L$
is at least $\mathcal R(x,y,L)$.  The number of such geodesics is
therefore at most $x/\mathcal R(x,y,L)$.  Together with the estimate
\begin{equation}\label{pl:mcshane-ratio}
 \frac{x}{\mathcal R(x,y,z)}
 \le100(1+x)\left(1+e^{(z-x-y)/2}\right)
\end{equation}
this proves \eqref{l:mcshane-end} by setting $z=L$.  The identity
and estimate are given in \cite[Theorem~2.11 and Lemma~2.12]{HSWX}.  The general identity applies
to $S_{0,5}$, so the bound holds for each fixed outer-boundary pair
on the core.  Zero-length endpoints in the integrations are understood
by limits and have measure zero.

Let $\mathcal J_{\rm pair}(L,b)$ be the weighted mass
defined as in \eqref{l:planar-kernel}, with $N^\star_{L,b}$ restricted to triples for
which at least one anchor--witness pair fills the whole core.
Put $T=t_1+\cdots+t_5$.  Such a pair has total length at most $2L$,
so the perimeter bound gives $T\le4L$.  Equations
\eqref{l:whole-pair-count} and \eqref{l:end-count} bound the number of
these triples on each metric by
\begin{equation}\label{pl:whole-pair-pointwise}
 C_\delta(1+L)^2\mathbf1_{\{T\le4L\}}
 e^{5L/2-(1-\delta)T/2}.
\end{equation}
The remaining mark is an end figure-eight.  Dropping its disjointness
and crossing requirements relative to the ordered filling pair only
enlarges the count.  The two possible choices of the filling anchor
contribute a factor two, absorbed into $C_\delta$.

The Weil--Petersson volume $V_{0,5}(\mathbf t)$ is a polynomial
of total degree four in the boundary lengths, so it is bounded by
$C(1+L)^4$ on $T\le4L$.  We use the normalization $V_{0,3}=1$,
as in \cite[Theorem~2.3]{HSWX}.  The estimate
$t\psi(t)=2\sinh(t/2)\le e^{t/2}$ bounds the product of the five
outer factors by $e^{T/2}$.  Combining this with the boundary factor
in \eqref{pl:whole-pair-pointwise} leaves $e^{\delta T/2}$.  Hence
\begin{align}
 \mathcal J_{\rm pair}(L,b)
 &\le C_\delta(1+L)^6e^{5L/2}
       \int_{\substack{t_a\ge0\\T\le4L}}e^{\delta T/2}\,d\mathbf t
 \notag\\
 &\le C_\delta(1+L)^{11}e^{(5/2+2\delta)L}.
 \label{pl:pair-mass}
\end{align}
Here $e^{\delta T/2}\le e^{2\delta L}$ and the five-dimensional
simplex has volume $(4L)^5/5!$.  The power eleven combines two powers
from counting, four from the volume polynomial, and five from the
outer integration.

\subsection{The end-witness integral}

Let $\mathcal J_{\rm end}(L,b)$ be the weighted mass
defined as in \eqref{l:planar-kernel}, with $N^\star_{L,b}$ restricted to triples whose
witness has an end-pants support.  It suffices to bound one outer-label
assignment in \eqref{l:marked-planar-decomposition}.  Write the four
outer lengths belonging to the marked pants as $a_1,a_2,c_1,c_2$, their
internal cuff lengths as $u,v$, and the fifth outer length as $e$.
The two marked figure-eight length constraints imply
\begin{equation}\label{pl:marked-domain}
 a_1,a_2,c_1,c_2,u,v\le L,\qquad
 a_1+a_2+u\le2L,\qquad c_1+c_2+v\le2L.
\end{equation}
Every counted witness enters both marked pants, so
\eqref{l:fifth-cuff} gives $0\le e\le u+v+L$.  For each metric, the
number of eligible end witnesses is at most
$C(1+L)(1+e^{L/2-b})$ by \eqref{l:end-count}.  This bound is uniform
in the marked twist parameters, so it may be applied before integrating
them.  Unfolding the two disjoint marked internal cuffs gives the
length--twist measure $uv\,du\,dv$, with one factor for each cuff.
Cutting along them leaves three pants, each with volume $V_{0,3}=1$.

The bound on the fifth cuff still involves $u$ and $v$. We retain this
dependence when integrating its length, so that the two marked
perimeter bounds can be applied together. Using $t\psi(t)\le e^{t/2}$
for each outer factor gives
\begin{equation}\label{pl:fifth-integration}
 \int_0^{u+v+L}e^{e/2}\,de\le2e^{(u+v+L)/2}.
\end{equation}
The marked perimeter inequalities in \eqref{pl:marked-domain} yield
\begin{equation}\label{pl:perimeter-cancellation}
 e^{(a_1+a_2+c_1+c_2)/2}e^{(u+v+L)/2}\le e^{5L/2}.
\end{equation}
After applying this bound, enlarge the remaining domain to $[0,L]^6$.
The four outer length integrations contribute $L^4$, while the two
internal length integrations give
\begin{equation}\label{pl:internal-length-integrals}
 \int_0^L\int_0^Luv\,du\,dv=\frac{L^4}{4}.
\end{equation}
These six integrations contribute $L^8/4$.  Together with the factor
$1+L$ in the witness count, they give
\begin{equation}\label{pl:end-mass}
 \mathcal J_{\rm end}(L,b)
 \le C(1+L)^9\left(e^{5L/2}+e^{3L-b}\right).
\end{equation}
If an end witness uses the fifth outer cuff, then $b\le e\le L$ as well.
Dropping this restriction gives an upper bound since the integrand is
nonnegative and \eqref{l:end-count} holds on every metric.
By the geometric cover in the proof of
Lemma~\ref{lem:planar-mass}, every counted triple belongs to at least
one of the two families.  Thus
$\mathcal J^\star_{05}(L,b)\le
\mathcal J_{\rm pair}(L,b)+\mathcal J_{\rm end}(L,b)$, whether or not
the families overlap.  Equations~\eqref{pl:pair-mass} and
\eqref{pl:end-mass} therefore prove \eqref{l:planar-mass}.

\subsection{Exterior marks and all complementary genera}

Fix a nonnegative integer $r$. The core $U\cong S_{0,5}$ has complexity
$c=3$ and five outer cuffs. The local count in \eqref{l:planar-kernel}
depends only on the metric of $U$, so it is independent of the outer
gluing twists. All $r$ additional marked pants are exterior, with
pairwise disjoint closed supports. Thus in Lemma~\ref{lem:attachment-bounds}
we have $p=k=d=0$ and $K=K^\star_{L,b}$, with
$\mathcal J(K)=\mathcal J^\star_{05}(L,b)$.
Equation~\eqref{att:all-complements-bound} sums over all complementary
topologies and gives \eqref{l:planar-exterior-incidence} for $L\ge1$,
$b\ge0$, and $g$ sufficiently large in terms of $r$.

At the critical cutoff $H=\log g$ and $L=H-\log H+S$, with $S$ fixed,
we have $L/H\to1$ and
\begin{equation}\label{pl:critical-scales}
 \mathcal A_g(L)=e^S\frac LH=O_S(1),\qquad
 \frac{e^{3L}}{Lg^3}=\frac{e^{3S}}{LH^3}=O_S(H^{-4}).
\end{equation}
Taking $b=10\log L$ and applying
\eqref{l:planar-critical-mass} therefore gives the raw bound
$O_{r,S}(H^{-4})$.

% END sections/appendix_c_planar.tex

% BEGIN sections/appendix_e_sources.tex
\section{Conventions in the imported estimates}\label{app:sources}

We record how the counting conventions and parameter ranges in the
cited sources enter the proofs. The final subsection derives the
identified four-holed-sphere estimate for ordered cuff triples from
the bounds for actual pairs of curves in Section~\ref{sec:tails}.

\subsection{Integration and boundary volumes}
For ambient genus $g\geq3$, \cite[Theorem~2.5 and Remark~2.7]{AMI}
gives the integration formula for ordered oriented multicurves.  Orient
each selected cuff with its marked pair of pants on the left.  The gluing
coefficient is then one, and the two orders of the winding cuffs give
a divisor of two for each actual figure-eight.  Complementary components
are indexed by their cuff-label subsets, so equal numerical lengths
create no additional symmetry divisor.

In the exact disjoint-pants formula of Section~\ref{sec:disjoint},
every complementary component has at least two cuffs.  Hence no
$V_{1,1}$ factor occurs, and the source's special half-volume convention
for the one-holed torus has no effect.  In the fixed-type upper bounds that
allow one-boundary components, the corresponding convention factors
are absorbed into the implicit constants.

The boundary-profile estimate \cite[Lemma~22]{NWX}, in the third
arXiv version listed in the bibliography, applies to positive genus and
at least one boundary, for all nonnegative boundary lengths.  The
finite polynomial bound in Appendix~\ref{app:genus} supplies the
corresponding estimate for the stable genus-zero types that occur here.
For each fixed number of boundaries,
\cite[Theorem~1.2]{MZ} gives the large-genus volume asymptotic with a
universal leading constant.  Its inverse square-root genus factor is
retained in the sums over complementary genera.

\subsection{Actual curves, source labels, and low complexity}
The geodesics counted in \cite[Section~2.5]{HSWX} are primitive and
unoriented.  Its first-moment formula, Proposition~3.1, counts actual
figure-eights, whereas its overlap estimates count ordered cuff triples
satisfying relaxed length conditions.  Each retained actual figure-eight
determines two such triples by ordering its winding cuffs.  An unordered
pair of such curves on distinct supports therefore gives exactly $2\cdot2\cdot2=8$
ordered source pairs.  As shown in
Lemma~\ref{lem:f-eight-transfer}, this passage from actual curves to
source triples preserves the length cutoff.

For the relaxed touching count $D^{\mathrm{rel}}(X,L)$ in
\eqref{f:eight-transfer}, let $b=10\log L$ and let $\Sigma_1,\Sigma_2$
be the representative-tuple sums in \cite[equation~(4.63)]{HSWX}.
With one shared cuff of length $t\in[b,L]$, the other two cuffs of
each pant have combined length at most $2L-b$, since its perimeter is
at most $2L$. With two shared cuffs, all four cuff lengths are at most
$L$. These are the respective length constraints in the two sums.
Each unordered support pair underlying the count $D^{\mathrm{rel}}$ has
such a representative, and the representative determines both supports.
The two support orders and up to $3!$ cuff labelings per support give
\begin{equation}
 D^{\mathrm{rel}}(X,L)\leq72\bigl(\Sigma_1(X,L)+\Sigma_2(X,L)\bigr).
 \label{f:touching-relabeling}
\end{equation}
The integration bounds \cite[equations~(4.64)--(4.65)]{HSWX} give,
for $L>1$ and large $g$,
\begin{equation}
 \mathbb E_g\Sigma_1\leq C\frac{e^{2L}}{g^2L^6},\qquad
 \mathbb E_g\Sigma_2\leq C\frac{L^4e^L}{g^2}.
 \label{f:touching-representative-estimates}
\end{equation}
Since $L^{10}e^{-L}$ is bounded for $L>1$, the two terms give the
relaxed-count estimate of \cite[Proposition~4.17]{HSWX},
\begin{equation}
 \mathbb E_gD^{\mathrm{rel}}(X,L)
 \leq C\frac{e^{2L}}{g^2L^6}.
 \label{f:touching-source-estimate}
\end{equation}

The formula in Section~\ref{sec:disjoint} sums all stable
complementary realizations allowed by the individual good-support
conditions.  It therefore does not use the connected-joint-complement
reduction in the proof of \cite[Proposition~4.4]{HSWX}.
For pair cores, Section~\ref{sec:foundations} uses the perimeter bound
$4L$ and chooses the subsurface parameter accordingly.  The stronger
$2L$ cutoff used in \cite[Lemma~4.5]{HSWX} is not needed.

In the auxiliary-genus proof of Lemma~\ref{w:low-local-mass},
\eqref{tail:other-small} bounds the actual pairs of genus-one and
identified four-holed-sphere type at the same length cutoff.
For the embedded four-holed sphere, Lemma~\ref{tail:planar} gives the
actual-pair estimate at large $L$, and \cite[Proposition~4.16]{HSWX}
covers bounded $L$.
We fix the length cutoff before applying Fatou's lemma as the auxiliary
genus tends to infinity. The resulting constant is uniform in the
cutoff, so the bound also applies at the critical length.

The planar calculation in Appendix~\ref{app:planar} uses the general
McShane--Mirzakhani identity and summand estimate in
\cite[Theorem~2.11 and Lemma~2.12]{HSWX}.  These apply directly to the
five-holed sphere.

\subsection{Fixed types and growing cores}
For a compact hyperbolic surface $Y$ with nonempty geodesic boundary,
\cite[Definition~17 and Theorem~18]{WXPGT} counts ordered tuples of
nontrivial, nonperipheral, primitive unoriented closed geodesics that fill
$Y$.  The length parameter $R>0$ bounds the sum of the curve lengths.
For fixed complexity $m=|\chi(Y)|\geq1$, fixed integer tuple size
$q\geq1$, and fixed $0<\delta<1/2$, the number of such tuples is
bounded by
\begin{equation}\label{s:tuple-bound}
C(1+R)^{q-1}\exp\bigl(R-(1-\delta)\ell(\partial Y)/2\bigr).
\end{equation}
The constant depends only on $m,q,\delta$ and is uniform in the
metric, boundary lengths, and $R$.  We use this estimate at fixed
complexity, and a separate argument when the core complexity grows.

The growing-core argument uses
\cite[Section~7.1, equation~(21), the following definition, and
Proposition~33]{WX}.  Fix $a>0$ and put $T=a\log g$.  Its family
consists of connected geodesic-boundary subsurfaces with perimeter at
most $2T$ and area at most $4T$.  Distinct intrinsic boundary circles
may have the same ambient geodesic image.  For every fixed integer
$K\geq1$ and sufficiently large $g$, Proposition~33, divided by $V_g$,
bounds the expected sum of $e^{-\ell(\partial Y)/4}$ over family
members with $|\chi(Y)|\geq K$ by $C_{a,K}Te^{7T/2}g^{-K}$, as
recorded in \eqref{f:subsurface-tail}.  The area condition gives the
upper cutoff $|\chi(Y)|\leq\lfloor4T/(2\pi)\rfloor$.
For completed images, the exposed perimeter counts the distinct outer
cuffs. Identification seams lie inside the completed image.

The doubling arguments in Sections~\ref{sec:foundations}
and~\ref{sec:higher} preserve primitivity, distinct unoriented classes,
and lengths. They use \cite[Theorem~2.8]{HSWX}: on every closed
genus-$G$ surface, $G\geq2$, the primitive unoriented count below
$L>0$ is at most $(G-1)e^{L+7}$. Additional marks are bounded
deterministically on the ambient surface before taking expectation.

\subsection{Point measures and the minimum}
Section~\ref{sec:probability} applies
\cite[Section~1, Theorem~1.1 and Lemma~1.2]{Kallenberg} to the
countable dissecting-covering semiring of rational half-open intervals.
Their boundaries are null for the diffuse target intensity; finite
joint count convergence and first-moment tightness are proved there.
Left half-lines require the separate tail truncation.

Buser's theorem~\cite[Theorem~4.2.4]{Buser2010} identifies the
non-simple systole with the shortest figure-eight length, with repeated
traversals of simple geodesics excluded as in
Section~\ref{sec:introduction}. The same statement is recalled
in~\cite[Introduction and Section~2.4]{HSWX}. The higher arrival results concern
the ordered figure-eight lengths only.

\subsection{Uniform tails and exponential moments}
For the negative tail, Section~\ref{sec:tails} uses
\cite[Proposition~3.1]{HSWX} with its error constant independent of
both $g\geq3$ and $L\geq2\operatorname{arccosh}3$. Below this length
threshold the figure-eight count vanishes. For the positive tail,
the deterministic upper bound in \cite[Lemma~2.7]{HSWX} controls
lengths beyond the observation window.

The small-core estimates keep the lower cuff threshold $b$ symbolic.
Lemma~\ref{tail:planar} uses the planar domains and integrals in
\cite[Lemmas~4.13--4.15 and equations~(4.57)--(4.61)]{HSWX}.
For the embedded genus-one core, the two witness families and the
exterior volume bound come from \cite[Lemmas~4.9--4.10 and
equation~(4.41)]{HSWX}.
For an identified four-holed-sphere core,
Lemmas~\ref{tail:witness} and~\ref{tail:small} prove the required
cutoff-dependent estimate by a completion-and-witness argument.
Lemma~\ref{s:identified-pair} applies this estimate to the ordered cuff
triples of \cite[Proposition~4.12]{HSWX}.

For touching supports, Lemma~\ref{tail:touch} applies
Proposition~\ref{prop:touching} with one marked curve and arbitrary
cuff cutoff $0\le b\le L$.

For the moving window, Section~\ref{sec:tails} uses $T=3H$,
$H=\log g$, in \cite[Propositions~33--34]{WX}, with fixed complexities
$3$--$14$, tail threshold $15$, and actual perimeter at most $4L$.
These bounds also apply with $L=a_g+t$ and cuff cutoff $b=L/5$
throughout $0\le t\le H/20$. The resulting positive-tail estimate
and the deterministic length bound give an exponential bound uniform
over all genera $g\ge3$.

\subsection{Role marks and a singular shape limit}
Lemma~\ref{mk1:lem:integration} proves the one-point integration
formula for nonnegative Borel weights in the ordered cuff roles, and
for signed or complex weights with finite absolute integral. The proof
uses a fair ordering of the two winding cuffs. The ordered integration
orbit has coefficient one, and its six cuff labelings represent three
actual figure-eights with two winding orders each. Averaging the two
orders gives the coefficient $1/2$. This argument extends the counting
calculation in~\cite[equation~(3.2)]{HSWX} to weights that need not be
invariant under cuff permutations.

Conditional on the surface, each curve receives an independent
winding order that is used at every cutoff. Section~\ref{sec:marked-factorials}
forms the integer-valued factorial counts before averaging the orders.
The first-intensity and fixed-order connected-complement formulas have
complementary genera $g-2$ and $g-2m$ with fixed boundary counts,
so their genera are eventually positive. Disconnected complements are
estimated separately.

For the marked vague limit, coordinate boxes with endpoints avoiding
$1/2$ have $\nu$-null relative boundaries: the first two shape
coordinates have continuous distributions and the third equals $1/2$. Their products
with rational half-open length intervals give the semiring used in
Section~\ref{sec:marked-process}; half-line counts and arrivals also
require its tail truncations. Fixed continuity-cell vectors converge
in total variation. The total-variation distance between the full process
laws equals one for every $g\ge3$ by Proposition~\ref{mk3:prop:tv-fails}.

\subsection{Proof of the identified four-holed-sphere estimate}

We deduce the identified-boundary estimate of \cite[Proposition~4.12]{HSWX}.  Its count $C^1_{0,4}(X,L)$ concerns
ordered pairs of ordered cuff triples from \cite[Definition~4.2]{HSWX},
whose intrinsic filled core is $S_{0,4}$ and whose completed image is
$S_{1,2}$.  Each triple $(c_1,c_2,c_3)$ bounds an embedded pair of pants
with three distinct ambient cuff geodesics and connected complement,
and satisfies
\[
 \ell(c_1)\le L,\qquad \ell(c_2)+\ell(c_3)\le L,\qquad
 \ell(c_j)\ge h:=10\log L\quad(1\le j\le3).
\]
These cuff inequalities define the relaxed count in the source.

\begin{lemma}[Identified four-holed-sphere pairs]\label{s:identified-pair}
There are absolute constants $C<\infty$ and $g_0\ge3$ such that, for
every integer $g\ge g_0$ and every $L>1$,
\begin{equation}\label{s:identified-pair-bound}
 \mathbb E_g C^1_{0,4}(X,L)\le Cg^{-2}e^{2L}.
\end{equation}
\end{lemma}

\begin{proof}
If $h>L$, the source count vanishes. Otherwise set
$L'=L+2\operatorname{arccosh}3$. For each triple let $\gamma$ be the
figure-eight with winding cuffs $c_2,c_3$, and write
$x=\ell(c_2)$, $y=\ell(c_3)$ and $z=\ell(c_1)$.
Since $z\le L$ and $|x-y|\le x+y\le L$, the length identity
\eqref{f:pants-length} gives
\[
 \cosh\frac{\ell(\gamma)}2
 =\cosh\frac z2+\cosh\frac{x+y}2+\cosh\frac{x-y}2
 \le 3\cosh\frac L2\le\cosh\frac{L'}2.
\]
Thus both curves have length at most $L'$ and retain their good
supports, cuff lower bound $h$ and intrinsic filled core with one
identified boundary pair. They form an ordered pair counted by
$Q_{04^1}(L',h)$. Each curve determines its supporting pants and its
unordered winding cuffs, so the two winding orders in each slot give
at most four preimages. Hence
\[
 C^1_{0,4}(X,L)\le 4Q_{04^1}(L',h).
\]

Apply the absolute estimate~\eqref{tail:other-small} of
Lemma~\ref{tail:small} with $0<h\le L'$. Its leading term is
$g^{-2}(L')^5e^{2L'-h/2}$. Since $e^{-h/2}=L^{-5}$ and $L'-L$ is
constant, this is $O(g^{-2}e^{2L})$ for $L>1$. Each remaining term
is a fixed power of $L'$ times $g^{-2}e^{\alpha L'}$ with
$\alpha\le3/2$, hence obeys the same bound. The genus threshold in
\eqref{tail:other-small} is independent of $L'$ and $h$.
\end{proof}

% END sections/appendix_e_sources.tex

% BEGIN references.tex

% END references.tex

\end{document}